\documentclass[preprint,12pt]{elsarticle}
\journal{ }
\date{}

\usepackage{graphicx, amsmath, amsfonts, amssymb}
\usepackage{hyperref}
\usepackage{subfigure}
\usepackage{algorithm,algorithmic}
\usepackage{mathtools,bm}
\usepackage{multirow,booktabs,makecell}
\usepackage{tikz,url}
\usetikzlibrary{shapes,arrows}
\tikzstyle{process} = [rectangle,minimum width=2cm,minimum height=1cm,text centered,text width =4cm,draw=black]

\newcommand{\deri}{\text{d}}

\newcommand{\calC}{\mathcal{C}}

\newcommand{\bbR}{\mathbb{R}}

\newcommand{\bbRplus}{\bbR^+}
\newcommand{\bbRbar}{\overline{\bbR}}

\newcommand{\bmx}{\bm{x}}

\newcommand{\bmv}{\bm{v}}
\newcommand{\bmb}{\bm{\beta}}
\newcommand{\bmm}{\bm{m}}
\newcommand{\bma}{\bm{\alpha}}

\newcommand{\bmn}{\bm{n}}

\newcommand{\bmO}{\bm{0}}
\newcommand{\Th}{\mathcal{T}_h}
\newcommand{\Ih}{\mathcal{I}_h}

\newcommand{\pol}{\mathcal{P}}

\DeclareMathOperator{\divg}{div}

\newtheorem{theorem}{Theorem}[section]

\newtheorem{remark}[theorem]{Remark}

\numberwithin{equation}{section}

\begin{document}
\begin{frontmatter}

\title{High order computation of optimal transport, mean field planning, and mean field games}
\author[1]{Guosheng Fu\fnref{fn1}}
\ead{gfu@nd.edu}

\author[2]{Siting Liu\fnref{fn2}}
\ead{siting6@math.ucla.edu}

\author[2]{Stanley Osher\fnref{fn2}}
\ead{sjo@math.ucla.edu}

\author[3]{Wuchen Li\fnref{fn3}}
\ead{wuchen@mailbox.sc.edu}

\cortext[cor]{Corresponding author}
\fntext[fn1]{G. Fu's work is supported in part by NSF DMS-2134168.}
\fntext[fn2]{S. Osher and S. Liu's work are supported in part by AFOSR MURI FP 9550-18-1-502, and ONR
grants:  N00014-20-1-2093, and N00014-20-1-2787.}

\fntext[fn3]{W. Li's work is supported by AFOSR MURI FP 9550-18-1-502, AFOSR YIP award 2023, and NSF RTG: 2038080.}

\affiliation[1]{organization={Department of Applied and Computational Mathematics and Statistics, University of Notre Dame},
                city={Notre Dame},
                postcode={IN 46556},
                country={USA}}
\affiliation[2]{organization={Department of Mathematics, University of California, Los Angeles},
                city={Los Angeles},
                postcode={CA 90095},
                country={USA}}     \affiliation[3]{organization={Department of Mathematics, University of South Carolina},
                city={Columbia},
                postcode={SC 29208},
                country={USA}}          

\begin{abstract}
Mean-field games (MFGs) have shown strong modeling capabilities for large systems in various fields, driving growth in computational methods for mean-field game problems. However, high order methods have not been thoroughly investigated. In this work, we explore applying general high-order numerical schemes with finite element methods in the space-time domain for computing the optimal transport (OT), mean-field planning (MFP), and MFG problems. We conduct several experiments to validate the convergence rate of the high order method numerically. Those numerical experiments also demonstrate the efficiency and effectiveness of our approach.
\end{abstract}

\end{frontmatter}


\section{Introduction}
\label{sec:intro}

Proposed by Lasry and Lions~\cite{lasry2007mean} and independently by Caines, Huang, and Malham{\'e}~\cite{huang2006large}, the mean-field game models an infinite number of identical agents' interactions in a mean-field manner, characterizing the equilibrium state of the system.
Thanks to its substantial descriptive ability, the MFG becomes an important approach to studying complex systems with large populations of interacting agents,  such as crowd dynamics, financial markets, power systems, pandemics, etc.\cite{carmona2015mean, casgrain2020mean,lachapelle2011mean,aurell2018mean,bagagiolo2014mean,kizilkale2019integral,lee2021controlling,lee2022mean}.
The mean-field planning is a class of MFGs where the distribution of agents at terminal time is imposed~\cite{porretta2013planning}.
On the other hand, the Benamou-Brenier dynamic formulation~\cite{benamou2000computational}  of the optimal transport problem connects with the variational form of the potential mean-field games.
It can be treated as a special case of the mean-field planning problem, which aims to find an efficient way of moving one probability distribution to another.
Along with the empirical success of MFG and OT in modeling and real-world applications, the study of mean-field game is also expanding.
From the PDE view, the mean-field game model can be described by a system of coupled partial differential equations: a forward-in-time Fokker-Planck (FP) equation governs the evolution of the population and a backward-in-time Hamilton-Jacobi-Bellman (HJB) equation for the value function that characterizes the control problem. For a review of MFG theory, we refer to\cite{LasryLions2007,gueantlasrylions11,CardaNotes,gomsaude'14}.

With such a wide range of applications, computational methods play a crucial role since most MFG and OT problems do not have analytical solutions.
While some recent computational approaches take advantage of machine learning methods and game theories~\cite{ruthotto2020machine,lin2021alternating,agrawal2022random, cui2021approximately,lauriere2022scalable,guo2022general,cardaliaguet2017learning,hadikhanloo2019finite}, classical numerical methods are mostly developed discretization using finite difference schemes or semi-Lagrangian schemes. 
In~\cite{achdou2010mean}, the MFG system is discretized using finite difference scheme and then solved by Newton's method. Semi-Lagrangian methods are studied in~\cite{carlini2014fully}. As for MFGs and OTs that can be written in a variational form, optimization methods, such as augmented Lagrangian, Primal-dual Hybrid Gradient, Alternating Direction Method of Multipliers, are applied to solve the discretized system\cite{benamou2000computational,benamou2014augmented, achdou2016mean,bencarsan'17, silva19,liu2021computational}.
Recently, computation of MFGs on mainfolds has been investigated in~\cite{yu2022computational}. 
For the survey of the numerical methods, we refer to~\cite{achdou2020mean,lauriere2021numerical}.
Within the augmented Lagrangian framework, the (low-order) finite element discretization has been used frequently; see, e.g.,  \cite{benamou2000computational,benamou2014augmented, Andreev17,Igbida18}.

Pioneering works on computational OT/MFG focus on first or second order methods; the general high order method is not well studied. Yet, high order methods generally have faster convergence rates in numerical analysis and provide more accurate solutions on a much coarse computational mesh than low order methods.  
Therefore, exploring high order computational methods for mean-field games and optimal transport problems is vital.

In this work, we propose a general high order numerical method for solving the optimal transportation problem and mean-field game (control) problems using the finite element method. More precisely,
\begin{enumerate}
    \item We discretize the augmented Lagrangian formulations of the  MFP and MFG systems  using high-order space-time finite elements. Considering derivation information used in the saddle point formulation, we approximate the value function (dual variable) $\phi$ using high order $H^1$-conforming  finite elements, while the density and momentum (primal variables) $\rho,\bmm$ are approximated via a high-order (discontinuous) integration rule space which only records values on the high-order (space-time) integration points. Our discrete saddle-point problem is then solved via the ALG2 algorithm, following~\cite{benamou2014augmented}.
    To the best of our knowledge,  this is the first time high-order schemes with more than second order accuracy being applied.
    \item 
   We present a series of comprehensive experiments to showcase the efficacy and efficiency of the proposed numerical algorithms. These experiments numerically validate the convergence rate of the algorithms as a function of mesh size and polynomial degree. In particular, we show a high-order method on a coarse mesh is more accurate than a low order method on a fine mesh with the same number of degrees of freedom.
   Furthermore, we apply the finite element scheme to a set of mean-field planning and mean-field game problems on non-rectangular domains (with obstacles) and computational graphics, demonstrating the validity and practicality of our method.
\end{enumerate}

This paper is organized as follows. Section 2 review the dynamic formulation of optimal transportation, mean-field planning, and mean-field games.
Section 3 presents the high-order schemes we designed for computing the above problems and the companion algorithm. Section 4 demonstrates the effectiveness of the high-order method with numerical experiments. Finally, we make some conclusions and remarks in Section 5.

\section{OT, MFP, and MFG}
\label{sec:model}

In this section, we briefly review dynamic MFP and MFG problems.

\subsection{Dynamic MFP}
Consider the model on time interval $[0,1]$ and space region $\Omega\subset\bbR^D$.
Let $\rho$ be the density of agents through $t\in[0,1]$, $\bmm$ be the flux of the density which models strategies (control) of the agents, and $(\rho,\bmm)\in\calC$:
\begin{equation}
\calC:=\left\{\begin{aligned}
(\rho,\bmm): 
&\rho:[0,1]\times\Omega\to\bbR^+,\|\rho\|_{L^1} <+\infty,
    \int_\Omega\rho(t,\bmx)\deri\bmx=1,\forall t\in[0,1],\\
&\bmm:[0,1]\times\Omega\to\bbR^D \text{ is Lebesgue measurable}, 
\end{aligned}\right\}.
\label{opt:cst}
\end{equation}
We are interested in $\rho$ with given initial and terminal density $\rho_0,\rho_1$ and
$(\rho,\bmm)$ satisfying zero boundary flux and mass conservation law, which satisfies the constraint set $\calC(\rho_0,\rho_1)$:
\begin{equation}
\calC(\rho_0,\rho_1):=\calC\cap\left\{\begin{aligned}
(\rho,\bmm):&\partial_t\rho + \divg_{\bmx} \bmm=0,\\
& \bmm\cdot\bmn = 0\text{ for }\bmx\in\partial\Omega, \rho(0,\cdot) = \rho_0, \rho(1,\cdot) = \rho_1,
\end{aligned}\right\}.
\label{mfp:cst}
\end{equation}
where the equation hold in the sense of distribution.

We denote $L:\bbRplus\times\bbR^D\to\bbRbar:=\bbR\cup\{\infty\}$ as the dynamic cost function and $A:\bbR\to\bbRbar$ as a function modeling interaction cost.
The goal of MFP is to minimize the total cost among all feasible $(\rho,\bmm)\in\calC(\rho_0,\rho_1).$ Therefore, the problem can be formulated as 
\begin{equation}
    \inf_{(\rho,\bmm)\in\calC(\rho_0,\rho_1)} \int_0^1\int_\Omega L(\rho(t,\bmx),\bmm(t,\bmx))
    +A(\rho(t,\bmx))\deri \bmx \deri t.
\label{mfp:problem}
\end{equation}

It is clear to see $\calC(\rho_0,\rho_1)$ is convex and compact. In addition, the mass conservation law $\partial_t\rho+\divg_{\bmx}\bmm=0$ and zero flux boundary condition $\bmm\cdot\bmn=0,\bmx\in\partial\Omega$ imply that $\calC(\rho_0,\rho_1)\neq\emptyset$ if and only if $\int_\Omega\rho_0=\int_\Omega\rho_1.$ 
Once $\calC(\rho_0,\rho_1)$ is non-empty, the existence and uniqueness of the optimizer depends on $L$ and $F$.
In this paper, we consider a typical dynamic cost function $L$ by \begin{equation}
    L(\beta_0,\bmb_1) :=  \begin{cases} 
    \frac{\|\bmb_1\|^2}{2\beta_0} &\text{ if } \beta_0>0\\
    0              &\text{ if } \beta_0=0,\bmb_1=\bmO\\
    +\infty        &\text{ if } \beta_0=0, \bmb_1\neq\bmO.
    \end{cases}.
\label{dynamic:L}
\end{equation}
Various choices of the interaction function will be given in the numerics section.

If the interaction cost function $A=0$, the MFP becomes the dynamic formulation of optimal transport problem:
\begin{equation}
  \text{(OT) }  \min_{\rho,\bmm\in\calC(\rho_0,\rho_1)} \int_0^1\int_\Omega L(\rho(t,\bmx),\bmm(t,\bmx)) \deri \bmx \deri t. 
\label{ot:problem}
\end{equation}
Since $\bmm=\rho\bmv$, this definition of $L$ makes sure that $\bmm=\bmO$ wherever $\rho=0$. 
OT can be viewed as a special case of MFP where masses move freely in $\Omega$ through $t\in[0,1]$.

To simplify notation, we denote 
an element in the set $\calC$ as 
$\bma:=(\alpha_0, \bma_1)\in \calC$.
Introducing the Lagrangian multiplier $\phi:[0,1]\times \Omega$ for the constraint \eqref{mfp:cst}, the 
MFP problem \eqref{mfp:problem} can be reformulated as the following saddle-point problem:
\begin{subequations}
\label{pd-mfp:problem}
\begin{equation}
    \inf_{\bma} 
    \sup_{\phi}
   F(\bma)
   -
   G(\phi)-\langle\bma, \nabla_{t,x}\phi\rangle,
\label{pd-mfp:problem1}
\end{equation}
where 
\begin{align}
\label{f-form}
   F(\bma):=&\;
       \int_0^1\int_\Omega 
   L(\alpha_0(t,\bmx),\bma_1(t,\bmx))
   +A(\alpha_0(t,\bmx))\deri \bmx \deri t ,\\
\label{g-form}
   G(\phi):=&\;
\int_\Omega -\phi(1,\bmx)\rho_1(\bmx)+\phi(0,\bmx)
\rho_0(\bmx)\deri \bmx,
\end{align}
\end{subequations}
$\nabla_{t,x}=(\partial_t, \mathrm{grad}_{\bmx})$ is the space-time gradient operator, and 
$\langle\bma,\bmb\rangle:=\int_0^1\int_{\Omega}
\bma\cdot\bmb\,\deri \bmx \deri t$ is the space-time integral.
The KKT system for this saddle-point system (with cost function $L$ in \eqref{dynamic:L}) is the following PDE system on the space-time domain $[0,1]\times \Omega$
\begin{subequations}
\label{pd-mfp:KKT}
\begin{align}
\partial_t\rho + \divg_{\bmx} \bmm=&\;0,\\
\frac{\bmm}{\rho} - \mathrm{grad}_{\bmx}\phi=&\;0,\\
\partial_t \phi + \frac{|\bmm|^2}{2\rho^2} =&\; A'(\rho),
\end{align}
with boundary conditions
\begin{align}
\label{mfp-bcm}
\bmm\cdot\bmn = &\;0, \quad \text{ on }[0,1]\times \partial\Omega,\\
\label{mfp-bcr}
\rho(0,\bmx) = &\; \rho_0(\bmx),\quad 
\rho(1,\bmx) = \; \rho_1(\bmx) \quad\text{ on }\Omega.
\end{align}
\end{subequations}

Denoting $[L+ A]^*(\bma^*)$ as
the convex conjugate (Legendre transformation) of $L(\alpha_0, \bma_1)+ A(\alpha_0)$ with $L$ given in \eqref{dynamic:L},
i.e., 
\begin{align*}
[L+ A]^*(\bma^*) = &\;
\sup_{\bma}\bma\cdot\bma^*-L(\alpha_0, \bma_1)-\lambda A(\alpha_0)\\
=&\;
\sup_{\alpha_0}\alpha_0\cdot(\alpha_0^*+|\bma_1^*|^2)-\underbrace{L(\alpha_0, \alpha_0\bma_1^*)}_{=\frac12\alpha_0|\bma_1^*|^2}- A(\alpha_0)\\
=&\;
\sup_{\alpha_0}\alpha_0\cdot(\alpha_0^*+\frac12|\bma_1^*|^2)- A(\alpha_0)\\
=&\;A^*(\alpha_0^*+\frac12|\bma_1^*|^2),
\end{align*}
where in the second equality we used the 
optimality condition $\bma_1 = \alpha_0\bma_1^*$.
By duality, we have 
\begin{align}
\label{duality}
L(\alpha_0, \bma_1)+ A(\alpha_0)
=
\sup_{\bma^*}\bma\cdot\bma^*- A^*(\alpha_0^*+\frac12|\bma_1^*|^2).
\end{align}
Using the above relation, 
we have the following dual formulation of the saddle-point problem \eqref{pd-mfp:problem1}:
\begin{align}
\label{pd-mfp:dual}
    \sup_{\bma} 
    \inf_{\phi, \bma^*}
   F^*(\bma^*)
+G(\phi)+\langle\bma, \nabla_{t,x}\phi-\bma^*\rangle,
\end{align}
where 
\[
   F^*(\bma^*)
=\int_0^1\int_{\Omega}A^*(\alpha_0^*+\frac12|\bma_1^*|^2)\,\deri \bmx \deri t.
\]

Introducing the augmented Lagrangian
\[
L_r(\phi, \bma, \bma^*):= 
   F^*(\bma^*)
+G(\phi)+\langle\bma, \nabla_{t,x}\phi-\bma^*\rangle
+\frac{r}{2}\langle\nabla_{t,x}\phi-\bma^*,\nabla_{t,x}\phi-\bma^*\rangle,
\]
where $r$ is a positive parameter, it is clear that
the corresponding saddle-point problem 
\begin{align}
\label{aug-mfp}
\sup_{\bma} 
    \inf_{\phi, \bma^*}L_r(\phi, \bma, \bma^*)
\end{align}
has the same solution as \eqref{pd-mfp:dual}.

\subsection{Dynamic MFG}
For MFG, the terminal density $\rho_1$ is not explicitly provided but it satisfies a given preference. 
The goal of MFG is to minimize the total cost among all feasible $(\rho,\bmm)\in\calC(\rho_0)$:
\begin{equation}
    \inf_{(\rho,\bmm)\in\calC(\rho_0)} 
F((\rho,\bmm))
+\underbrace{\int_\Omega \Gamma(\rho(1,\bmx))\deri\bmx}_{:=R(\rho(1,\cdot))},
\label{mfg:problem}
\end{equation}
where $\Gamma:\bbR\rightarrow\bbRbar$ is the terminal cost, and the constraint set $\calC(\rho_0)$ is similar to $\calC(\rho_0,\rho_1):$
\begin{equation}
    \calC(\rho_0) 
    := \calC\cap\left\{\begin{aligned}
        (\rho,\bmm):
        &\partial_t\rho + \divg_{\bmx} \bmm=0,\\
        & \bmm\cdot\bmn = 0\text{ for }\bmx\in\partial\Omega, \rho(0,\cdot) = \rho_0,
    \end{aligned}\right\}.
\label{mfg:cst}
\end{equation}

Similar to MFP, we reformulate the problem \eqref{mfg:problem} into a saddle-point problem:
\begin{equation}
    \inf_{\bma, \rho_1} 
    \sup_{\phi}
   F(\bma)
   +R(\rho_1)+
 (\rho_1, \phi(1,\cdot))-
 (\rho_0, \phi(0,\cdot))-
 \langle\bma, \nabla_{t,x}\phi\rangle,
\label{pd-mfg:problem}
\end{equation}
in which 
$(\alpha, \beta):=\int_{\Omega}\alpha\beta\,\deri\bmx$
is the spatial integration.
Here the KKT system of the saddle-point problem \eqref{pd-mfg:problem} is simply the MFP system \eqref{pd-mfp:KKT}
with boundary condition \eqref{mfp-bcr}
replaced by the following:
\[
\rho(0,\bmx) = \; \rho_0(\bmx),\quad 
\phi(1,\bmx) = \; -\Gamma'(\rho_1(\bmx)) \quad\text{ on }\Omega.
\]

Introducing the dual variables $\bma^*$ and $\rho_1^*$
for $\bma$ and $\rho_1$, respectively, 
we get the following equivalent saddle-point problem:
\begin{align}
\label{pd-mfg:dual}
    \sup_{\bma, \rho_1} 
    \inf_{\phi, \bma^*, \rho_1^*}
&\;   F^*(\bma^*)
+\langle\bma, \nabla_{t,x}\phi-\bma^*\rangle\nonumber\\
&\;+R^*(\rho_1^*)-(\rho_1, \phi(1,\cdot)+\rho_1^*)
+(\rho_0, \phi(0,\cdot)),
\end{align}
where 
$
R^*(\rho_1^*):=\int_{\Omega}\,\Gamma^*(\rho_1(\bmx))\deri \bmx,
$
with $\Gamma^*$ being the convex conjugate of $\Gamma$.

The augmented Lagrangian reformulation of \eqref{pd-mfg:dual} is the following:
\begin{align}
\label{aug-mfg}
    \sup_{\bma, \rho_1} 
    \inf_{\phi, \bma^*, \rho_1^*}
&\;   F^*(\bma^*)+R^*(\rho_1^*)+(\rho_0, \phi(0,\cdot))
\nonumber\\
&\;+\langle\bma, \nabla_{t,x}\phi-\bma^*\rangle
+\frac{r_1}{2}
\langle\nabla_{t,x}\phi-\bma^*,\nabla_{t,x}\phi-\bma^*\rangle\nonumber\\
&\;-(\rho_1, \phi(1,\cdot)+\rho_1^*)
+\frac{r_2}{2}
(\phi(1,\cdot)+\rho_1^*,\phi(1,\cdot)+\rho_1^*),
\end{align}
where $r_1,r_2$ are two positive parameters.

\begin{remark}
Following the seminal works in \cite{benamou1999numerical,benamou2014augmented},
we propose our high-order schemes for MFP and MFG based on the augmented Lagrangian formulations \eqref{aug-mfp}.
and \eqref{aug-mfg}.
The discrete saddle-point problem is then solved using the ALG2 algorithm \cite{FortinBook}.
The major novelty of our scheme is the use of high-order 
space-time finite elements for the discretization of the variables in \eqref{aug-mfp} and \eqref{aug-mfg}.
This is the first time high-order schemes with more than second order accuracy being applied
to such problems.
\end{remark}

\section{High-order schemes for OT, MFP and MFG}
\label{sec:alg}
In this section, we discretize the augmented Lagrangian problems \eqref{aug-mfp} and \eqref{aug-mfg} using high-order space-time finite element spaces.
We start with notation including the mesh and definition of finite element spaces to be used.
We then formulate the discrete saddle-point problems using these finite element spaces, which is solved iteratively using the ALG2 algorithm \cite{FortinBook}.
Throughout this section, we restrict the discussion to $D=2$ spatial dimensions.

Since space/time derivative information is needed for $\phi$, we approximate it using (high-order) $H^1$-conforming finite elements. On the other hand, since no derivative information appear for $\bma$, $\bma^*$, 
(and $\rho_1$ and $\rho_1^*$ for MFG), it is natural to  approximate these variables only on the (high-order) integration points.
\subsection{The finite element spaces and notation}
Let $\Ih=\{I_j\}_{j=1}^N$ be a triangulation of the time domain $[0, 1]$ with $I_j=[x_{j-1}, x_{j}]$, and 
$0=x_0<x_1<\cdots<x_N=1$.
Let $\Th=\{T_\ell\}_{\ell=1}^M$ be a conforming triangulation of the 
spatial domain $\Omega$, where we assume the element $T_\ell:=\Phi_{T_\ell}(\widehat{T})$ is obtained from a polynomial mapping $\Phi_{T_\ell}$ from the reference element $\widehat{T}$, which, is a unit triangle or unit square. 
We obtain the space-time mesh for $\Omega_T:=[0,1]\times \Omega$ using tensor product of the spatial and temporal meshes:
\[
\Ih\otimes \Th:=\{I_j\otimes T_\ell: \forall j\in \{1,\cdots, N\}, \text{ and } \ell\in \{1,\cdots, M\}\}.
\]

We denote $\pol^k(I)$ as the polynomial space of degree no greater than $k$ on the interval $I$, and $\pol^k(\widehat{T})$ as the polynomial space of degree no greater than $k$ if $\widehat{T}$
is a unit triangle, or the 
tensor-product polynomial space of degree no greater than $k$ in each direction if $\widehat{T}$
is a unit square, for $k\ge 1$.
The mapped polynomial space on a spatial physical element
$T\in\Th$ 
is denoted as 
\[
\pol^k(T):=
\{\widehat{v}\circ(\Phi_T)^{-1}:\; \forall 
\widehat{v}\in\pol^k(\widehat{T}))
\}.
\]

We denote $\{\widehat{\bm \xi}_i\}_{i=1}^{N_k}$
as a set of quadrature points with positive weights 
$\{\widehat{\omega}_i\}_{i=1}^{N_k}$ that is accurate for polynomials of degree up to $2k+1$ on the reference element $\widehat T$, i.e.,
\begin{align}
\label{quad-2k}
\int_{\widehat T}\widehat f\, \deri\bmx
=\sum_{i=1}^{N_k}\widehat{\omega}_i\widehat f(\widehat{\bm \xi}_i), \quad\forall \widehat f\in\pol^{2k+1}(\widehat T).
\end{align}
Note that when $\widehat T$ is a reference square, we simply use the Gauss-Legendre quadrature rule with $N_k = (k+1)^2$, which is optimal.
On the other hand, when $\widehat T$ is a reference triangle, the optimal choice of quadrature rule is more complicated; see, e.g., 
\cite{Zhang09,Witherden15} and references cited therein.
For example, the number $N_k$ for $0\le k\le 6$ of the symmetric quadrature rules on a triangle provided in \cite{Zhang09} are given in Table \ref{tab:quad}.
\begin{table}[ht!]
\begin{tabular}{cccccccc}
 & $k=0$ & $k=1$ & $k=2$ & $k=3$ & $k=4$& $k=5$ & $k=6$\\\hline
$N_k$ on Triangle 
& 1 & 6 & 7 & 15 & 19&28&37\\
\end{tabular}
\caption{Number of quadrature points $N_k$ for the quadrature rule 
on a triangle that is accurate up to degree $2k+1$ for $0\le k\le 6$.
}
\label{tab:quad}
\end{table}
The integration points and weights on a physical element $T_\ell$ are simply obtained via mapping: $\{\bm \xi_i^{\ell}:=\Phi_{T_\ell}(\widehat{\bm \xi}_i)\}_{i=1}^{N_k}$, and 
$\{\omega_i^{\ell}:=|\mathrm{grad}_{\bmx}\Phi_{T_\ell}(\widehat{\bm \xi}_i)|\widehat \omega_i\}_{i=1}^{N_k}$.
Moreover, we denote $\{\eta_i^j\}_{i=1}^{k+1}$
as the set of ($k+1$) Gauss-Legendre quadrature points 
on the interval $I_j$ with corresponding weights 
$\{\zeta_i^j\}_{i=1}^{k+1}$.
To simplify the notation, 
we denote the set of physical integration points and weights
\begin{subequations}
\label{quad-rule}
\begin{align}
\Xi_h^k:=&\;\{\bm\xi_i^{\ell}:\;\;1\le i\le N_k, \, 1\le \ell\le M\},\\
\Omega_h^k:=&\;\{\omega_i^{\ell}:\;\;1\le i\le N_k, \, 1\le \ell\le M\},\\
H_h^k:=&\;\{\eta_i^{j}:\;\;1\le i\le k+1, \, 1\le j\le N\},\\
Z_h^k:=&\;\{\zeta_i^{j}:\;\;1\le i\le k+1, \, 1\le j\le N\}.
\end{align}
\end{subequations}
Moreover, we denote $(\cdot, \cdot)_h$ as the discrete inner-product on the mesh $\Th$ using the quadrature points $\Xi_h^k$ and weights $\Omega_h^k$:
\[
(\alpha, \beta)_h:=\sum_{\ell=1}^M\sum_{i=1}^{N_k}\alpha(\bm\xi_i^{\ell})
\beta(\bm\xi_i^{\ell})\omega_i^{\ell},
\]
and 
$\langle\cdot, \cdot\rangle_h$ as the discrete inner-product on the space-time mesh $\Ih\otimes \Th$ using the quadrature points $\Xi_h^k$, $H_h^k$ and weights $\Omega_h^k$, $Z_h^k$:
\[
\langle\alpha, \beta\rangle_h:=\sum_{\ell=1}^M\sum_{i_s=1}^{N_k}
\sum_{j=1}^N\sum_{i_t=1}^{k+1}
\alpha(\eta_{i_t}^j, \bm\xi_{i_s}^{\ell})
\beta(\eta_{i_t}^j, \bm\xi_{i_s}^{\ell})\omega_{i_s}^{\ell}
\zeta_{i_t}^{\ell}.
\]

We are now ready to present our finite element spaces: 
\begin{align}
\label{fes-V}
V_h^k:=&\;\{v\in H^1(\Omega_T): \;\; v|_{I_j\times T_\ell}
\in\pol^k(I_j)\otimes\pol^k(T_\ell)\;\;
\forall j, \ell\},\\
\label{fes-W}
W_h^k:=&\;\{w\in L^2(\Omega_T): \;\; w|_{I_j\times T_\ell}
\in \pol^{k}(I_j)\otimes W^k(T_\ell)\;\;
\forall j, \ell\},\\
\label{fes-M}
M_h^k:=&\;\{\mu\in L^2(\Omega):\;\;\; \;\; \mu|_{T_\ell}
\in W^k(T_\ell)\;\;\forall \ell\},
\end{align}
where $V_h^k$ is an $H^1$-conforming space on the 
space-time mesh $\Ih\otimes\Th$, 
$W_h^k$ is an $L^2$-conforming space on the 
space-time mesh $\Ih\otimes\Th$, 
and $M_h^k$ is an $L^2$-conforming space on the 
spacial mesh $\Th$, 
in which the local space 
\[
W^k(T_\ell):=\pol^{k}(T_\ell)\oplus \delta W_k(T_\ell),
\]
is associated with the integration rule in \eqref{quad-2k}
such that $\dim W^k(T_\ell) = N_k$, and the nodal conditions \begin{align}
\label{collocation}
\varphi_i^{\ell}({\bm \xi}_j^{\ell})=
\delta_{ij},\quad \forall 1\le j\le N_k,
\end{align}
in which  $\delta_{ij}$ is the Kronecker delta function determines a unique solution $\varphi_i^\ell\in W^k(T_\ell)$. 
This implies that $\{\varphi_i^\ell\}_{i=1}^{N_k}$ is a set of nodal bases for the space ${W}^k(T_\ell)$, i.e., 
\begin{align}
\label{nodal-basis}
{W}^k(T_\ell)=\mathrm{span}_{1\le i\le N_k}\{\varphi_i^\ell\}.
\end{align}
When $T^\ell$ is mapped from a reference square, 
we have $N_k=(k+1)^2$, hence $W^k(T_\ell)$ is simply the (mapped) tensor product polynomial space $\pol^{k}(T_\ell)$.
Moreover, we emphasize that the explicit expression of the basis function $\phi_i^\ell$ does not matter in our construction, as only their nodal degrees of freedom (DOFs) on the quadrature nodes will enter into the numerical integration.
Furthermore, let $\{\psi_i^j(t)\}_{i=1}^{k+1}$ be 
the set of basis functions for $\pol^{k}(I_j)$
corresponding to the Gauss-Legendre quadrature nodes 
$\{\eta_i^j\}_{i=1}^{k+1}$, i.e., $\psi_i^{j}\in \pol^{k}(I_j)$ satisfies 
\begin{align*}
\psi_i^{j}({\eta}_l^{j})=
\delta_{il},\quad \forall 1\le l\le k+1.
\end{align*}
With this notation by hand, we have 
\begin{align}
\label{W-basis}
W_h^k = \mathrm{span}
\left\{\psi_{i_t}^j(t)\varphi_{i_s}^\ell(\bmx):\;\;
\begin{tabular}{l}
$1\le i_t\le k+1$,
$1\le i_s\le N_k$,\\ 
$1\le j\le N$, 
$1\le \ell\le M$\\
\end{tabular}
\right\}
\end{align}
and 
\begin{align}
\label{M-basis}
M_h^k = \mathrm{span}
\left\{\varphi_{i_s}^\ell(\bmx):\;\;
1\le i_s\le N_k, 1\le \ell\le M
\right\}
\end{align}

We approximate the dual variable $\phi$ using the $H^1$-conforming finite element space $V_h^{k+1}$, 
each components of $\bma$ and $\bma^*$ using the 
integration rule space $W_h^k$, and
the variables $\rho_1$ and $\rho_1^*$ (for MFG) using the 
integration rule space $M_h^k$.

\subsection{High-order FEM for MFP and MFG}
The discrete scheme for MFP \eqref{aug-mfp} reads as follows:
given a space-time mesh $\Ih\otimes \Th$ and a polynomial degree $k\ge 0$, find $\bma_h, \bma_h^*\in [W_h^k]^{3}$, and $\phi_h\in V_h^{k+1}$ such that 
\begin{align}
\label{aug-mfp-h}
\sup_{\bma_h\in [W_h^k]^{3}} \;\;
    \inf_{\phi_h\in V_h^{k+1}, \bma_h^*\in[W_h^k]^{3}}L_{r,h}(\phi_h, \bma_h, \bma_h^*),
\end{align}
where  the discrete augmented Lagrangian is 
\begin{align}
\label{Lrh}
L_{r,h}:= &\;
   F_h^*(\bma_h^*)
+G_h(\phi_h)+\langle\bma_h, \nabla_{t,x}\phi_h-\bma_h^*\rangle_h\nonumber\\
&\;+\frac{r}{2}\langle\nabla_{t,x}\phi_h,
\nabla_{t,x}\phi_h\rangle-{r}\langle\nabla_{t,x}\phi_h,
\bma_h^*\rangle_h
+\frac{r}{2}\langle\bma_h^*,\bma_h^*\rangle_h,
\end{align}
in which 
\begin{align}
\label{Fh}
   F_h^*(\bma_h^*)
:=&\;\langle A^*(\alpha_{0,h}^*+\frac12|\bma_{1,h}^*|^2),1\rangle_h,\\
\label{Gh}
   G_h^*(\phi_h)
:=&\;-(\phi_h(1,\bmx),\rho_1(\bmx))_h
+(\phi_h(0,\bmx),\rho_0(\bmx))_h.
\end{align}
Note that all terms in the discrete augmented Lagrangian \eqref{Lrh} are integrated using numerical integration
$(\cdot,\cdot)_h$ or $\langle\cdot,\cdot\rangle_h$, except the space-time Laplacian term in the second row of \eqref{Lrh}, which is integrated using exact integration $\langle\cdot, \cdot\rangle$ to avoid a singular matrix
for the Laplacian.

Similarly, the discrete scheme for MFG \eqref{aug-mfg}
reads as follows:
given a space-time mesh $\Ih\otimes \Th$ and a polynomial degree $k\ge 0$, find $\bma_h, \bma_h^*\in [W_h^k]^{3}$,
$\rho_{1,h}, \rho_{1,h}^*\in M_h^k$,
and $\phi_h\in V_h^{k+1}$ such that 
\begin{align}
\label{aug-mfg-h}
\sup_{\bma_h\in [W_h^k]^{3},\rho_{1,h}\in M_h^k} \;
    \inf_{\phi_h\in V_h^{k+1}, \bma_h^*\in[W_h^k]^{3}, \rho_{1,h}^*\in M_h^k
    }L_{r,h}^{MFG}(\phi_h, \bma_h, \rho_{1,h},\bma_h^*,  \rho_{1,h}^*),
\end{align}
where  the discrete augmented Lagrangian is 
\begin{align}
L_{r,h}^{MFG}
=&
 F_h^*(\bma_h^*)+R_h^*(\rho_{1,h}^*)+(\rho_0, \phi_h(0,\cdot))_h
\nonumber\\
&\;+\langle\bma_h, \nabla_{t,x}\phi-\bma_h^*\rangle_h
\;-(\rho_{1,h}, \phi_h(1,\cdot)+\rho_{1,h}^*)_h\nonumber\\
&\;+\frac{r_1}{2}\langle\nabla_{t,x}\phi_h,
\nabla_{t,x}\phi_h\rangle-{r_1}\langle\nabla_{t,x}\phi_h,
\bma_h^*\rangle_h
+\frac{r_1}{2}\langle\bma_h^*,\bma_h^*\rangle_h\nonumber\\
&
+\frac{r_2}{2}
(\phi_h(1,\cdot),\phi_h(1,\cdot))
+{r_2}
(\phi_h(1,\cdot),\rho_{1,h}^*)_h
+\frac{r_2}{2}
(\rho_{1,h}^*,\rho_{1,h}^*)_h,
\end{align}
in which 
\begin{align}
\label{Rh}
   R_h^*(\rho_{1,h}^*)
:=&\;(\Gamma^*(\rho_{1,h}^*),1)_h.
\end{align}

\subsection{The ALG2 algorithm}
The discrete saddle-point problems \eqref{aug-mfp-h}
and \eqref{aug-mfg-h} can be solved efficiently using the ALG2 algorithm \cite{FortinBook}, where minimization of 
$\phi_h$, $\bma_h^*$, and $\rho_{1,h}^*$
are decoupled.
For simplicity, we only illustrate the main steps for the 
discrete MFG problem \eqref{aug-mfg-h}; see also \cite{benamou1999numerical, benamou2014augmented}.
One iteration of ALG2 contains the following three steps.

\subsubsection*{Step A: update $\phi_h$}
Minimize $L_{r,h}^{MFG}$ with respect to the first component by solving the elliptic problem:
Find $\phi_h^{m+1}\in V_h^{k+1}$ such that
it is the solution to 
\begin{align*}
 \inf_{\phi_h\in V_h^{k+1}}L_{r,h}^{MFG}(\phi_h, \bma_h^{m}, \rho_{1,h}^{m},\bma_h^{*,m}, \rho_{1,h}^{*,m}).
\end{align*}
This is simply a linear, constant-coefficient,  space-time diffusion problem: Find $\phi_h^{m+1} \in V_h^{k+1}$
such that
\begin{align}
\label{phi-solve}
&r_1\langle\nabla_{t,x}\phi_h^{m+1},
\nabla_{t,x}\psi_h\rangle
+r_2
(\phi_h^{m+1}(1,\cdot),\psi_h(1,\cdot))
\\
&=\langle r_1\bma_h^{*,m}-\bma_h^{m},
\nabla_{t,x}\psi_h\rangle_h
-
(r_2\rho_{1,h}^{*,m}-\rho_{1,h}^{m},\psi_h(1,\cdot))_h
-(\rho_{0},\psi_h(0,\cdot))_h,\nonumber
\end{align}
for all $\psi_h\in V_h^{k+1}$.

\subsubsection*{Step B: update $\bma_h^{*}$ and $\rho_{1,h}^{*}$}
Minimize $L_{r,h}^{MFG}$ with respect to the last two components by solving the nonlinear problem:
Find $\bma_h^{*,m+1}\in [W_h^{k}]^3$ 
and $\rho_{1,h}^{*,m+1}\in M_h^k$ such that
they are the solutions to 
\begin{align*}
 \inf_{\bma_h^{*}\in [W_h^{k}]^3,
 \rho_{1,h}^{*}\in M_h^{k}}L_{r,h}^{MFG}(\phi_h^{m+1}, \bma_h^{m}, \rho_{1,h}^{m},\bma_h^{*}, \rho_{1,h}^{*}).
\end{align*}

Using the basis functions in \eqref{W-basis} and \eqref{M-basis},
we write 
\begin{align*}
\bma_h=&\;\sum_{\ell=1}^M\sum_{i_s=1}^{N_k}\sum_{j=1}^N\sum_{i_t=1}^{k+1}{\sf {\bm a}}_{\ell, i_s,j,i_t}\psi_{i_t}^j(t)\varphi_{i_s}^\ell(\bmx),
\quad \rho_{1,h}=&\;\sum_{\ell=1}^M\sum_{i_s=1}^{N_k}{\sf r}_{\ell, i_s}\varphi_{i_s}^\ell(\bmx),\\
\bma_h^*=&\;\sum_{\ell=1}^M\sum_{i_s=1}^{N_k}\sum_{j=1}^N\sum_{i_t=1}^{k+1}{\sf {\bm a}}_{\ell, i_s,j,i_t}^*\psi_{i_t}^j(t)\varphi_{i_s}^\ell(\bmx),\quad\rho_{1,h}^*=&\;\sum_{\ell=1}^M\sum_{i_s=1}^{N_k}{\sf r}_{\ell, i_s}\varphi_{i_s}^\ell(\bmx),
\end{align*}
with ${\sf {\bm a}}_{\ell, i_s,j,i_t}, {\sf {\bm a}}_{\ell, i_s,j,i_t}^*$,
${\sf r}_{\ell, i_s}$ and ${\sf r}_{\ell, i_s}^*$.

By the choice of the numerical integration and the nodal bases
for $W_h^k$ and $M_h^k$, we observe that this optimization problem is decoupled for each DOF of $\bma_{h}^{*,m+1}$
and $\rho_{1,h}^*$, hence can be efficiently solved pointwisely:
for each $\ell, i_s, j, i_t$, find 
$\sf {\bm a}_{\ell, i_s, j, i_t}^{*,m+1}\in \bbR^3$ such that it solves
\begin{align}
\label{alpha-newton}
 \inf_{{\sf {\bm a}^{*}}=({\sf a_0^*}, {\sf\bm a_1^*})\in \bbR^3
 }&A^*({\sf a_0^*}+\frac12|{\sf\bm a_1^*}|^2)
+\frac{r_1}{2}|{\sf\bm a}^*|^2\nonumber\\
&-({\sf\bm a}_{\ell, i_s, j, i_t}^{m}+r_1\nabla_{t,x}\phi_h^{m+1}(\eta_{i_t}^j, \bm\xi_{i_s}^\ell))\cdot{\sf\bm a}^*,
\end{align}
and find $\sf r_{\ell, i_s}^{*,m+1}\in \bbR$
such that it solves
\begin{align}
\label{rho-newton}
 \inf_{{\sf r^{*}}\in \bbR^+}\Gamma^*({\sf r}^*)
+\frac{r_2}{2}|{\sf r}^*|^2-
({\sf r}_{\ell, i_s}^{m}-r_2\phi_h^{m+1}(1, \bm\xi_{i_s}^\ell))\cdot{\sf r}^*.
\end{align}
Both optimization problems can be efficiently solved in parallel using the Newton's method.

\subsubsection*{Step C: update $\bma_h$ and $\rho_{1,h}$}
This is a simple pointwise update for the DOFs of 
the Lagrange multipliers $\bma_h$ and $\rho_{1,h}$:
\begin{align}
\label{stepC}
{\sf {\bm a}}_{\ell, i_s, j, i_t}^{m+1}
= &\;
{\sf {\bm a}}_{\ell, i_s, j, i_t}^{m}+r_1(\nabla_{t,x}\phi_h^{m+1}(\eta_{i_t}^j, \bm\xi_{i_s}^\ell))-{\sf {\bm a}}_{\ell, i_s, j, i_t}^{m+1}),\\
{\sf r}_{\ell, i_s}^{m+1}
= &\;
{\sf r}_{\ell, i_s}^{m}-r_2(\phi_h^{m+1}(1, \bm\xi_{i_s}^\ell))+{\sf r}_{\ell, i_s}^{m+1}).
\end{align}

We use the $\ell_{\infty}$-errors in the Lagrange multipliers
\begin{align}
\label{err}
err_m^{a} := &\max_{\ell, i_s, j, i_t}|{\sf {\bm a}}_{\ell, i_s, j, i_t}^{m+1}-{\sf {\bm a}}_{\ell, i_s, j, i_t}^{m}|,\\
err_m^{r} :=& \max_{\ell, i_s}|{\sf r}_{\ell, i_s}^{m+1}-{\sf r}_{\ell, i_s}^{m}|,
\end{align}
to monitor the convergence of the ALG2 algorithm.

\begin{remark}
We specifically note that the use of the integration rule space $W_h^k$ and numerical integration is crucial for the efficient implementation of Step B in the ALG2 algorithm, which leads to a pointwise update per integration point. If this space and numerical integration were not chosen carefully, additional unnecessary degrees of freedom coupling maybe introduced, which slows down the overall algorithm. 
\end{remark}

\section{Numerical experiments}
\label{sec:num}
In this section, we conduct comprehensive experiments to show the efficiency and effectiveness of the proposed numerical algorithms. We restrict ourself to structured (hyper-)rectangular meshes. The case with unstructured meshes will be considered elsewhere.
We first numerical verify the convergence of rate of the algorithm related to the mesh size and polynomial degree. Throughout, we take the augmented Lagrangian parameters to be $r=r_1=r_2=1$.
Our numerical simulations are performed using the open-source finite-element software
{\sf NGSolve} \cite{Schoberl16}, \url{https://ngsolve.org/}.

\subsection{Convergence rates}
\label{ex1}
We first consider OT problems with known exact solutions.
Specifically, we take the domain $\Omega=\bbR^d$ with $d=1$ or $d=2$, cost $A(\rho)=0$ in \eqref{mfp:problem} with initial and terminal densities:
\[
\rho_0(\bmx)=\exp(-50|\bmx-\bmx_0|^2), \quad
\rho_1(\bmx)=\exp(-50|\bmx-\bmx_1|^2),
\]
where $\bmx_0=0.25, \bmx_1=0.75$ when spatial dimension $d=1$, and 
$\bmx_0=(0.25, 0.25), \bmx_1=(0.75,0.75)$ when spatial dimension $d=2$.
The exact solution is simply a traveling wave solution:
\begin{align*}
\rho_{ex}(t,\bmx)=&\;\exp(-50|\bmx-(1+2t)\bmx_0|^2), \\
m_{ex,i}(t,\bmx)=&\;0.5\exp(-50|\bmx-(1+2t)\bmx_0|^2), \quad
\forall 1\le i\le d,
\end{align*}
where $\bmm_{ex}=(m_{ex,1}, \cdots., m_{ex,d})$.
We truncate the domain $\Omega$ to be a unit box $[0,1]^d$, and 
replace the homogeneous boundary condition \eqref{mfp-bcm}
with a boundary source term 
\[
\bmm\cdot\bmn = \bmm_{ex}\cdot\bmn,\quad \text{ on }
[0,1]\times \partial\Omega.
\]
With this modification, the $G$-term in \eqref{pd-mfp:problem1}
contains an additional boundary source term:
\[
G(\phi):=
\int_{\Omega}-\phi(1,\bmx)\rho_1(\bmx)+
\phi(0,\bmx)\rho_0(\bmx)\,\deri\bmx+
\int_0^1\int_{\partial\Omega}\phi(t,\bmx) \bmm_{ex}\cdot\bmn\,\deri s\deri t.
\]

\begin{figure}[ht]
  \centering
  \subfigure[$k=0$: 4$\times$4 grid]{
\begin{tikzpicture}
\draw[step=0.8cm,black,thick] (-1.61, -1.61) grid (1.6,1.6);
\foreach \x in {-2,-1,...,2}{
\foreach \y in {-2,-1,...,2}{
        \node[draw,circle,inner sep=2pt,fill, red] at (0.8*\x,0.8*\y) {};
      }
}
\foreach \x in {-1.5,-0.5,...,1.5}{
\foreach \y in {-1.5,-0.5,...,1.5}{
        \node[draw,inner sep=2pt, fill, black] at (0.8*\x,0.8*\y) {};
      }
}
\end{tikzpicture}  
}
   \subfigure[$k=1$: 2$\times$2 grid]{
\begin{tikzpicture}
\draw[step=0.8cm,black,dotted, thin] (-1.61, -1.61) grid (1.6,1.6);
\draw[step=1.6cm,black,thick] (-1.61, -1.61) grid (1.6,1.6);
\foreach \x in {-2,-1,...,2}{
\foreach \y in {-2,-1,...,2}{
        \node[draw,circle,inner sep=2pt,fill, red] at (0.8*\x,0.8*\y) {};
      }
}
\foreach \x in {-1.57735026919,-0.42264973081
,0.42264973081,1.57735026919}{
\foreach \y in {-1.57735026919,-0.42264973081
,0.42264973081,1.57735026919}{
        \node[draw,inner sep=2pt, fill, black] at (0.8*\x,0.8*\y) {};
      }
}
\end{tikzpicture}  
}
   \subfigure[$k=3$: 1$\times$1 grid]{
\begin{tikzpicture}
\draw[step=0.8cm,black,dotted, thin] (-1.61, -1.61) grid (1.6,1.6);
\draw[black,thick] (-1.61, -1.61) ->
(-1.61,1.6)->(1.6,1.6)->(1.6,-1.61)->(-1.61,-1.61);
\foreach \x in {-2,-1,...,2}{
\foreach \y in {-2,-1,...,2}{
        \node[draw,circle,inner sep=2pt,fill, red] at (0.8*\x,0.8*\y) {};
      }
}
\foreach \x in {-1.722272,-0.679962,0.679962,1.722272}{
\foreach \y in {-1.722272,-0.679962,0.679962,1.722272}{
        \node[draw,inner sep=2pt, fill, black] at (0.8*\x,0.8*\y) {};
      }
}
\end{tikzpicture}  
}
  \caption{Coarse mesh DOFs.
  Circles: DOFs for $V_h^{k+1}$; 
  Squares: DOFs for $W_h^{k}$.
  The coarse mesh is $4\times 4$ for $k=0$, 
$2\times 2$ for $k=1$, 
and $1\times 1$ for $k=3$.
  }
  \label{figure:dofs}
\end{figure}
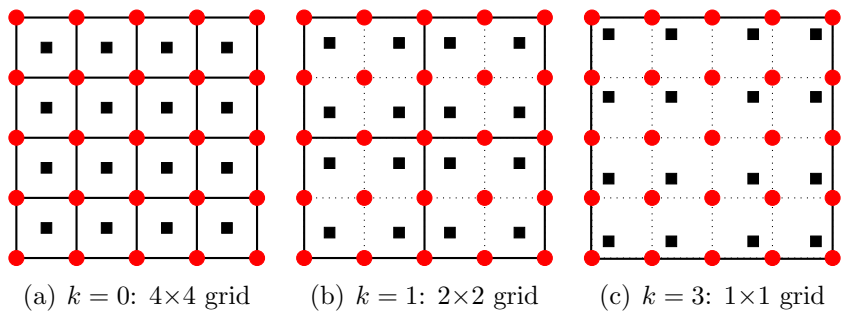

We apply the scheme \eqref{aug-mfp-h}
with polynomial degree $k=0, 1, 3$ on a sequence of uniform hypercubic meshes with $2^{s+2}/(k+1)$ cells in each direction
for $s=0,1,2,3$. The total number of DOFs on the $s$-level meshes is the same for each polynomial degree, which is 
$2^{(s+2)(d+1)}$ for $W_h^k$, and 
$(2^{(s+2)}+1)^{d+1}$ for $V_h^{k+1}$.
So their computational costs are similar.
We apply the ALG2 algorithm to \eqref{aug-mfp-h}
with a stopping tolerance $err_m^a<10^{-10}$ where $err_m^a$ is given in \eqref{err}. We take the parameter $r=1$.
The DOFs on the coarsest meshes for $d=1$ are shown in 
Figure~\ref{figure:dofs}.

We record the $L^2(\Omega_T)$-convergence rates of $\rho_h$
and $\bmm_h$, 
along with the 
convergence rate of the distance
\[
W_2^2 = \int_0^1\int_\Omega \frac{|\bmm|^2}{2\rho}\,\deri\bmx\deri t
\]
in Table~\ref{tab:conv1D} for $d=1$, and 
Table~\ref{tab:conv2D} for $d=2$.
We find that the convergence behavior for $d=1$ and $d=2$ are similar, in particular, (nearly) optimal $L^2$-convergence rates of $k+1$ are observed on the finest mesh for each case, and the average convergence rates for the distance $W_2^2$ is between $2k+2$ and $2k+4$.
Moreover, the advantage of higher order scheme is clearly observed on the fine meshes where the $k=3$ case 
on the $8^{d+1}$ mesh produces $L^2$-errors that are  50 times smaller, 
and $W_2^2$ error that is three orders of magnitude smaller
than the $k=0$ case on the $32^{d+1}$ mesh, although the same number of DOFs are used.
\begin{table}[tb]
\centering
\caption{Convergence rates of scheme~\ref{aug-mfp-h} applied to 1D OT problem.}
\begin{tabular}{ll|ll|ll|ll}
\toprule
 $k$&  mesh   & $L^2$-err in $\rho$  & order  &$L^2$-err in $\bmm$  & order  & $W_2^2$ error     & order  \\ \midrule
0  & $4^2$   
 & 2.068E-01& & 1.097E-01 & &2.834E-03&
\\
0 & $8^2$   & 
1.159E-01&0.84& 5.985E-02&0.87 & 5.472E-04&2.37
\\
0  & $16^2$   &  
  6.007E-02 & 0.95 & 2.970E-02 & 1.01 & 5.788E-05 & 3.24
  \\
0  & $32^2$   &  
  3.002E-02&1.00& 1.497E-02&0.99 & 4.196E-06
  &3.79
\\\midrule
1  & $2^2$   &  
  1.868E-01&& 1.110E-01 &&1.127E-02&
  \\
1 & $4^2$   & 
7.496E-02 & 1.32 &
3.863E-02 & 1.52 &
4.625E-04 & 4.61
\\
1 & $8^2$   & 
2.169E-02 & 1.78 &
1.077E-02 & 1.84 &
9.523E-06 & 5.60 
\\
1  & $16^2$   &  
5.683E-03 & 1.93 & 2.844E-03 & 1.92 & 1.611E-07 & 5.89
\\\midrule
3  & $1^2$   &  
2.148E-01 &&1.301E-01 &&3.337E-02&
\\
3  & $2^2$   &  
6.602E-02 & 1.70 &
3.548E-02 & 1.87 &
5.390E-04 & 5.95
  \\
3 & $4^2$   & 
7.234E-03 &3.19&3.595E-03 &3.30&5.044E-07&10.1
\\
3 & $8^2$   & 5.079E-04 & 3.83 & 2.542E-04 & 3.82 & 4.521E-09 & 6.80
\\
\bottomrule
\end{tabular}
\label{tab:conv1D}
\end{table}
\begin{table}[tb]
\centering
\caption{Convergence rates of scheme~\ref{aug-mfp-h} applied to 2D OT problem.}
\begin{tabular}{ll|ll|ll|ll}
\toprule
 $k$&  mesh   & $L^2$-err in $\rho$  & order  &$L^2$-err in $\bmm$  & order  & $W_2^2$ error     & order  \\ \midrule
0  & $4^3$   
&1.172E-01&& 8.385E-02&& 1.602E-03\\
0 & $8^3$   & 
6.832E-02 &0.78 & 4.879E-02 & 0.78 &
1.646E-04 & 3.28
\\
0  & $16^3$   &  
3.559E-02 & 0.94 & 2.505E-02 & 0.96 & 2.693E-05 & 2.61
  \\
0  & $32^3$   &  
1.787E-02 &0.99 &
1.262E-02 & 0.99& 2.391E-06&3.49
\\\midrule
1  & $2^3$   &  
1.113E-01 && 8.196E-02 && 7.008E-03
  \\
1 & $4^3$   & 4.540E-02 &1.29 &3.260E-02 &1.33&3.882E-05&7.50
\\
1 & $8^2$   & 1.326E-02 &1.78 & 9.354E-03 & 1.80 & 3.563E-06&3.45
\\
1  & $16^2$   &  3.474E-03 & 1.93 & 2.457E-03 & 1.93 & 3.854E-08 & 6.53
\\\midrule
3  & $1^2$   &  
1.432E-01 &&1.109E-01 &&2.278E-02
\\
3  & $2^2$   &  
3.873E-02 & 1.89& 2.804E-02 &1.98&2.795E-04&6.35
  \\
3 & $4^2$   & 
4.353E-03 & 3.15 & 3.072E-03 & 3.19 & 
2.004e-07 & 10.4
\\
3 & $8^2$   & 
3.068E-04 &3.83 & 2.170E-04 &3.82 & 3.977E-09
&5.65
\\
\bottomrule
\end{tabular}
\label{tab:conv2D}
\end{table}

\subsection{MFP with obstacles}
\label{ex2}
We consider a similar MFP problem used in \cite{Buttazzo09}, in which the spatial domain is 
a square excluding some obstacles that mass can not cross:
\[
\Omega=[-1,1]^2\backslash\{\Omega_1\cup\Omega_2\cup\Omega_3\cup\Omega_4\},
\]
where the obstacles $\Omega_1=[-0.2,0.2]\times[-1.0,-0.7]$, 
$\Omega_2=[-0.2,0.2]\times[-0.5,-0.1]$, 
$\Omega_3=[-0.2,0.2]\times[0.1,0.5]$, 
$\Omega_4=[-0.2,0.2]\times[0.7,1.0]$.
We take initial and terminal densities as
two Gaussians
\[
\rho_0(\bmx) = \frac{1}{2\pi\sigma^2}\exp(-\frac{1}{2\sigma^2}|\bmx-\bmx_0|^2), \quad 
\rho_1(\bmx) = \frac{1}{2\pi\sigma^2}\exp(-\frac{1}{2\sigma^2}|\bmx-\bmx_1|^2), 
\]
where the standard deviation $\sigma=0.1$, and 
$\bmx_0=(-0.65,0)$, $\bmx_1=(0.65,0)$.
We take the following 5 choices of interaction cost functions
in the MFP problem \eqref{mfp:problem}, whose convex conjugate are also recorded for completeness:
\begin{align*}
 \begin{cases}
 \text{Case 1}: A(\rho)=0,\quad 
A^*(\rho^*)=\begin{cases}
0 & \text{ if } \rho^*\le 0,\\
+\infty &\text{ if } \rho^*> 0.
\end{cases}, \\[.2ex]
 \text{Case 2}: A(\rho)=c\rho^2,
 \quad 
A^*(\rho^*)=\begin{cases}
0 & \text{ if } \rho^*\le 0,\\
(\rho^*)^2/(4c) &\text{ if } \rho^*> 0.
\end{cases}, \\[.2ex]
 \text{Case 3}: A(\rho)=c\rho\log(\rho),
 \quad
 A^*(\rho^*)=\exp(\rho^*/c-1),
 \\[.2ex]
 \text{Case 4}: A(\rho)=c/\rho,
\quad A^*(\rho^*)=\begin{cases}
-2\sqrt{-c\rho^*}& \text{ if } \rho^*\le 0,\\
+\infty &\text{ if } \rho^*> 0.
\end{cases}
 \\[.2ex]
 \text{Case 5}: A(\rho)=
 \begin{cases}
0& \text{ if } 0\le \rho\le \rho_{\max},\\
+\infty &\text{ else}.
\end{cases}
\quad
A^*(\rho^*)=\rho_{\max}(\rho^*)_+
\end{cases}
\end{align*}
where we take the scaling constant $c=0.1$ in Cases 2--4, 
and maximum density $\rho_{\max}=\frac{1}{2\pi\sigma^2}$
in Case 5.

\begin{figure}[tb]
\centering
\includegraphics[width=0.3\textwidth]{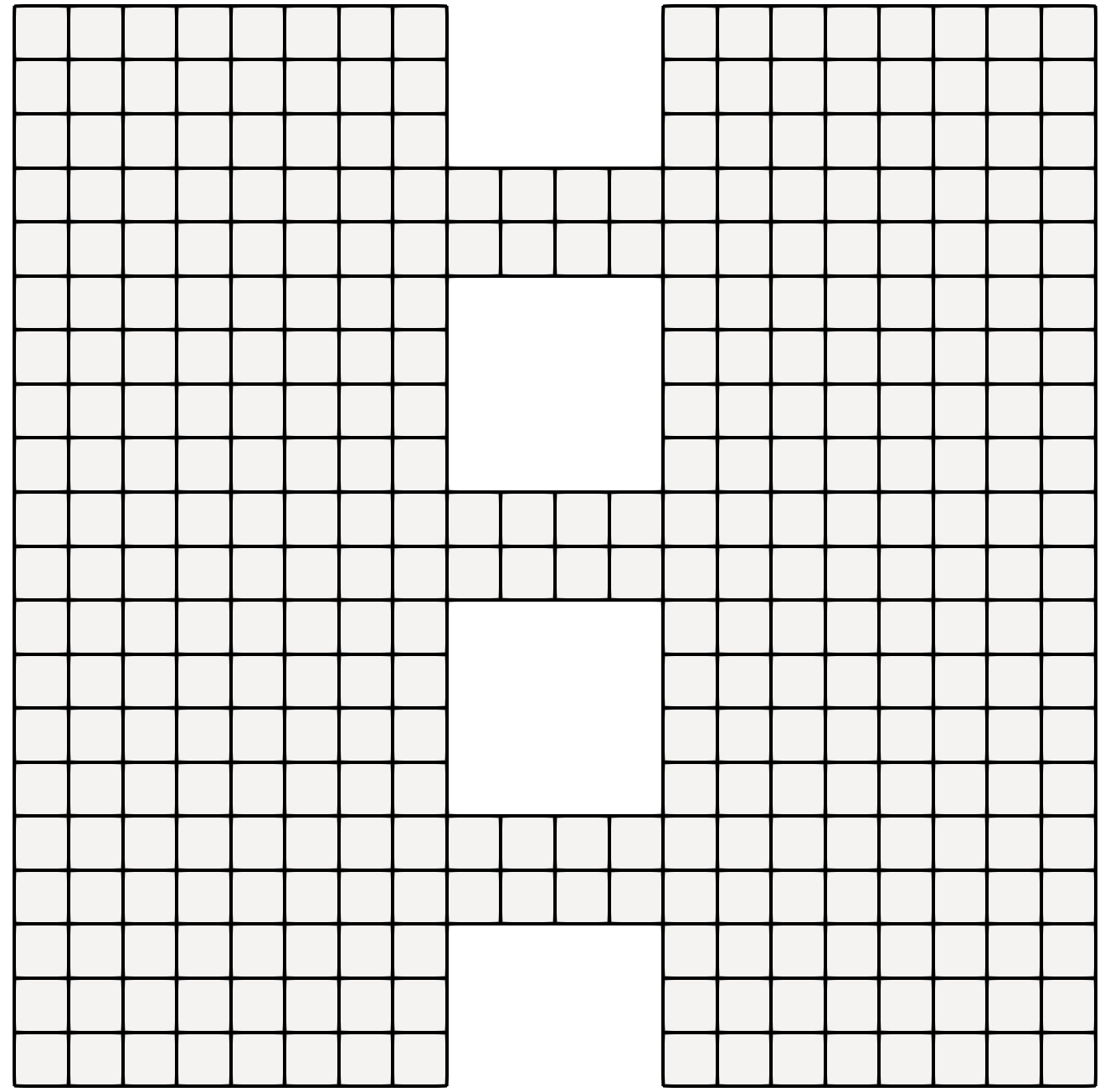}
\caption{A uniform rectangular mesh with $\Delta t=0.1$ for the spatial domain $\Omega$. }
\label{fig:mesh-case2}
\end{figure}

We apply the scheme \eqref{aug-mfp-h} with polynomial degree $k=3$ on a structured hexahedral mesh obtained from
tensor product of a uniform spatial rectangular mesh with mesh size $\Delta x = 0.1$ and uniform temporal mesh with  
$\Delta t = 0.1$. 
The spatial mesh for $\Omega$ is shown in Figure~\ref{fig:mesh-case2}.
We terminate the ALG2 algorithm when the error $err_m^a$
is less than $0.01$.
The number of iterations needed for convergence for the 5 cases are recorded in Table~\ref{tab:iter-case2}, where we find Case 2 has the smallest number of iterations.
\begin{table}[ht!]
\centering
\begin{tabular}{cccccccc}
 & Case 1 & Case 2 & Case 3 & Case 4 & Case 5&\\\hline
iterations 
&780&72&245&503&552\\
\end{tabular}
\caption{Example \ref{ex2}. Number of ALG2 iterations for each case.
}
\label{tab:iter-case2}
\end{table}

\begin{figure}[tb]
\centering
\subfigure[Case 1]{
\label{fig:1x}
\includegraphics[width=0.192\textwidth]{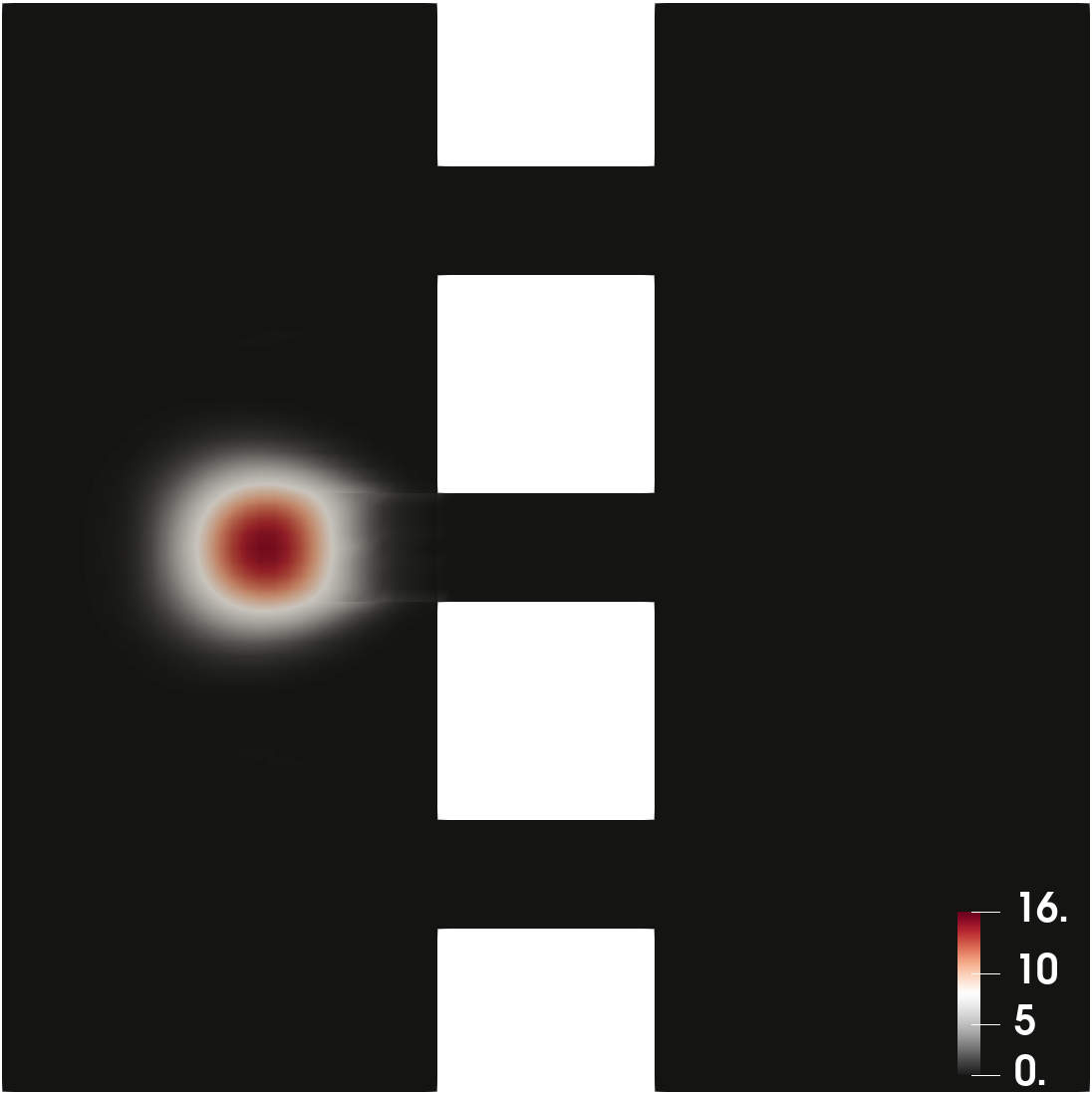}
\includegraphics[width=0.192\textwidth]{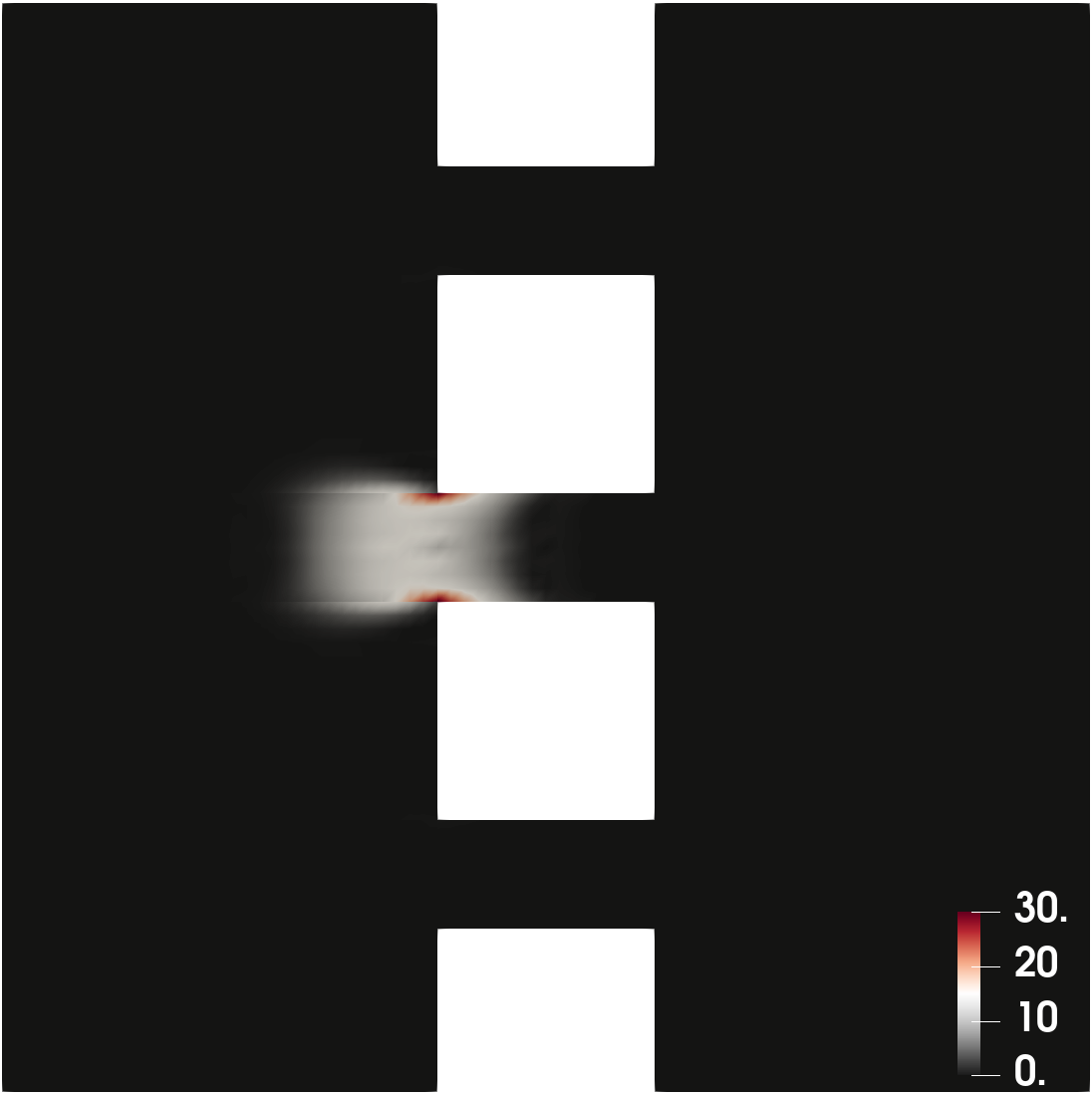}
\includegraphics[width=0.192\textwidth]{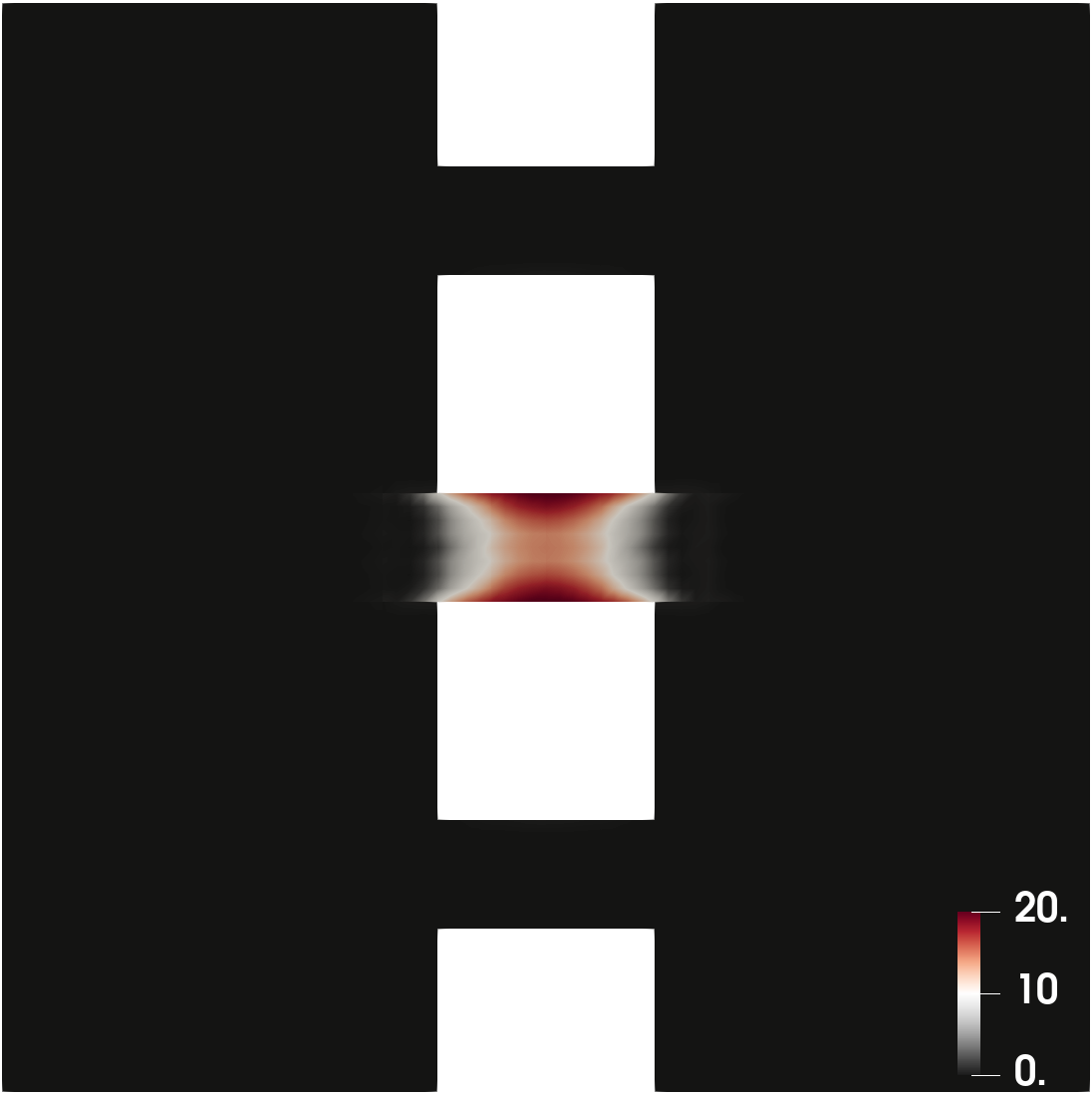}
\includegraphics[width=0.192\textwidth]{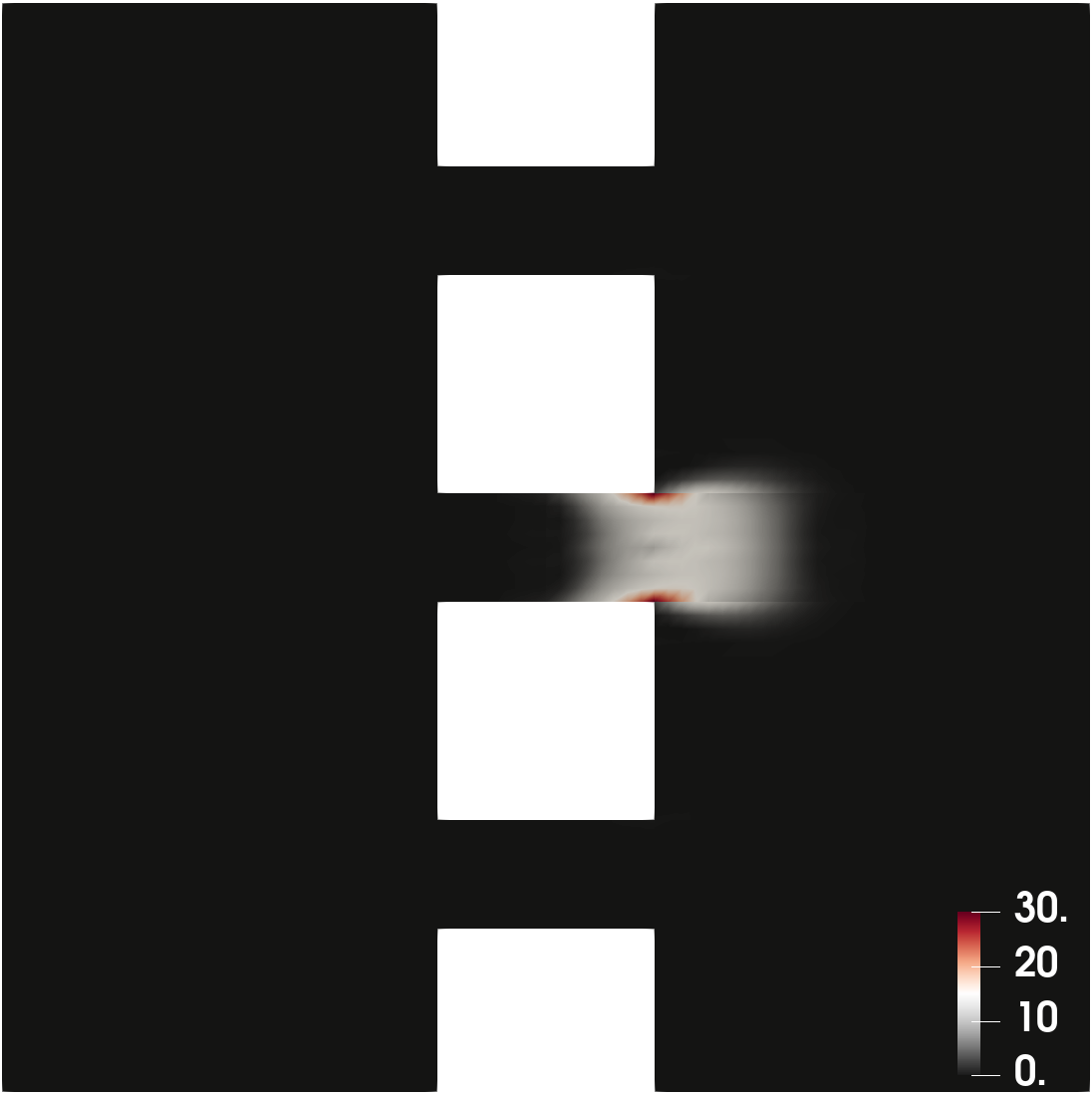}
\includegraphics[width=0.192\textwidth]{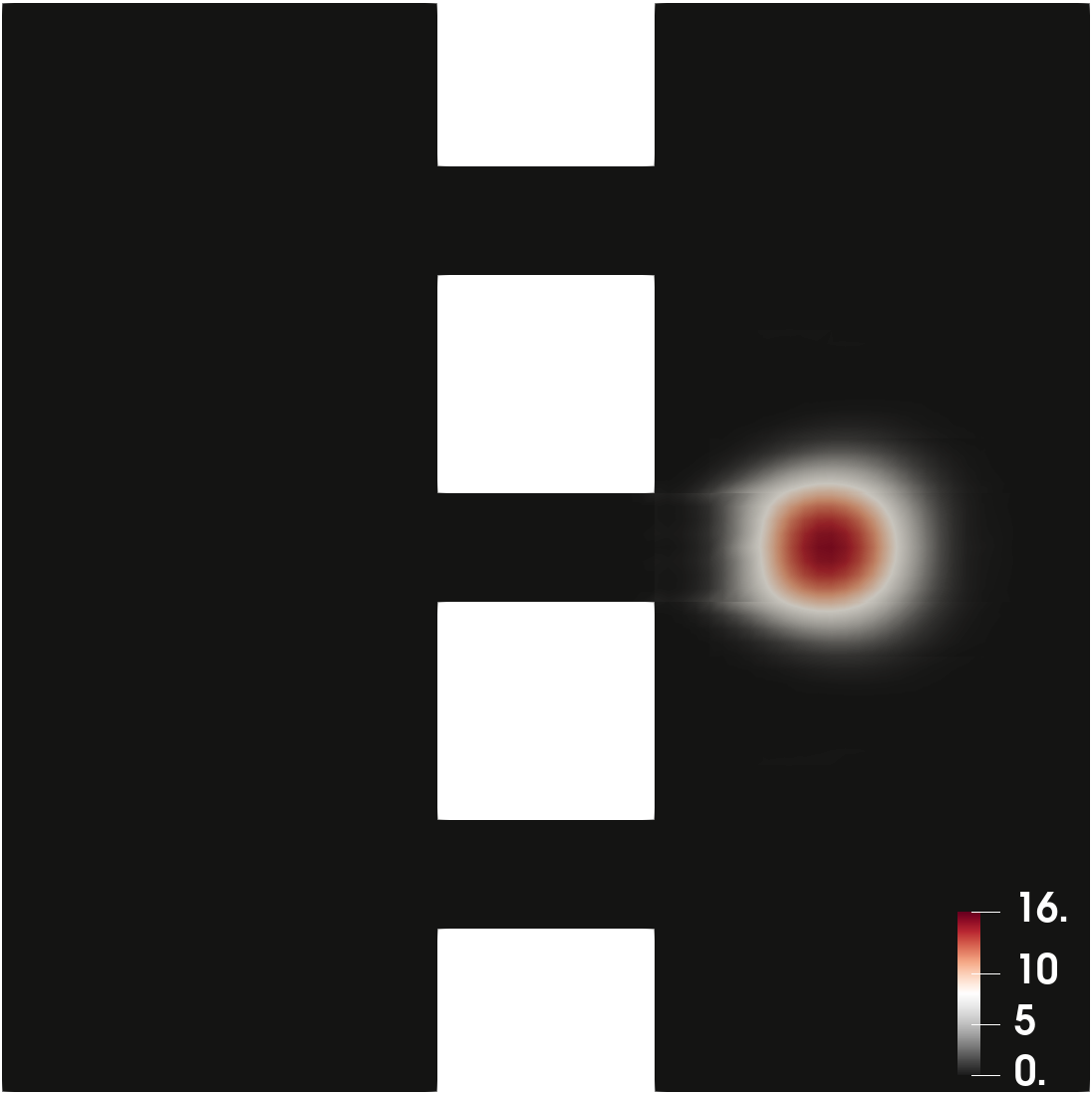}
}
\subfigure[Case 2]{
\label{fig:2x}
\includegraphics[width=0.192\textwidth]{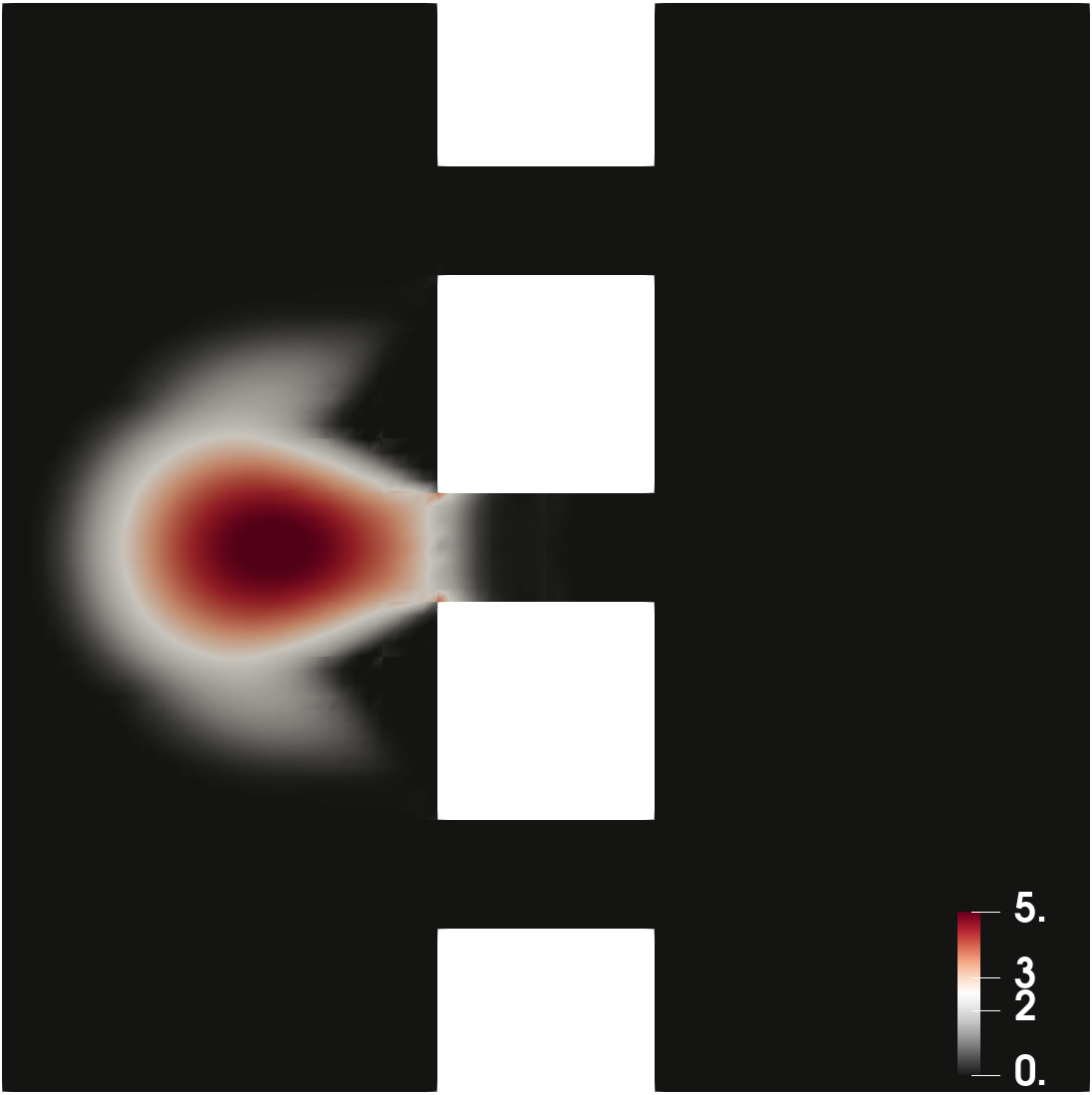}
\includegraphics[width=0.192\textwidth]{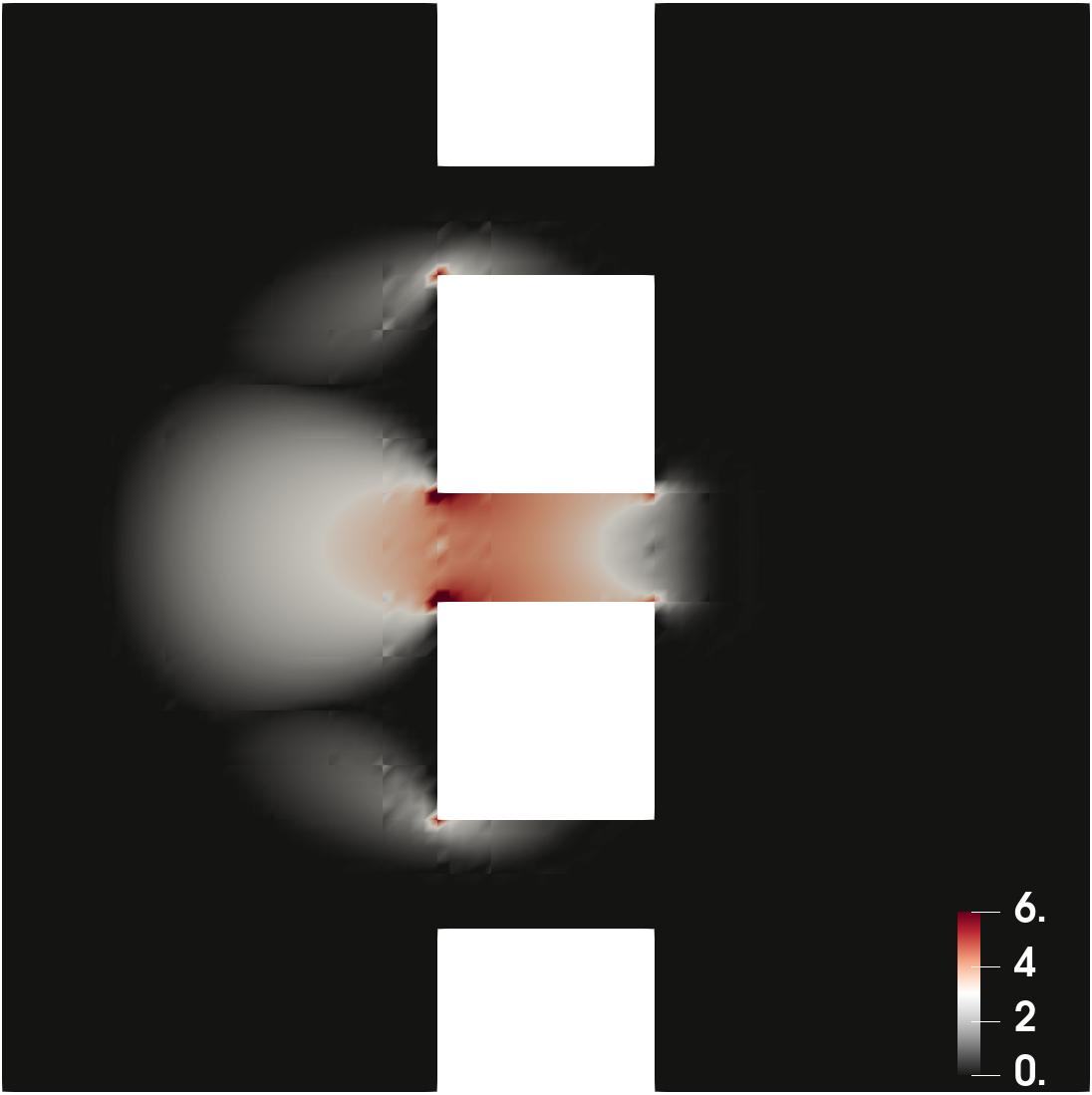}
\includegraphics[width=0.192\textwidth]{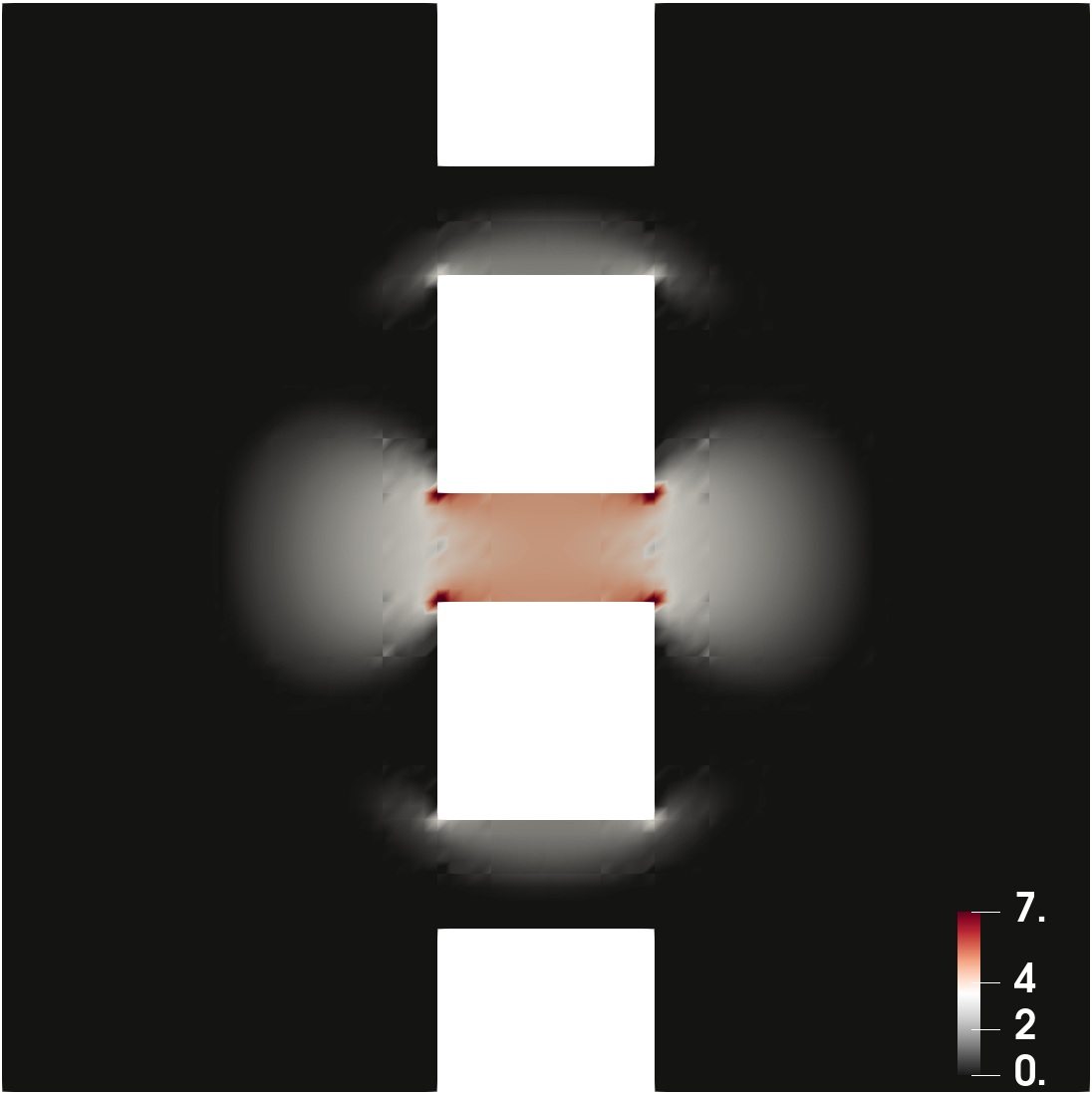}
\includegraphics[width=0.192\textwidth]{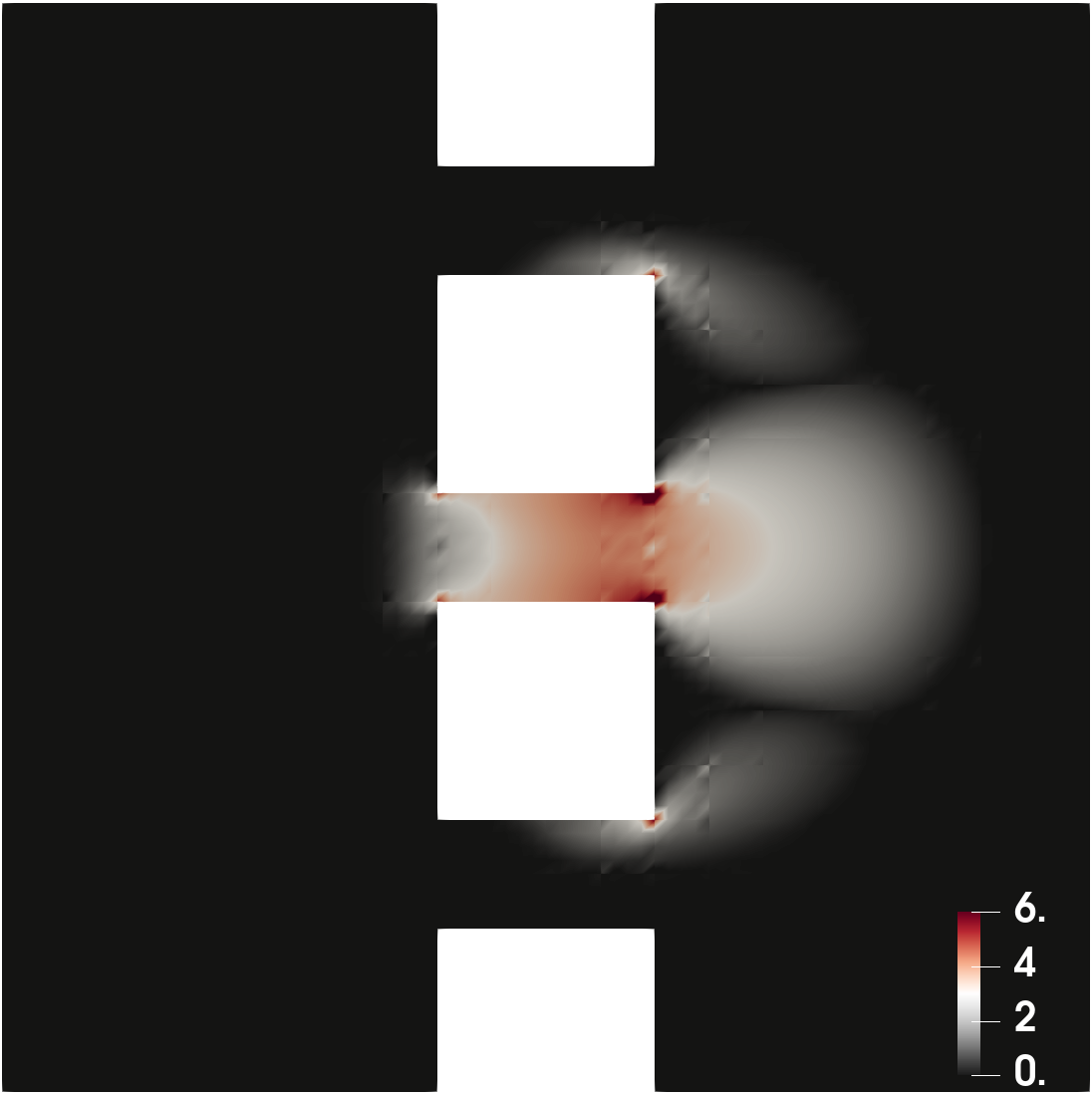}
\includegraphics[width=0.192\textwidth]{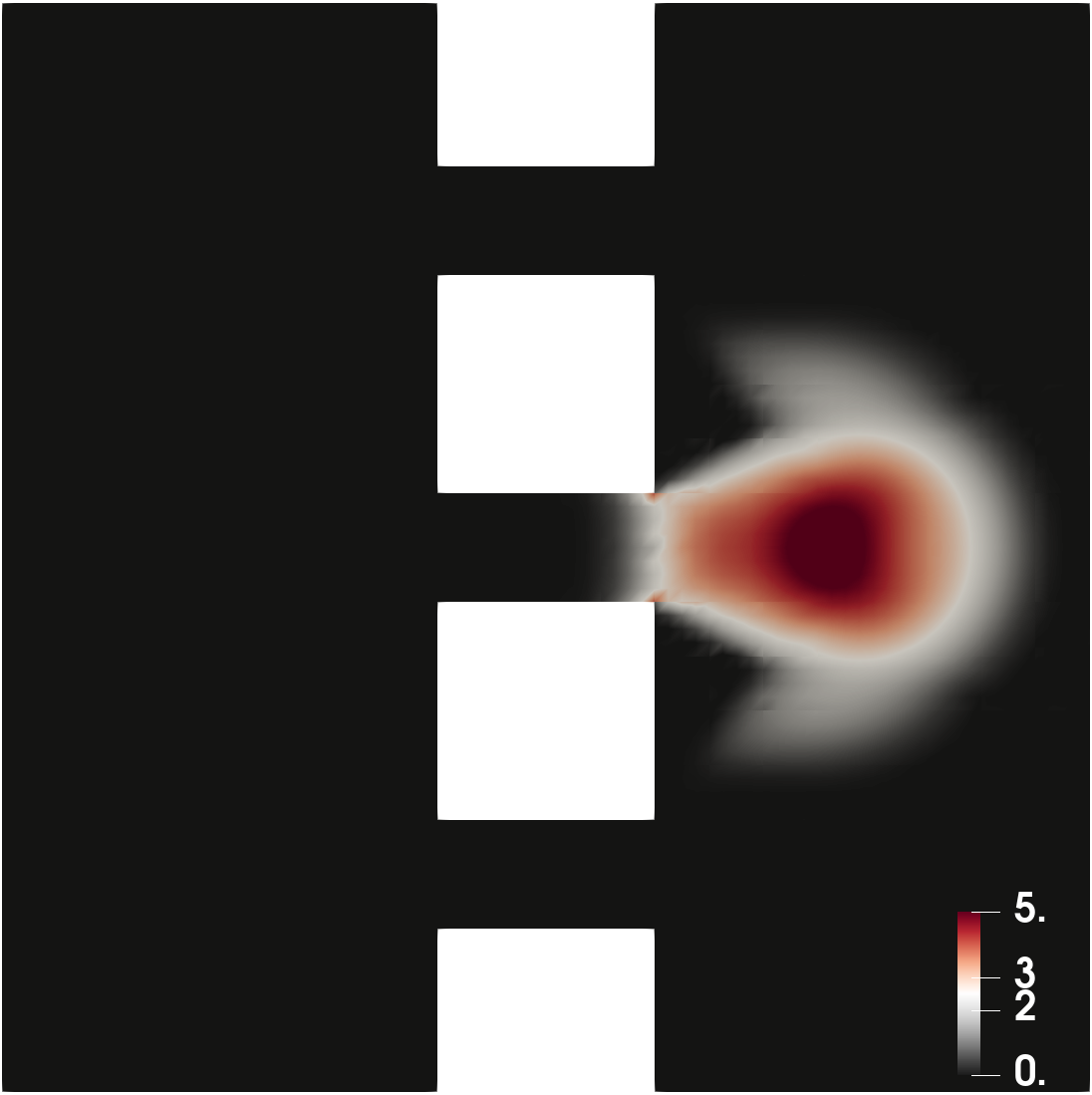}
}
\subfigure[Case 3]{
\label{fig:3x}
\includegraphics[width=0.192\textwidth]{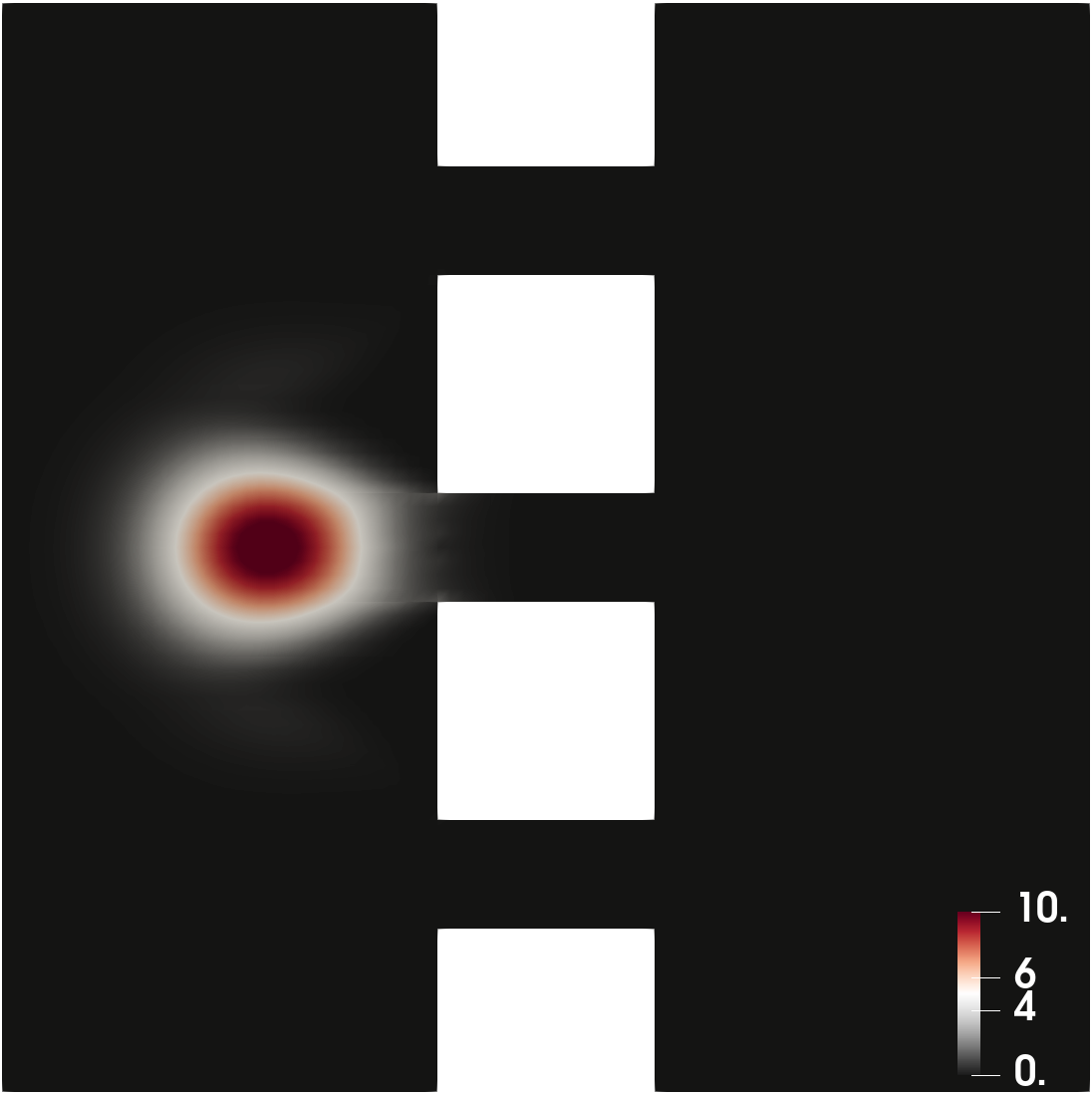}
\includegraphics[width=0.192\textwidth]{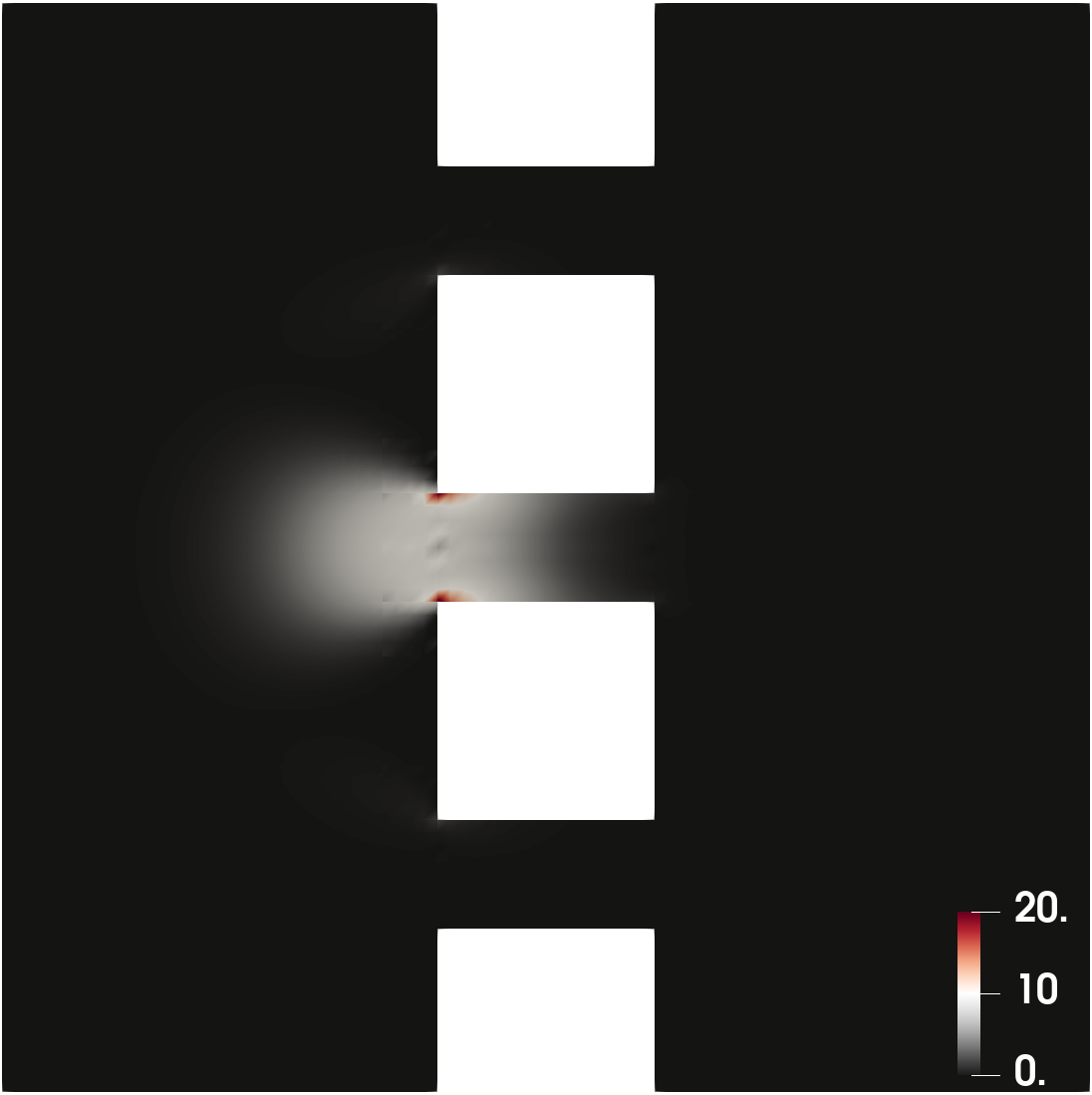}
\includegraphics[width=0.192\textwidth]{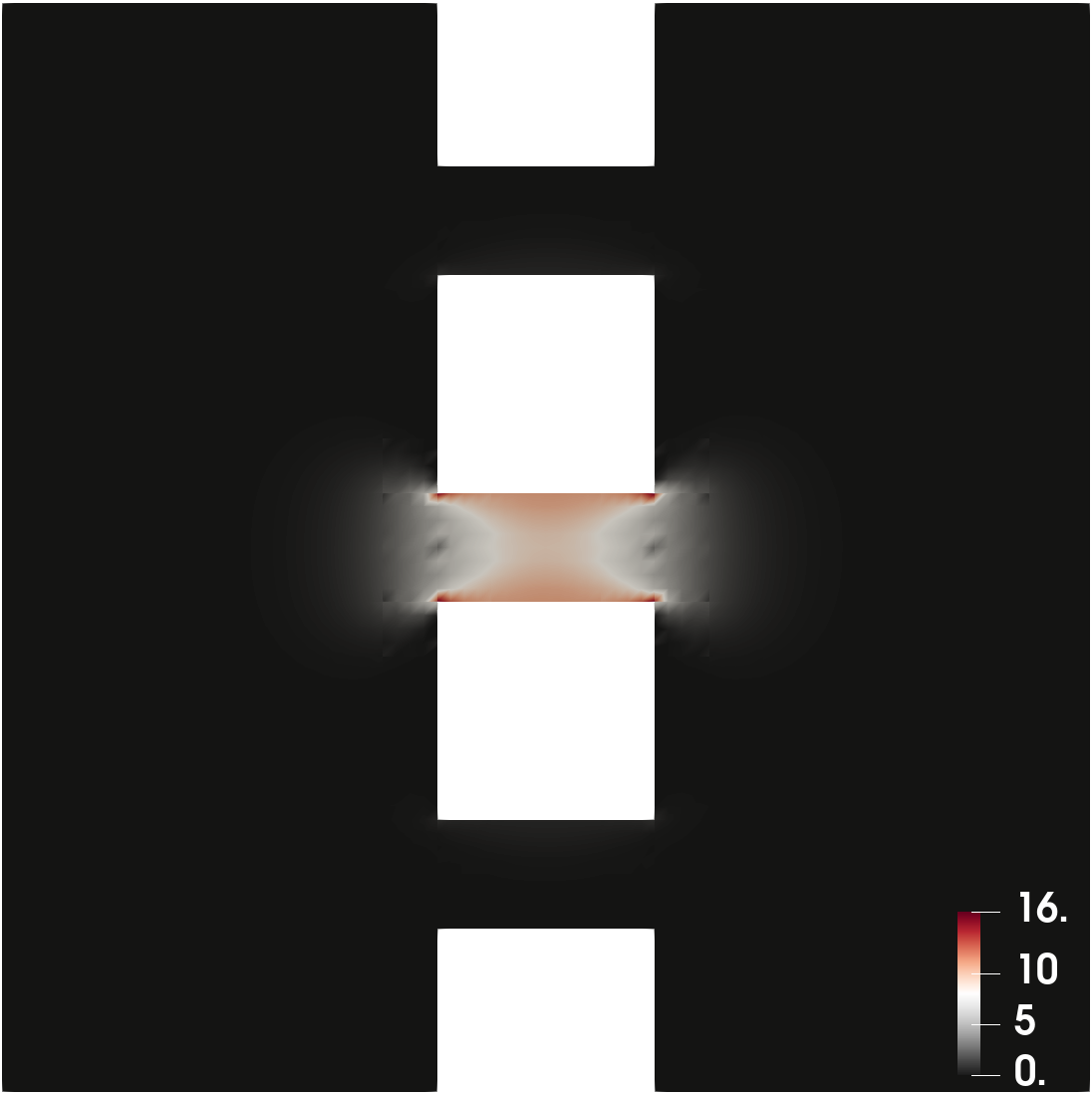}
\includegraphics[width=0.192\textwidth]{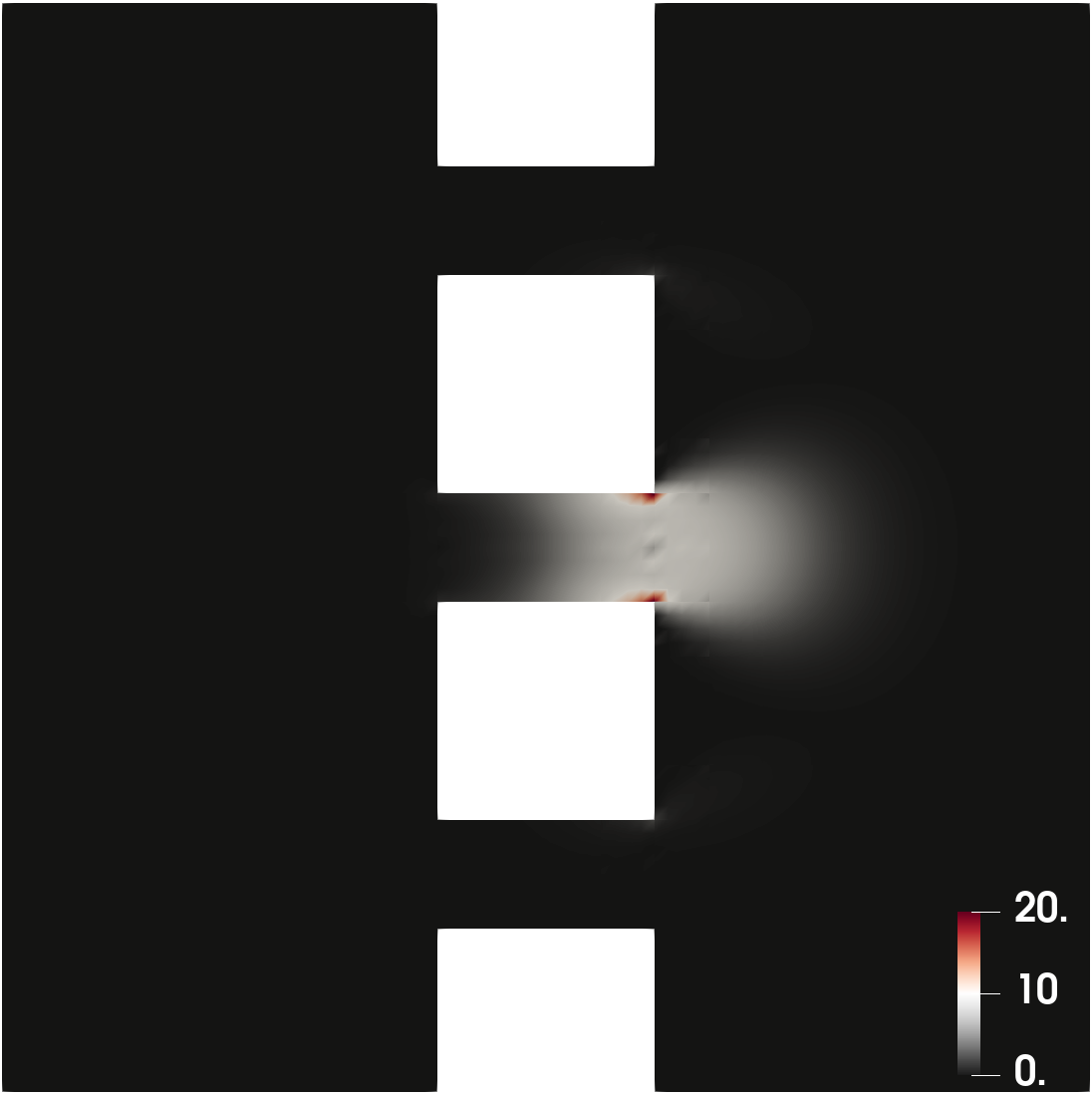}
\includegraphics[width=0.192\textwidth]{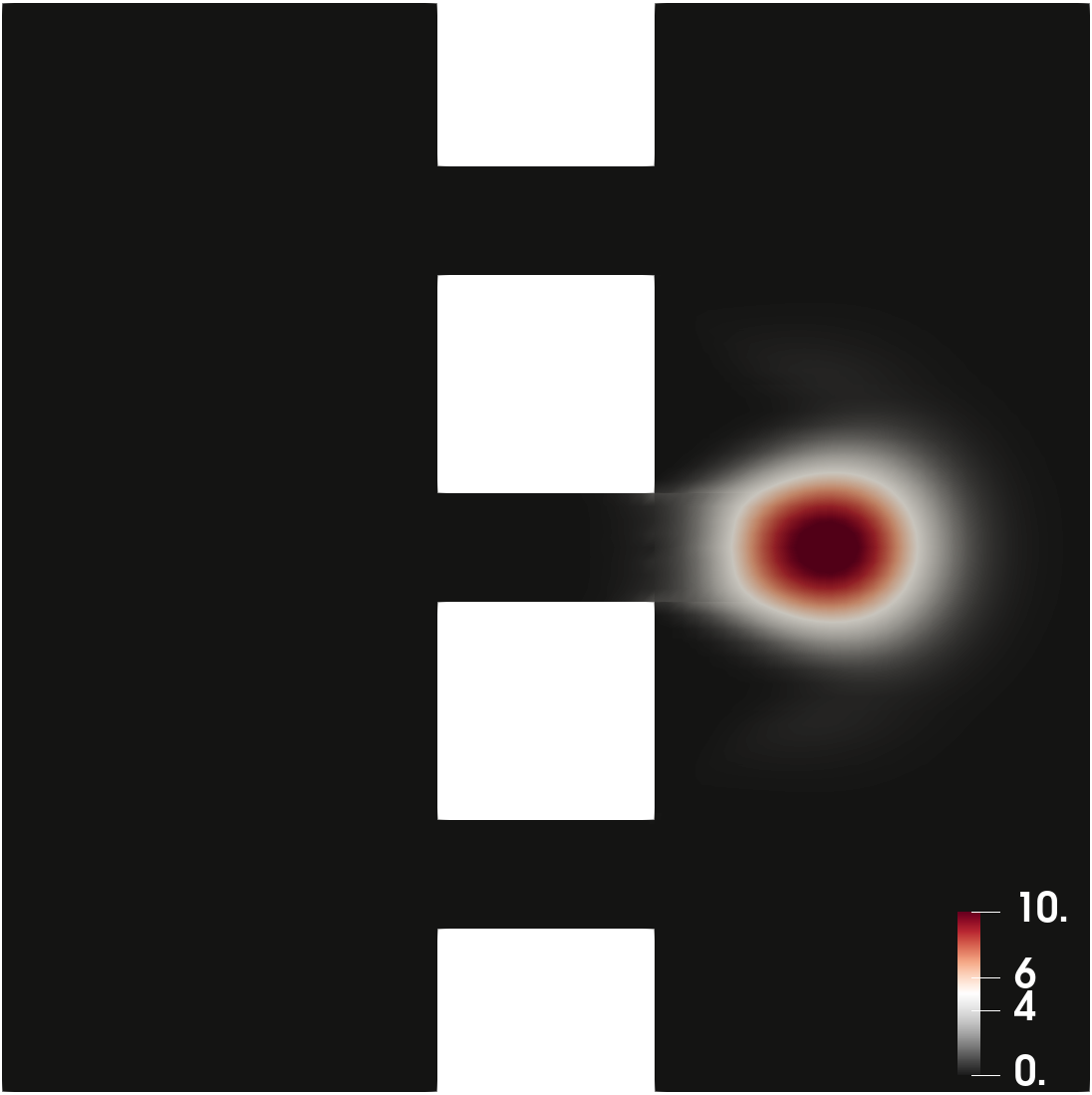}
}
\subfigure[Case 4]{
\label{fig:4x}
\includegraphics[width=0.192\textwidth]{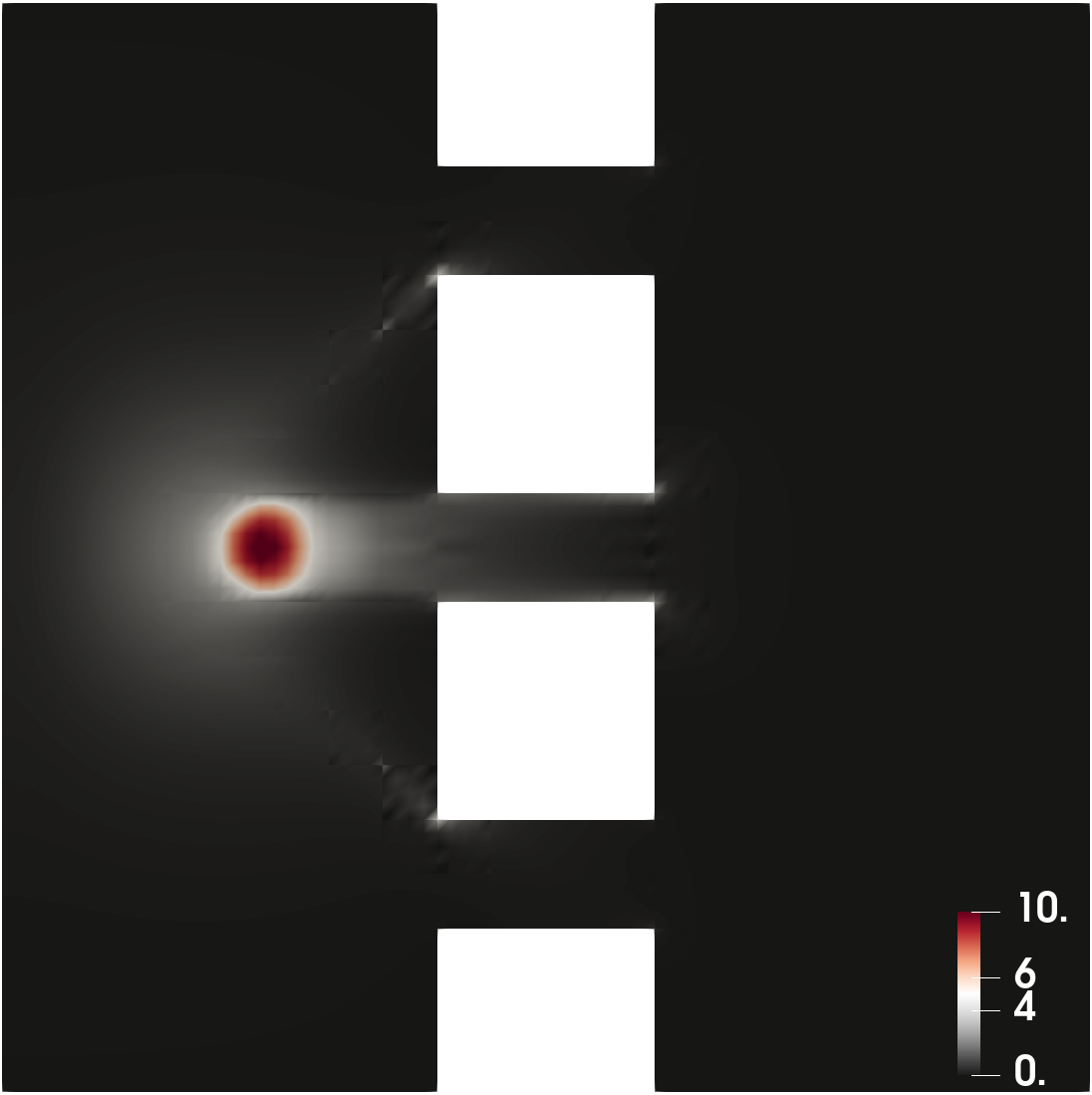}
\includegraphics[width=0.192\textwidth]{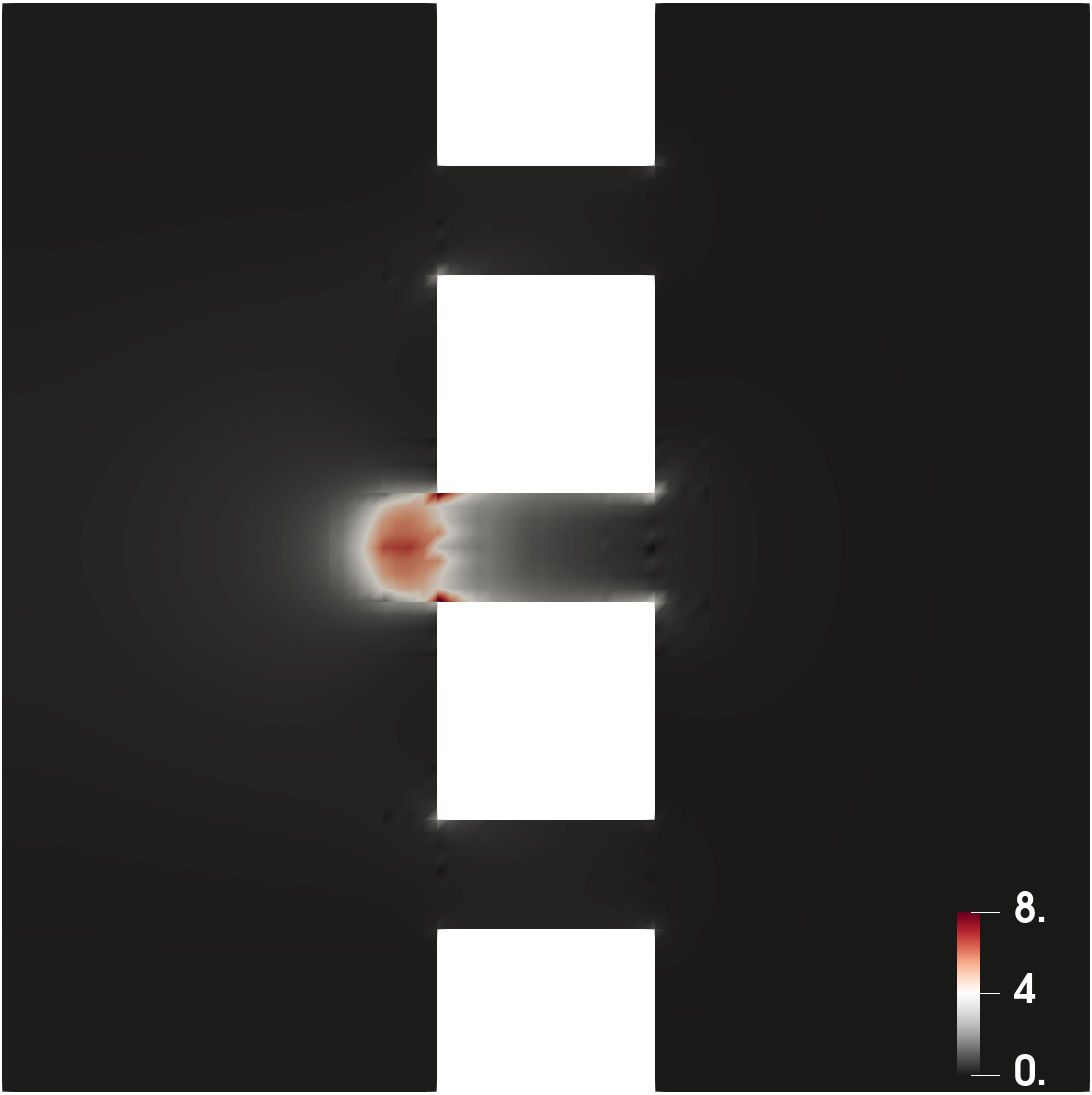}
\includegraphics[width=0.192\textwidth]{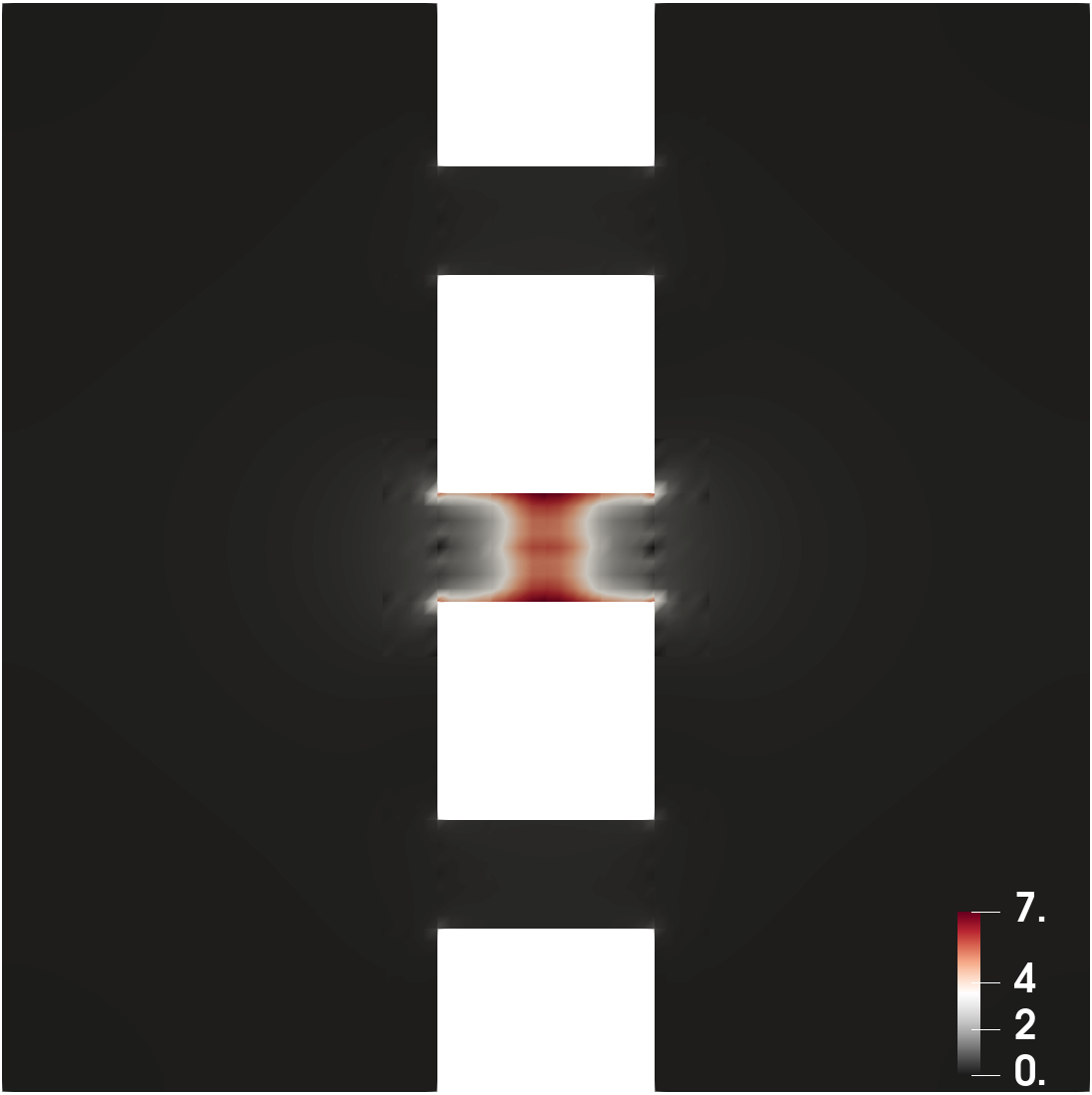}
\includegraphics[width=0.192\textwidth]{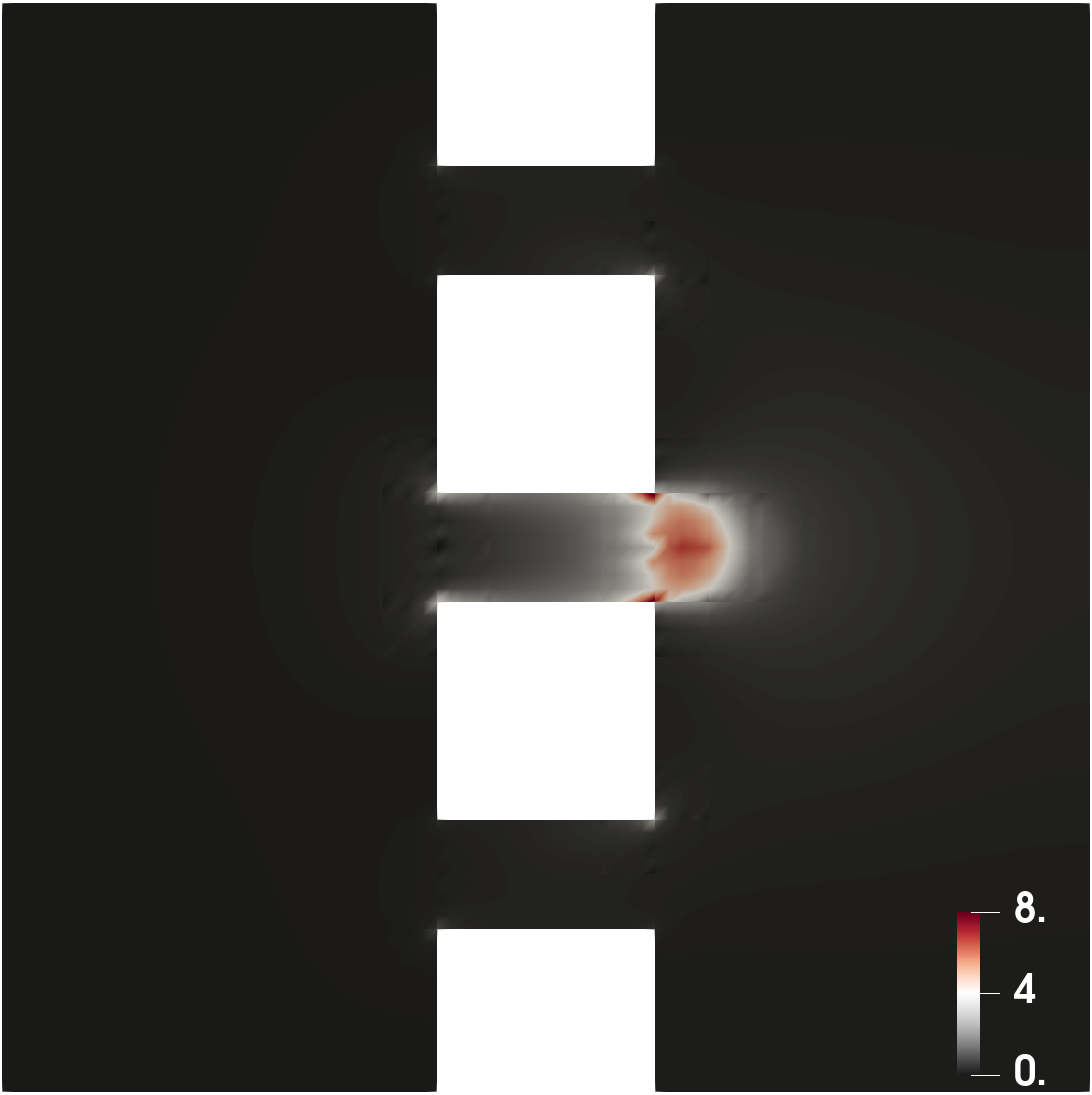}
\includegraphics[width=0.192\textwidth]{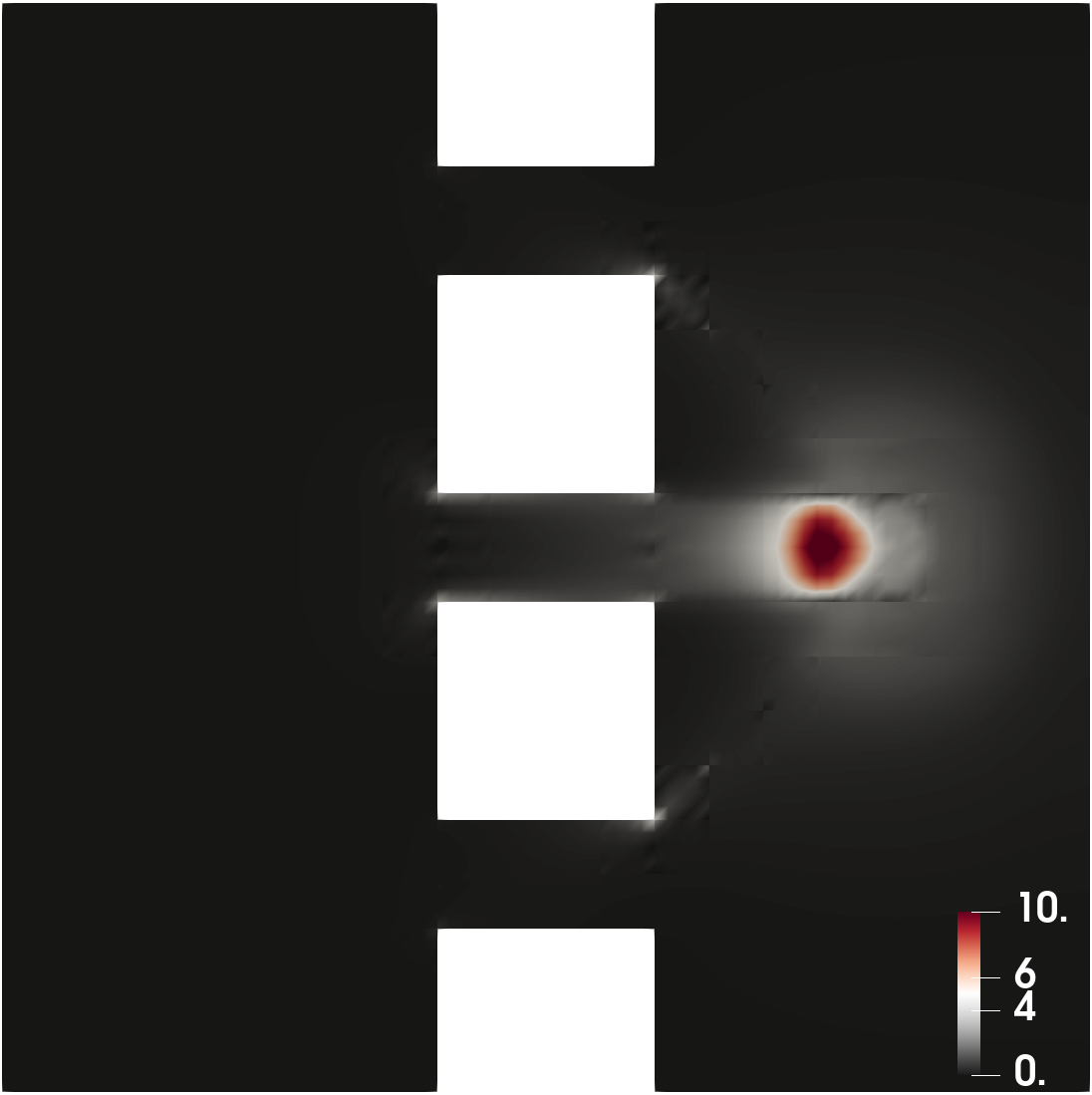}
}
\subfigure[Case 5]{
\label{fig:5x}
\includegraphics[width=0.192\textwidth]{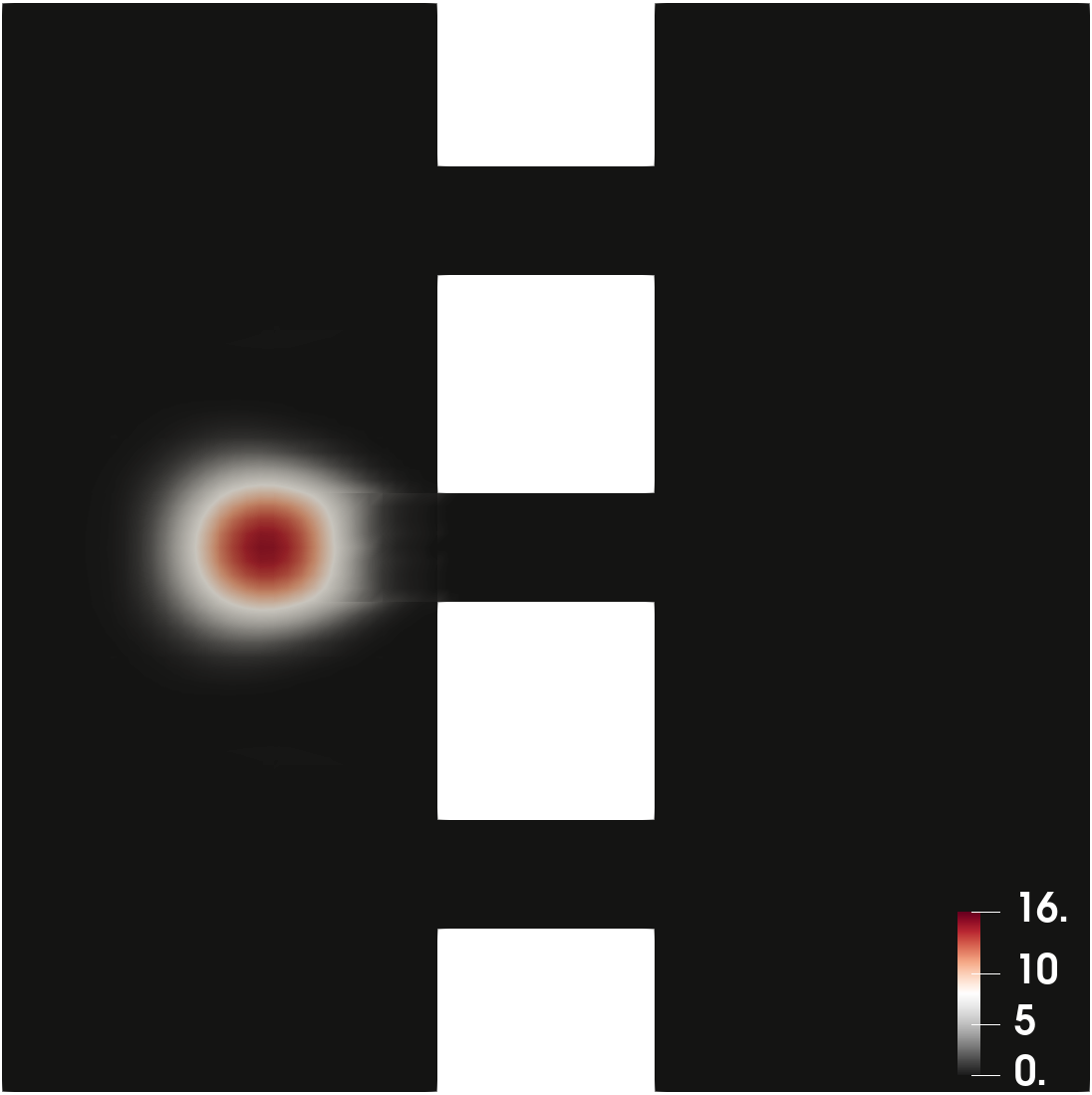}
\includegraphics[width=0.192\textwidth]{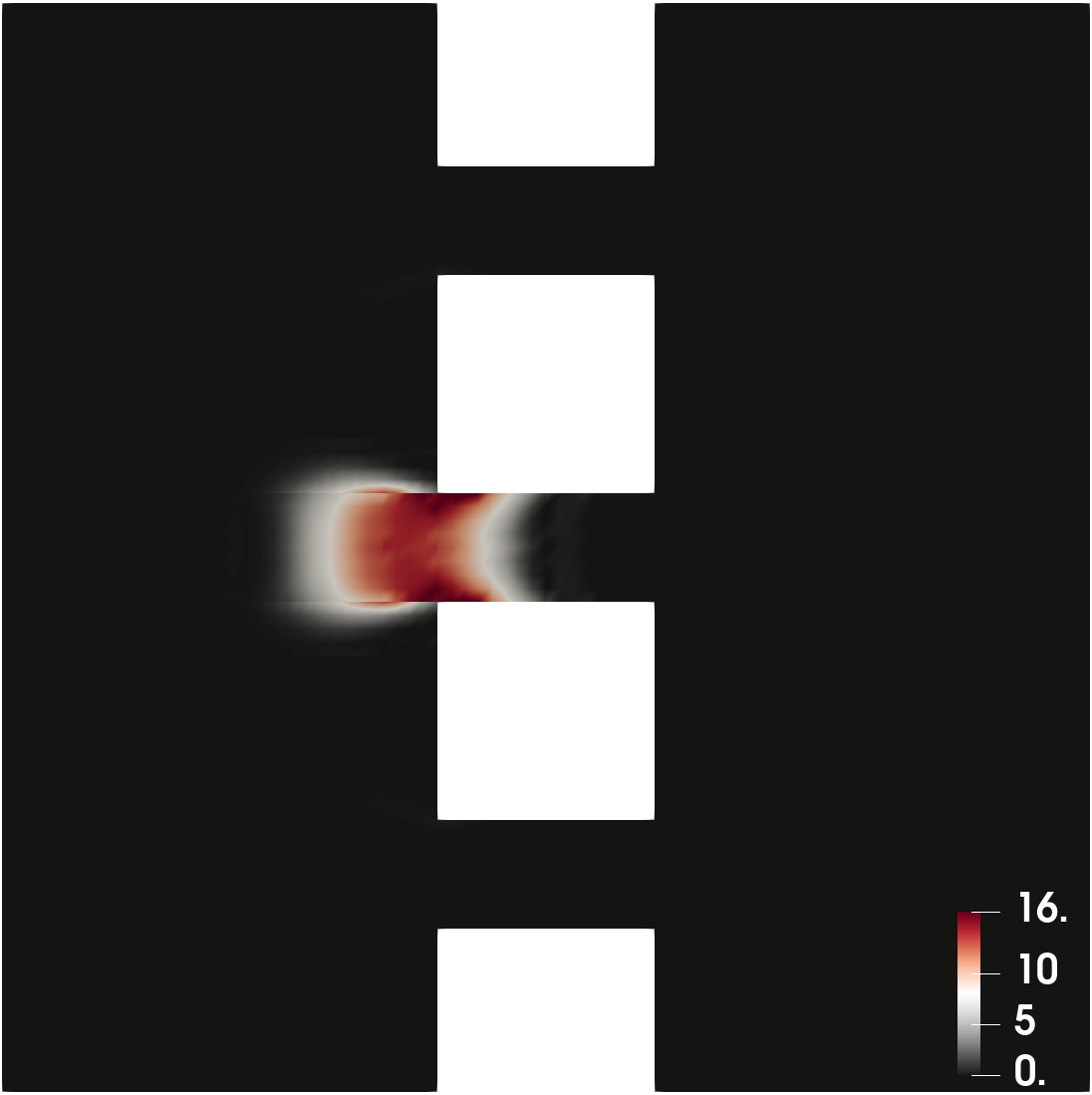}
\includegraphics[width=0.192\textwidth]{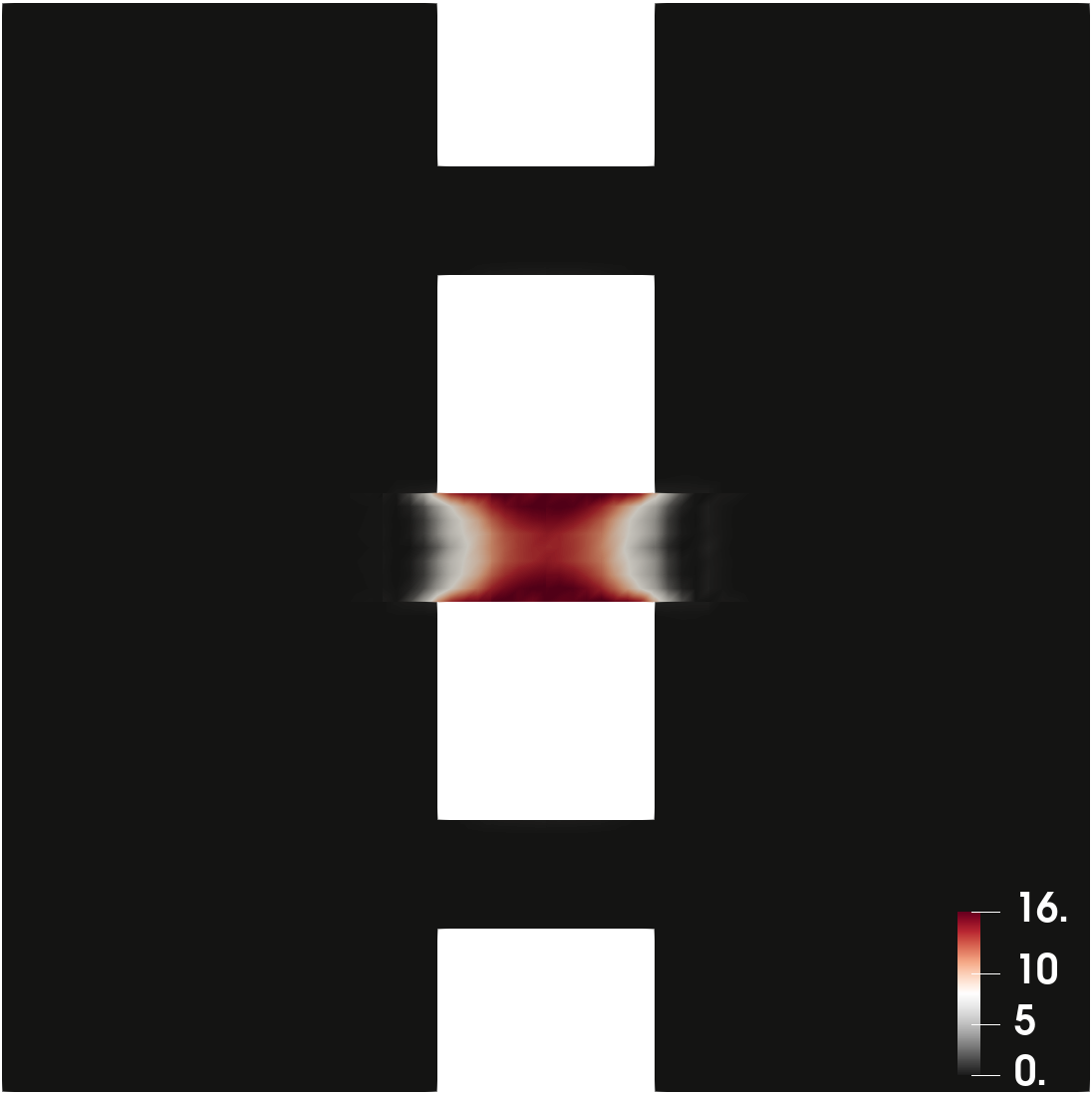}
\includegraphics[width=0.192\textwidth]{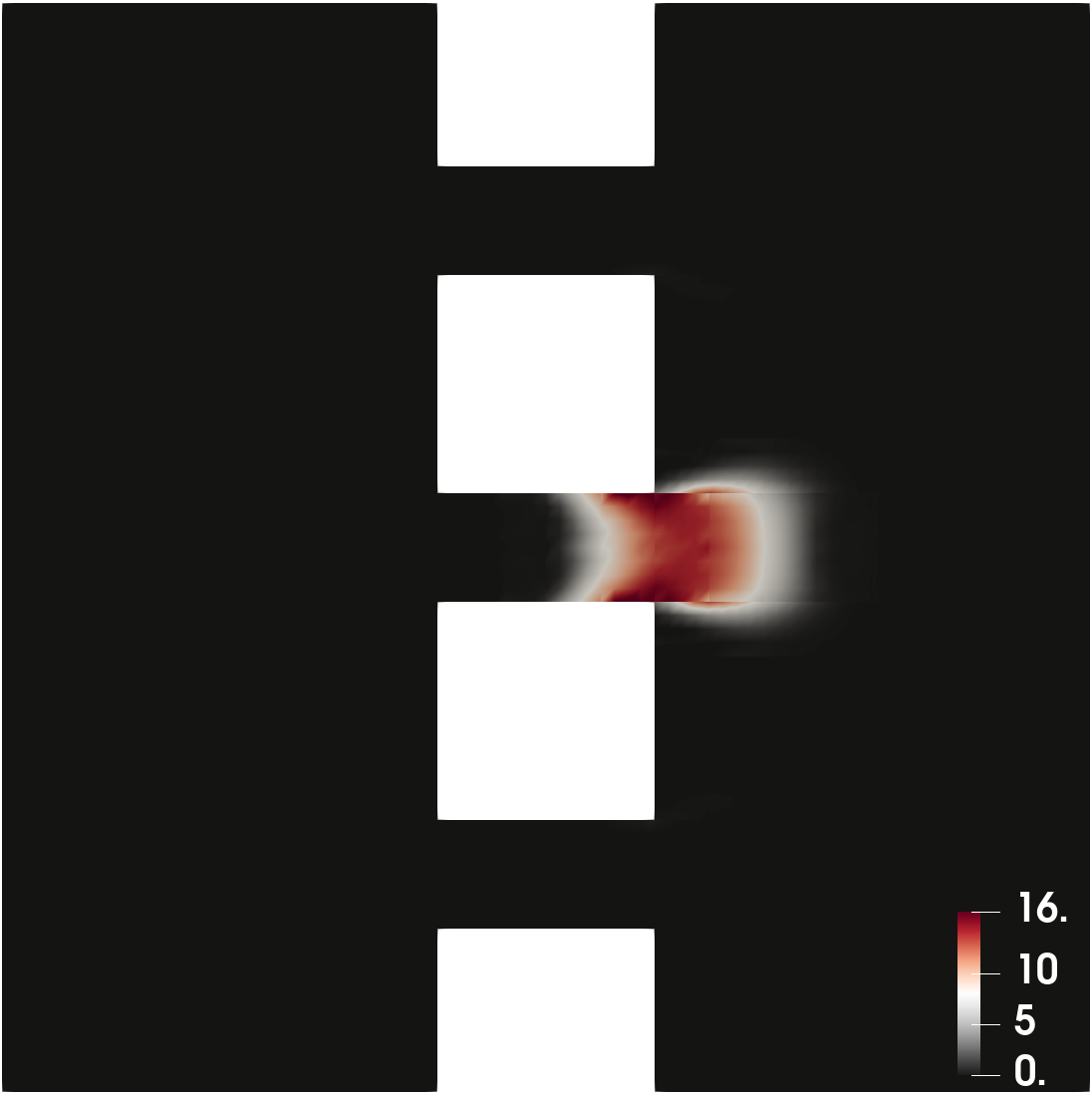}
\includegraphics[width=0.192\textwidth]{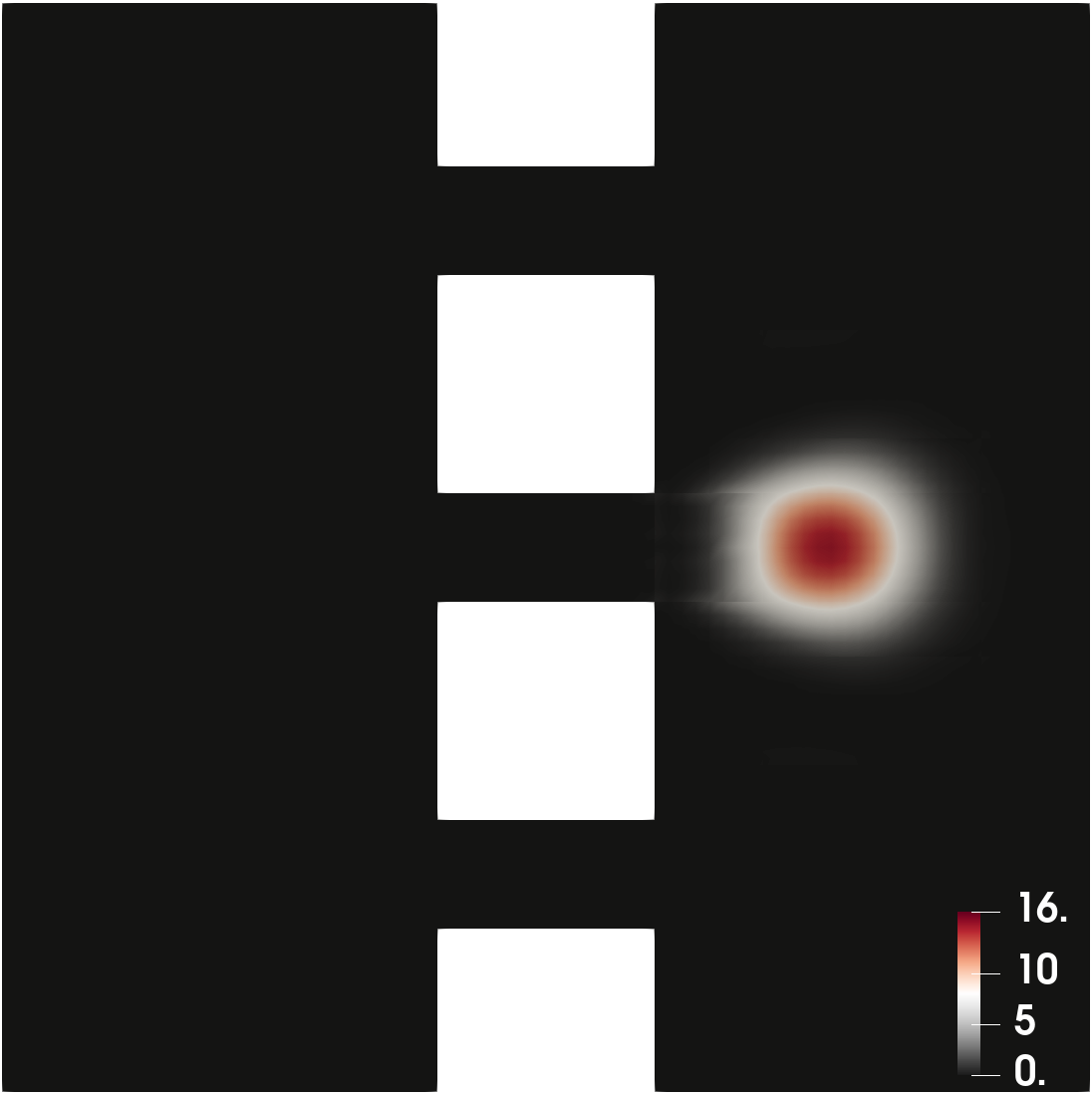}
}
\caption{Example \ref{ex2}. Snapshots of $\rho$ at 
$t=0.1,0.3,0.5,0.7,0.9$ 
(left to right).}
\label{fig:den-case2}
\end{figure}

Snapshots of the density contour at different times
are shown in Figure~\ref{fig:den-case2}.
The effects of different interaction 
cost functions on the density profile are clearly observed.

\subsection{MFG with obstacles}
\label{ex3}
We consider a similar setting as in Example \ref{ex2}, where we consider a MFG problem with terminal cost 
\[
\Gamma(\rho):=\begin{cases}
\frac12(\rho-\rho_T)^2 & \text{ if }\rho\ge0,\\
+\infty &\text{ otherwise,}
\end{cases}
\]
where the target density
\[
\rho_T:=
\frac{1}{2\pi\sigma^2}\left(\exp(-\frac{1}{2\sigma^2}|\bmx-(0.65,0.3)|^2)+
\exp(-\frac{1}{2\sigma^2}|\bmx-(0.65,-0.3)|^2)\right)
\]
with $\sigma=0.1$.
Note that we allow $\rho_T$ and $\rho_0$ to have different total masses here.

We apply the scheme \eqref{aug-mfg-h} with polynomial degree $k=3$ on the same mesh as in Example \ref{ex2}, and use the same stopping criterion.
The number of iterations needed for convergence for the 5 cases are recorded in Table~\ref{tab:iter-case3}, where again we find Case 2 has the smallest number of iterations.
\begin{table}[ht!]
\centering
\begin{tabular}{cccccccc}
 & Case 1 & Case 2 & Case 3 & Case 4 & Case 5&\\\hline
iterations 
&3510&82&476&503&798\\
\end{tabular}
\caption{Example \ref{ex3}. Number of ALG2 iterations for each case.
}
\label{tab:iter-case3}
\end{table}

\begin{figure}[tb]
\centering
\subfigure[Case 1]{
\label{fig:1x}
\includegraphics[width=0.192\textwidth]{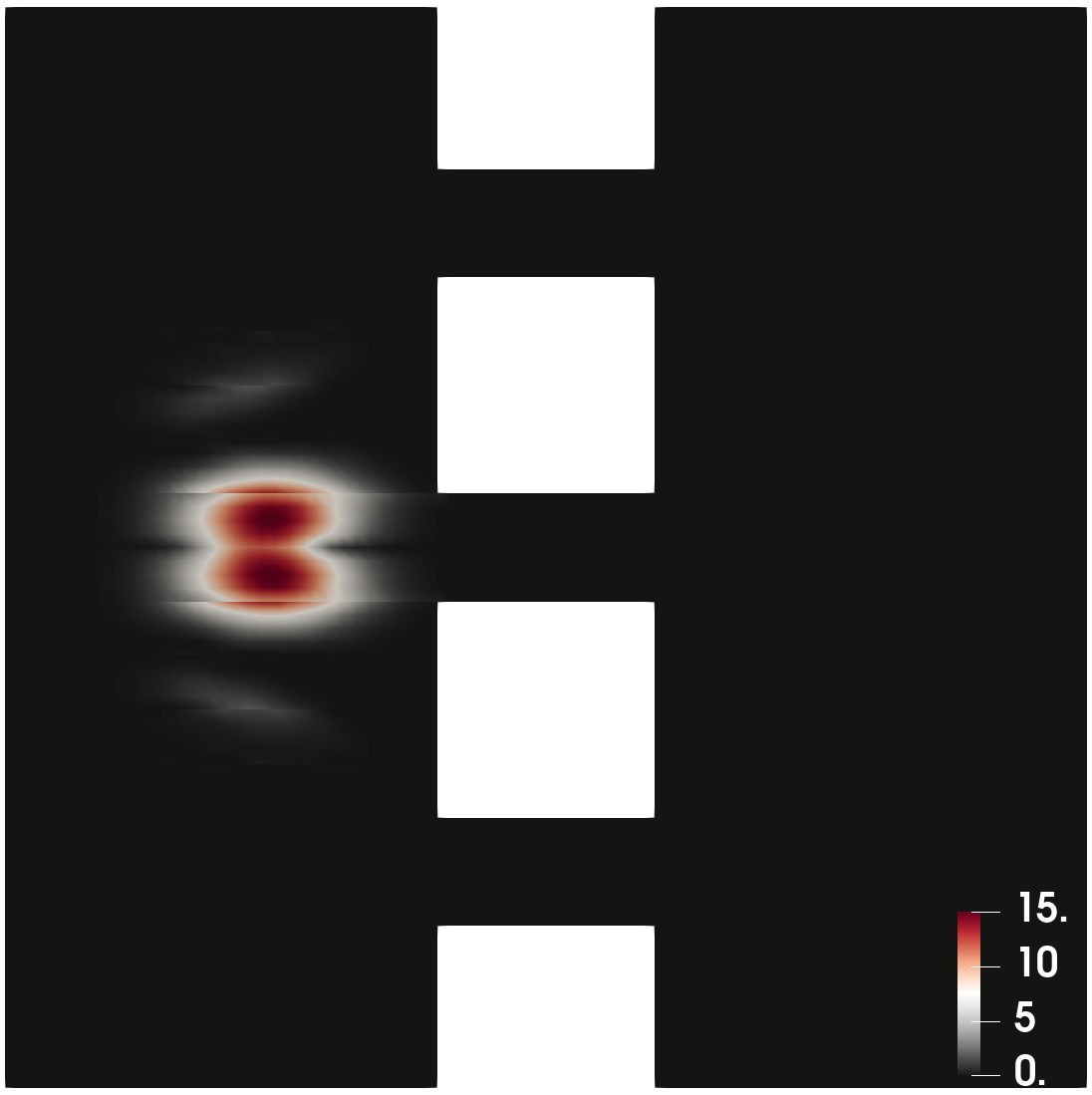}
\includegraphics[width=0.192\textwidth]{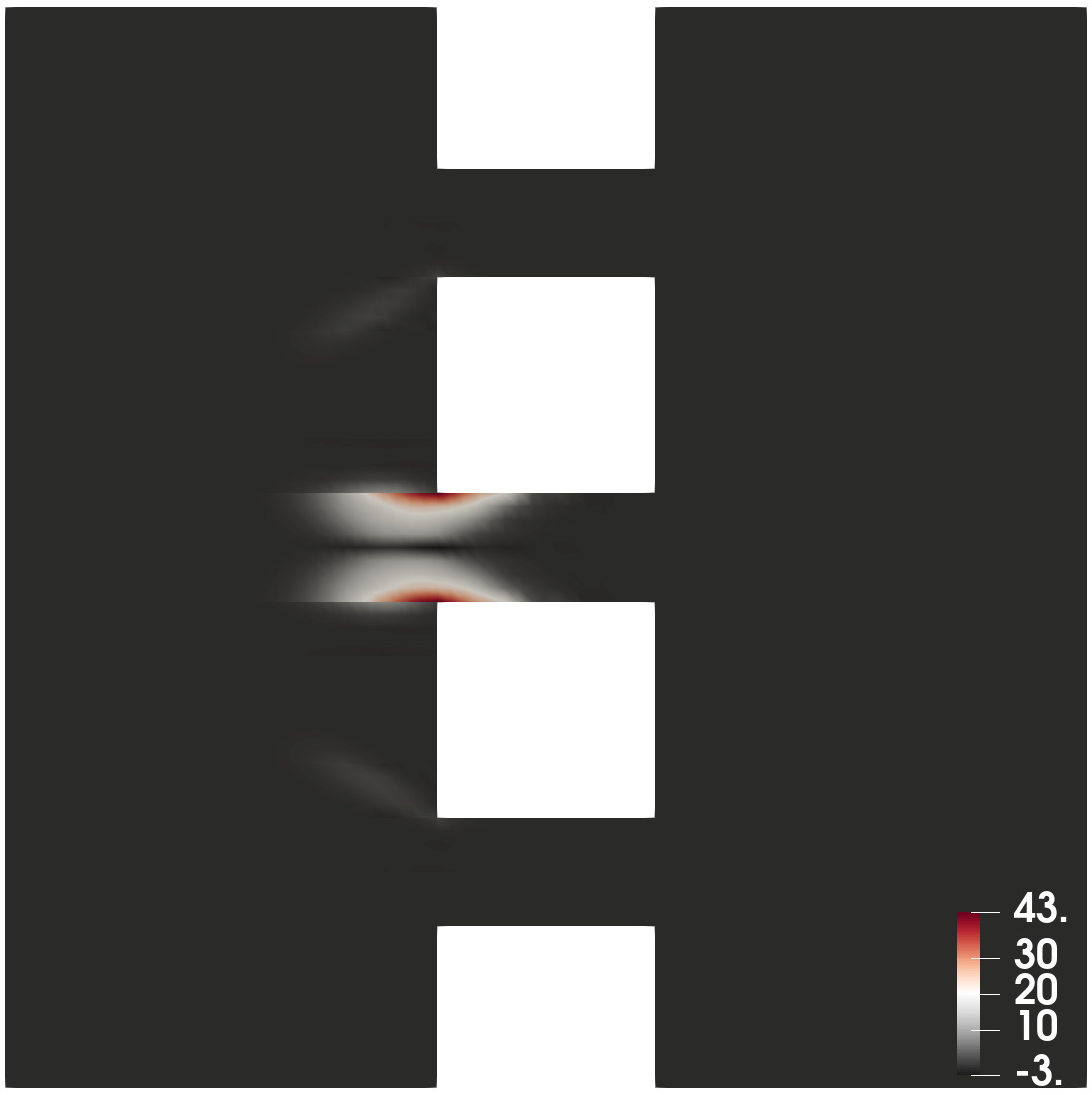}
\includegraphics[width=0.192\textwidth]{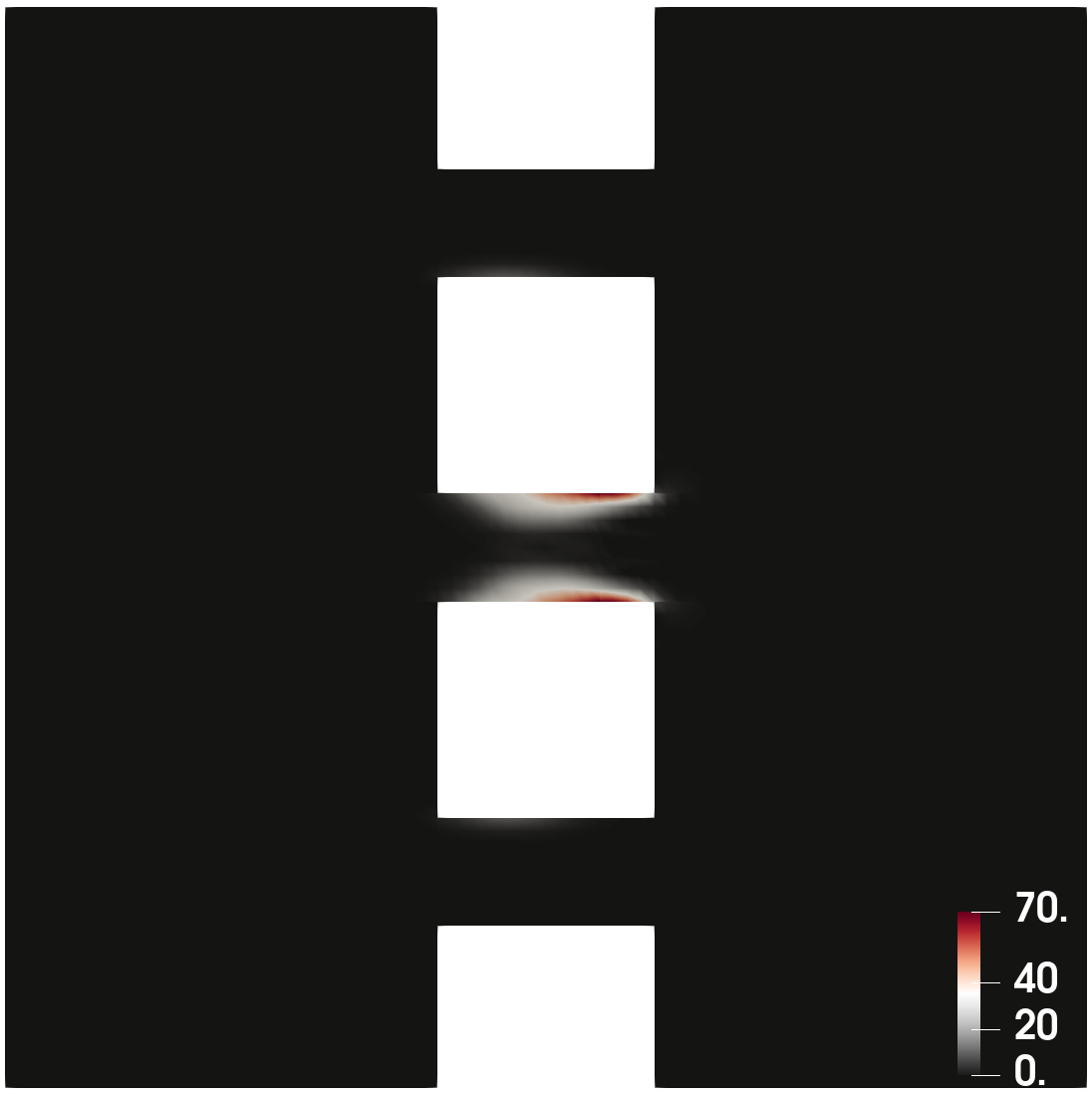}
\includegraphics[width=0.192\textwidth]{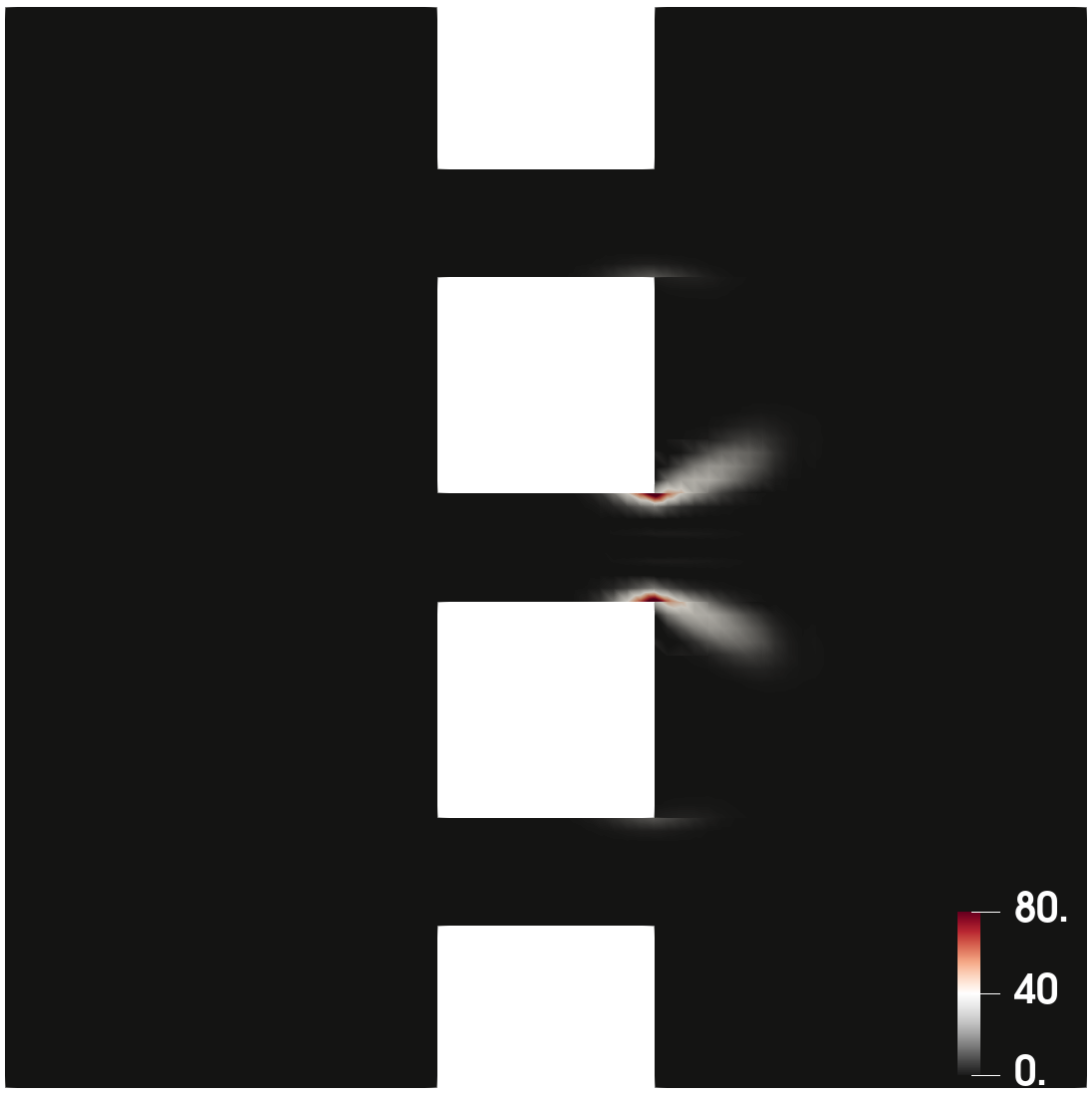}
\includegraphics[width=0.192\textwidth]{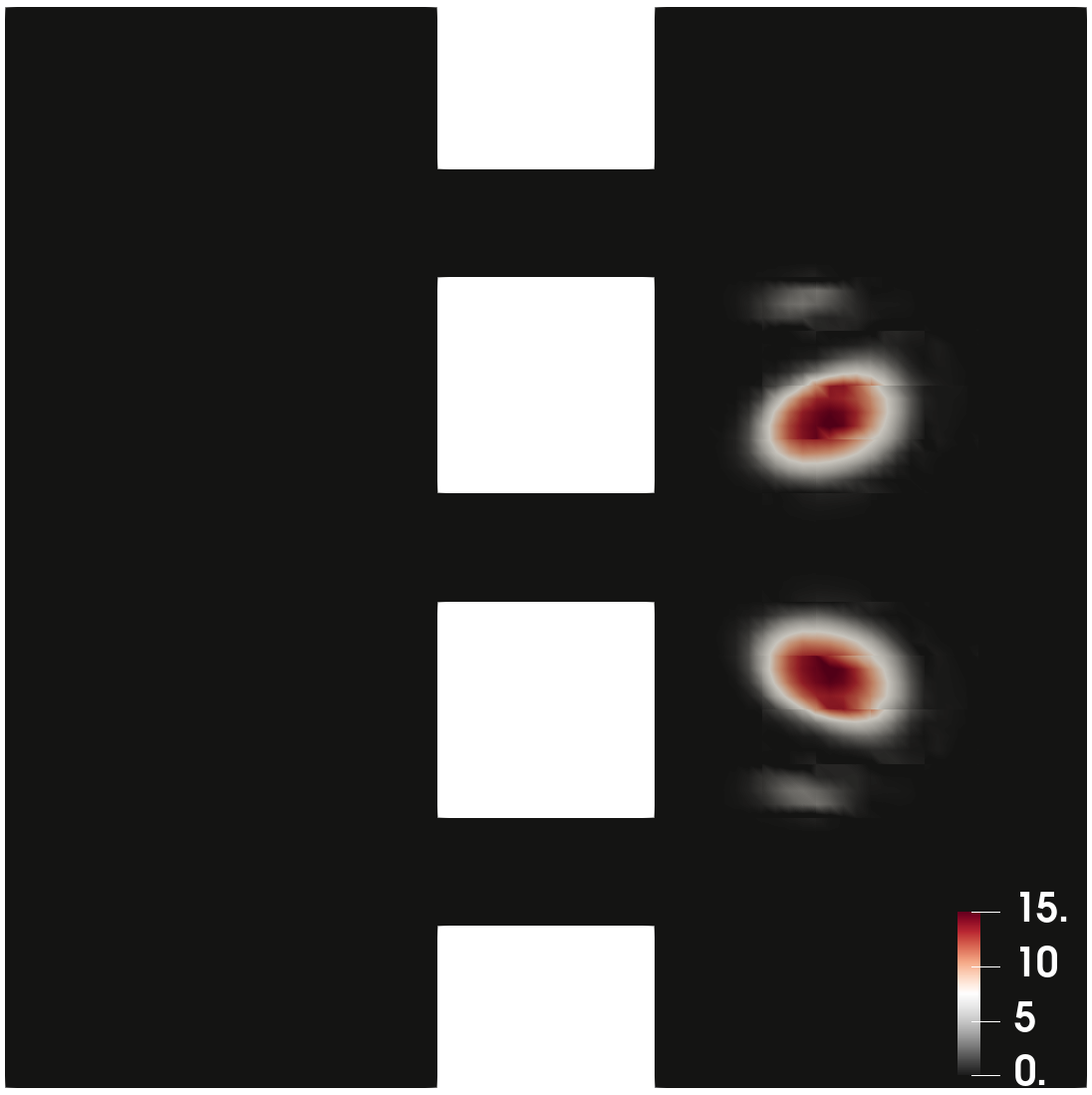}
}
\subfigure[Case 2]{
\label{fig:2x}
\includegraphics[width=0.192\textwidth]{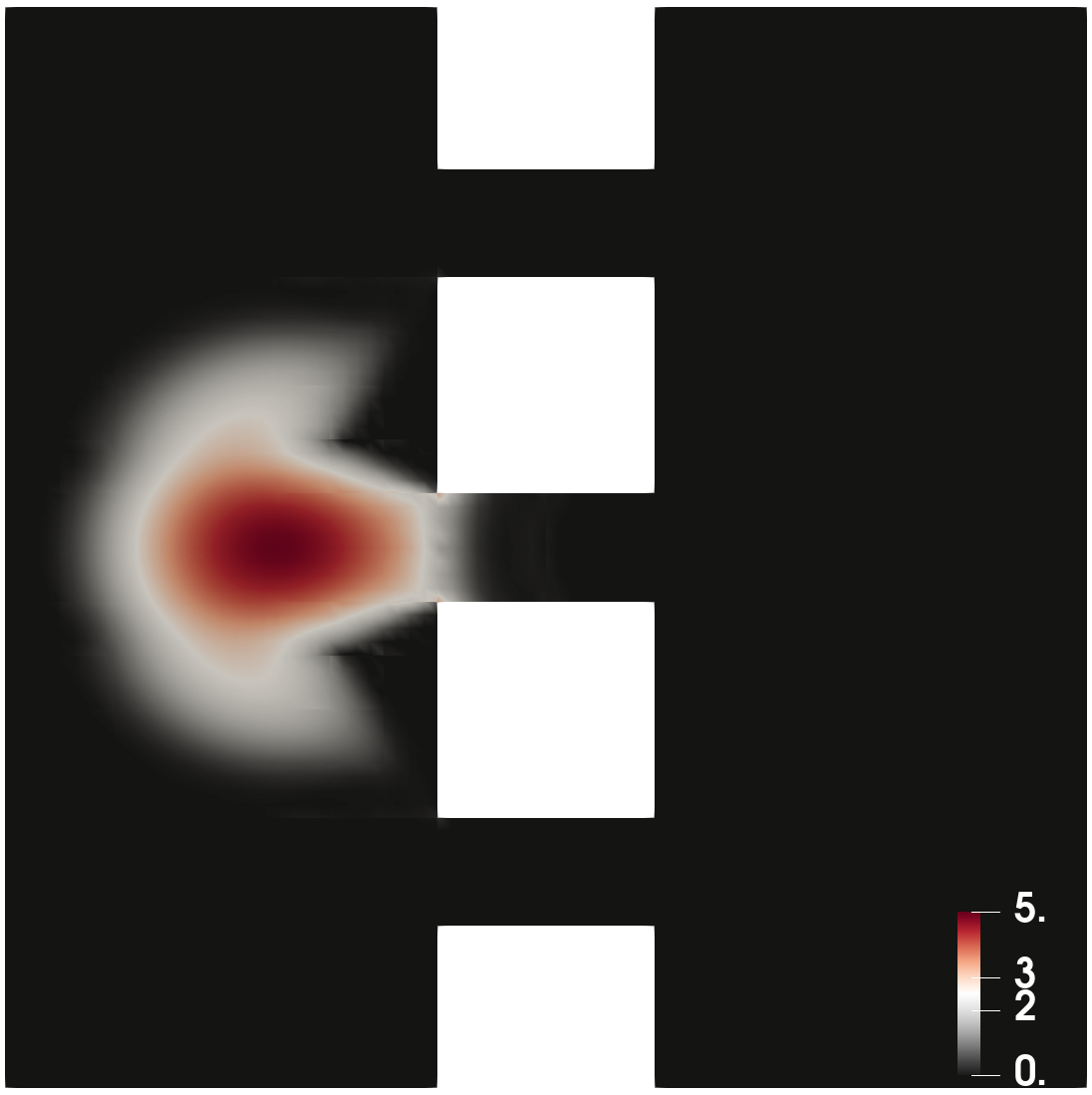}
\includegraphics[width=0.192\textwidth]{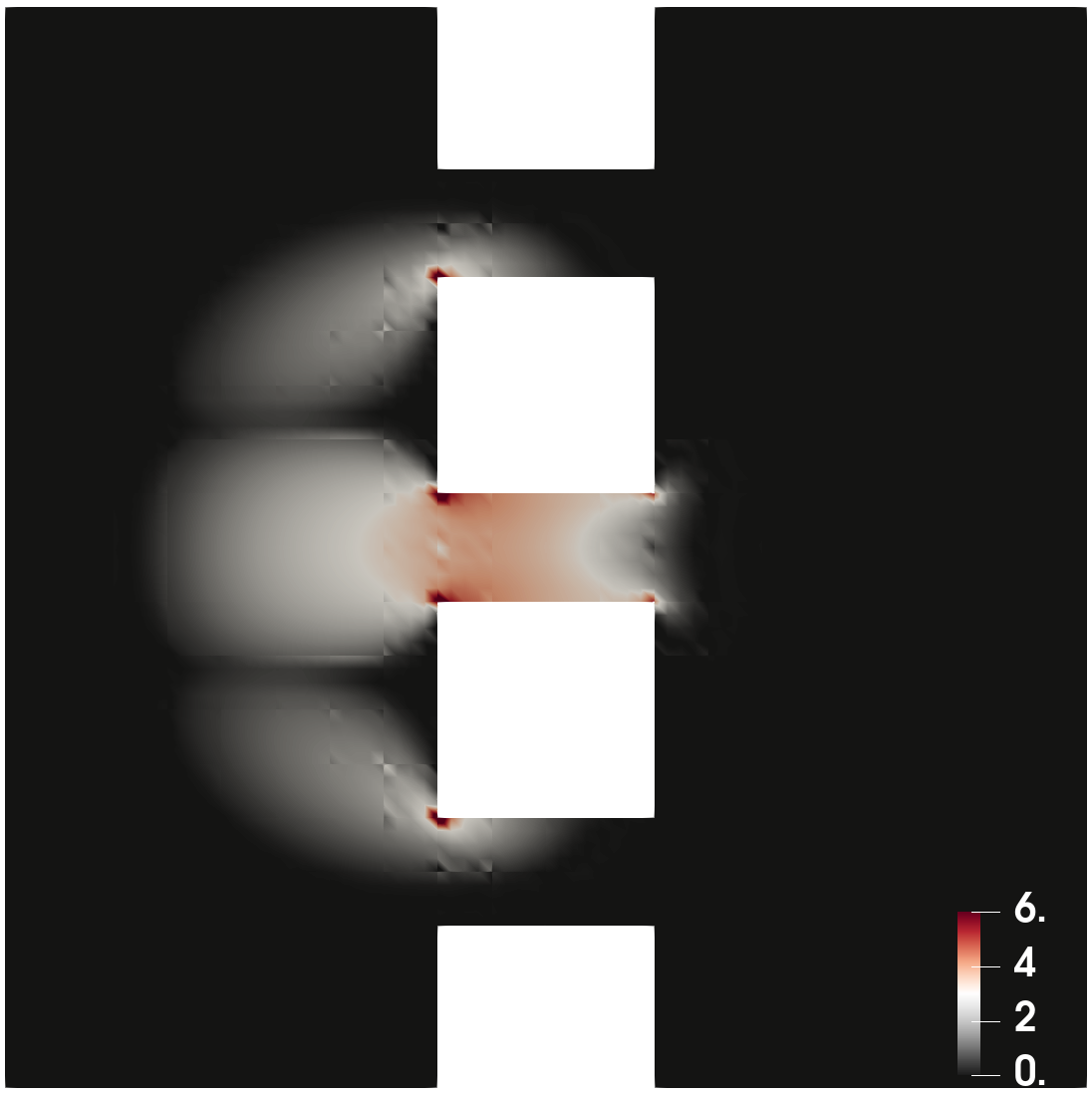}
\includegraphics[width=0.192\textwidth]{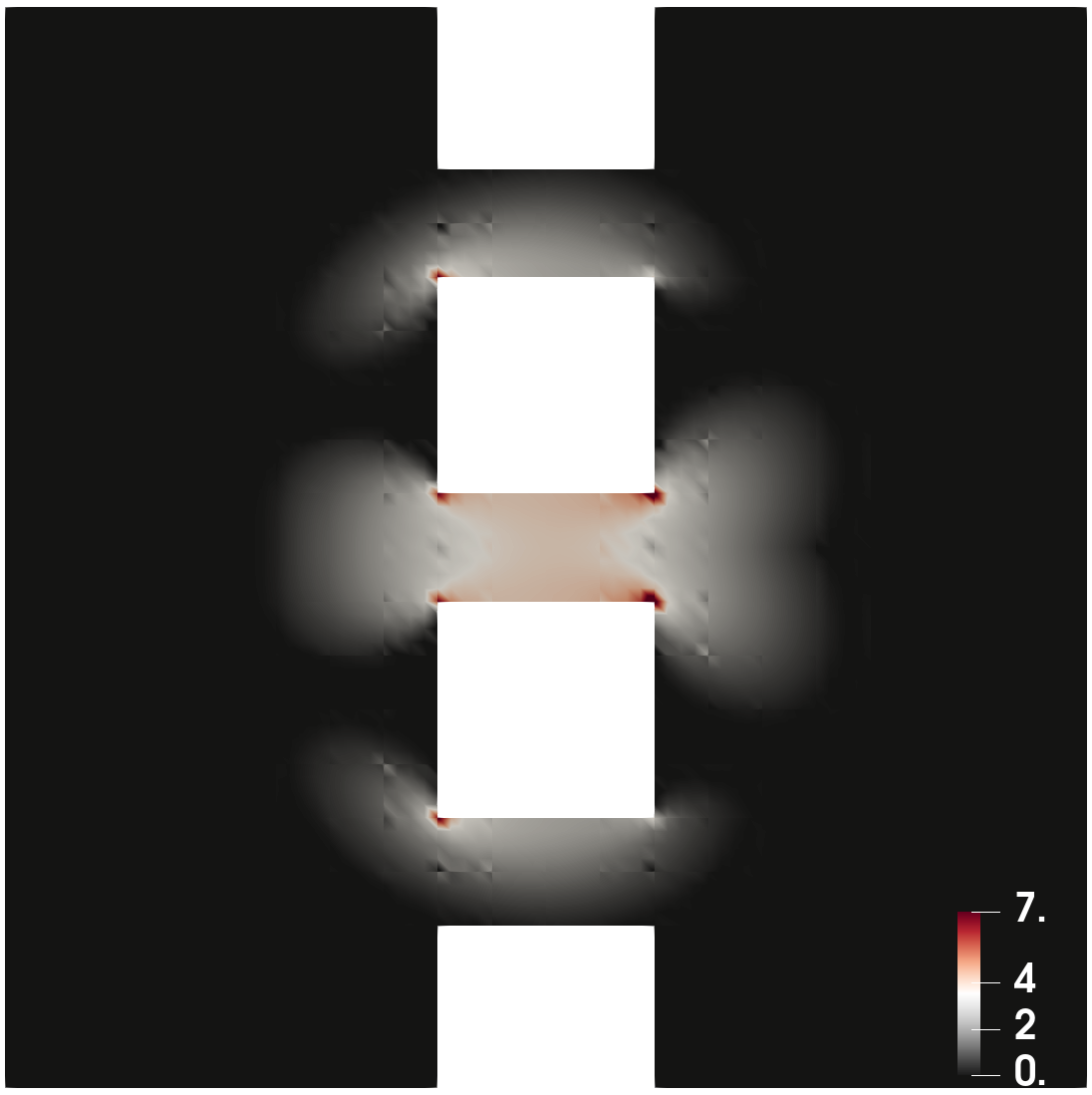}
\includegraphics[width=0.192\textwidth]{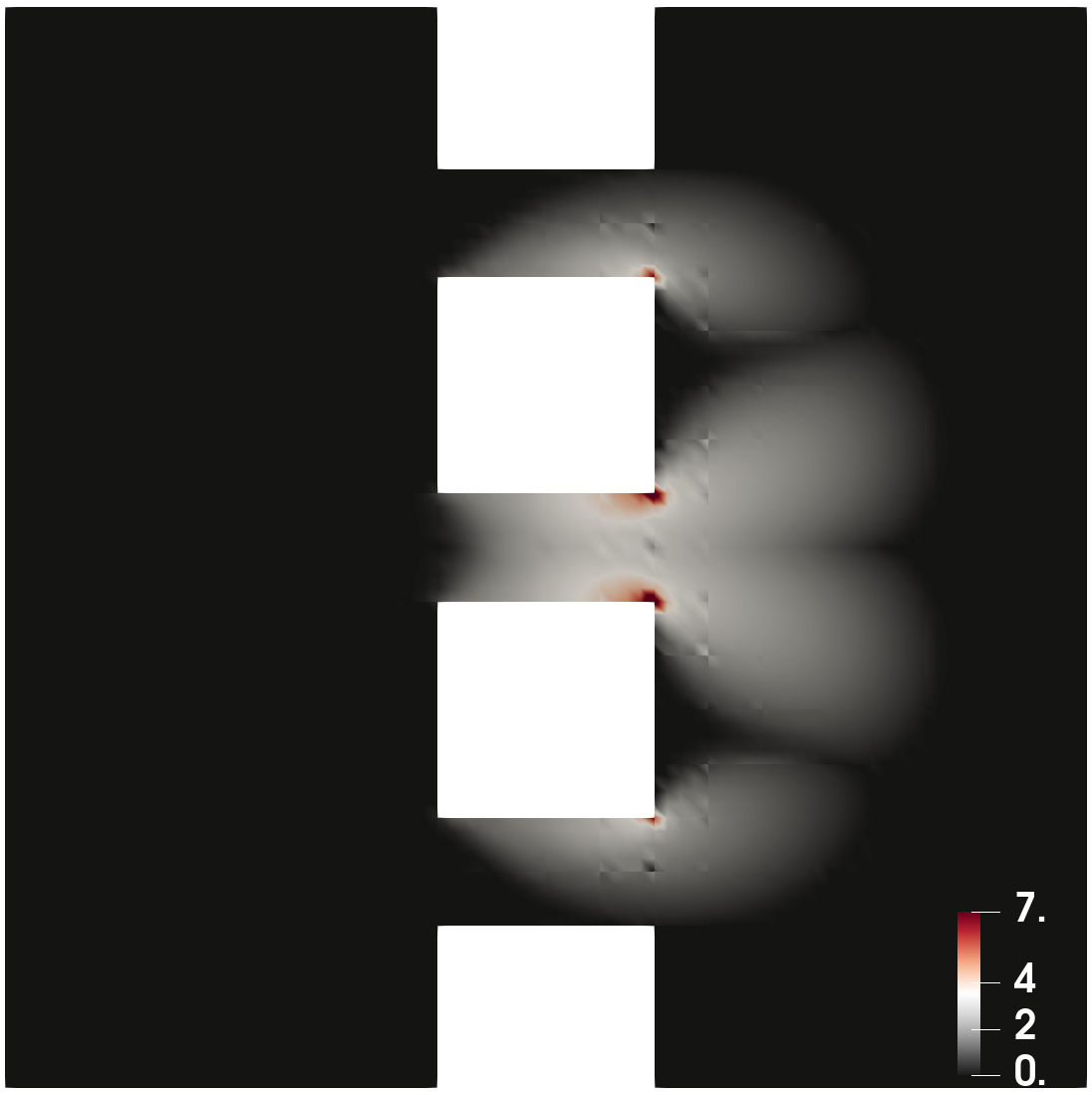}
\includegraphics[width=0.192\textwidth]{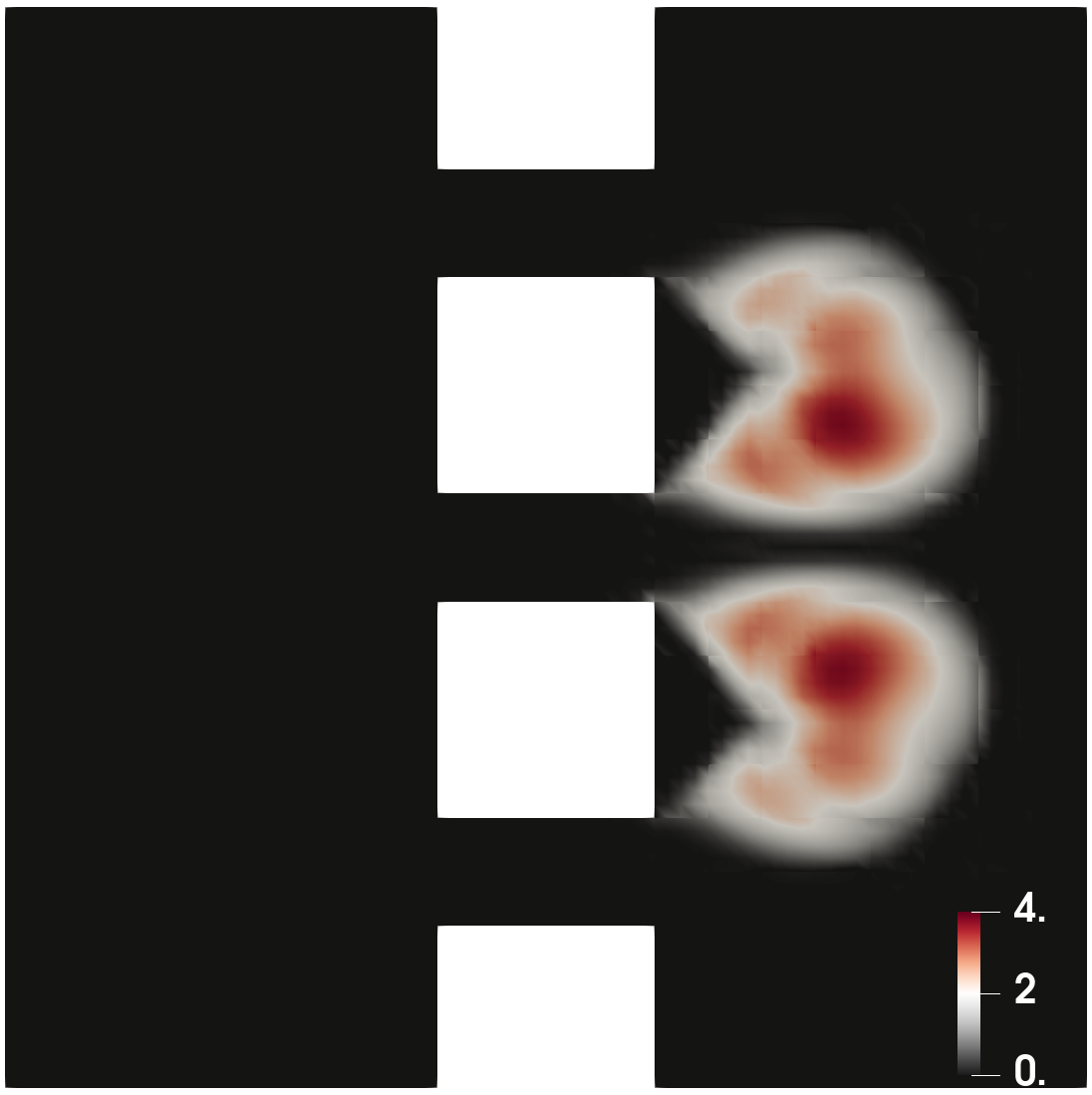}
}
\subfigure[Case 3]{
\label{fig:3x}
\includegraphics[width=0.192\textwidth]{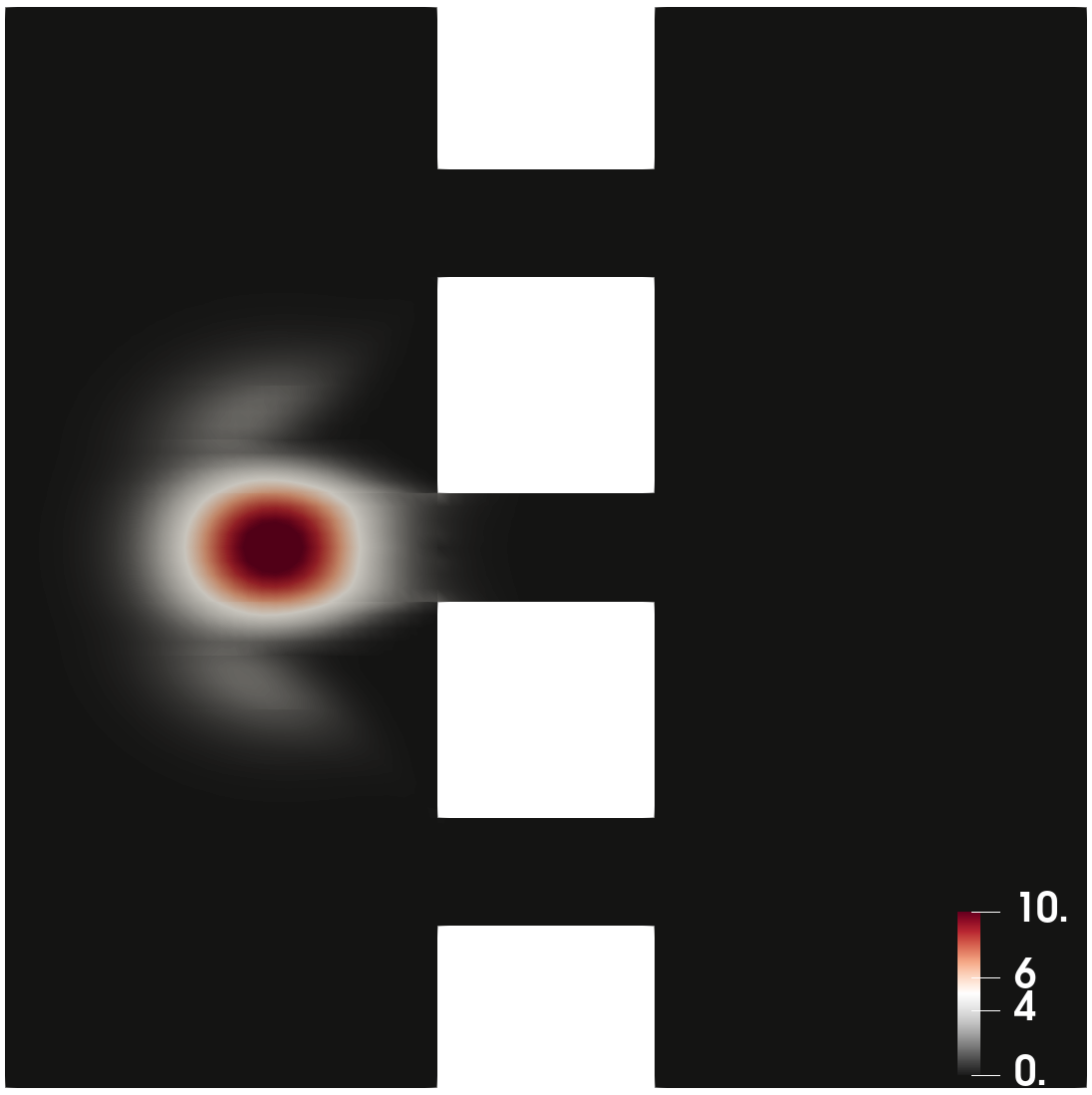}
\includegraphics[width=0.192\textwidth]{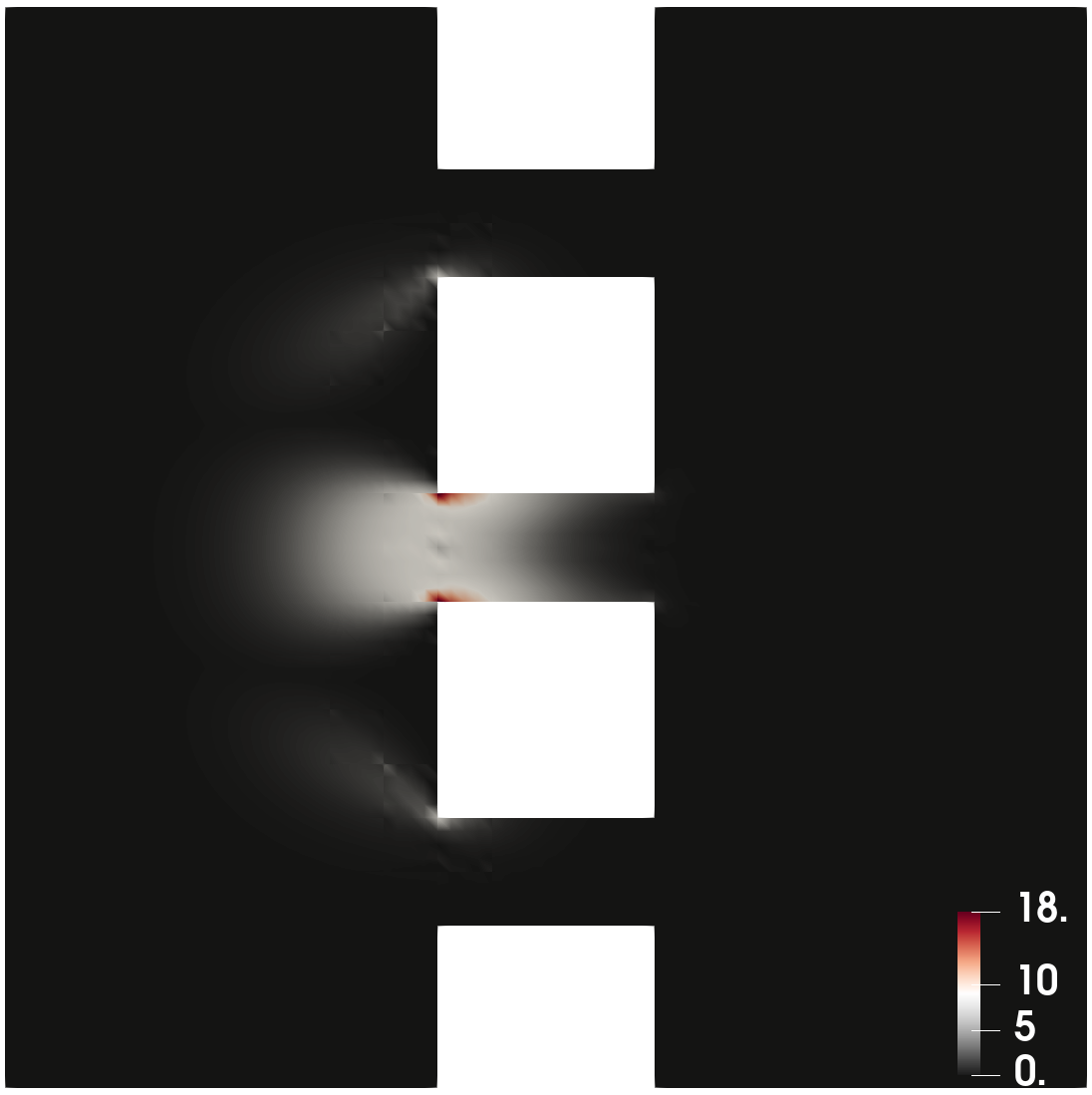}
\includegraphics[width=0.192\textwidth]{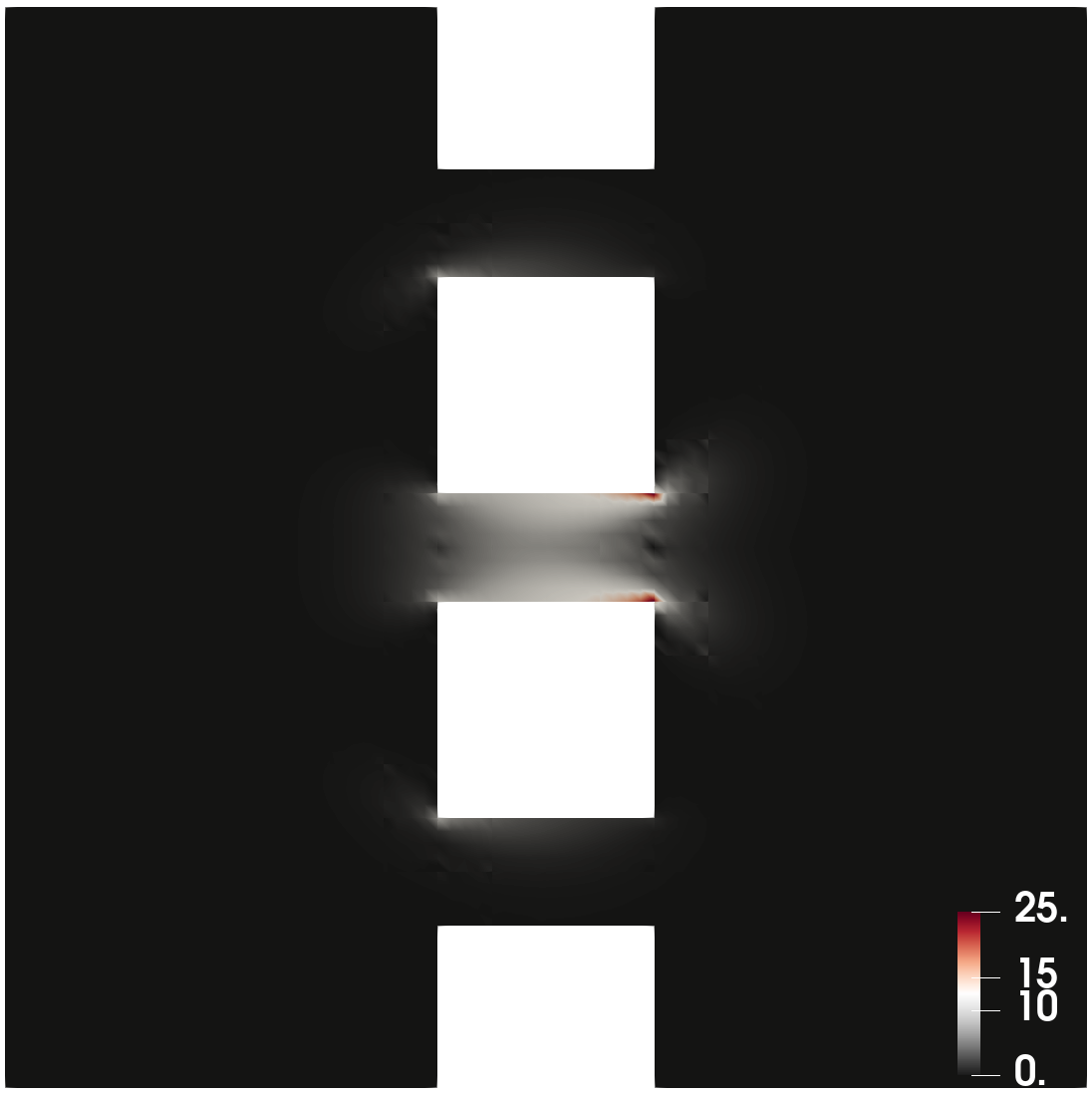}
\includegraphics[width=0.192\textwidth]{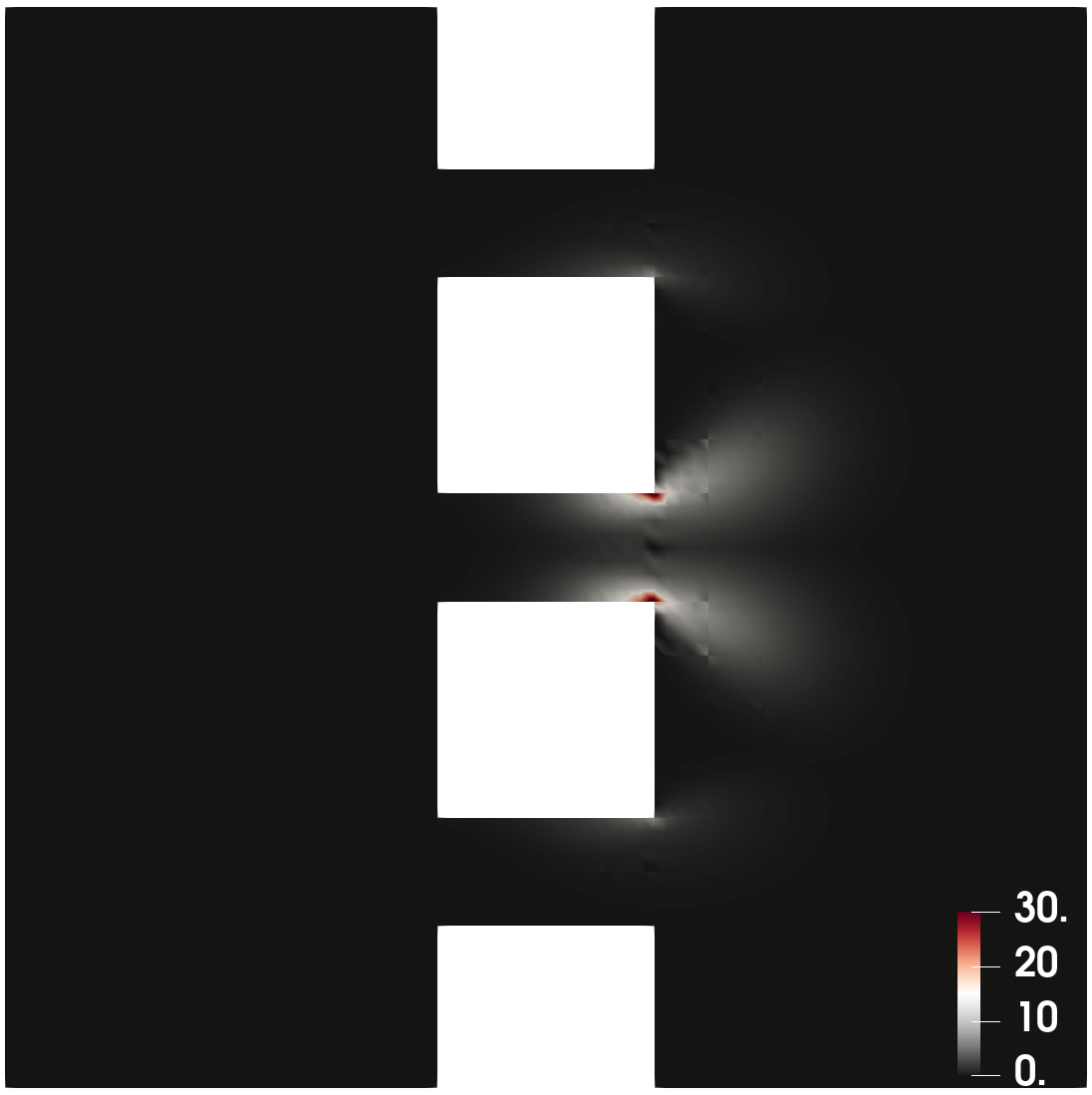}
\includegraphics[width=0.192\textwidth]{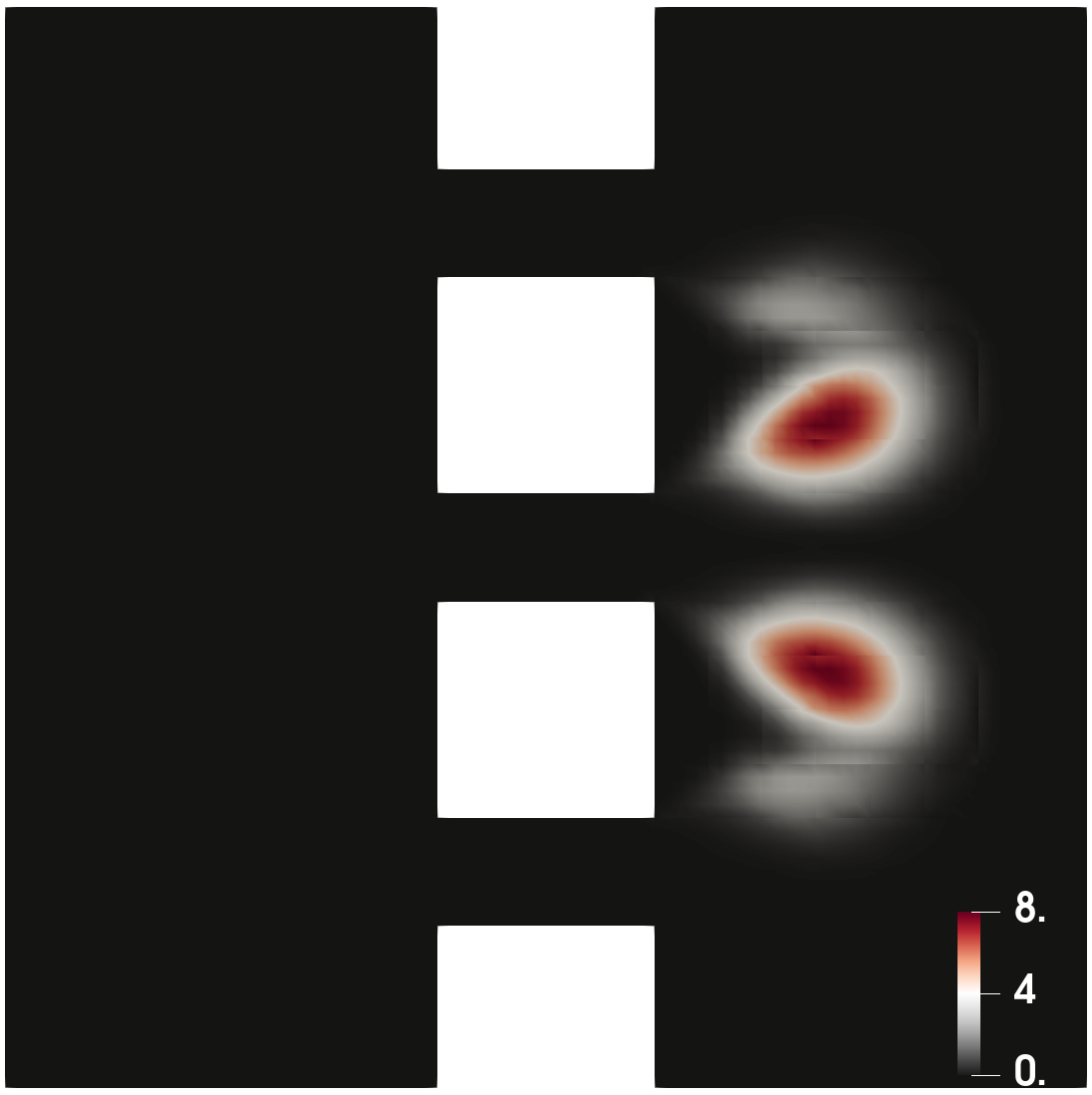}
}
\subfigure[Case 4]{
\label{fig:4x}
\includegraphics[width=0.192\textwidth]{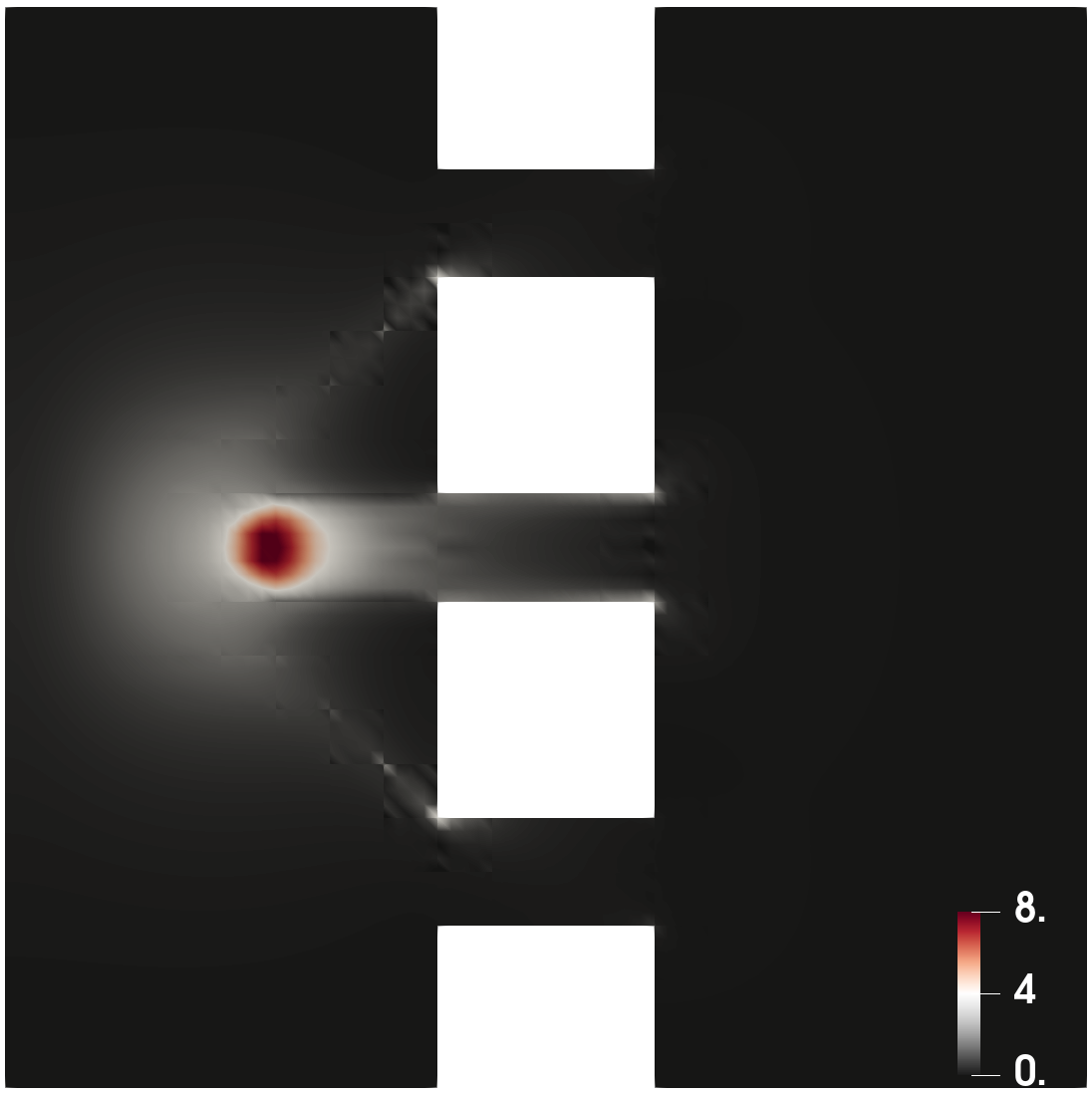}
\includegraphics[width=0.192\textwidth]{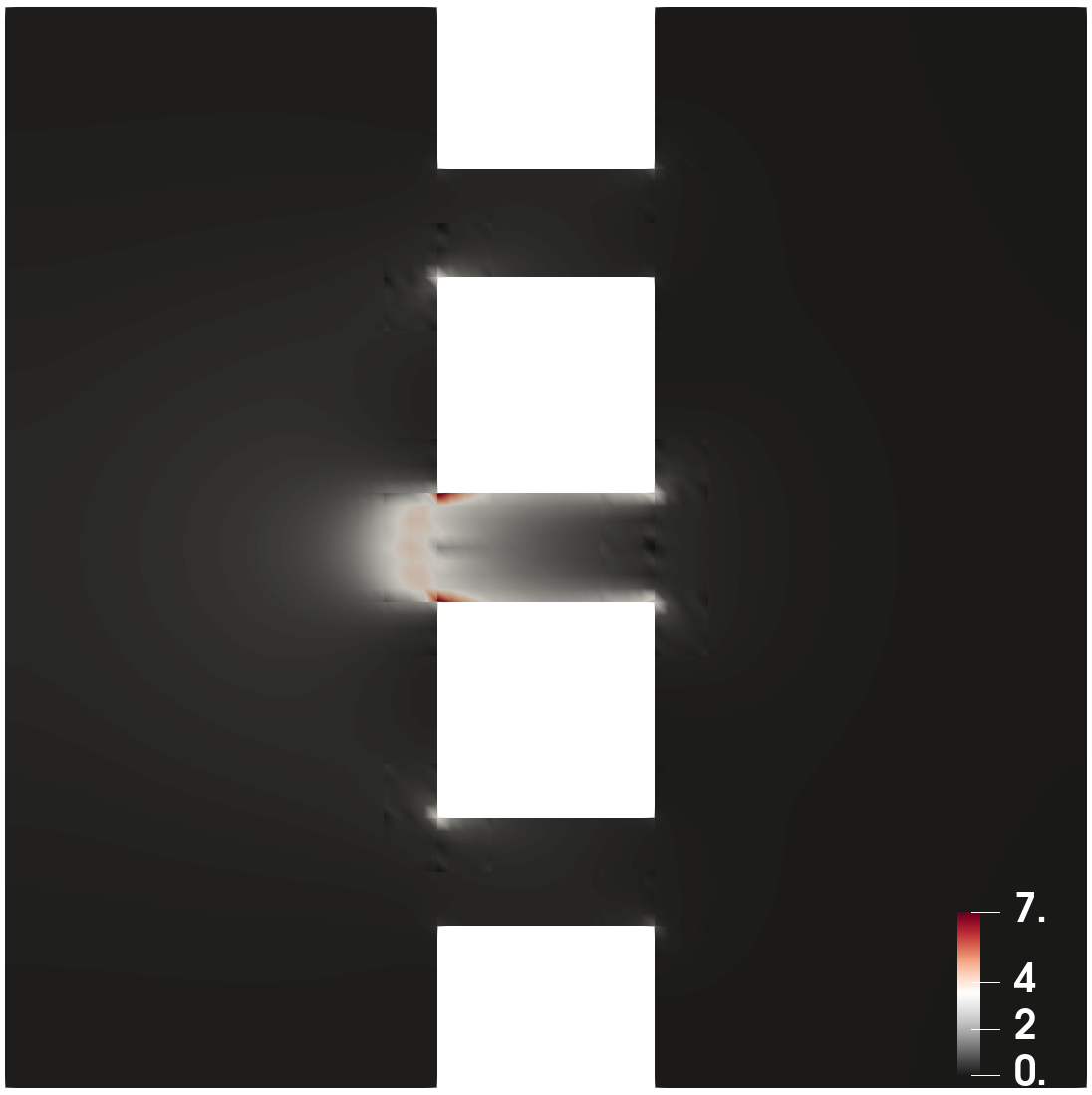}
\includegraphics[width=0.192\textwidth]{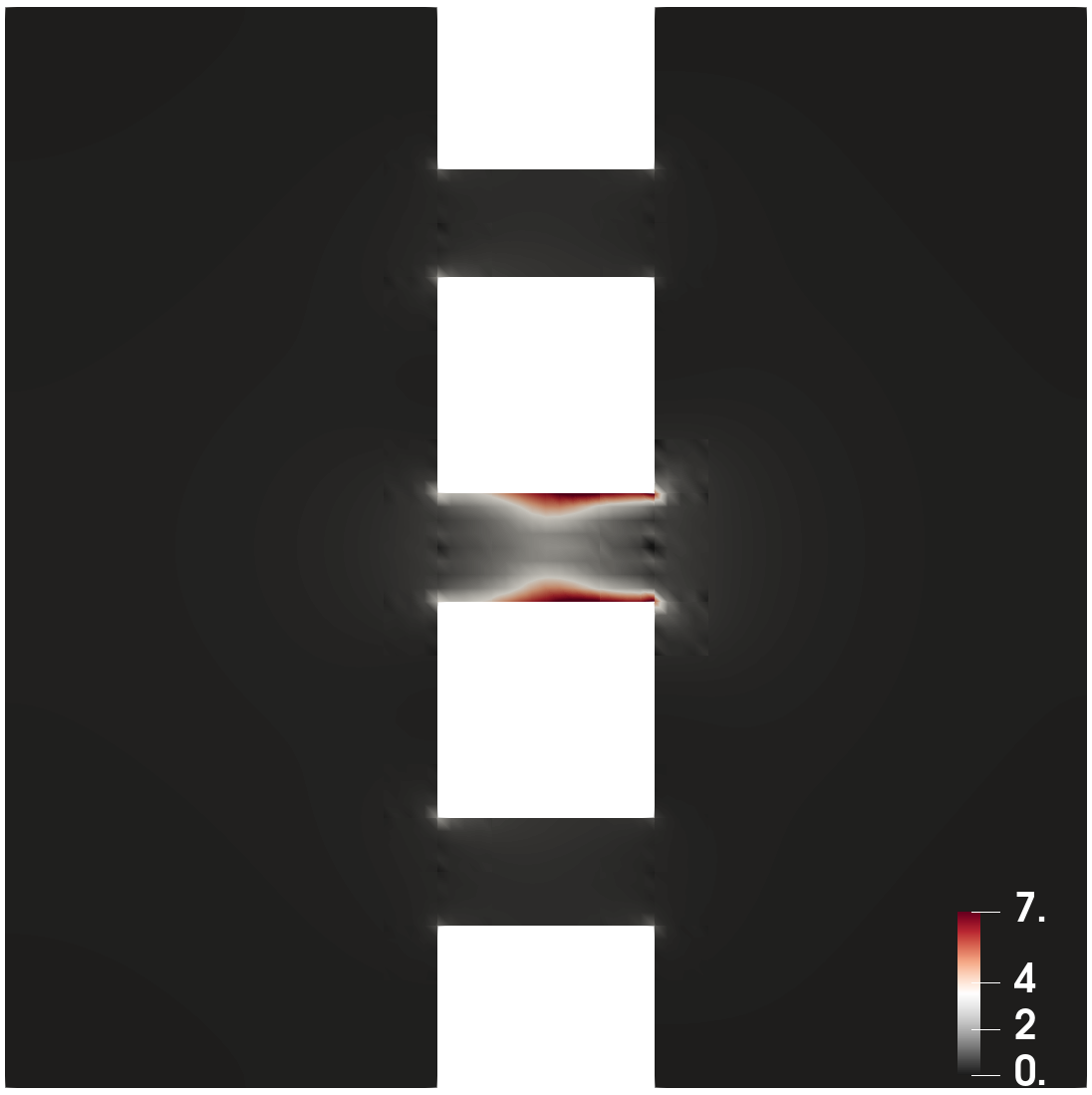}
\includegraphics[width=0.192\textwidth]{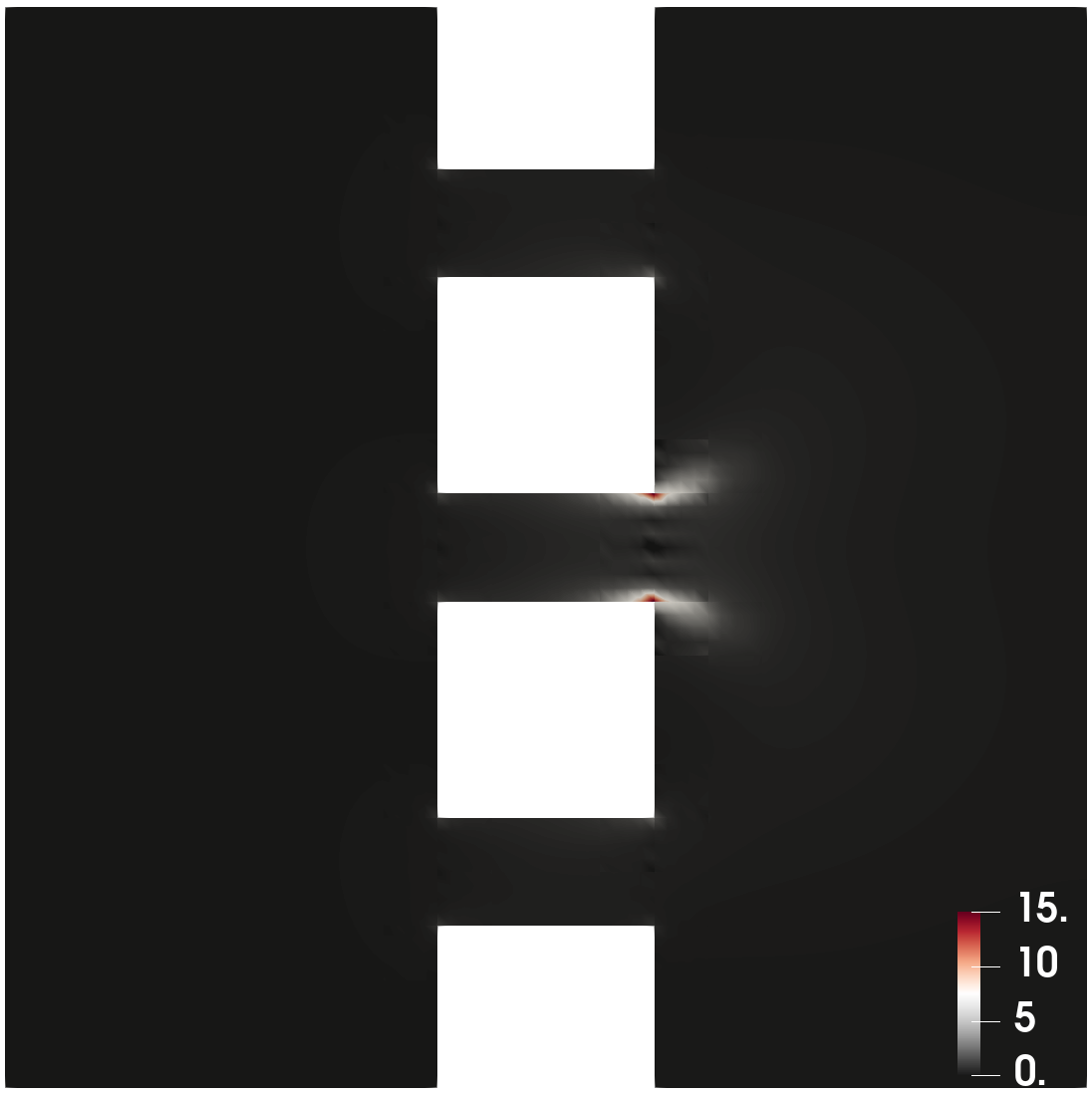}
\includegraphics[width=0.192\textwidth]{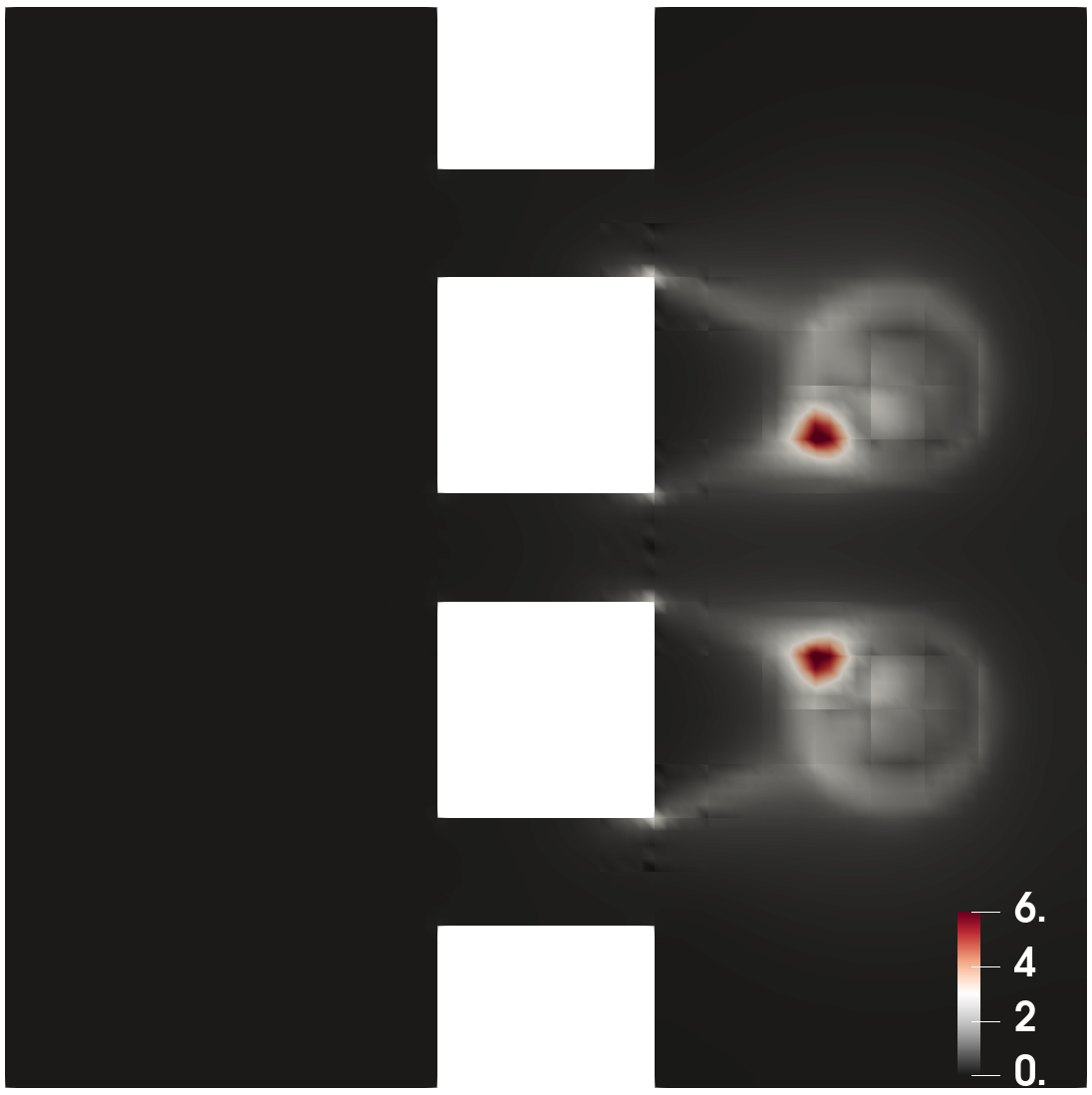}
}
\subfigure[Case 5]{
\label{fig:5x}
\includegraphics[width=0.192\textwidth]{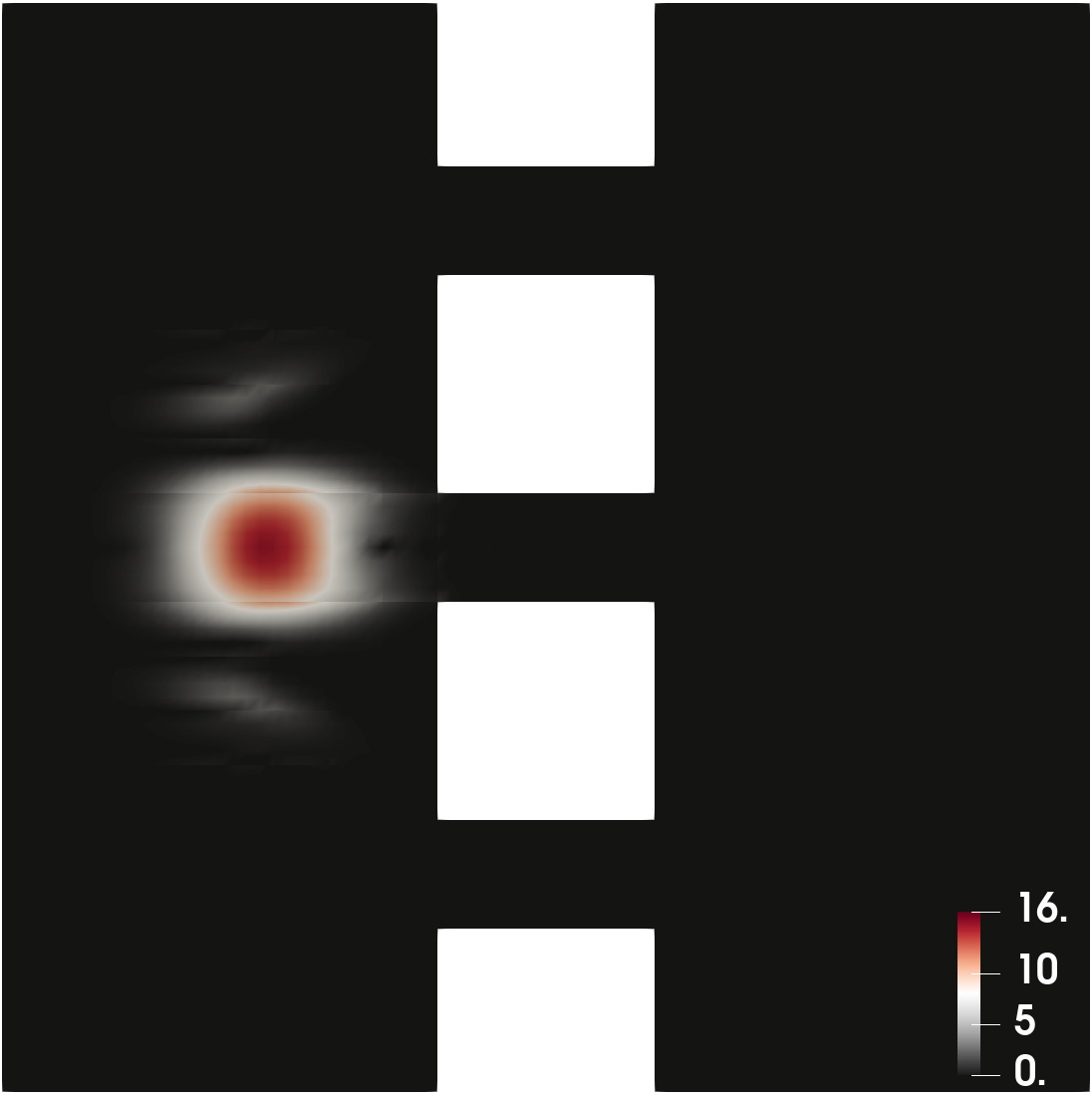}
\includegraphics[width=0.192\textwidth]{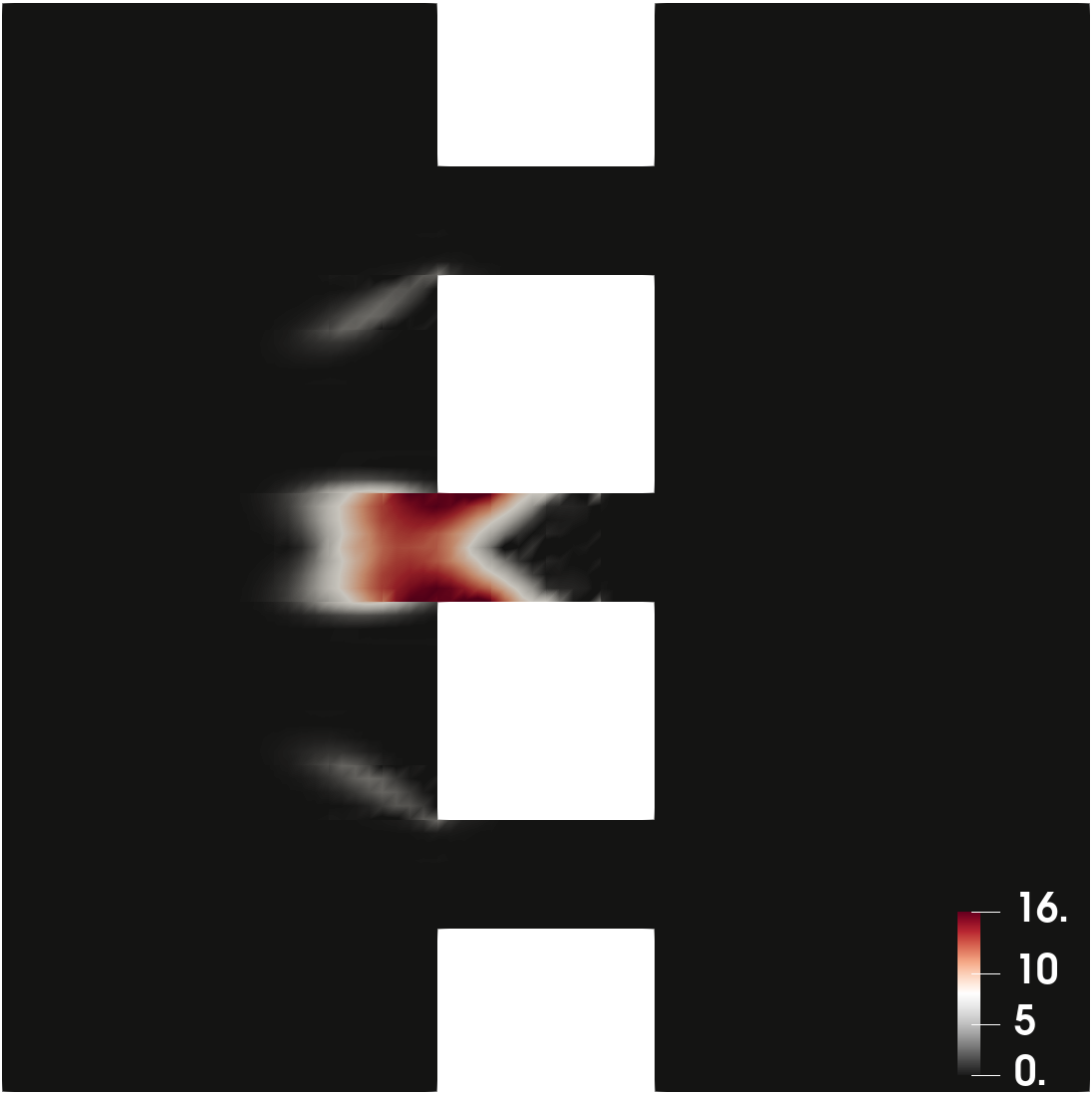}
\includegraphics[width=0.192\textwidth]{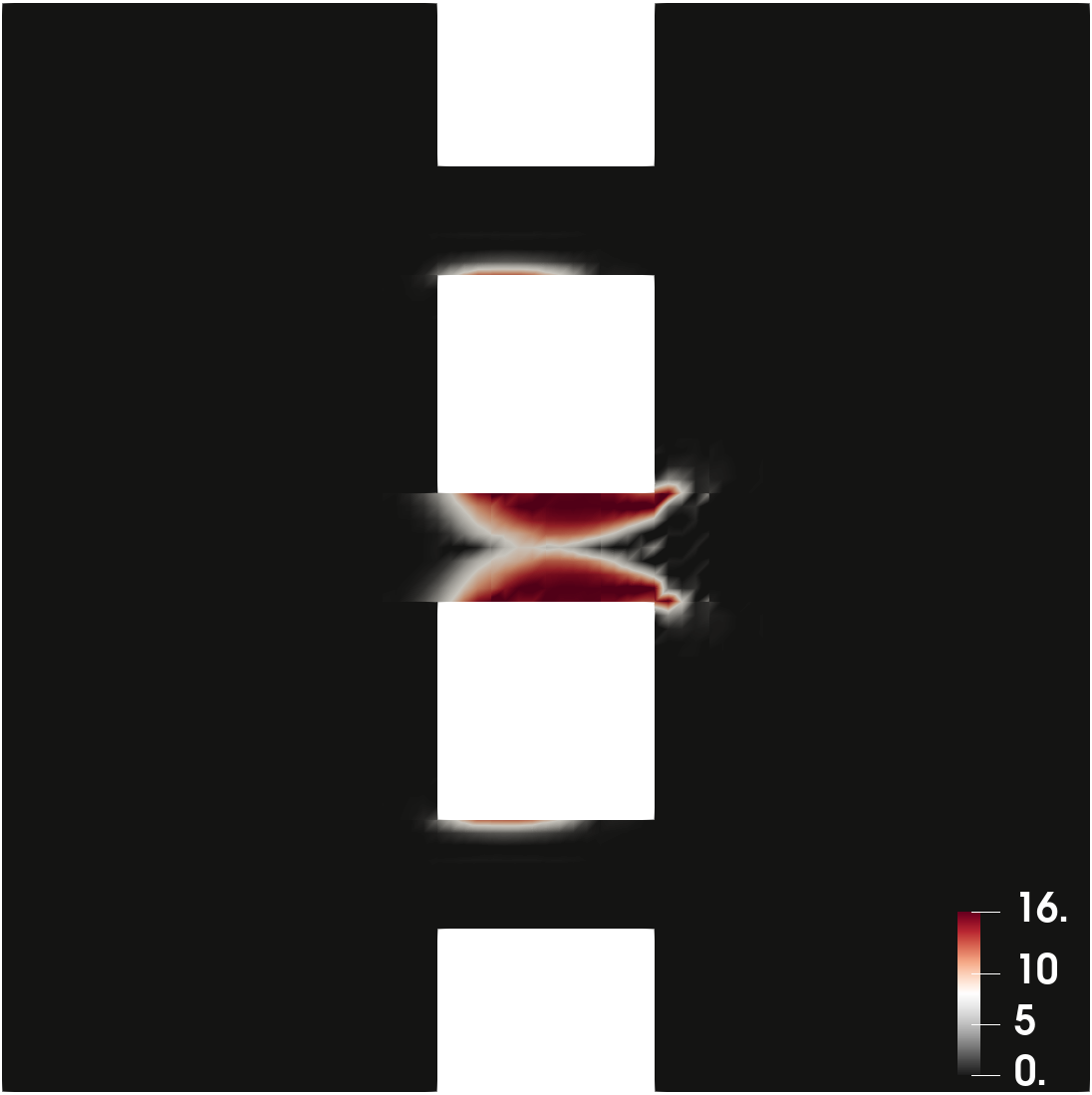}
\includegraphics[width=0.192\textwidth]{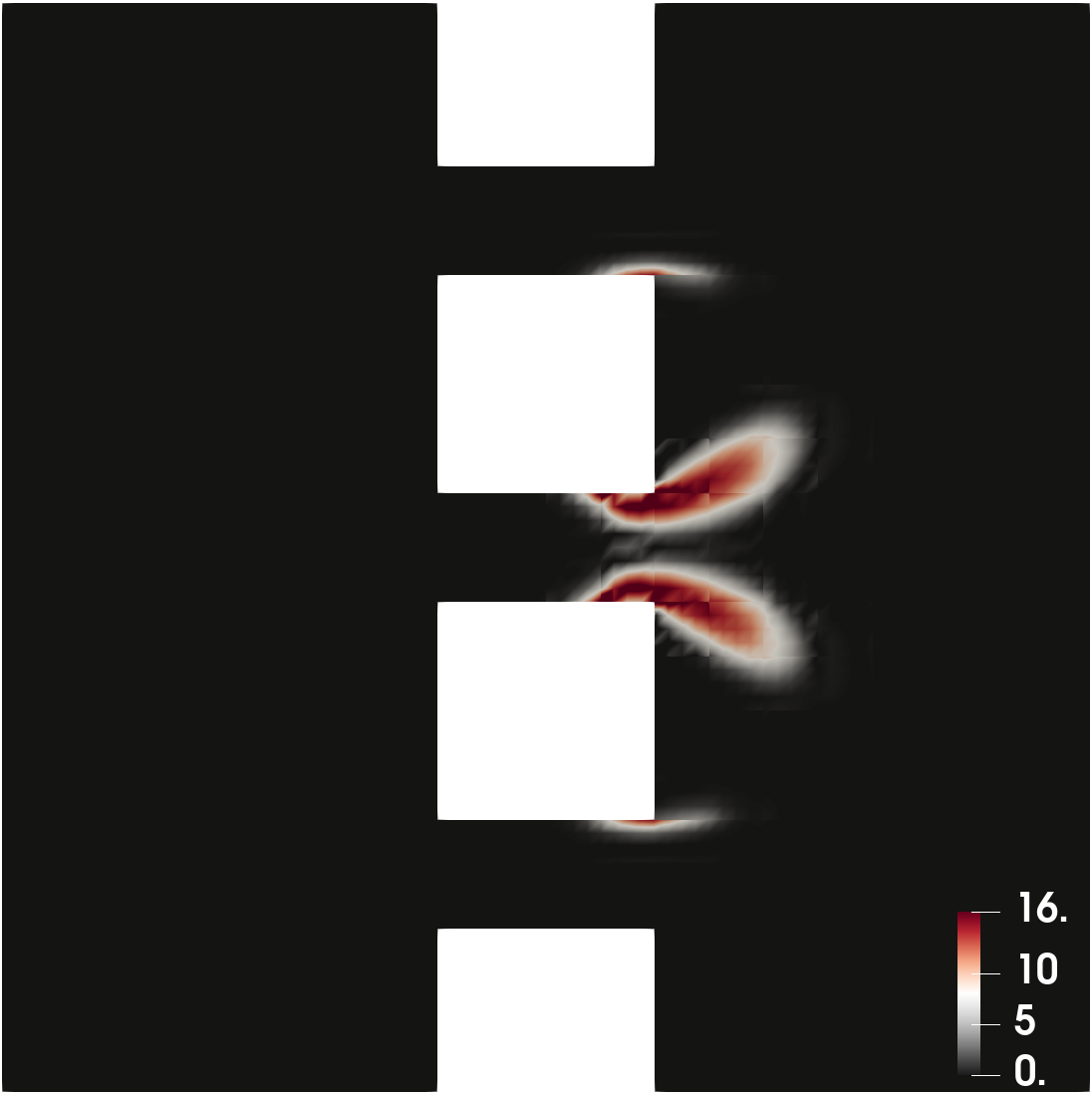}
\includegraphics[width=0.192\textwidth]{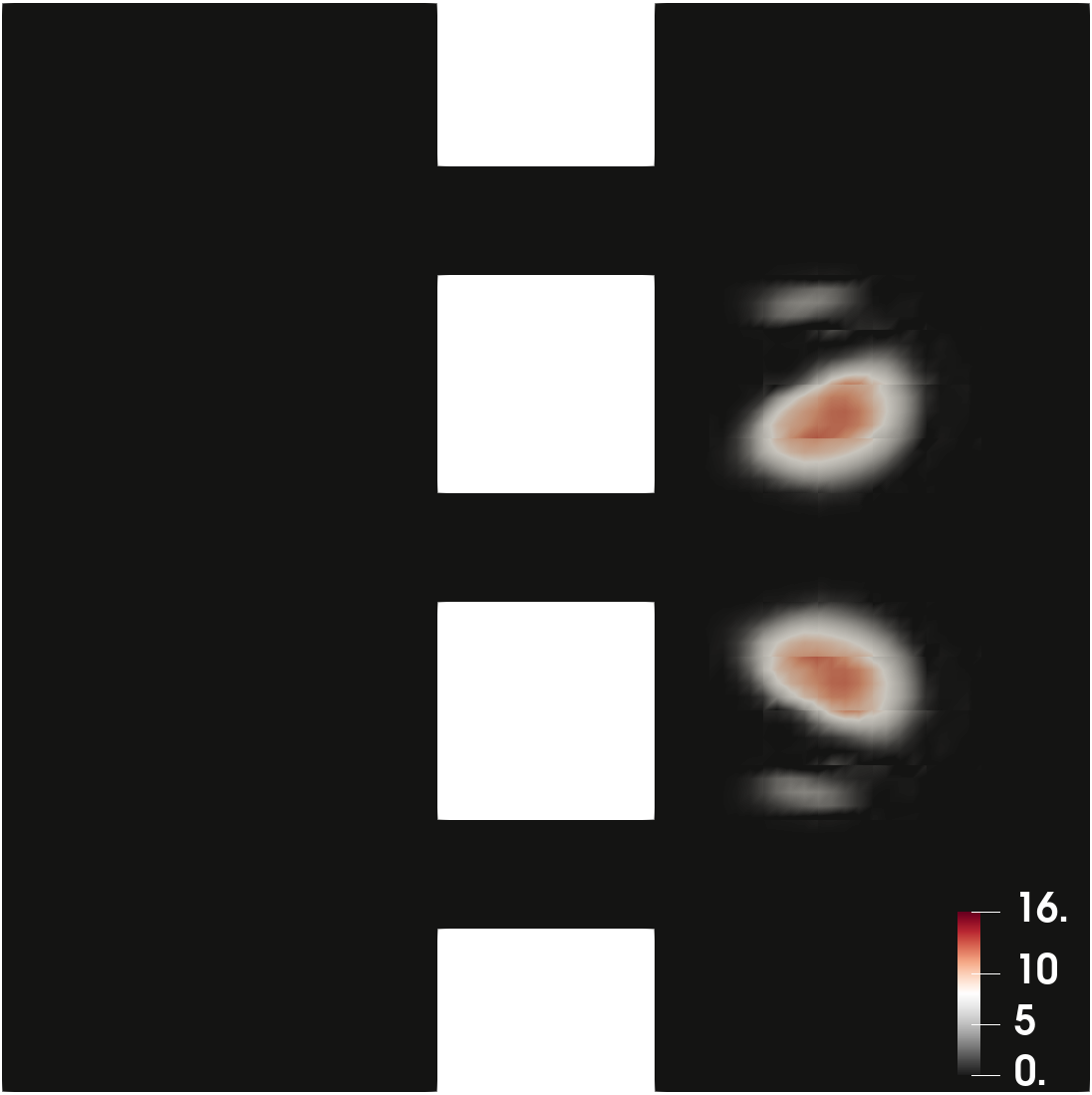}
}
\caption{Example \ref{ex3}. Snapshots of $\rho$ at 
$t=0.1,0.3,0.5,0.7,0.9$ 
(left to right).}
\label{fig:den-case3}
\end{figure}

Snapshots of the density contour at different times
are shown in Figure~\ref{fig:den-case3}.
The results are similar to those in Example \ref{ex2}, where different interaction cost function leads to different  density evolution.

\subsection{MFP between mascot images}
\label{ex4}
\begin{figure}[tb]
\centering
\subfigure[ND (Leprechaun)]
{
\includegraphics[width=0.31\textwidth]{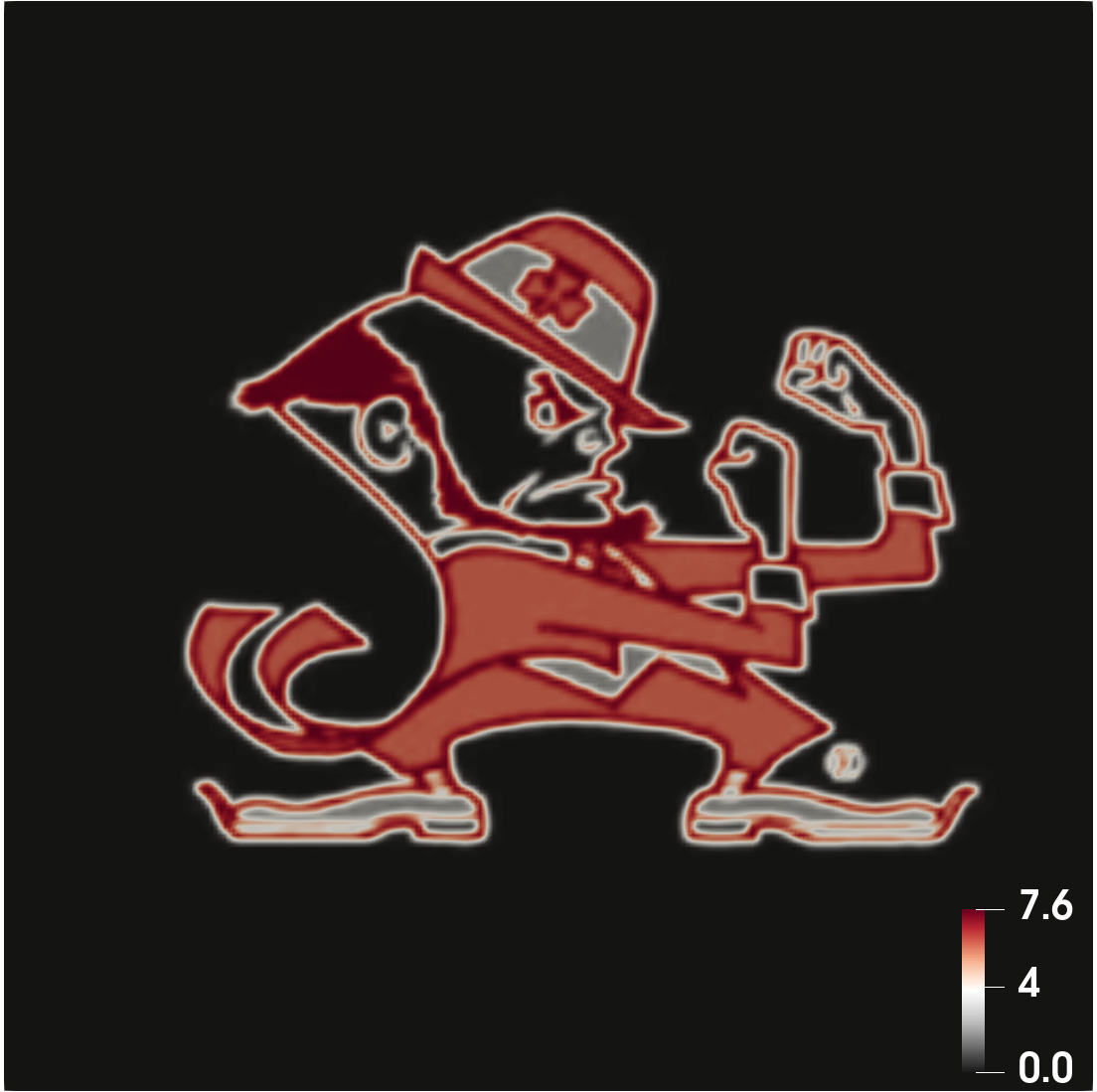}
}
\subfigure[UCLA (Bruins)  ]{\includegraphics[width=0.31\textwidth]{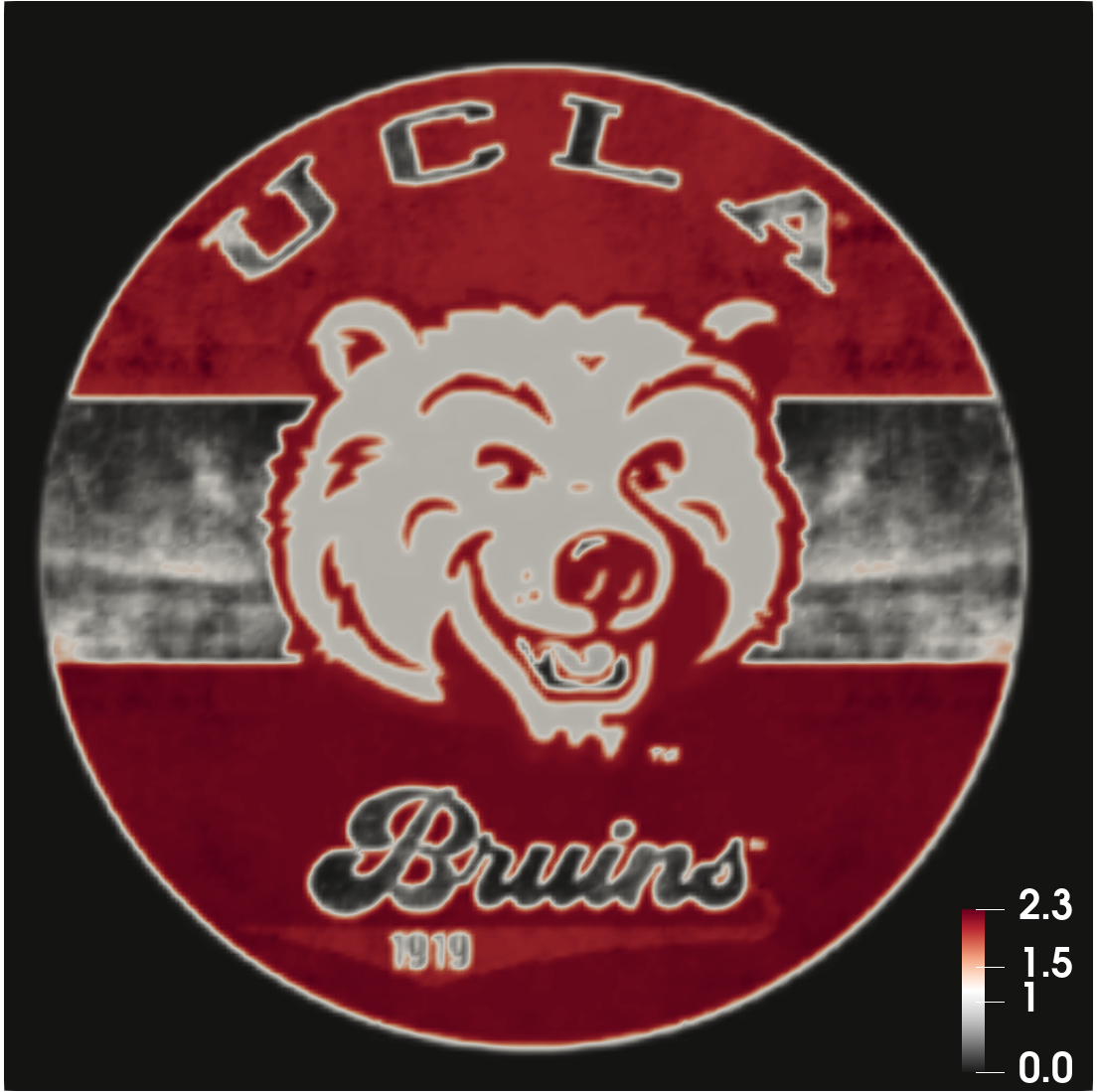}}
\subfigure[USC (Gamecocks)]{
\includegraphics[width=0.31\textwidth]{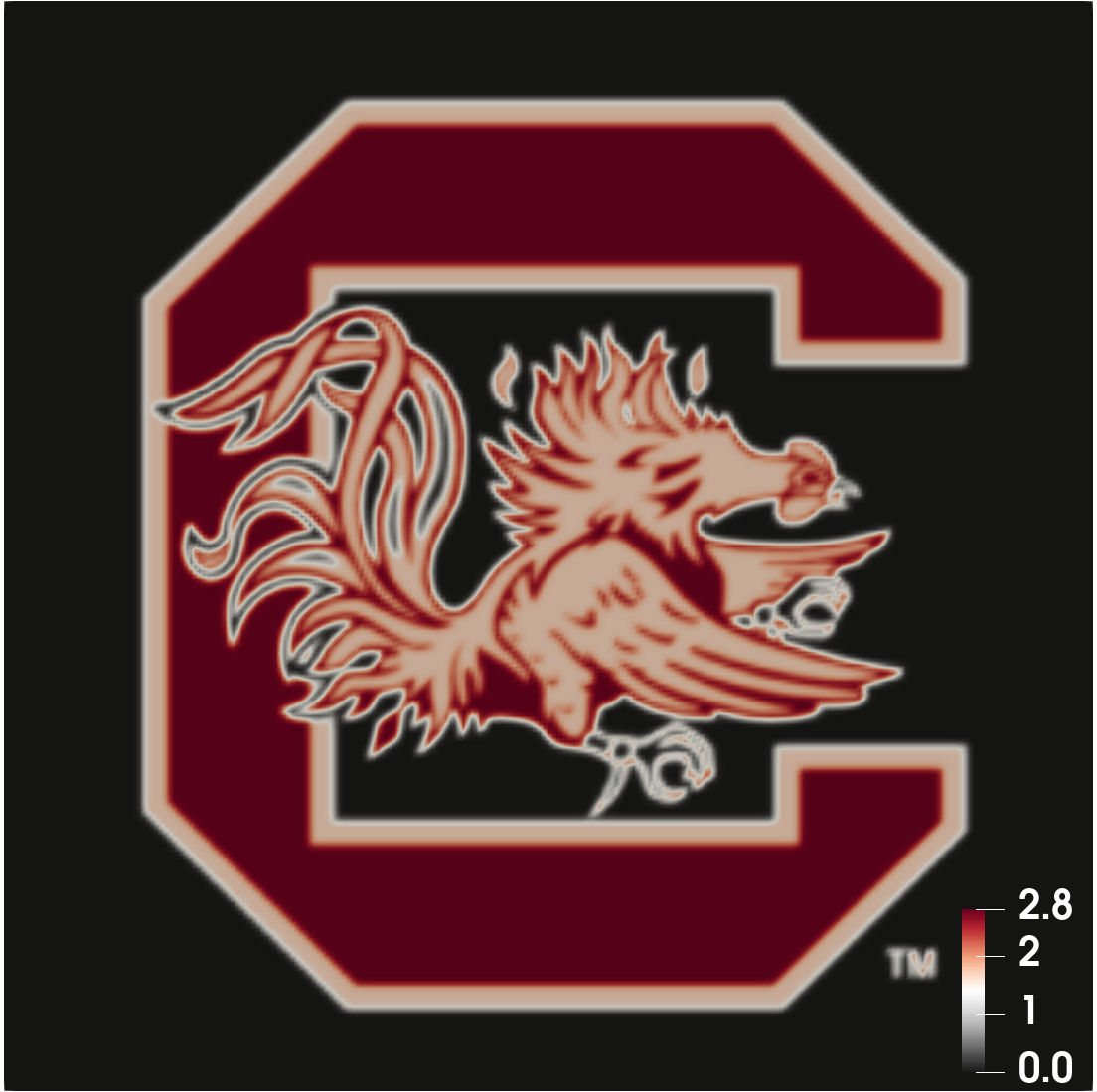}
}
\caption{Example \ref{ex4}. Initial/final densities. 
}
\label{fig:img}
\end{figure}

Our last example concerns with OT and MFP \eqref{mfp:problem} between images. 
The initial or terminal densities are 
normalized images of athletics mascots
from University of Notre Dame (Leprechaun), 
UCLA (Brunins), and University of South Carilina (Gamecocks); see Figure~\ref{fig:img}.
The spatial domain is a unit square $\Omega=[0,1]\times [0,1]$, and the initial/terminal densities are normalized to have unit mass.

We apply the scheme \eqref{aug-mfp-h} with polynomial degree $k=3$ on a structured hexahedral mesh of size $64\times 64\times 16$, where the time step size is 
$\Delta t = 1/16$. 
Three set of initial/terminal density pairs are considered: (i) ND$\rightarrow$ UCLA where initial density is the ND image and terminal density is the UCLA image, (ii) UCLA$\rightarrow$ USC where initial density is the UCLA image and terminal density is the USC image,
and (iii) USC$\rightarrow$ ND where initial density is the USC image and terminal density is the ND image. 
For each pair of data, we consider three choices of interaction cost, namely,  Case 1: $A(\rho) = 0$ (OT), Case 2: $A(\rho) = 0.01\rho\log(\rho)$, and Case 3: $A(\rho) = 0.01/\rho$.
The ALG2 algorithm is terminated when $err_m^a$ is less than 0.001. 
The number of iterations needed for convergence are recorded in Table~\ref{tab:iter-case4}.
\begin{table}[ht!]
\centering
\begin{tabular}{lccccccc}
 & Case 1 & Case 2 & Case 3\\\hline
ND$\rightarrow$UCLA & 2440&471&892\\ 
UCLA$\rightarrow$USC & 1511&211&244\\
USC$\rightarrow$ND & 3577&496&907\\
\end{tabular}
\caption{Example \ref{ex4}. Number of ALG2 iterations for each case.
}
\label{tab:iter-case4}
\end{table}

\begin{figure}[tb]
\centering
\subfigure[Case 1: $A(\rho) = 0$. ND $\rightarrow$ UCLA]
{
\includegraphics[width=0.192\textwidth]{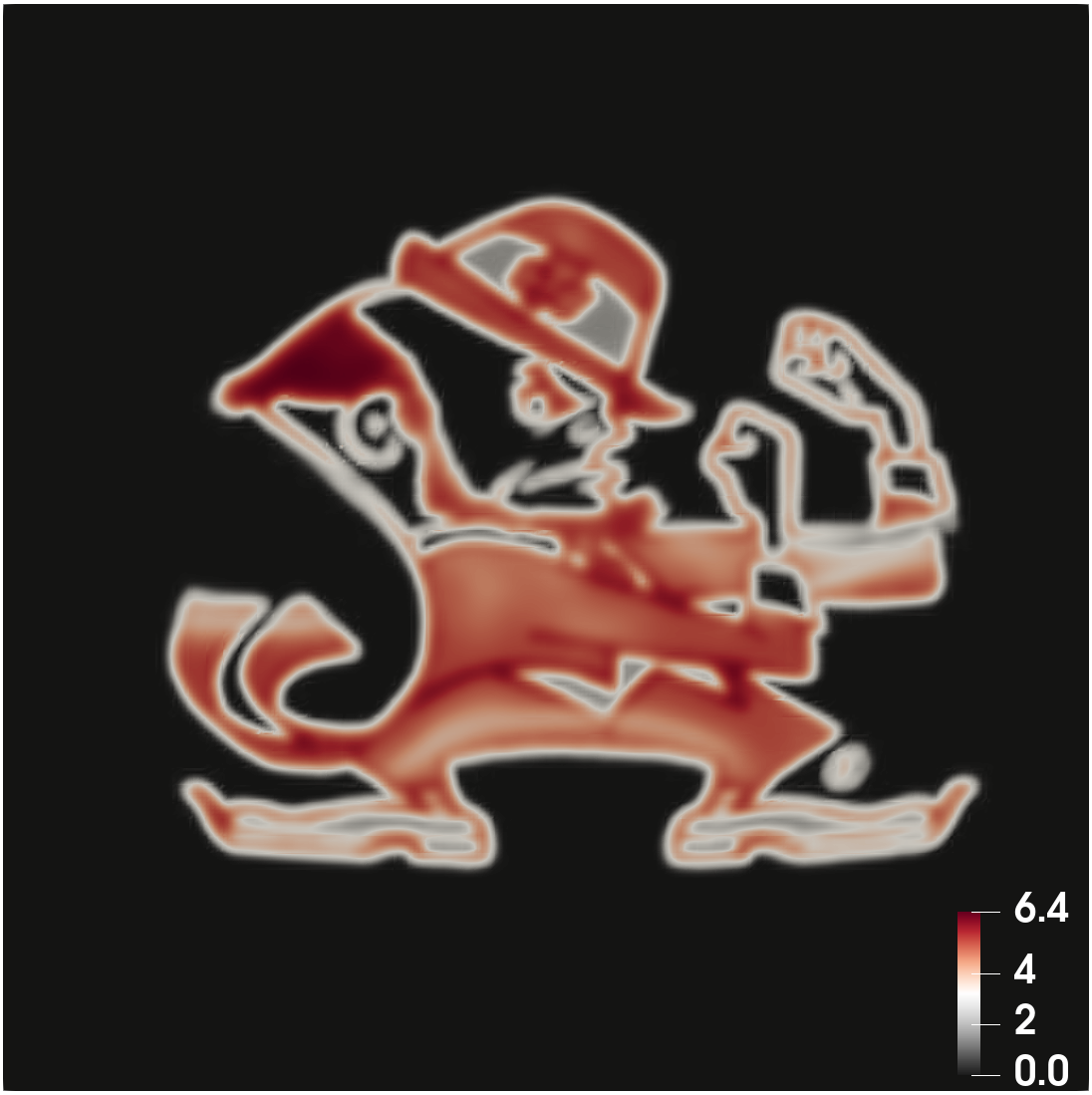}
\includegraphics[width=0.192\textwidth]{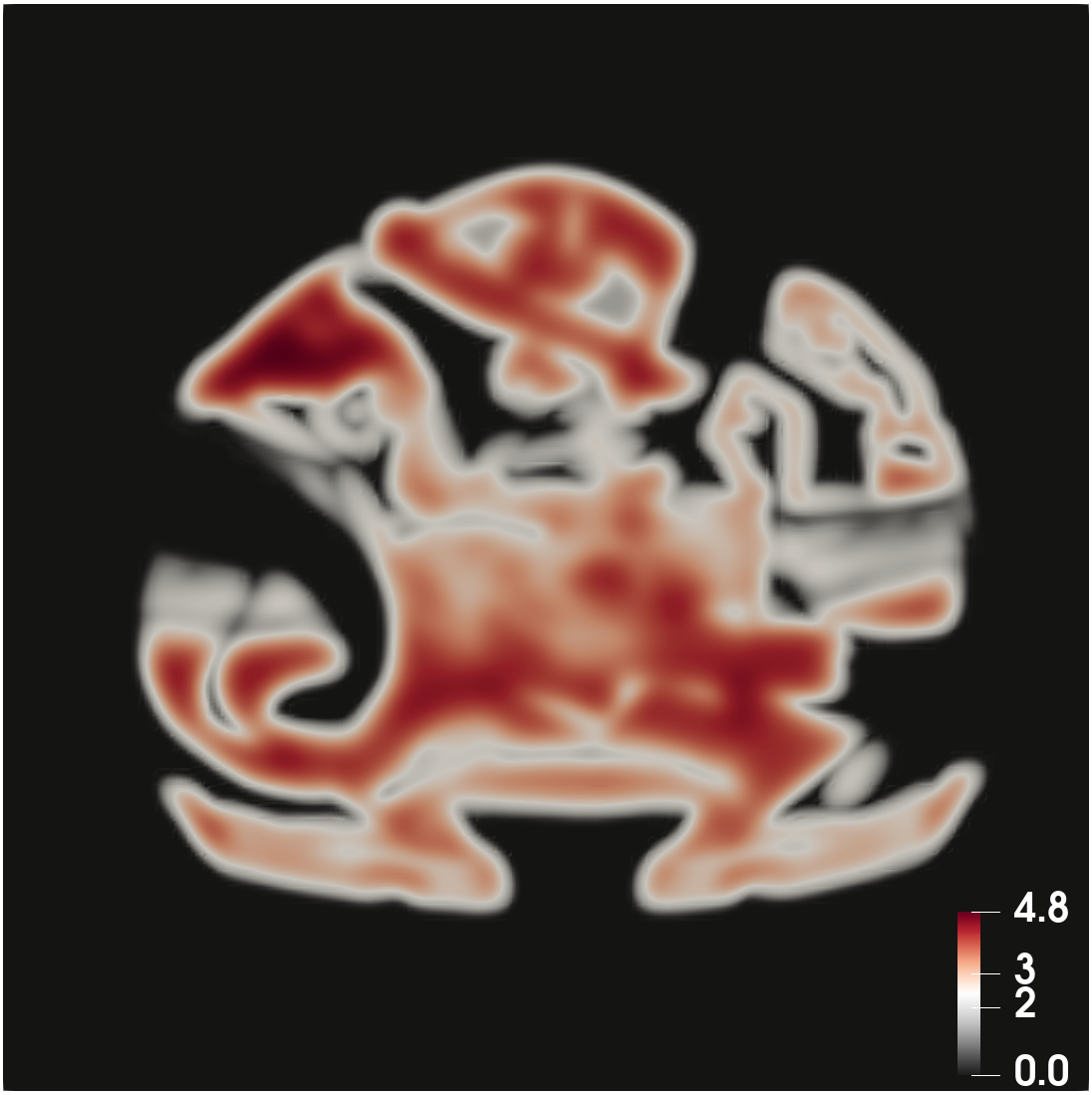}
\includegraphics[width=0.192\textwidth]{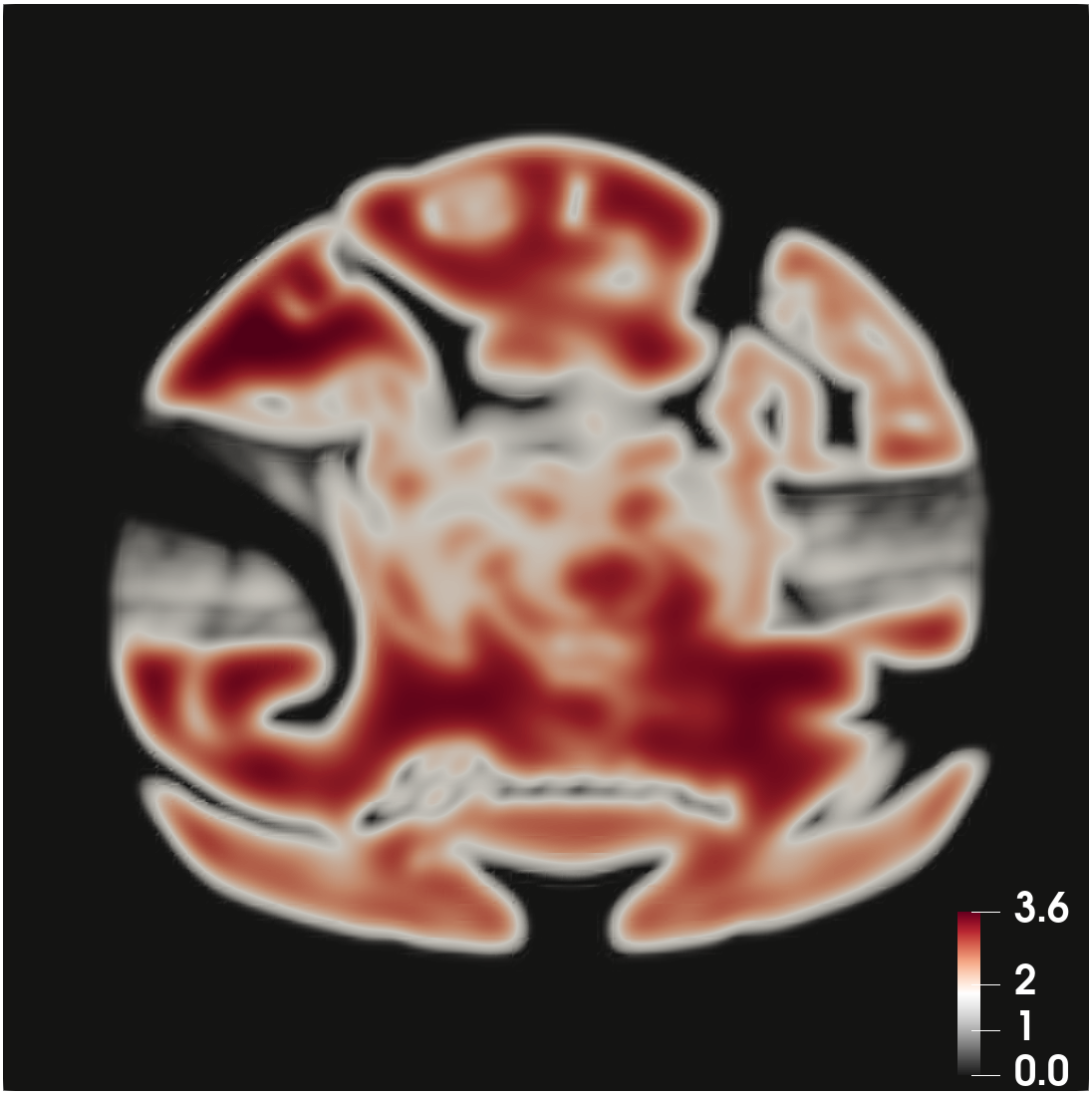}
\includegraphics[width=0.192\textwidth]{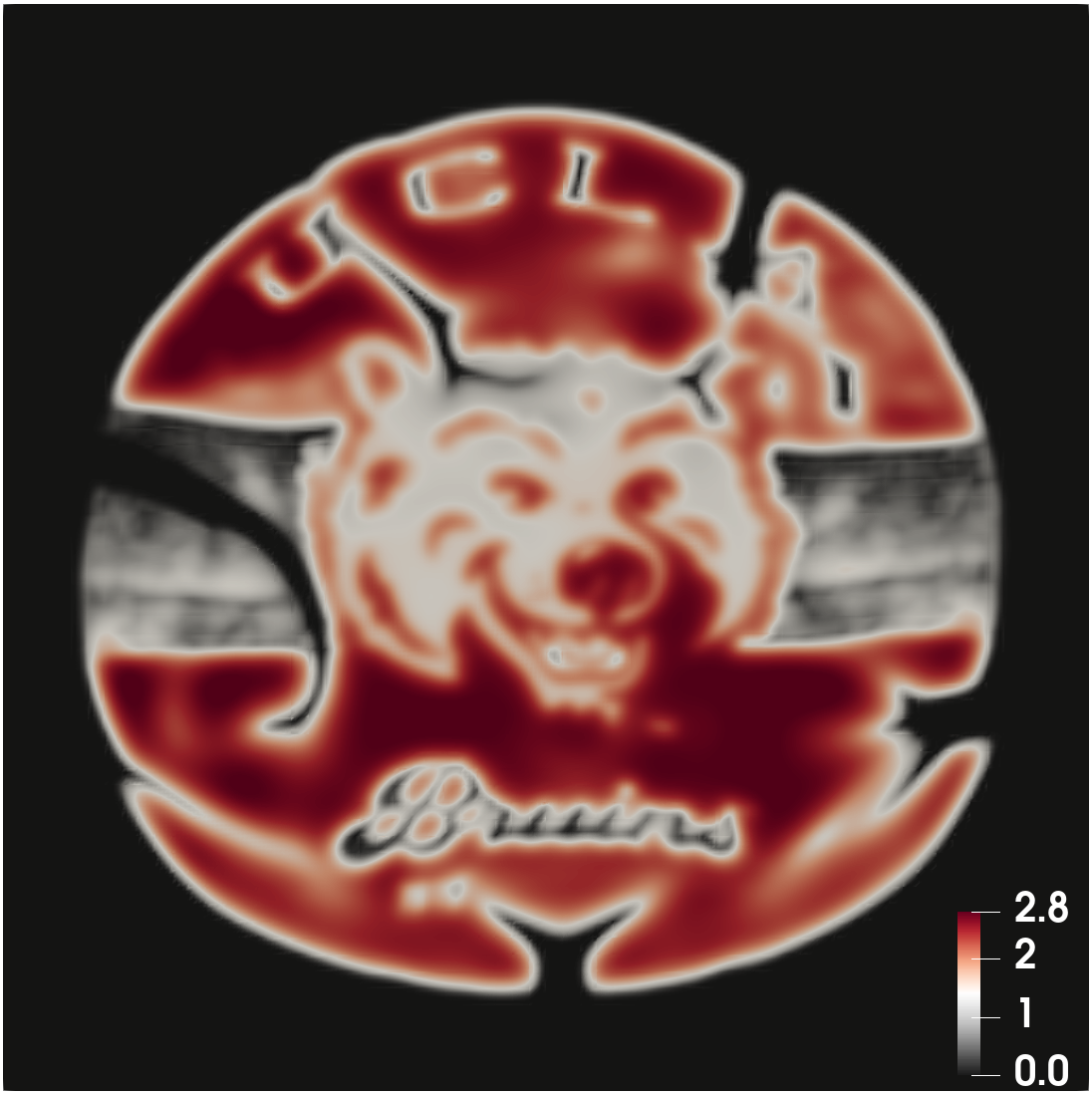}
\includegraphics[width=0.192\textwidth]{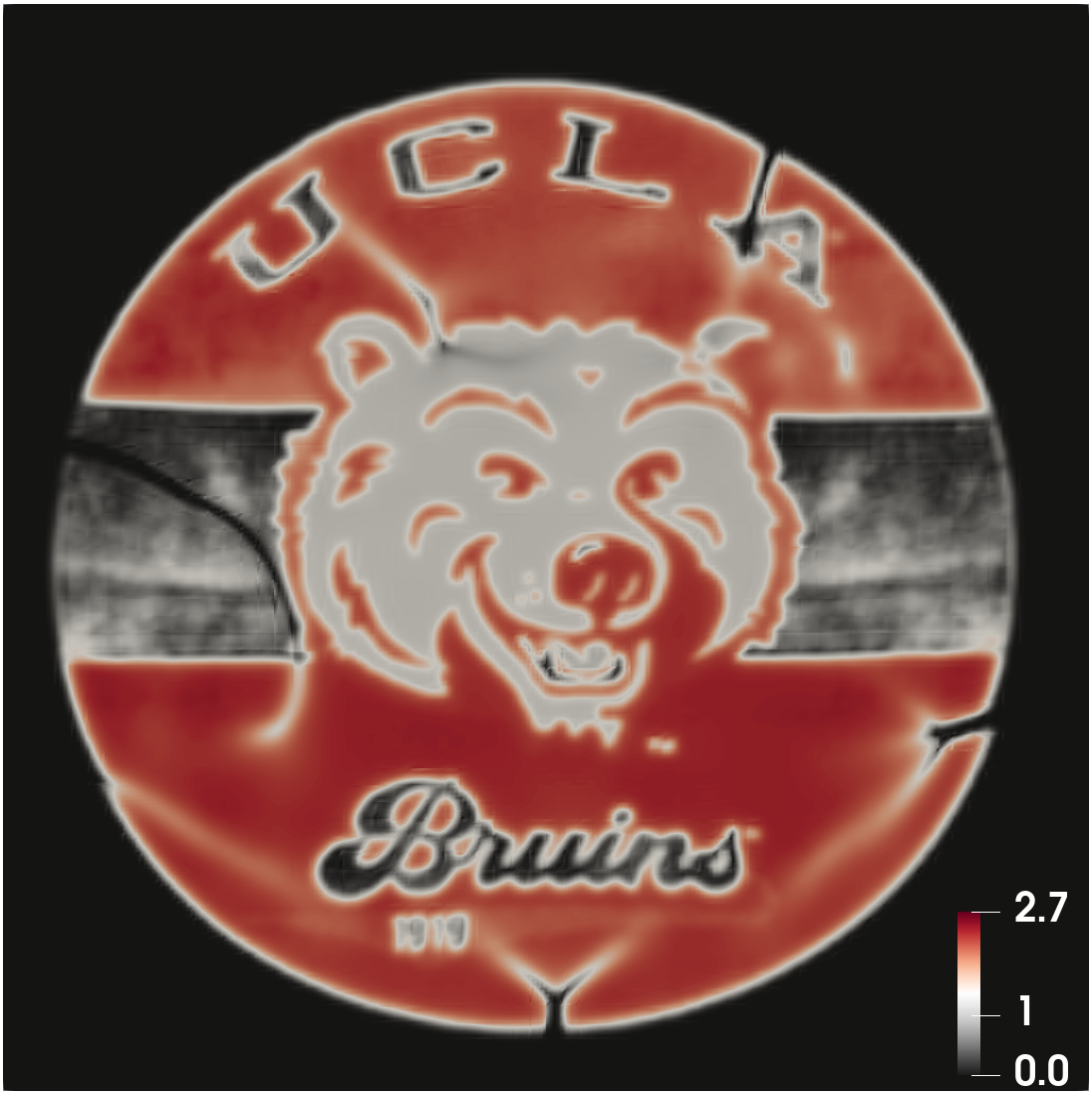}
}
\subfigure[Case 2: $A(\rho) = 0.01\rho\log(\rho)$.  ND $\rightarrow$ UCLA]
{
\includegraphics[width=0.192\textwidth]{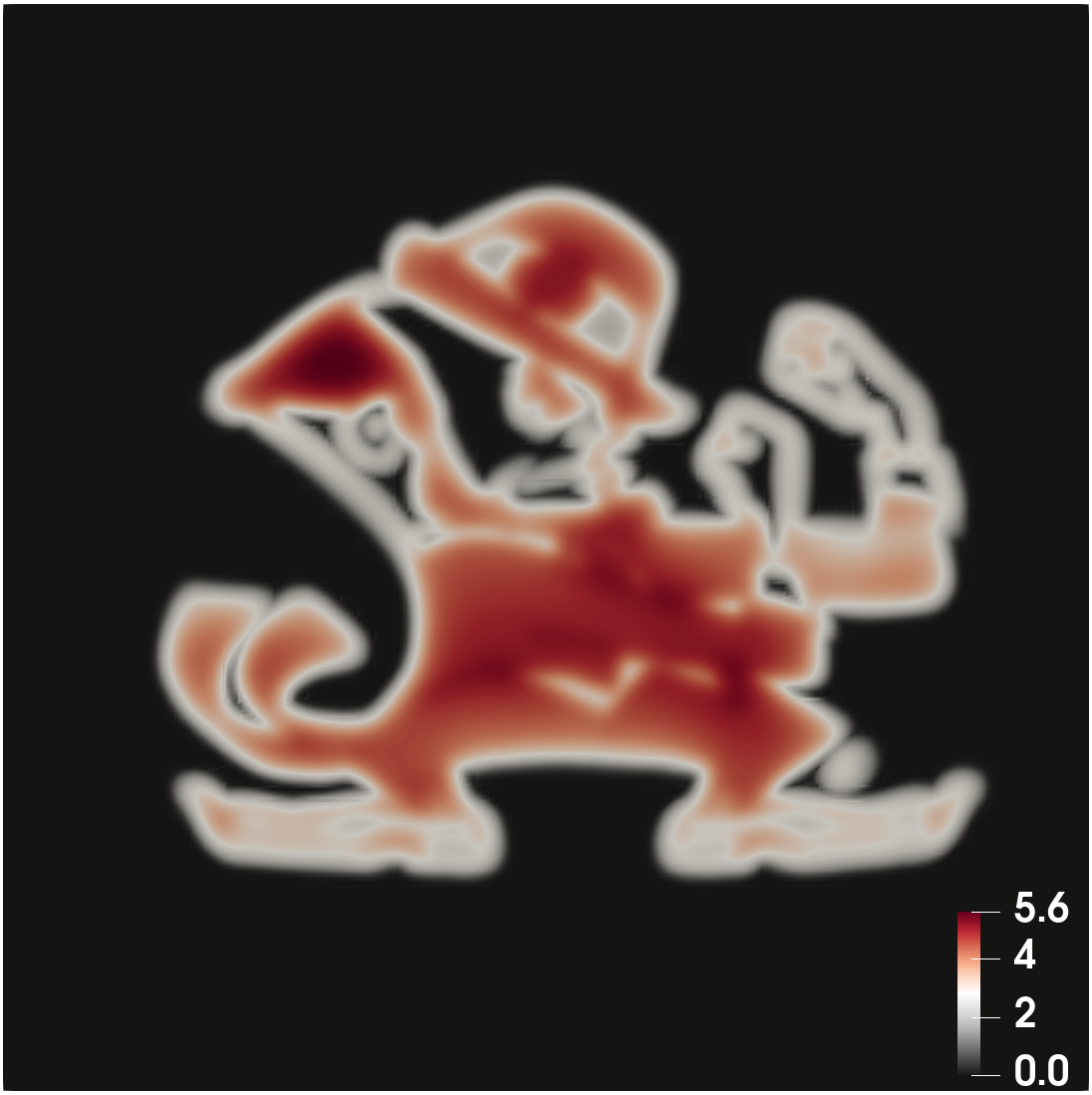}
\includegraphics[width=0.192\textwidth]{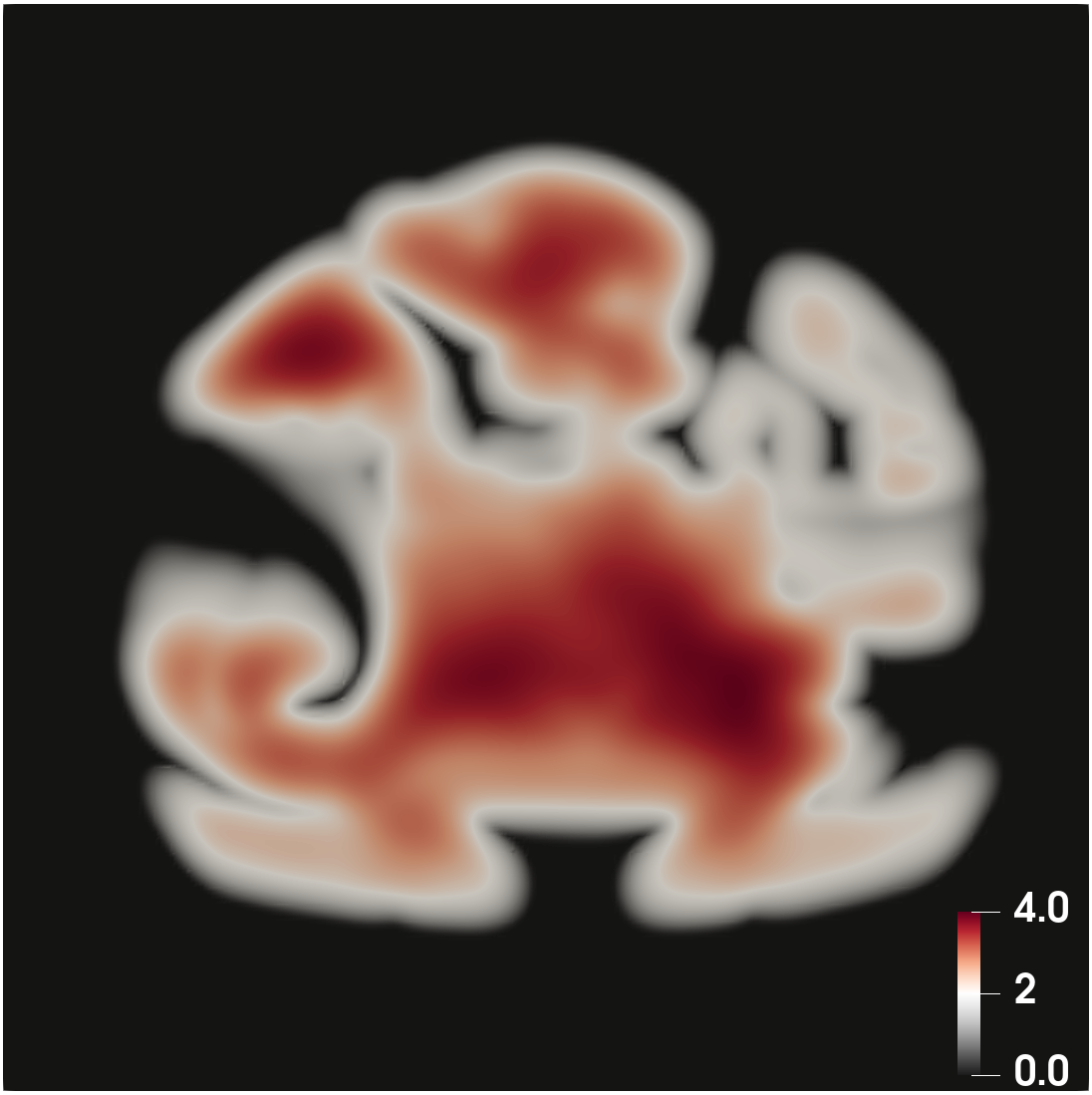}
\includegraphics[width=0.192\textwidth]{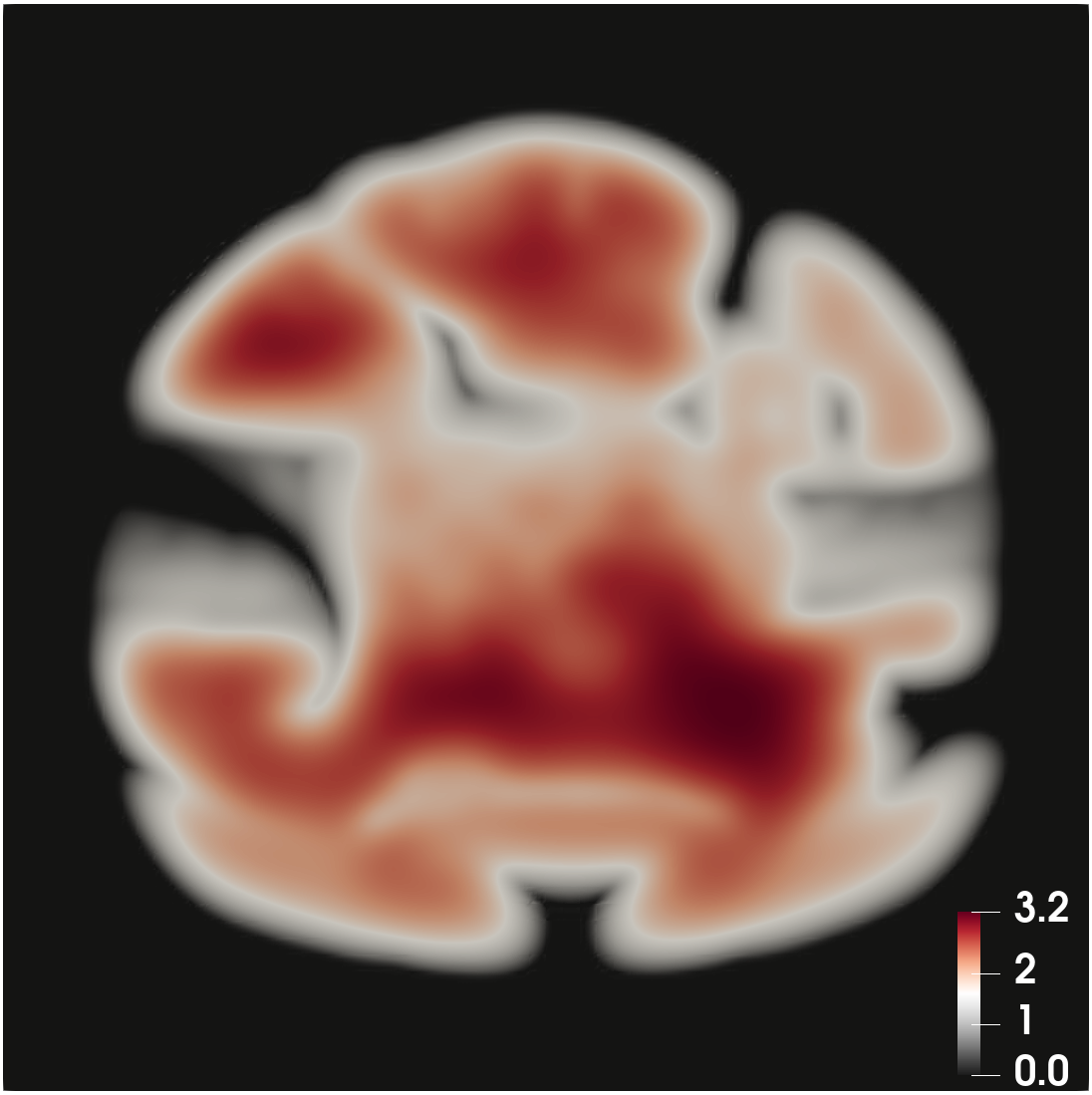}
\includegraphics[width=0.192\textwidth]{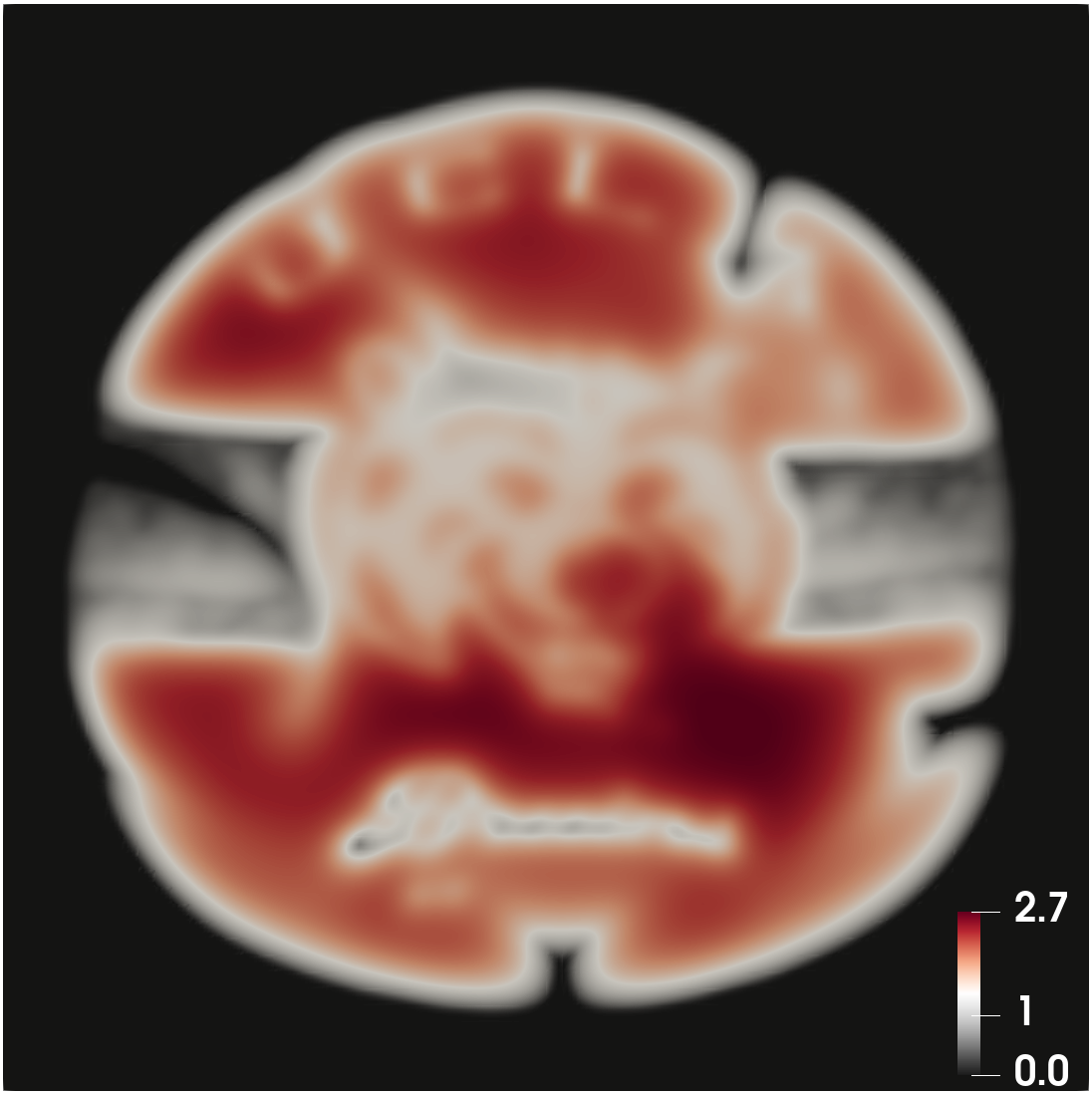}
\includegraphics[width=0.192\textwidth]{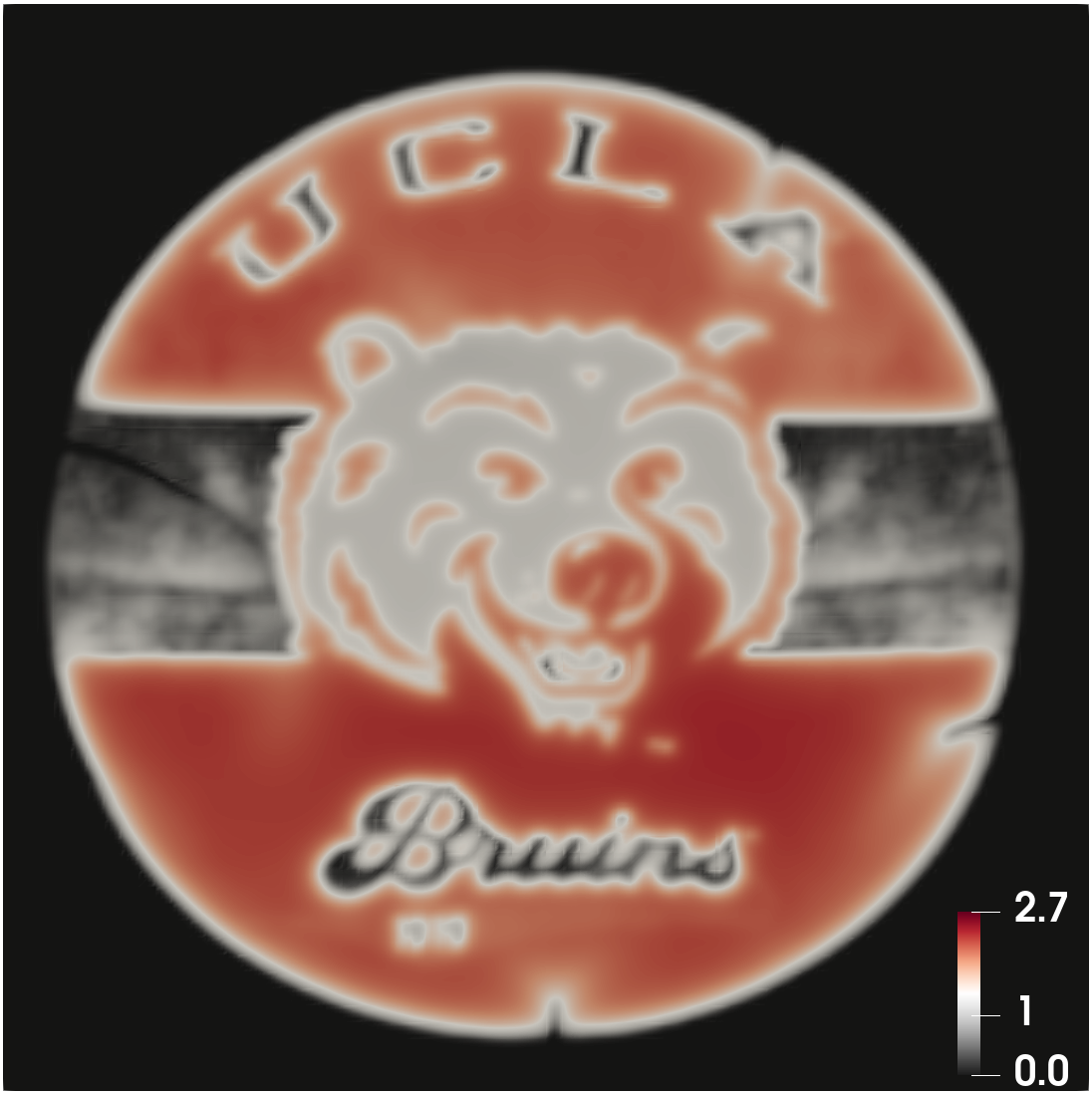}
}

\subfigure[Case 3: $A(\rho) = 0.01/\rho$. ND $\rightarrow$ UCLA]
{
\includegraphics[width=0.192\textwidth]{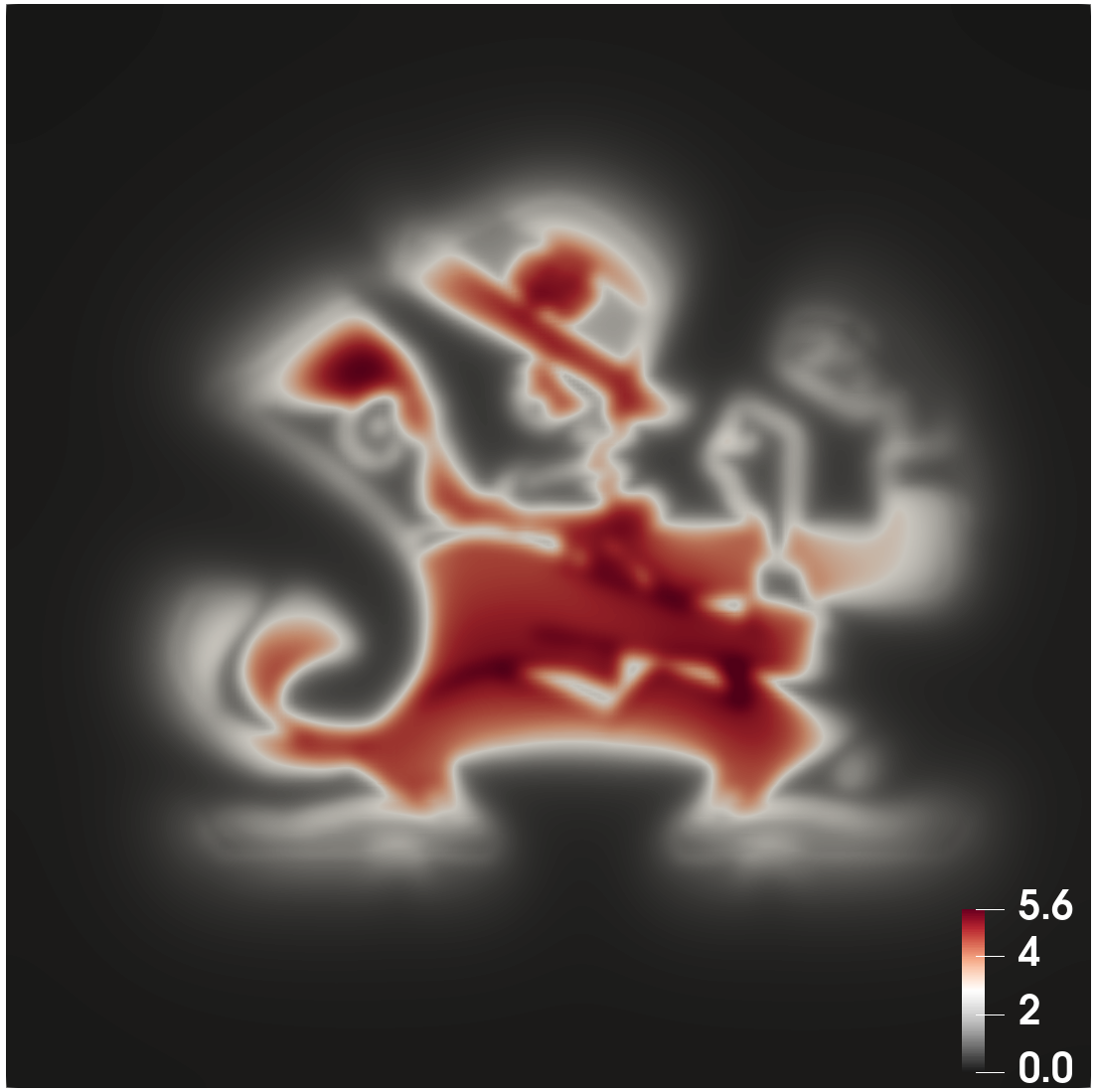}
\includegraphics[width=0.192\textwidth]{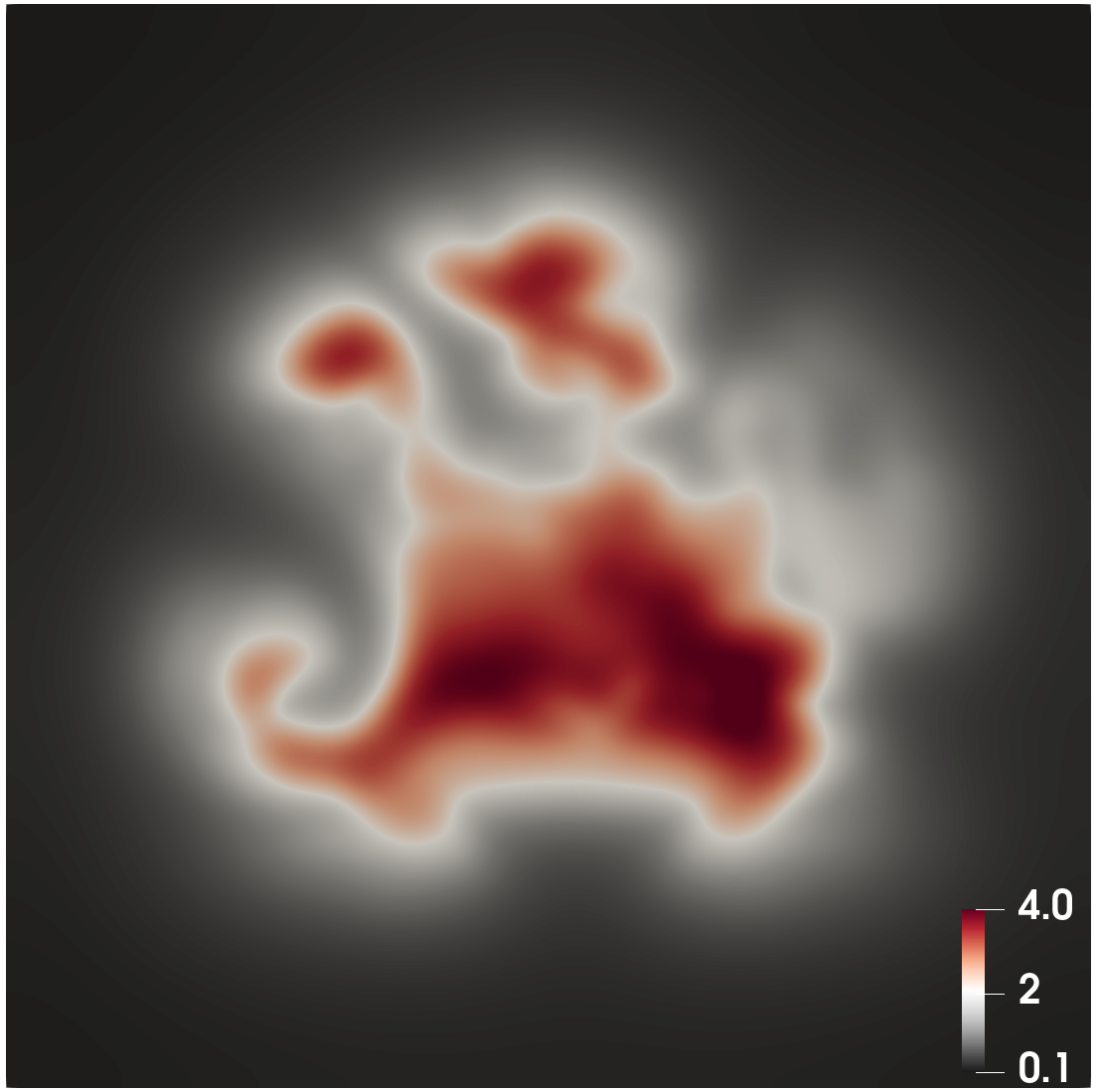}
\includegraphics[width=0.192\textwidth]{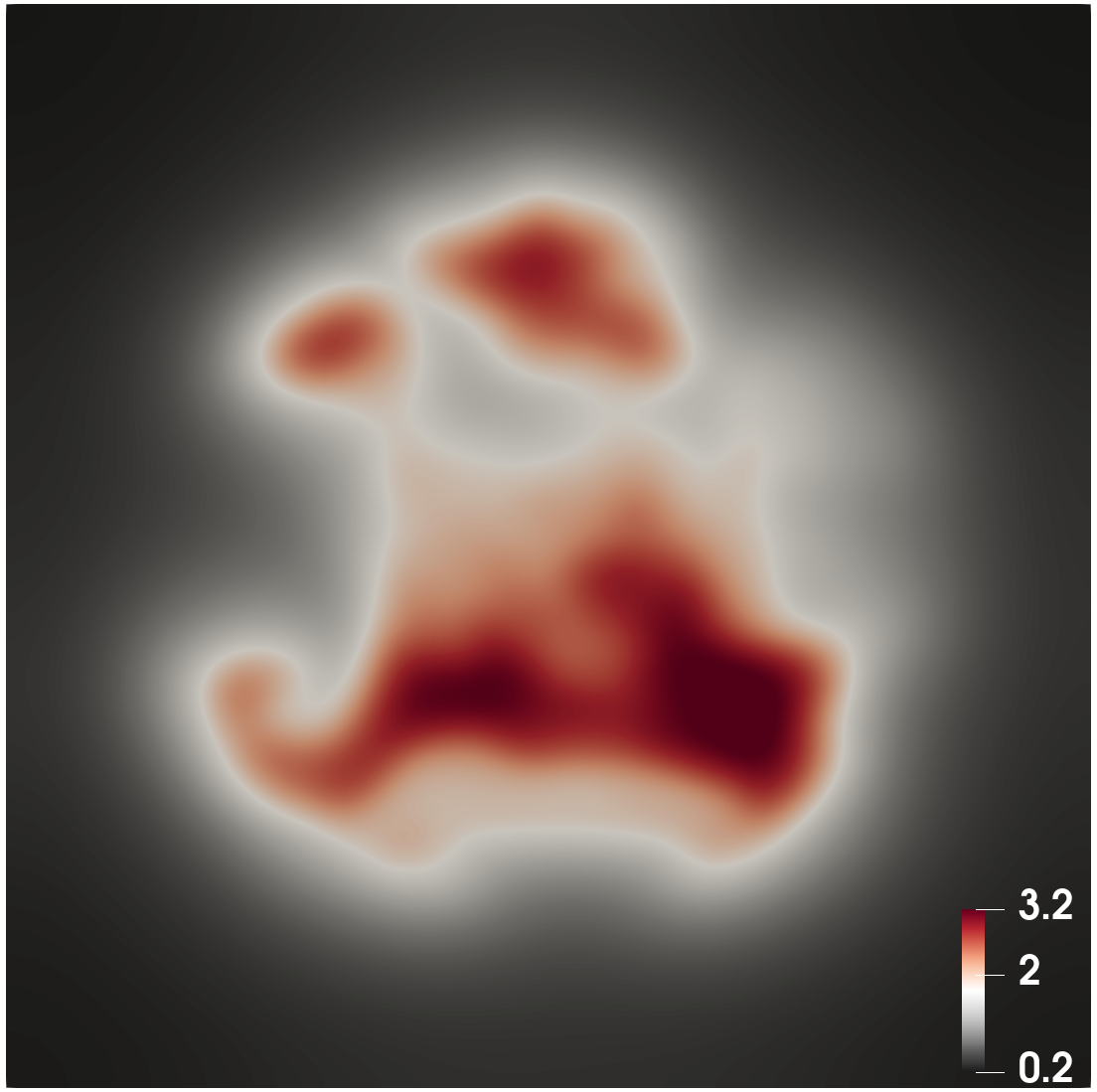}
\includegraphics[width=0.192\textwidth]{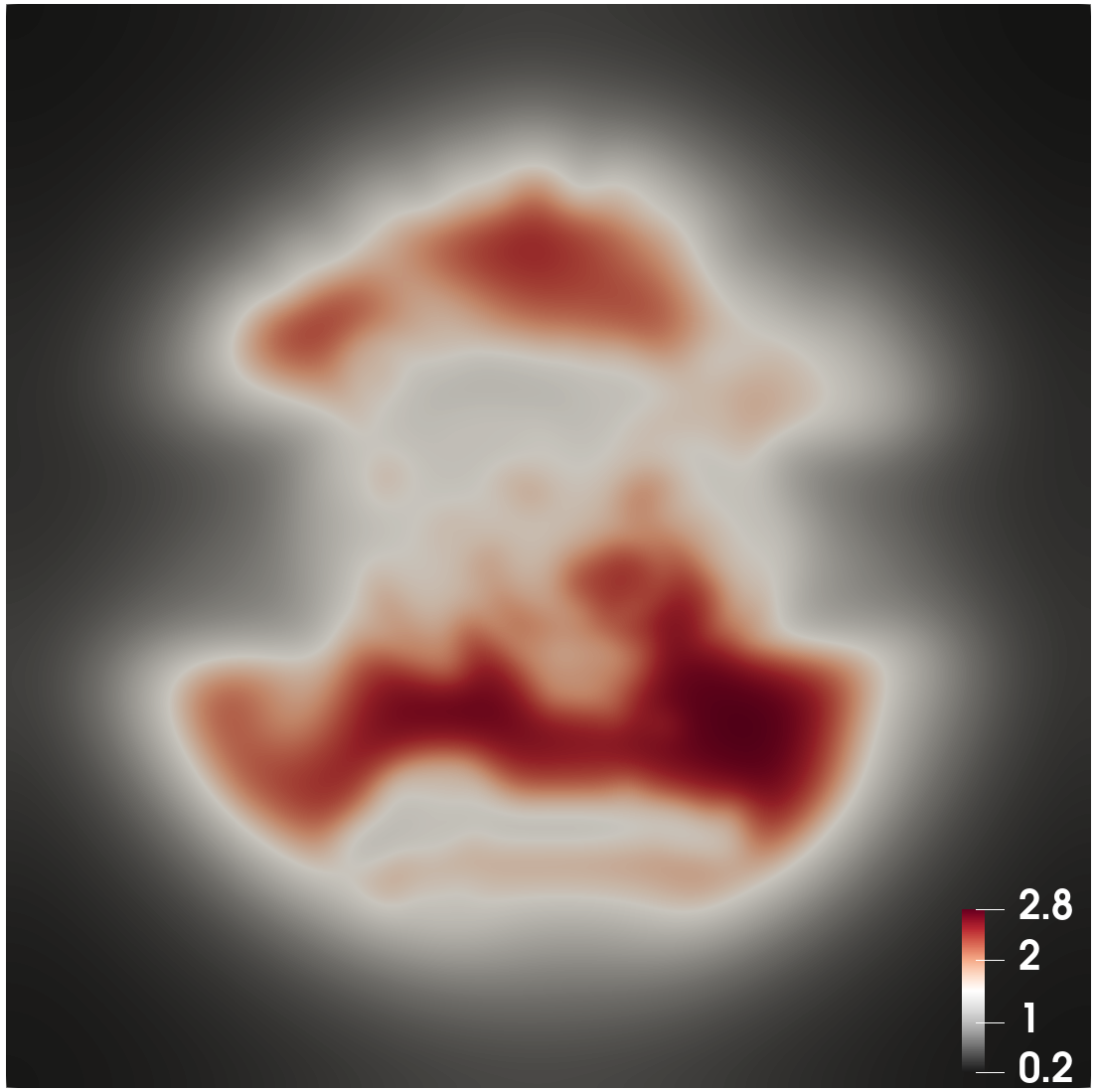}
\includegraphics[width=0.192\textwidth]{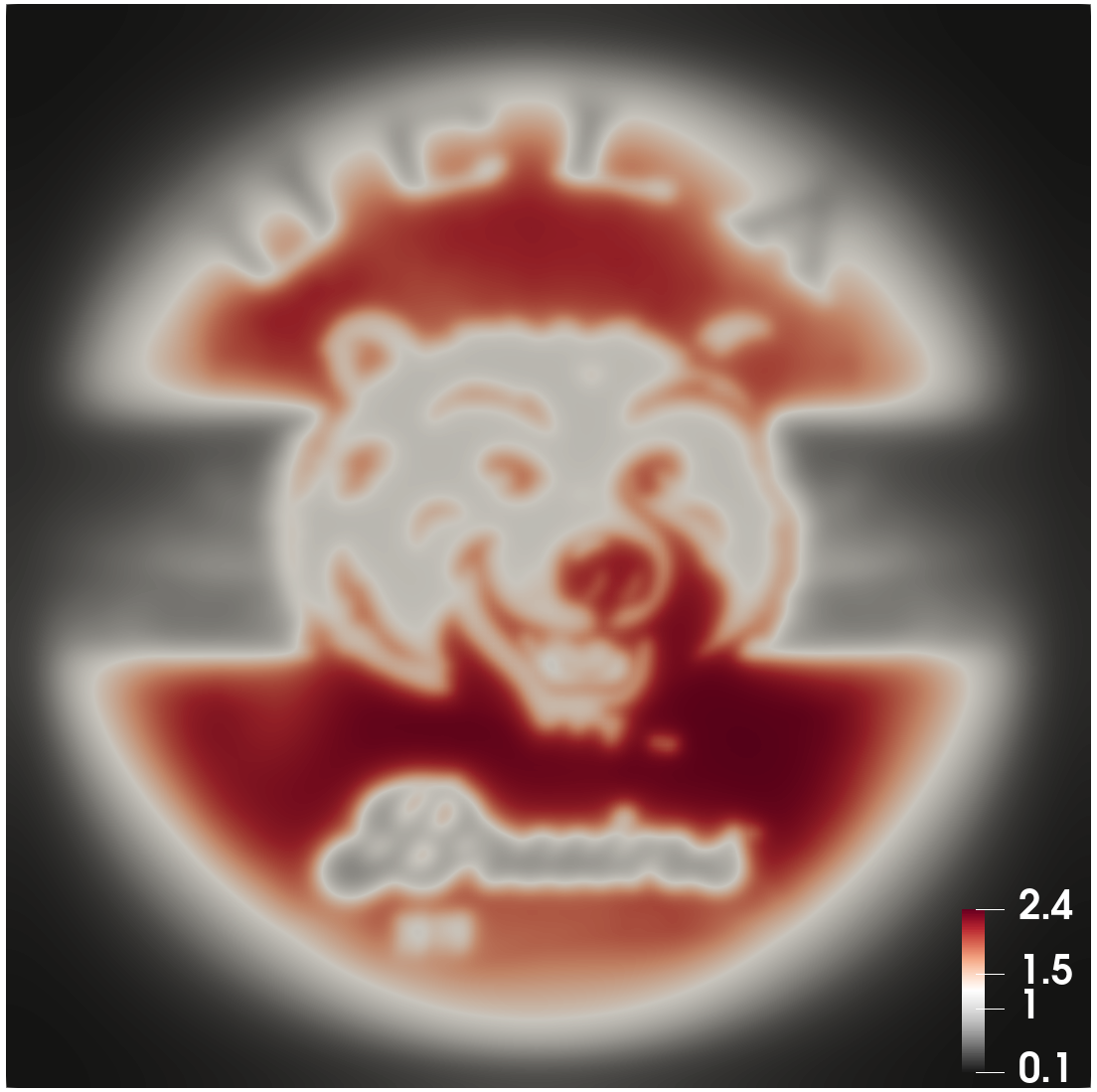}
}

\caption{Example \ref{ex4}. Initial density: ND. Terminal density: UCLA.
Snapshots of $\rho$ at 
$t=$ 0.1,0.3,0.5,0.7,0.9 (left to right).
}
\label{fig:den-case4A}
\end{figure}

\begin{figure}[tb]
\centering
\subfigure[Case 1: $A(\rho) = 0$. UCLA $\rightarrow$ USC]
{
\includegraphics[width=0.192\textwidth]{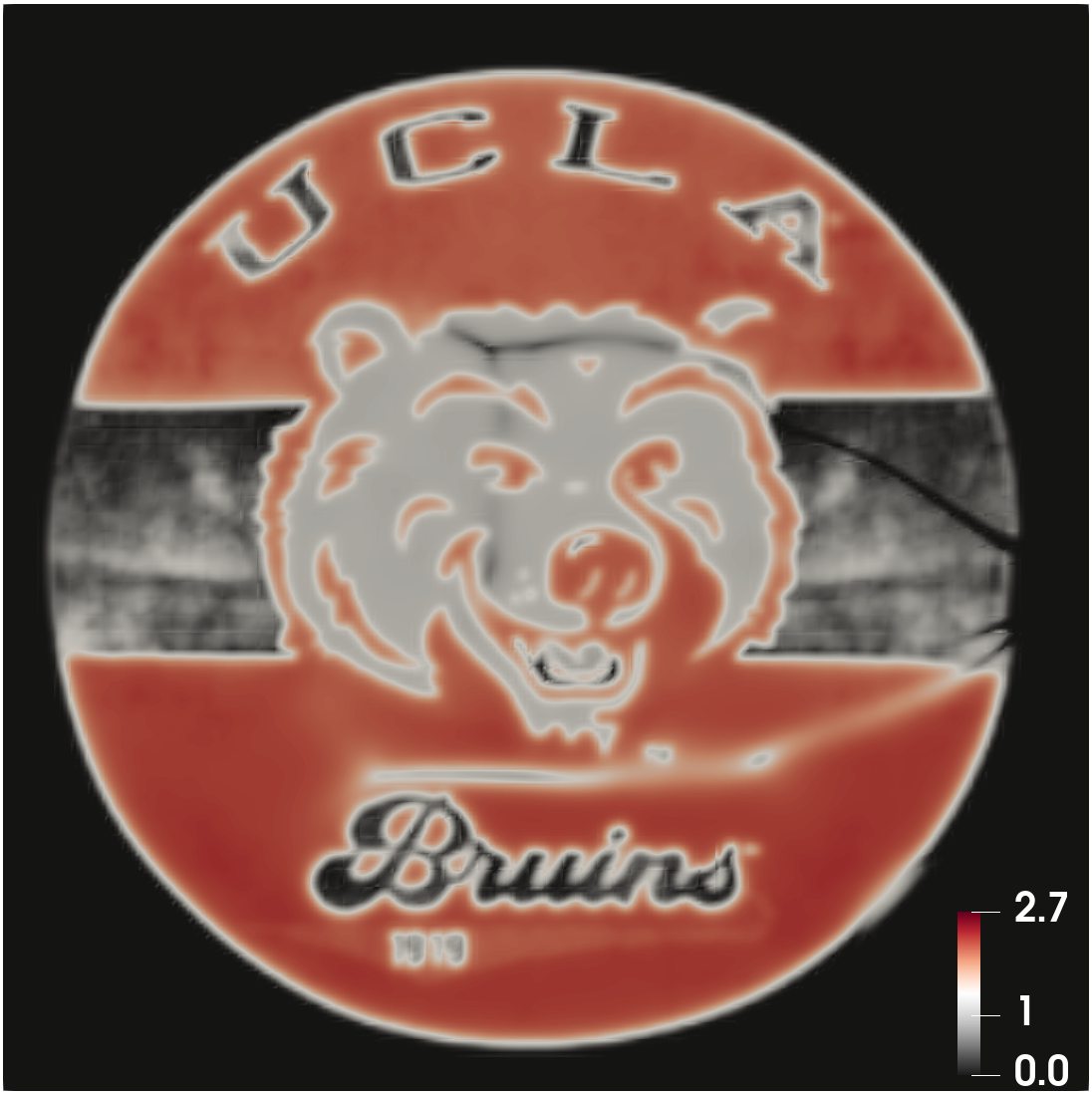}
\includegraphics[width=0.192\textwidth]{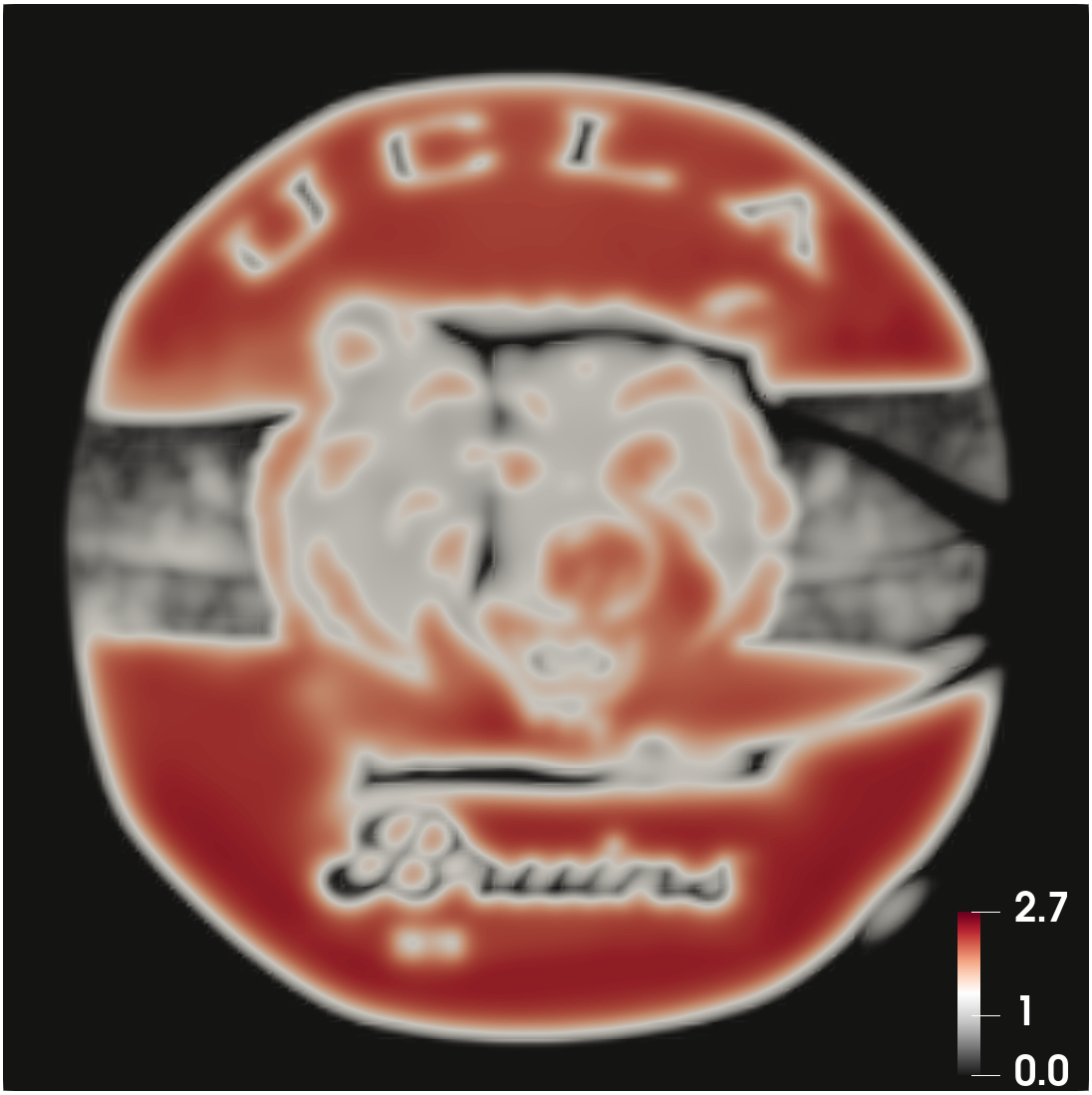}
\includegraphics[width=0.192\textwidth]{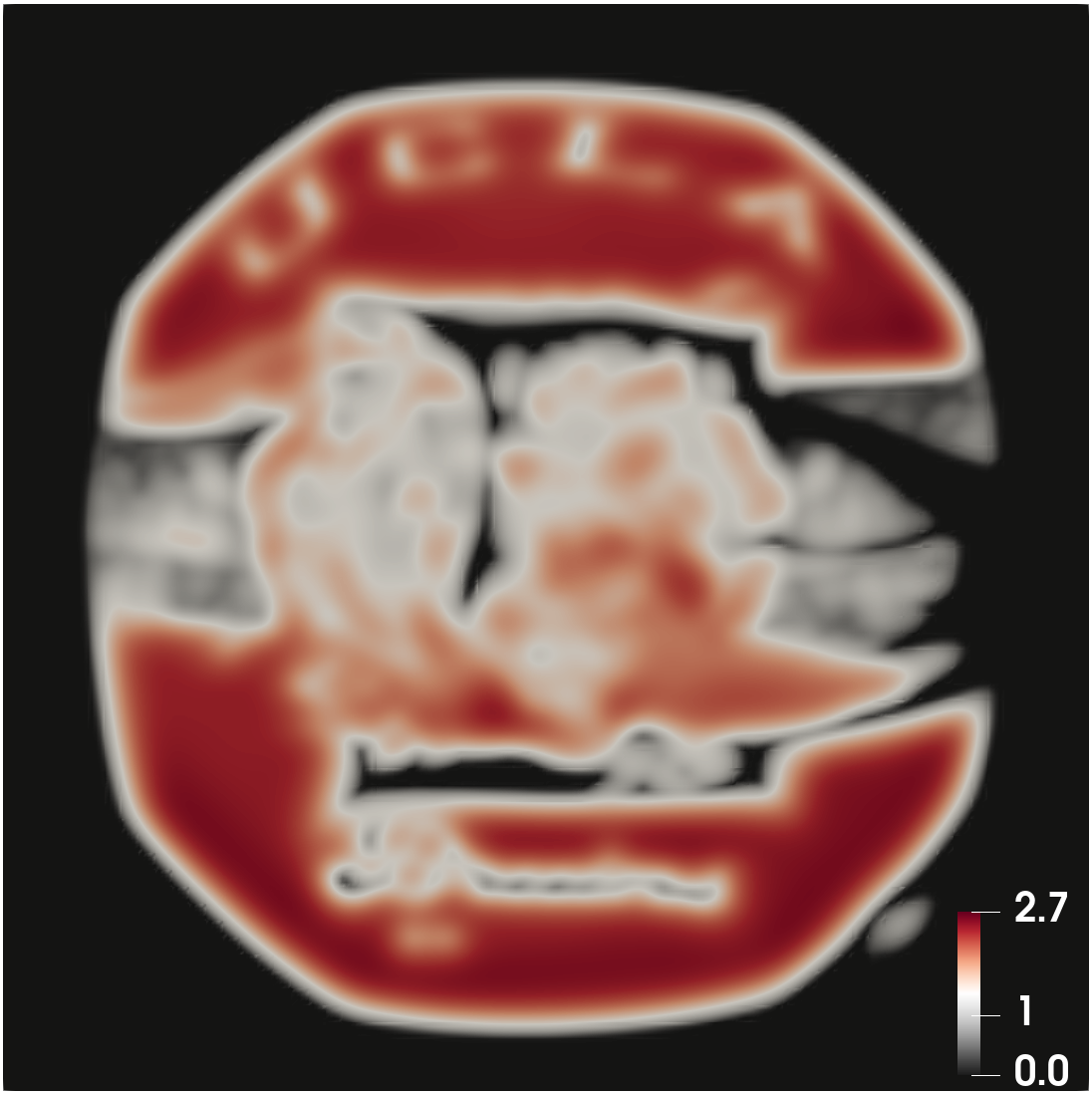}
\includegraphics[width=0.192\textwidth]{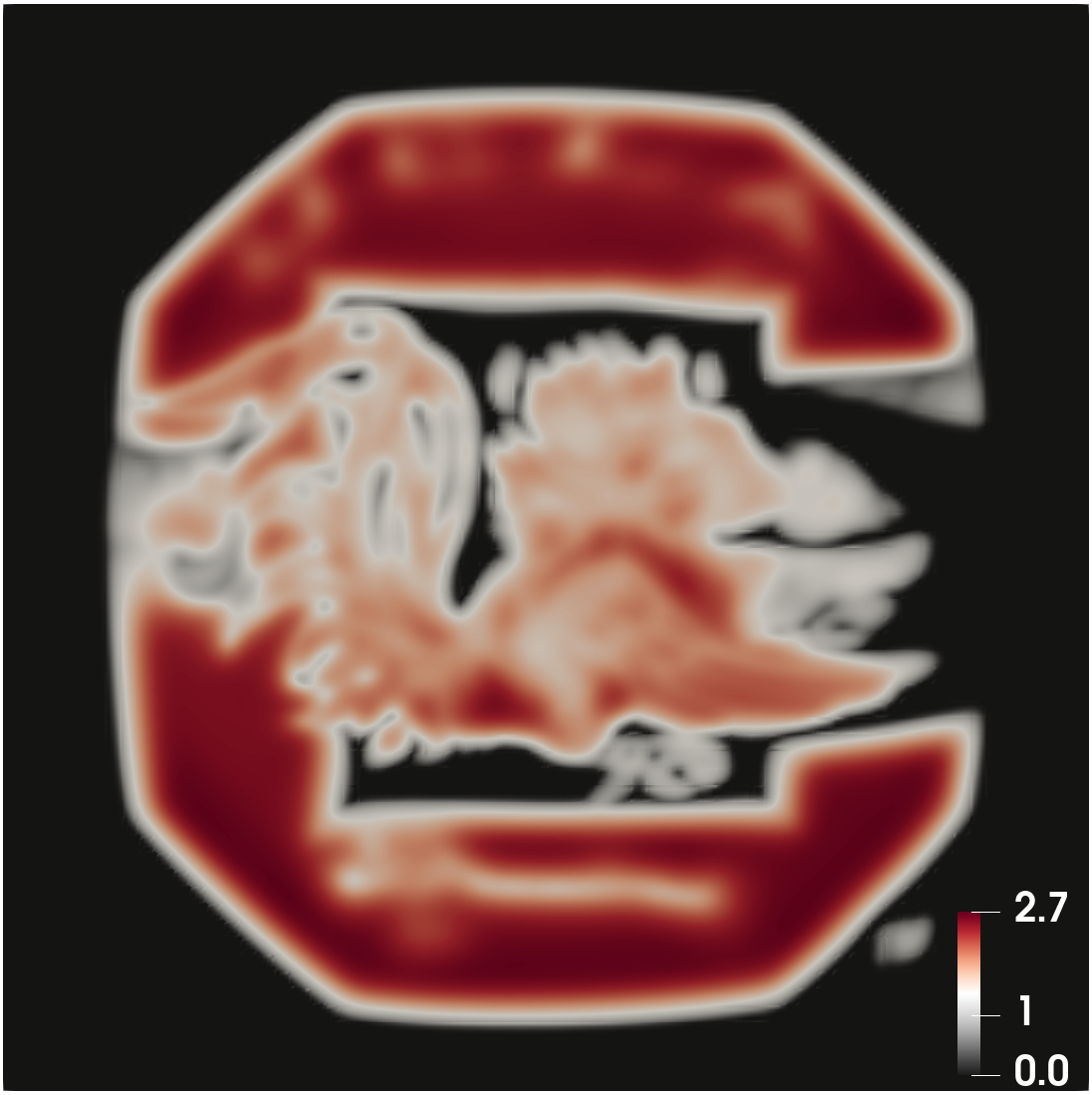}
\includegraphics[width=0.192\textwidth]{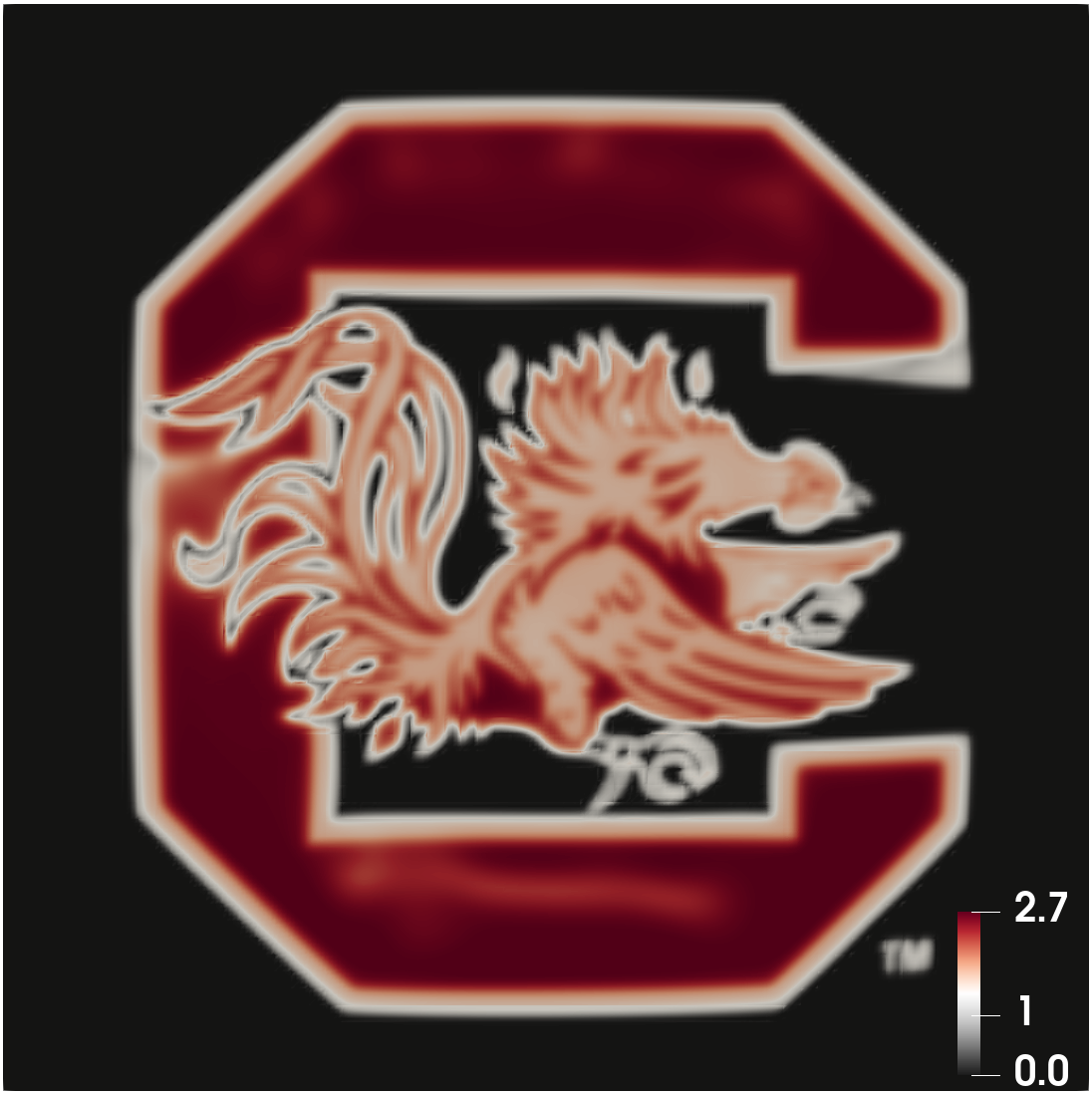}
}
\subfigure[Case 2: $A(\rho) = 0.01\rho\log(\rho)$.  UCLA $\rightarrow$ USC]
{
\includegraphics[width=0.192\textwidth]{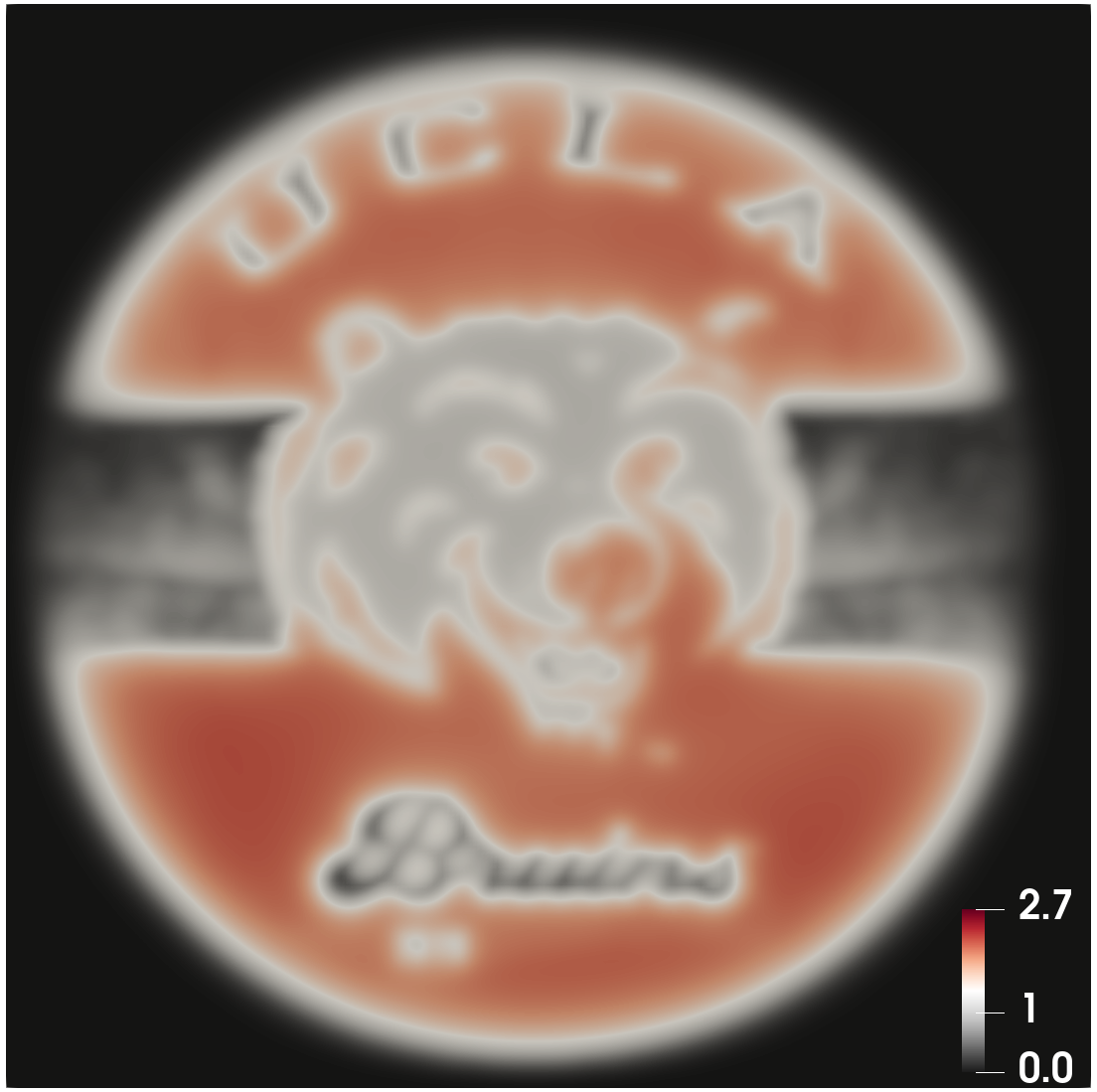}
\includegraphics[width=0.192\textwidth]{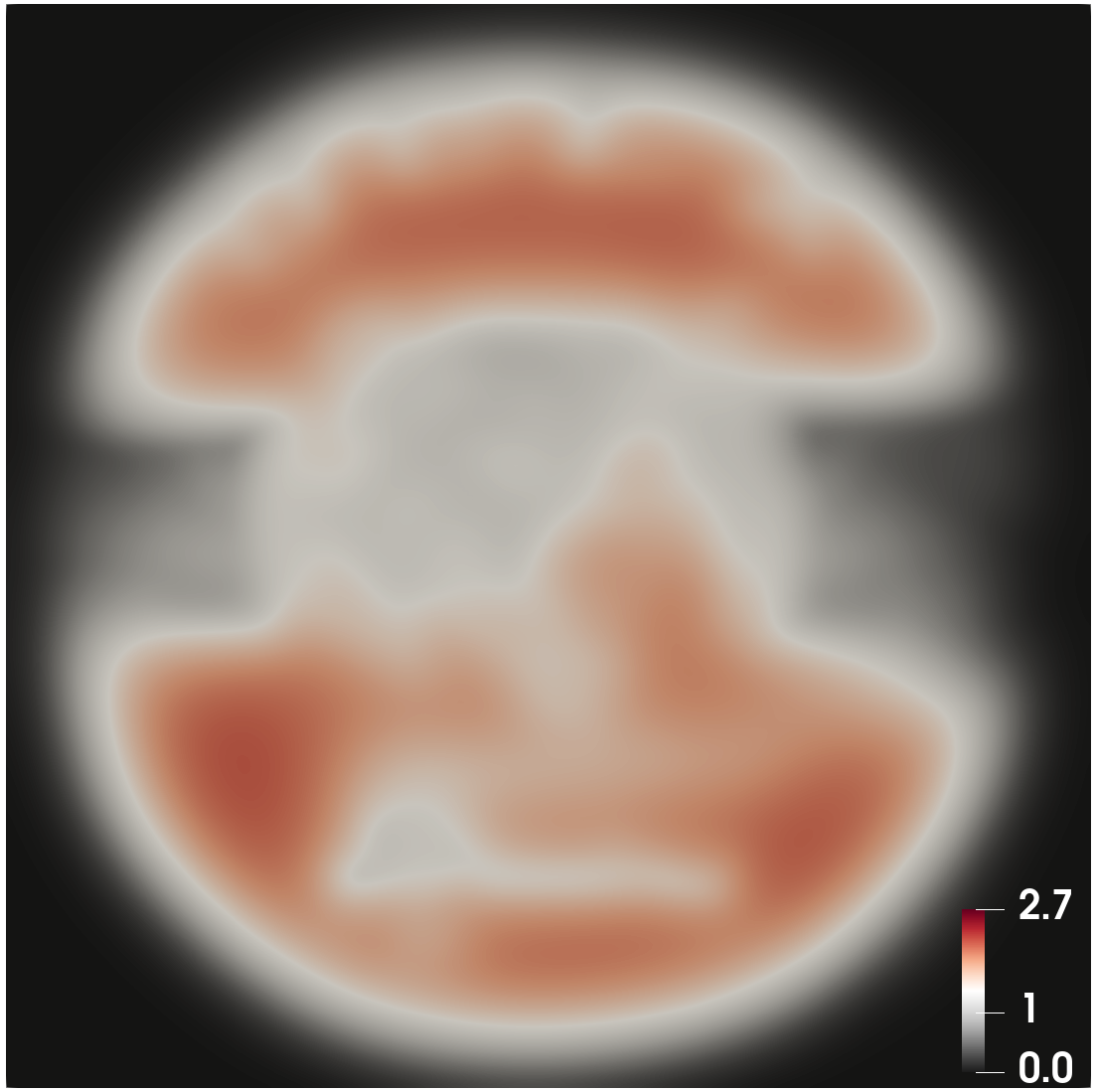}
\includegraphics[width=0.192\textwidth]{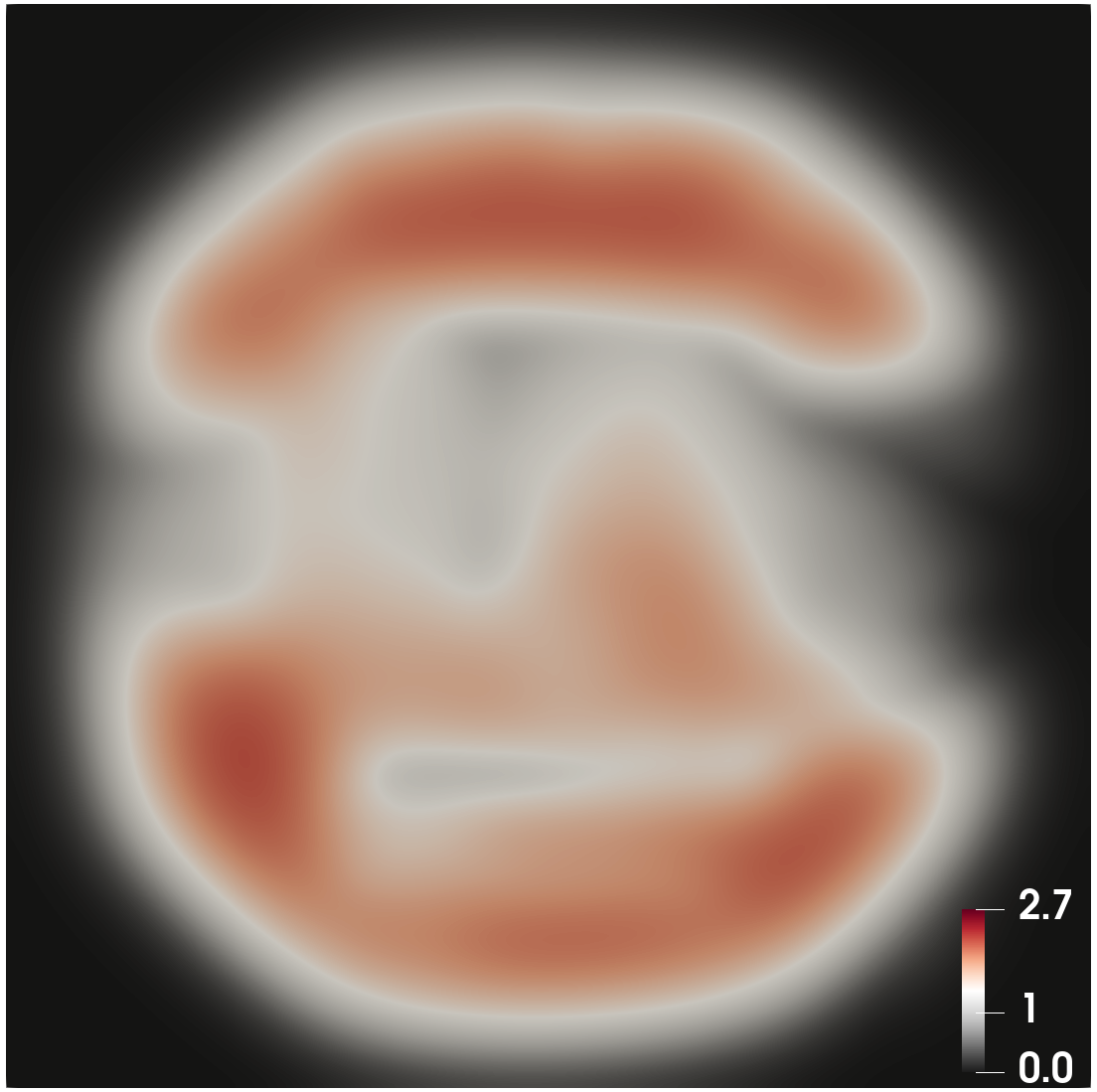}
\includegraphics[width=0.192\textwidth]{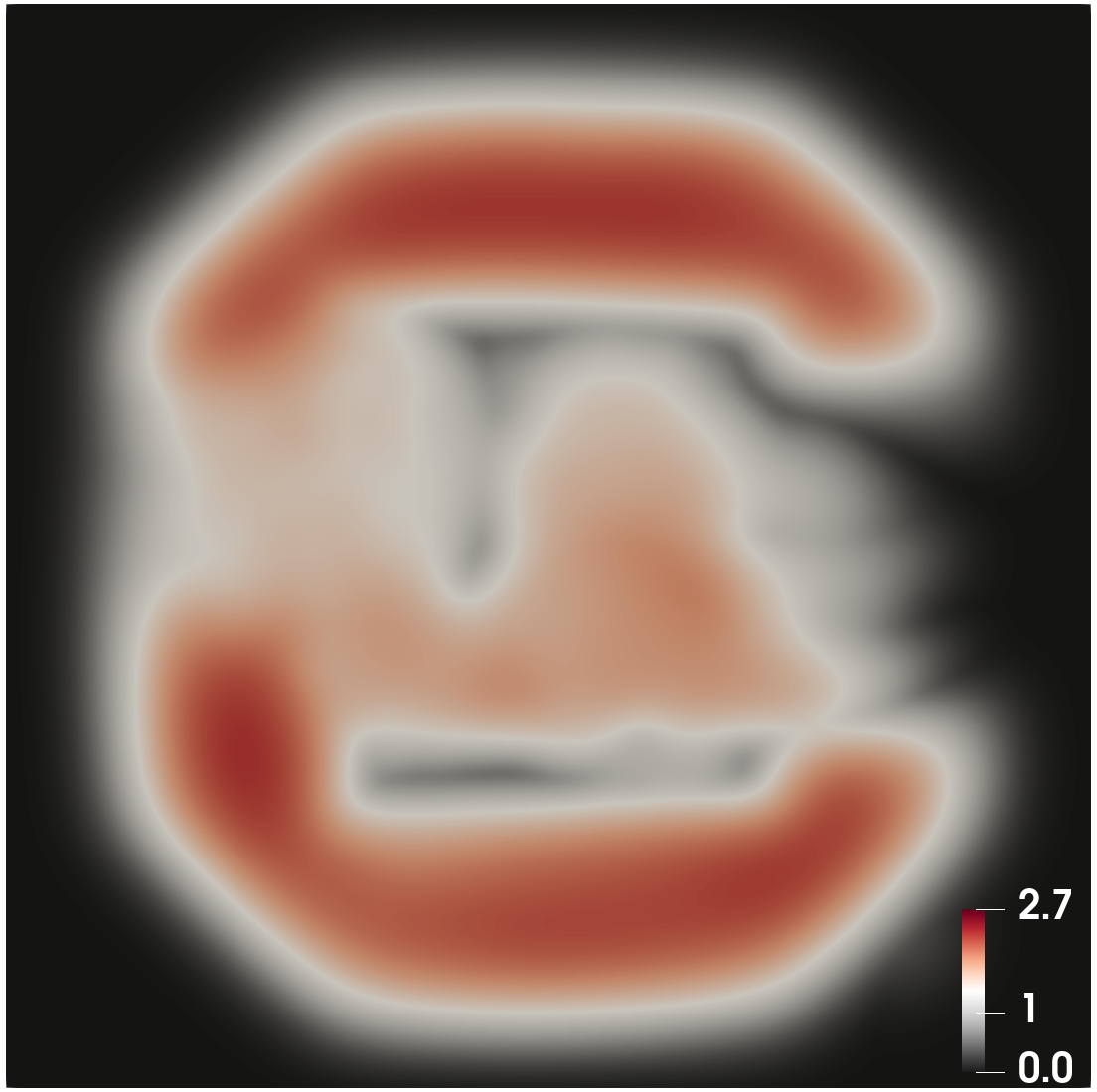}
\includegraphics[width=0.192\textwidth]{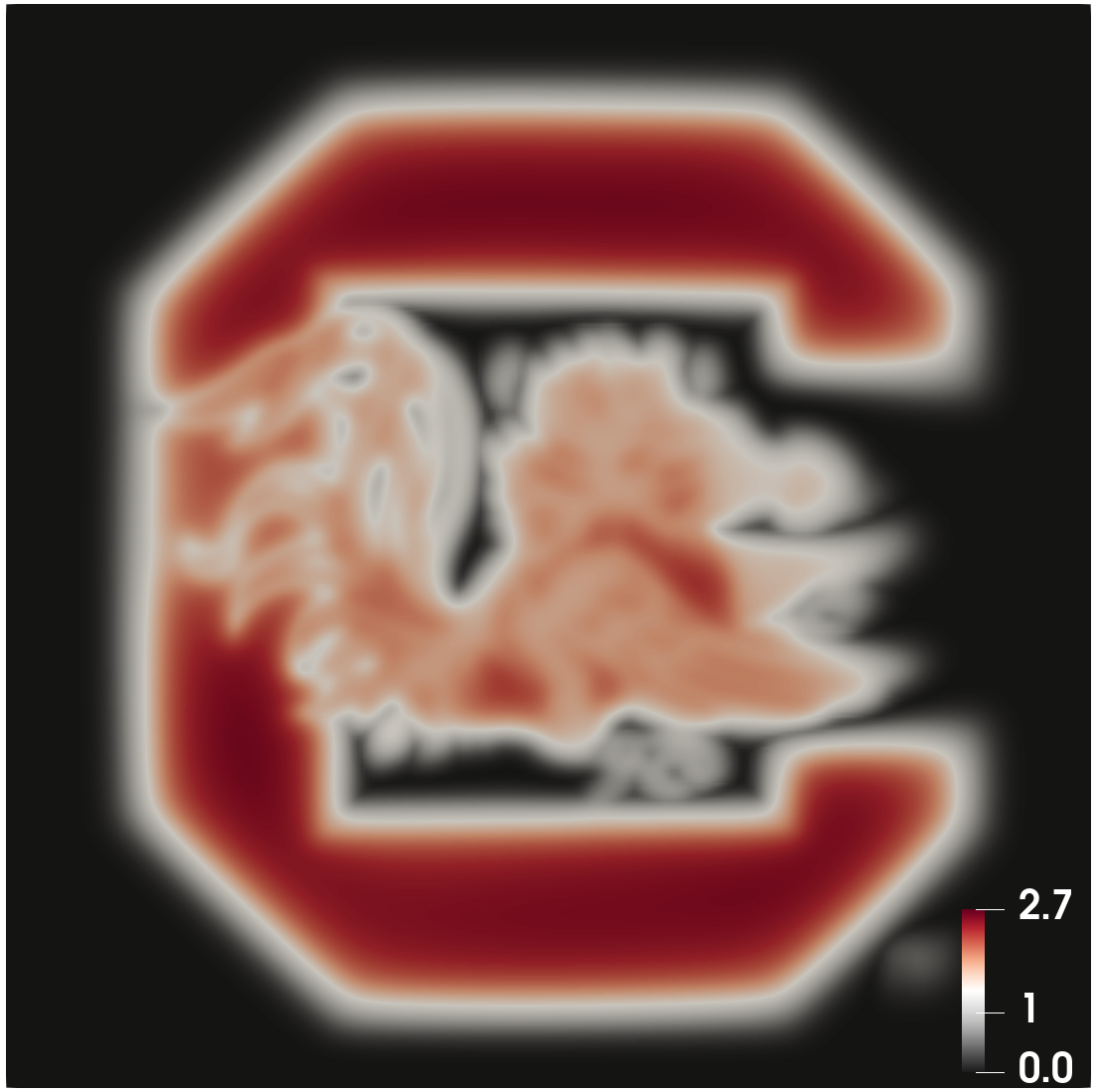}
}

\subfigure[Case 3: $A(\rho) = 0.01/\rho$. UCLA $\rightarrow$ USC]
{
\includegraphics[width=0.192\textwidth]{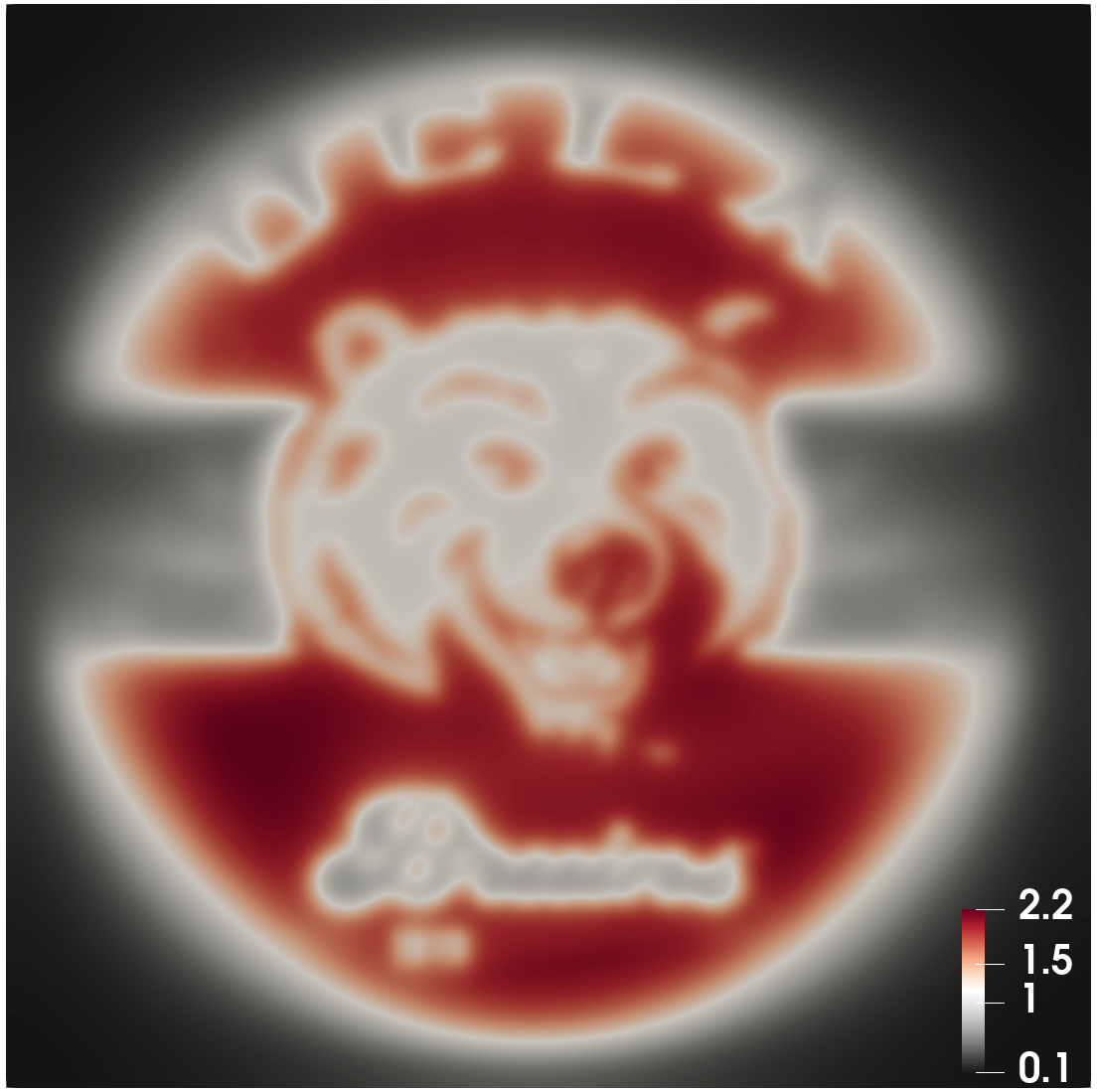}
\includegraphics[width=0.192\textwidth]{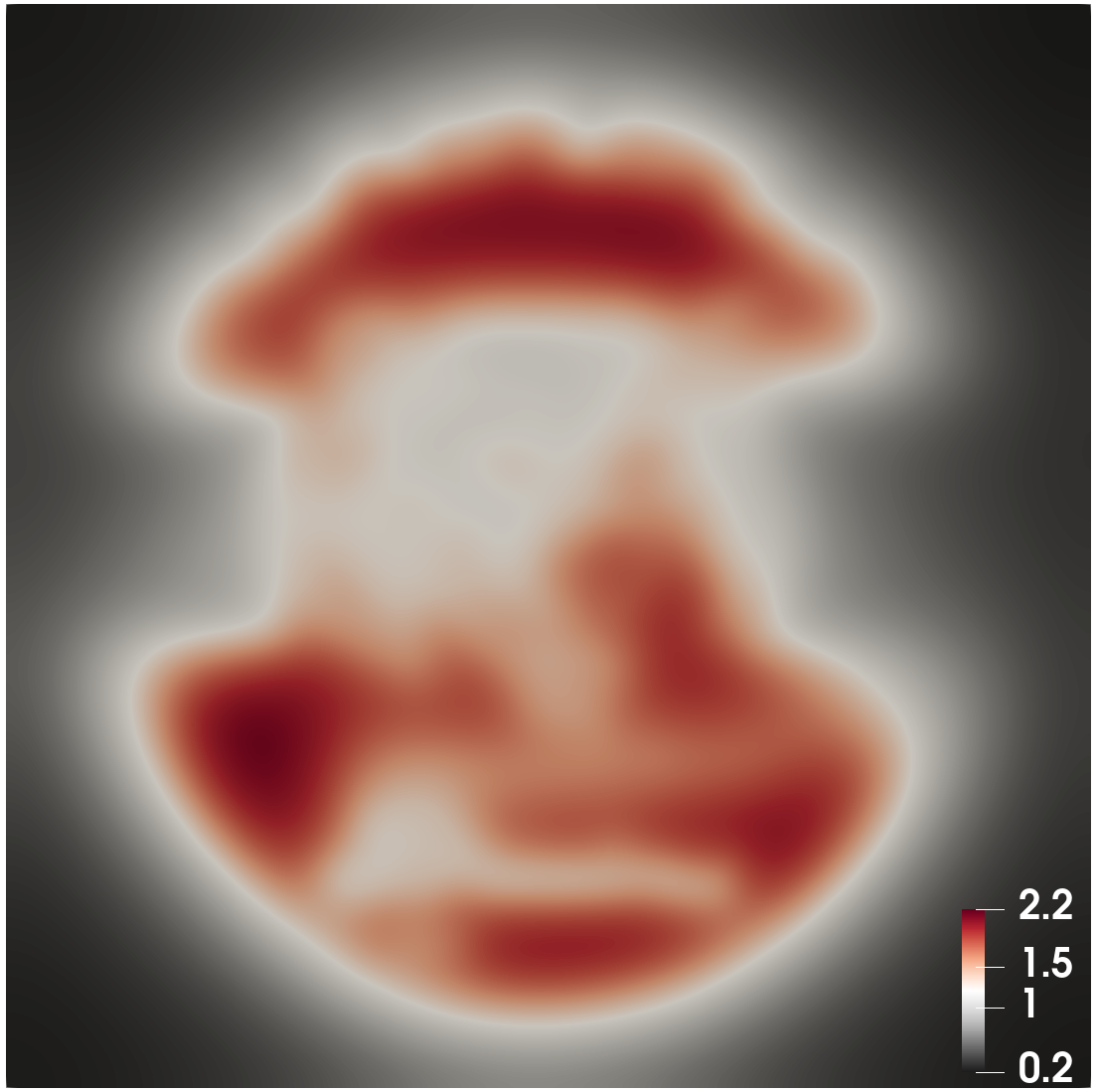}
\includegraphics[width=0.192\textwidth]{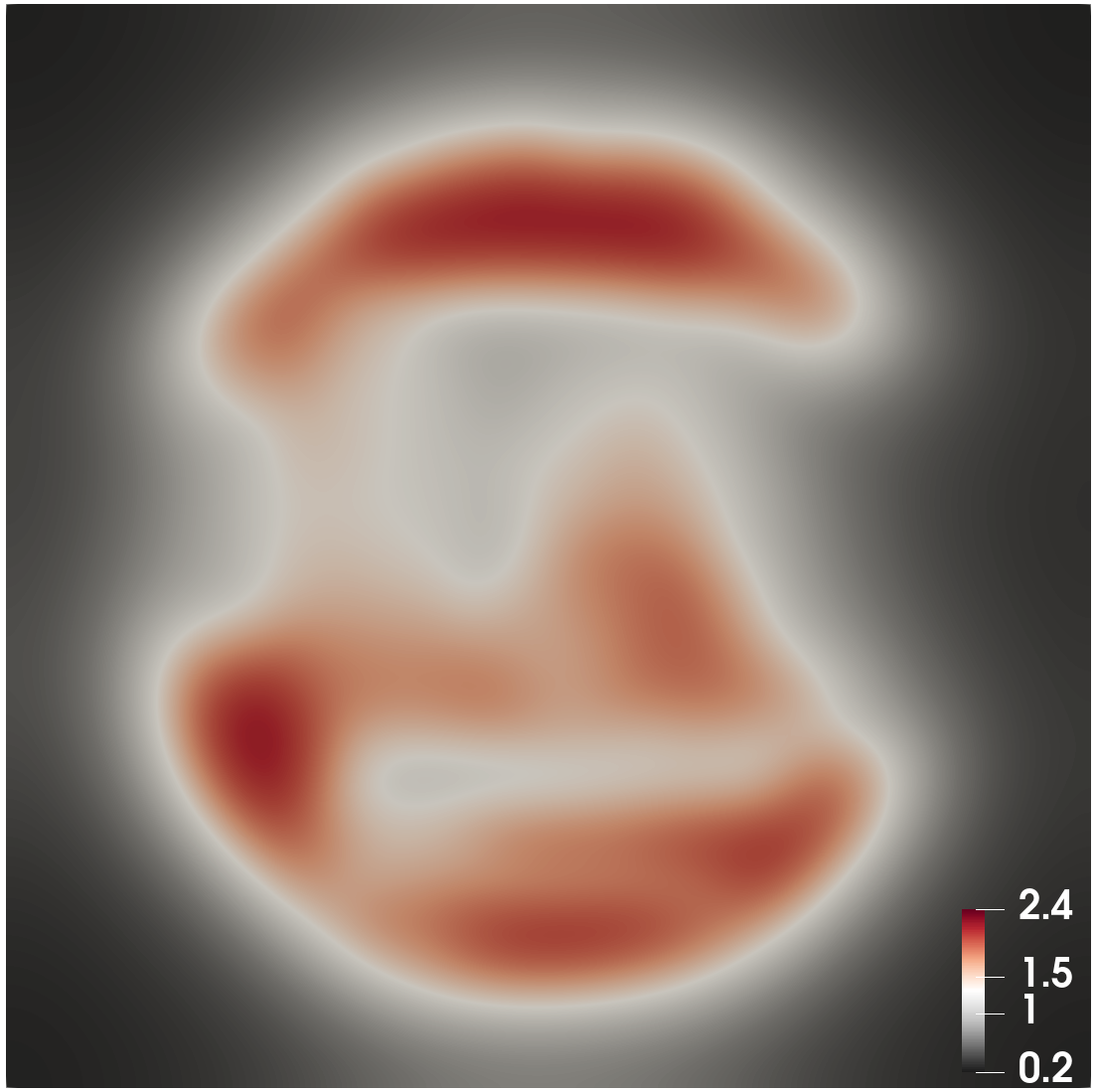}
\includegraphics[width=0.192\textwidth]{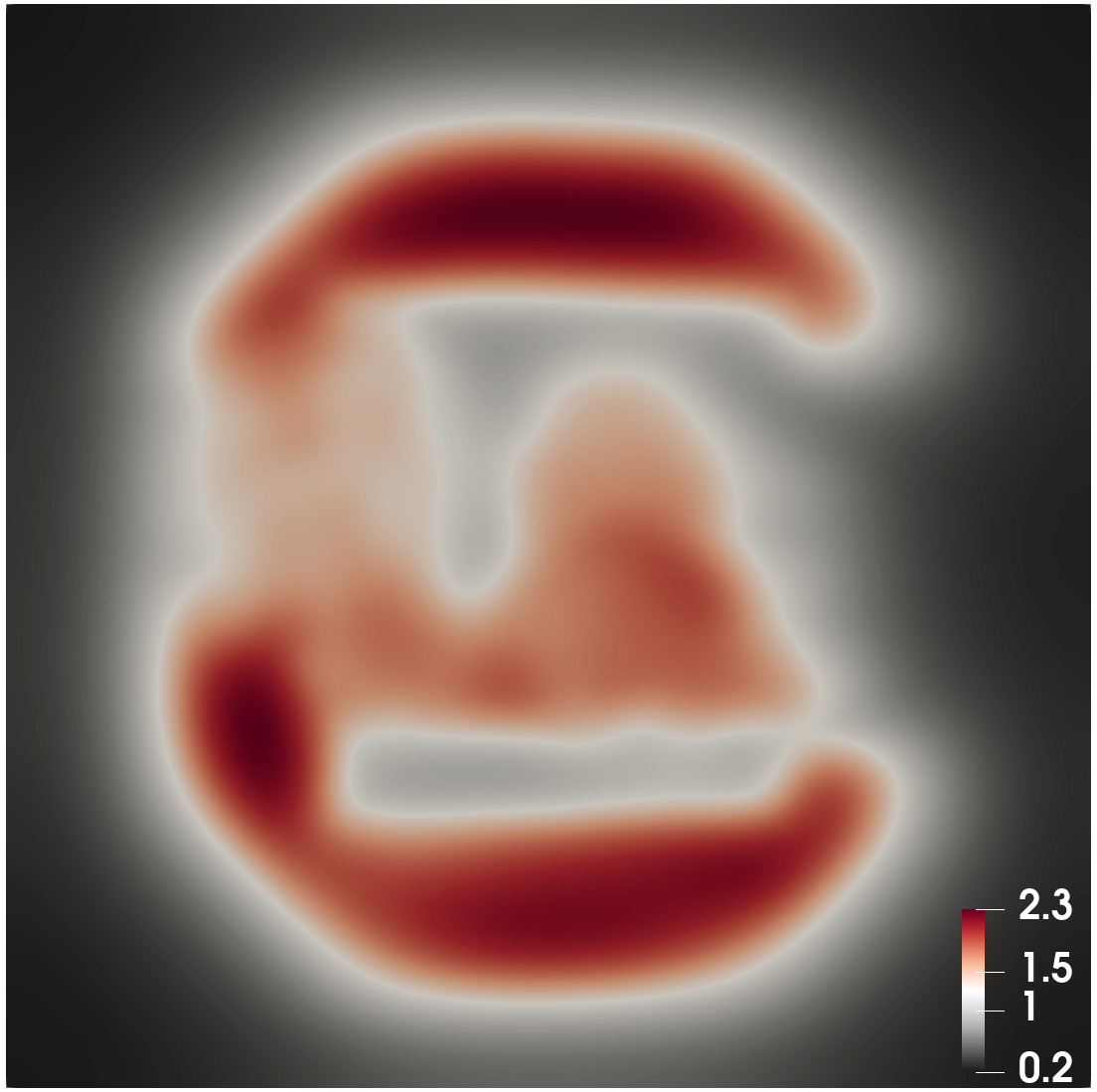}
\includegraphics[width=0.192\textwidth]{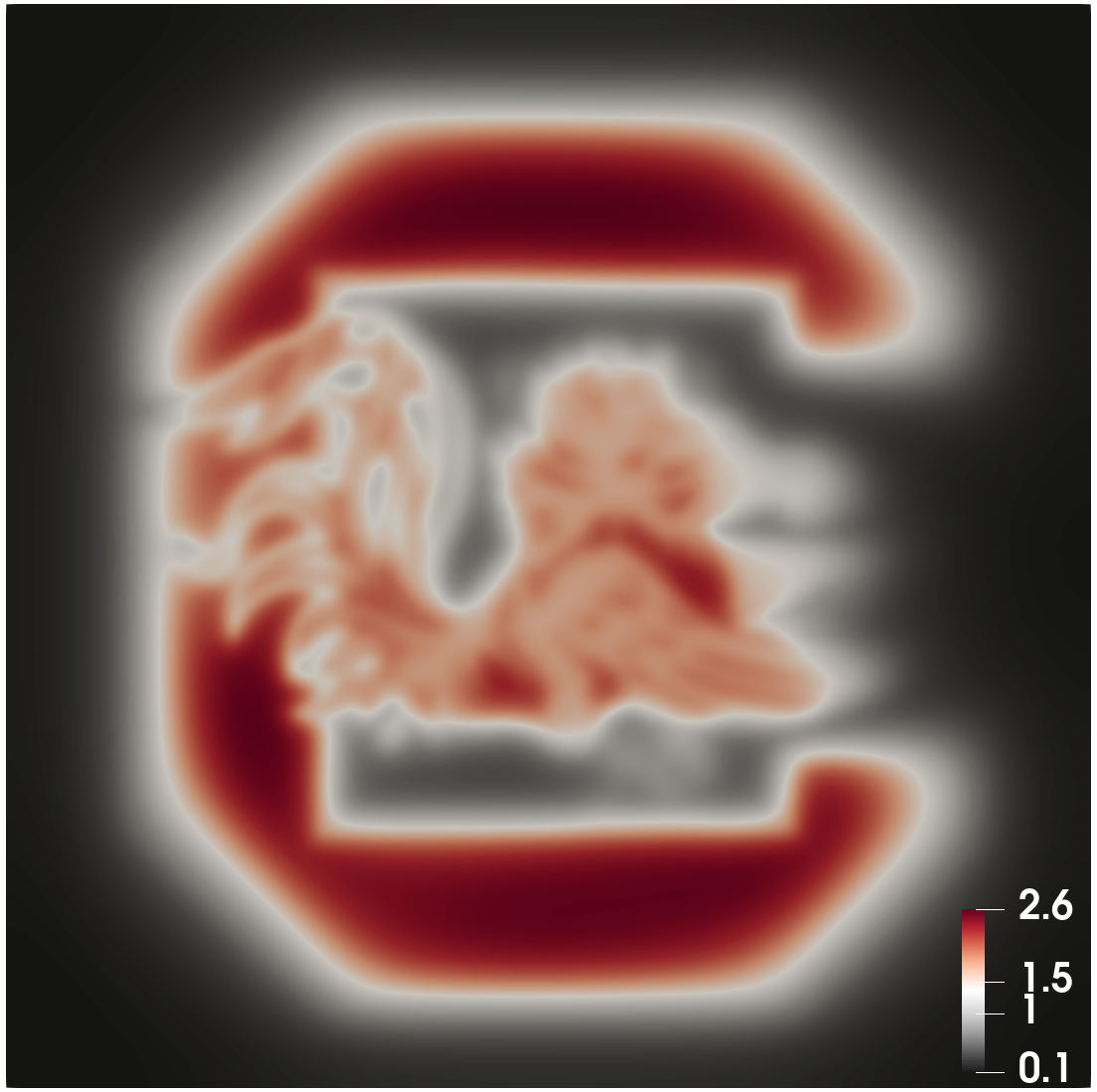}
}

\caption{Example \ref{ex4}. Initial density: UCLA. Terminal density: USC.
Snapshots of $\rho$ at 
$t=$ 0.1,0.3,0.5,0.7,0.9 (left to right).
}
\label{fig:den-case4B}
\end{figure}

\begin{figure}[tb]
\centering
\subfigure[Case 1: $A(\rho) = 0$. USC $\rightarrow$ ND]
{
\includegraphics[width=0.192\textwidth]{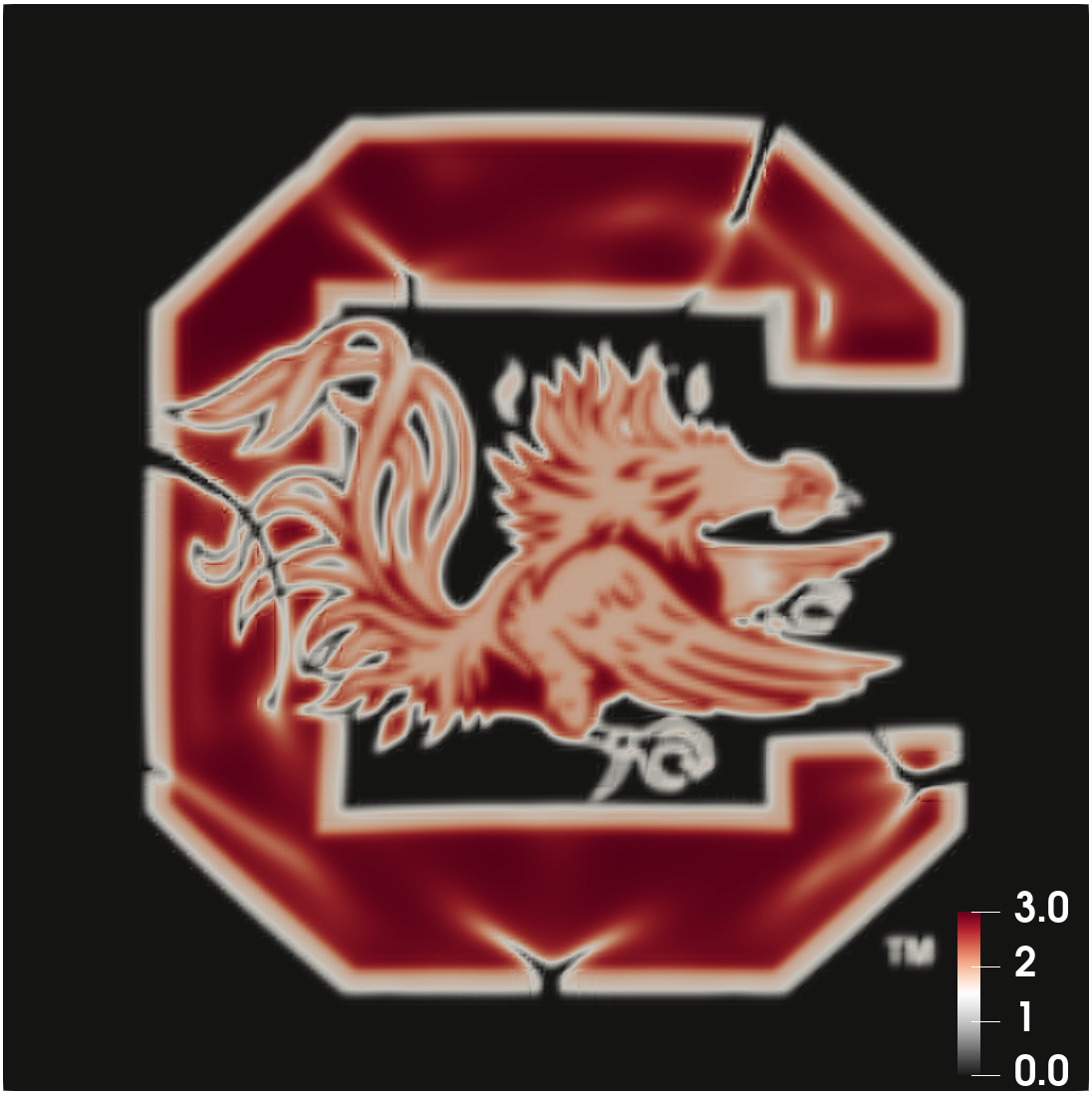}
\includegraphics[width=0.192\textwidth]{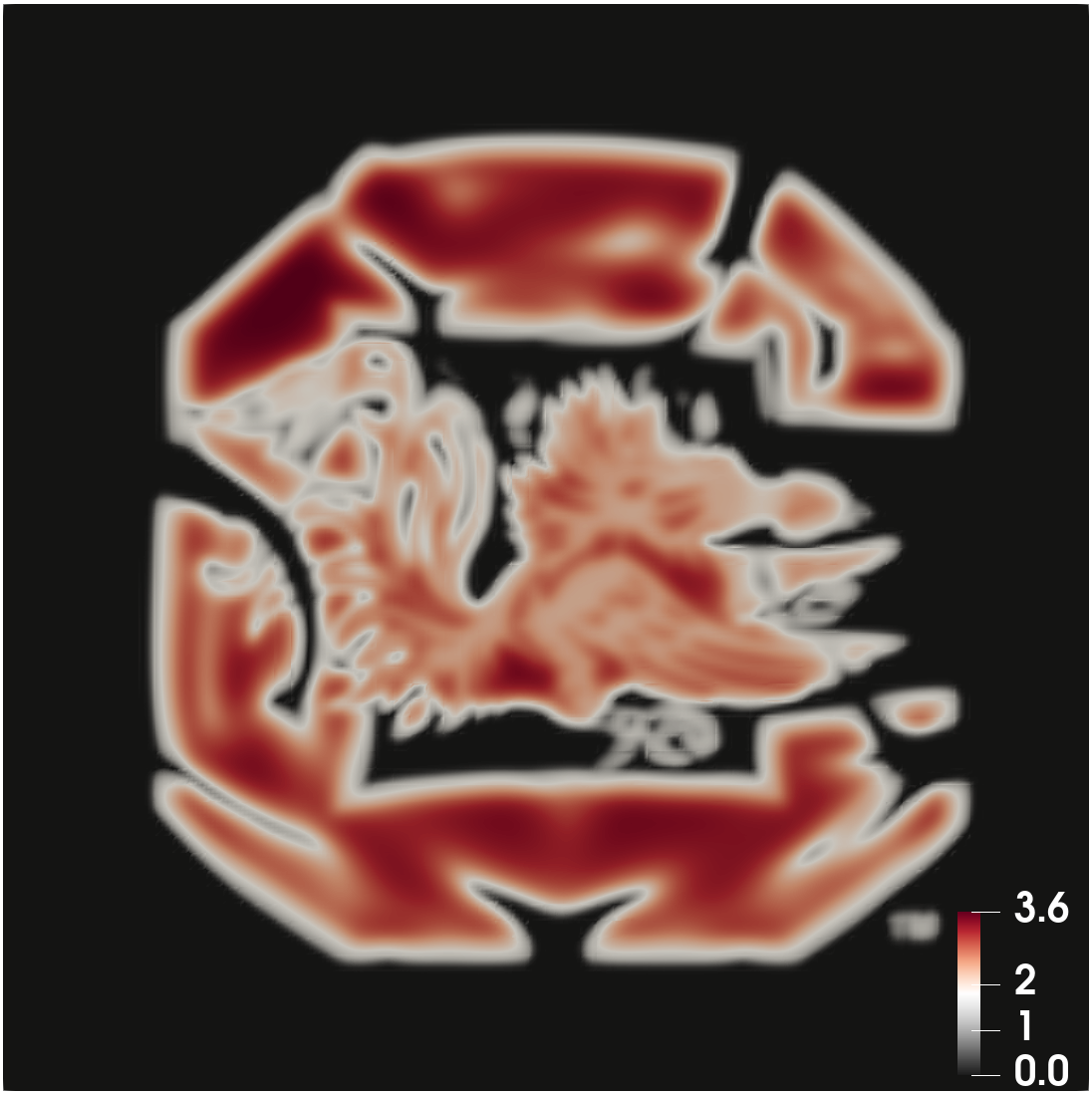}
\includegraphics[width=0.192\textwidth]{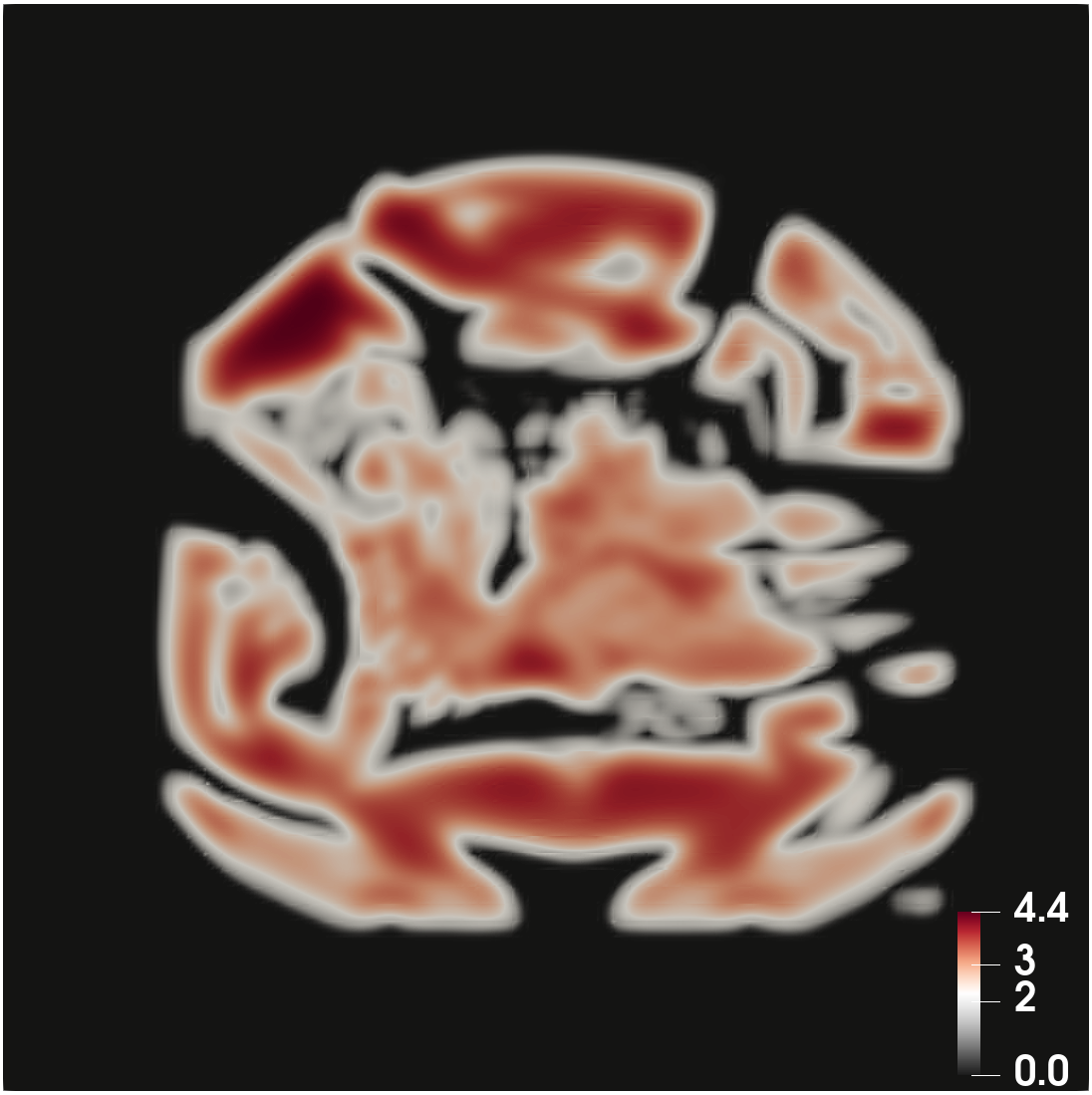}
\includegraphics[width=0.192\textwidth]{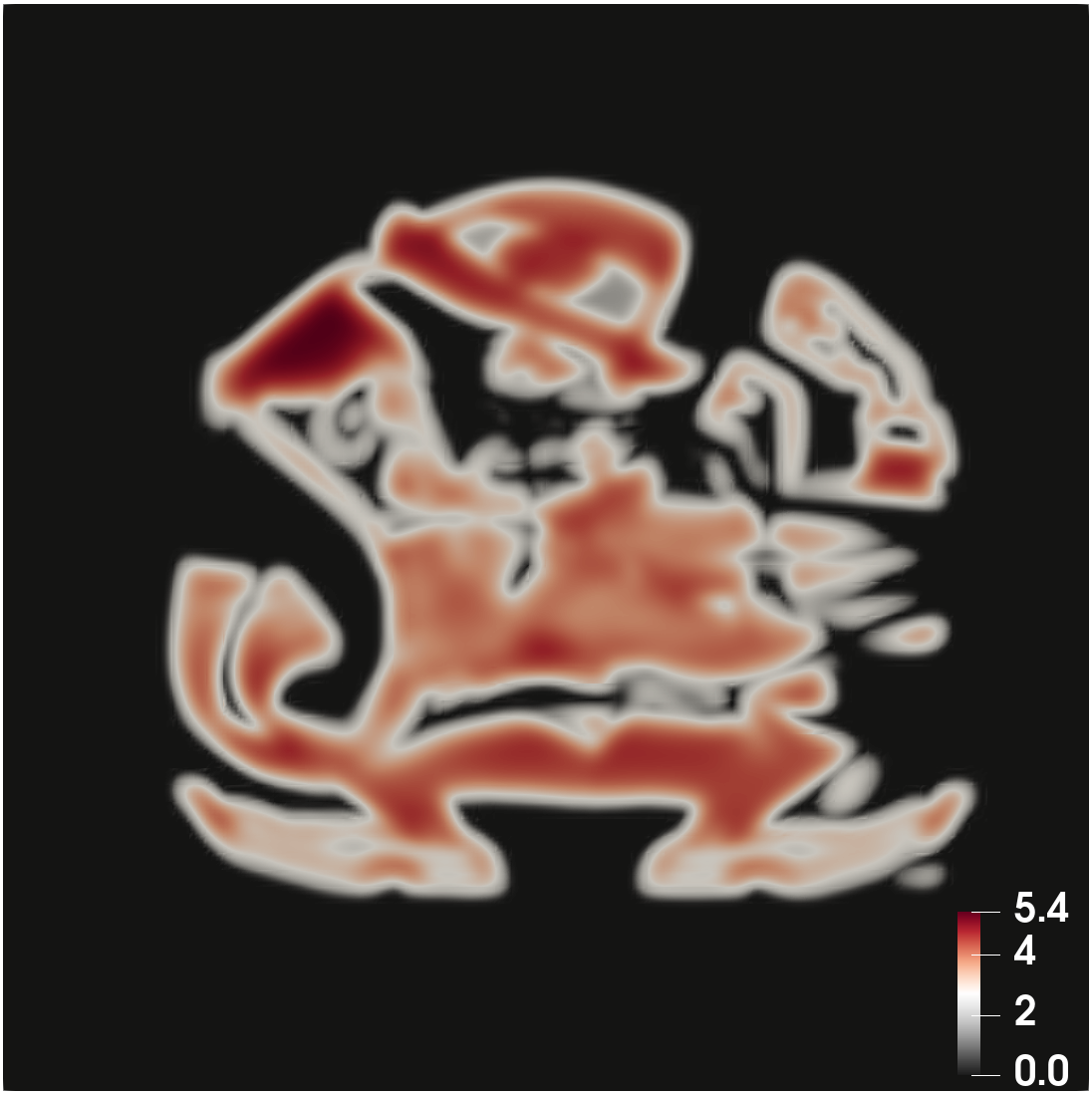}
\includegraphics[width=0.192\textwidth]{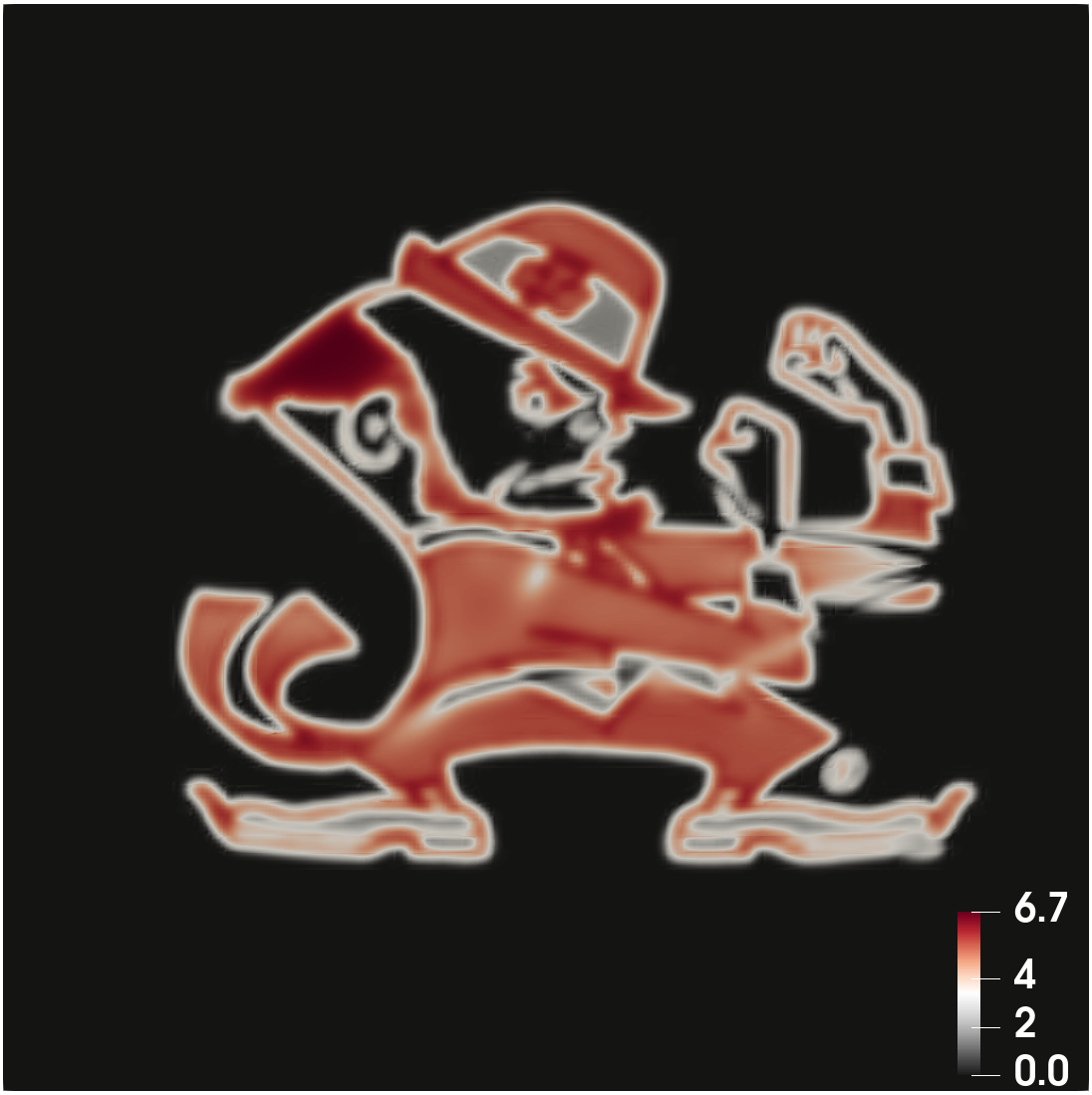}
}
\subfigure[Case 2: $A(\rho) = 0.01\rho\log(\rho)$.  USC $\rightarrow$ ND]
{
\includegraphics[width=0.192\textwidth]{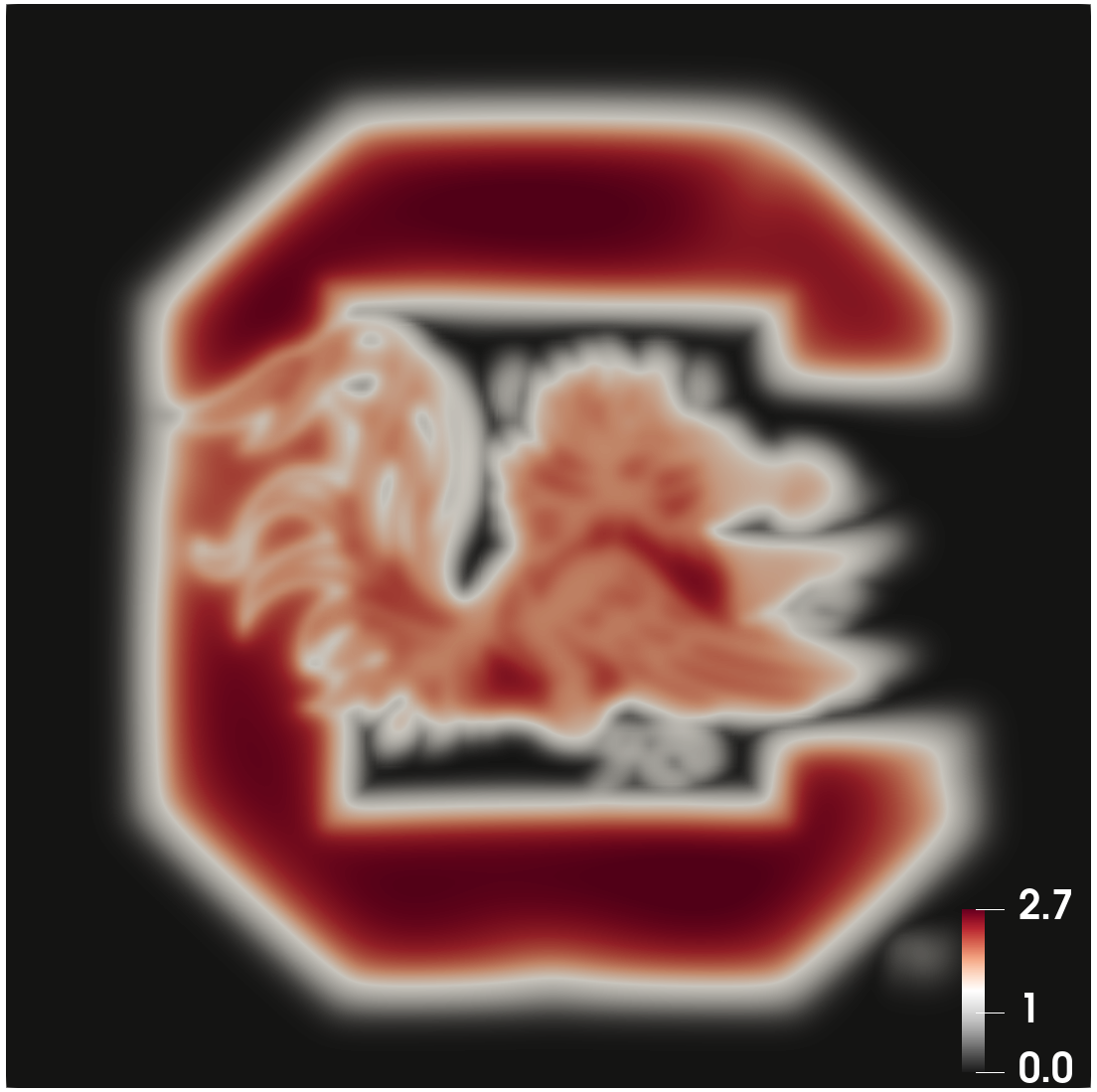}
\includegraphics[width=0.192\textwidth]{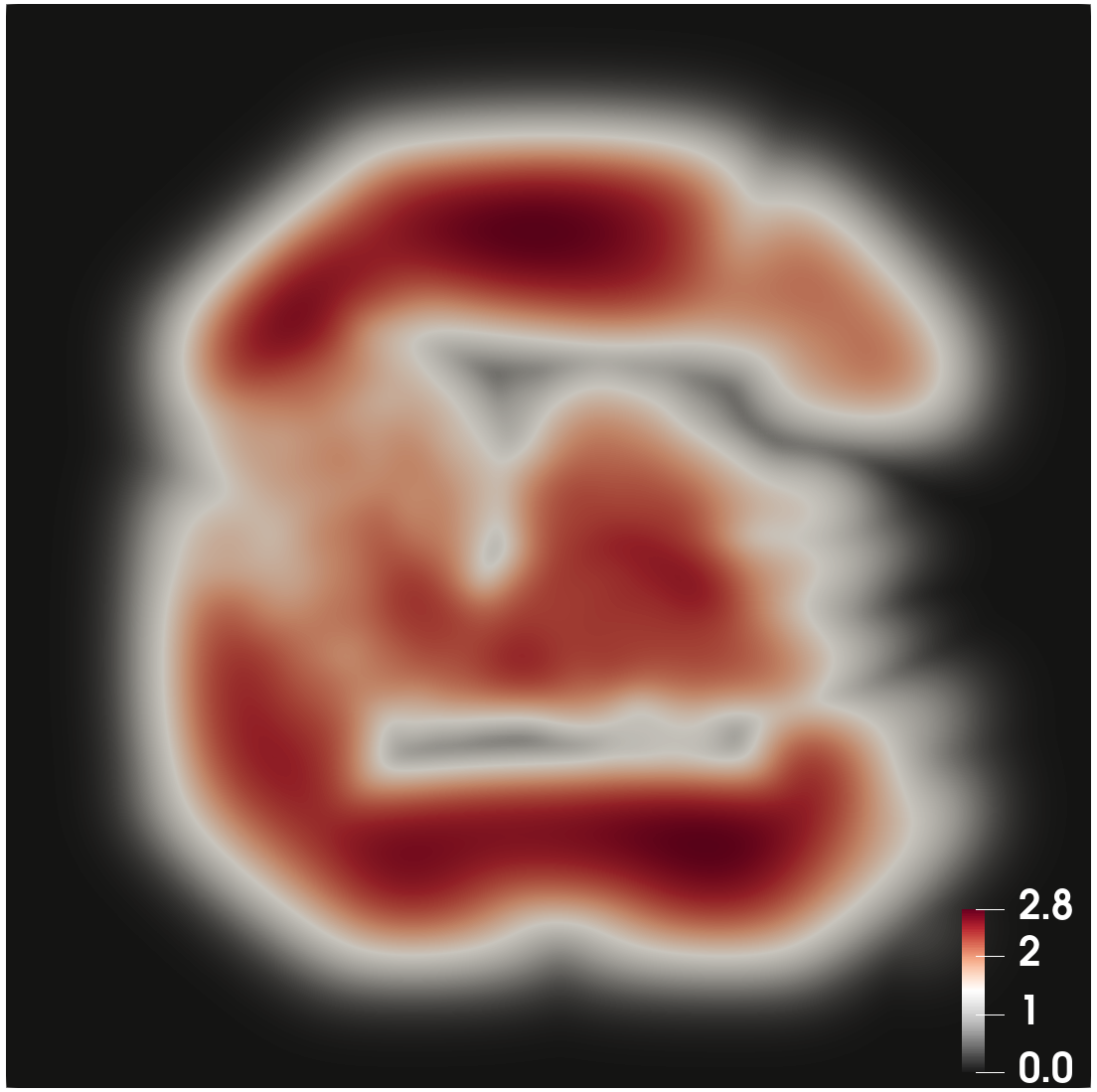}
\includegraphics[width=0.192\textwidth]{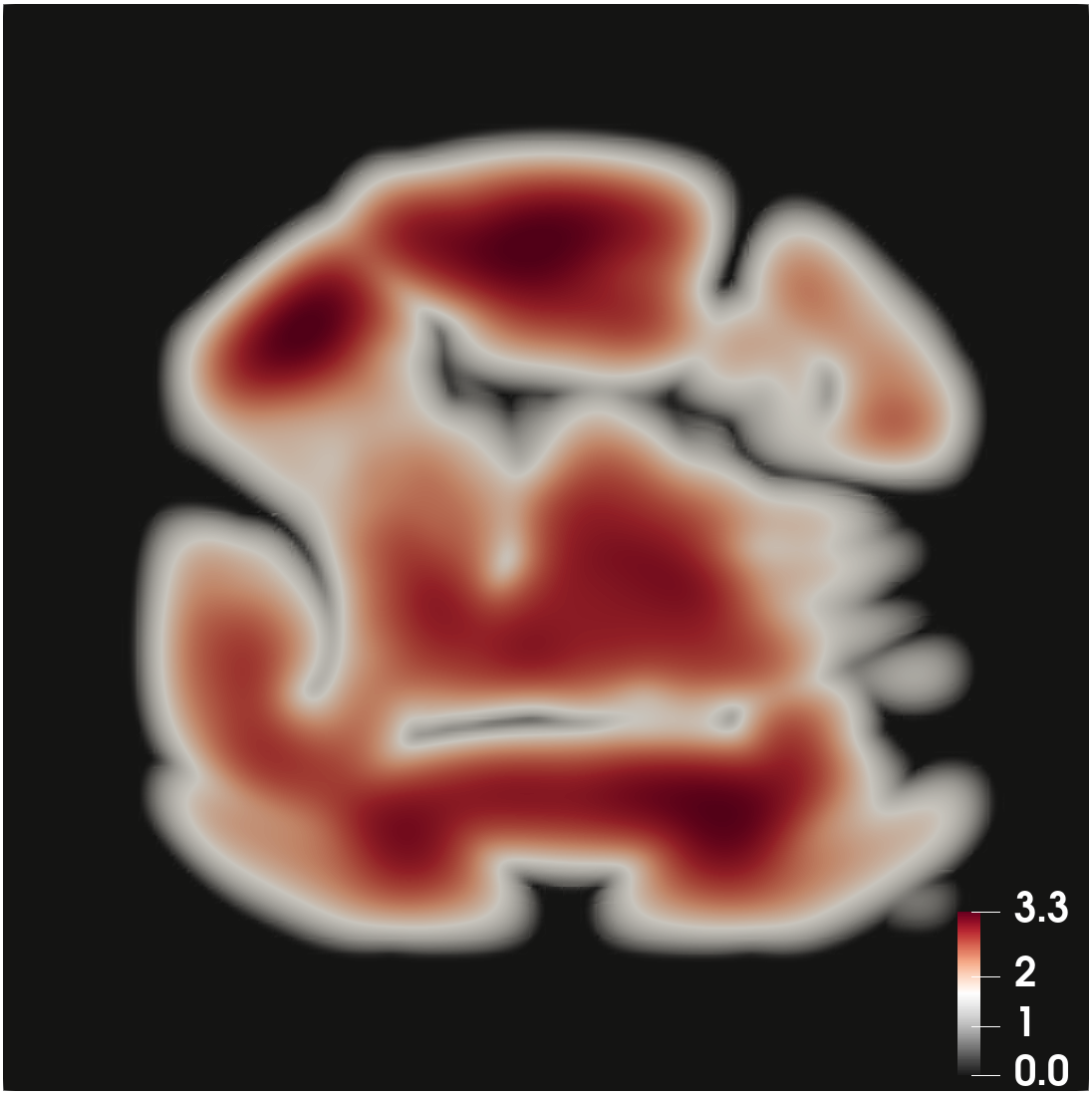}
\includegraphics[width=0.192\textwidth]{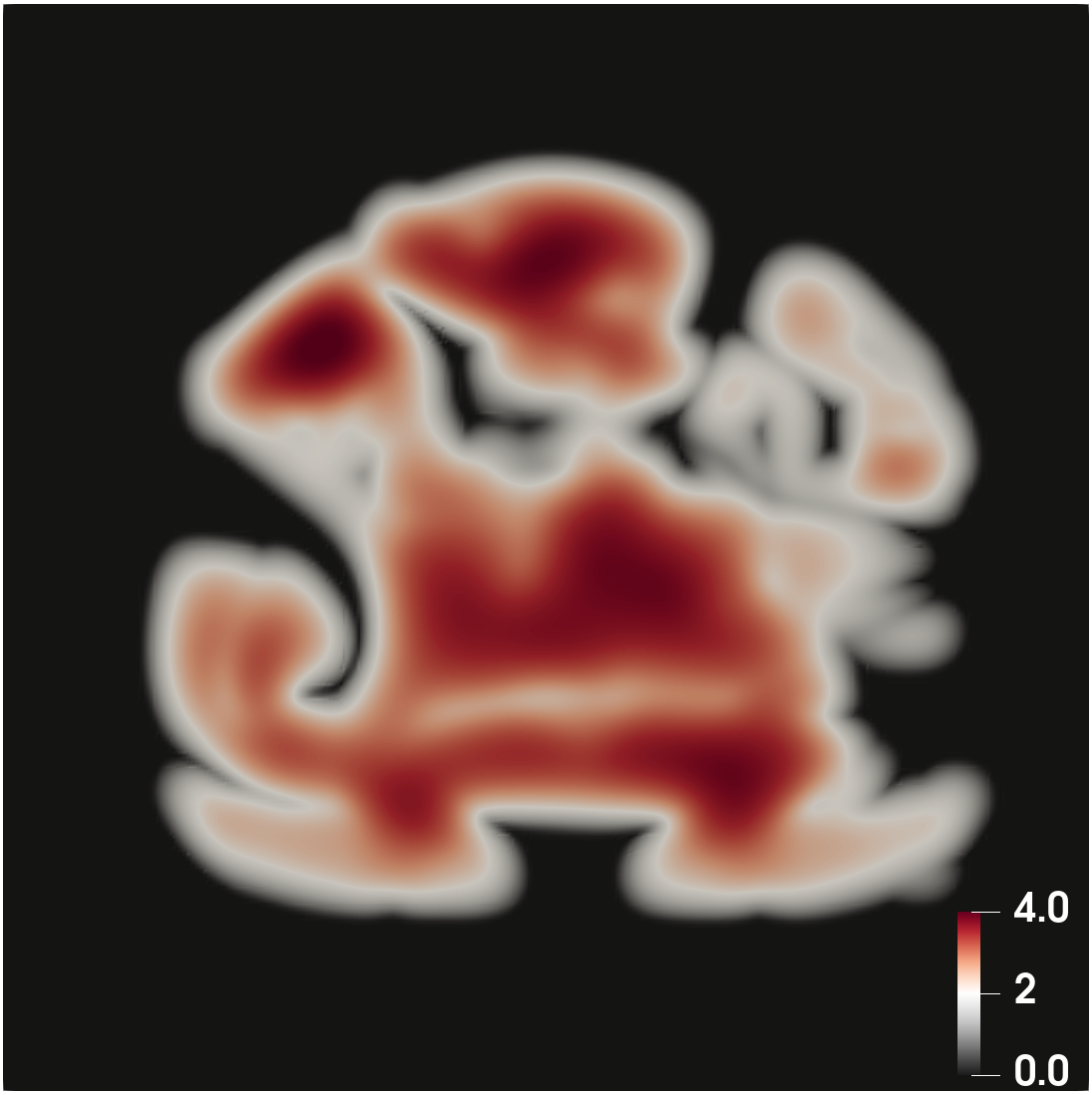}
\includegraphics[width=0.192\textwidth]{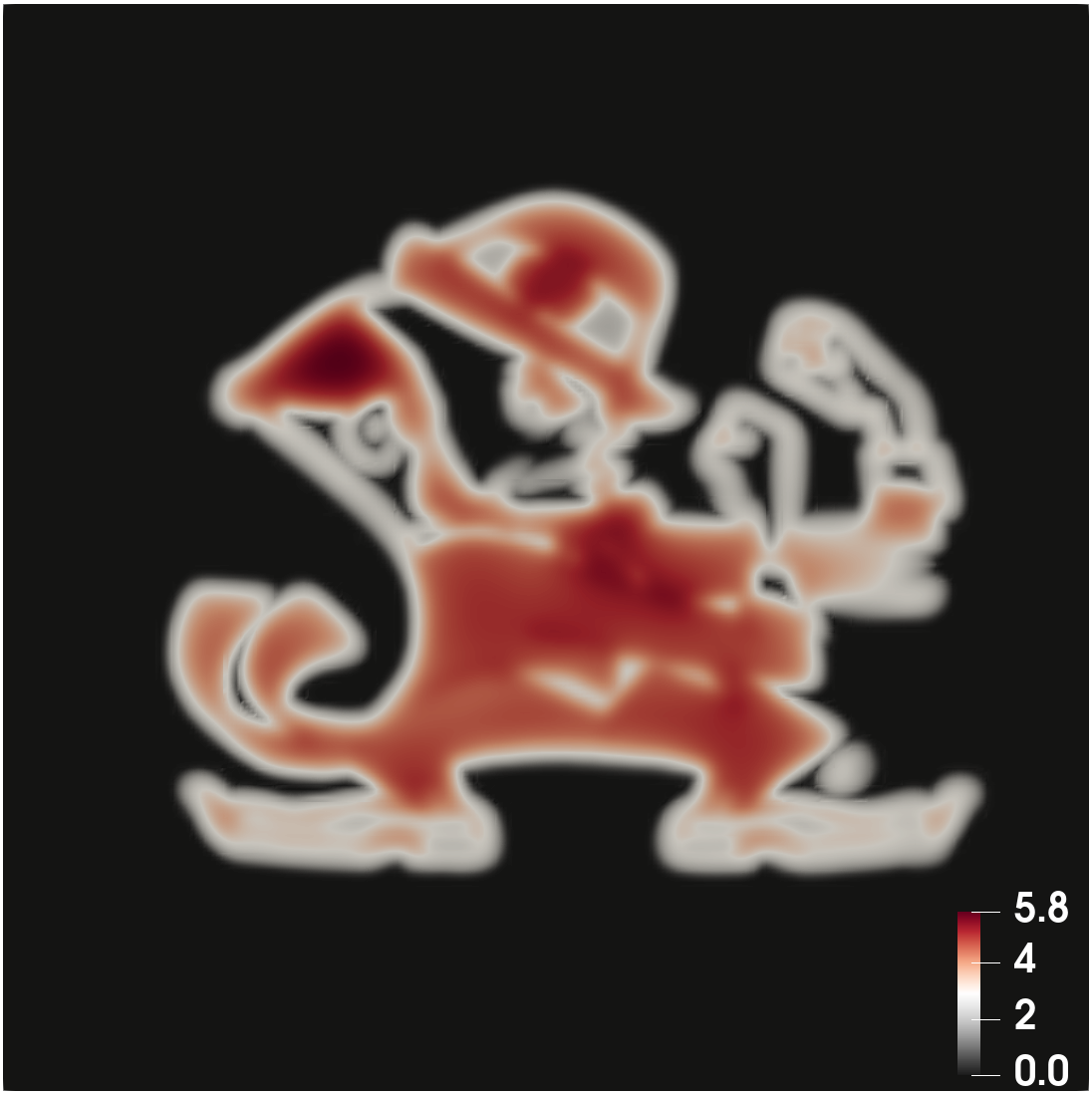}
}

\subfigure[Case 3: $A(\rho) = 0.01/\rho$. USC $\rightarrow$ ND]
{
\includegraphics[width=0.192\textwidth]{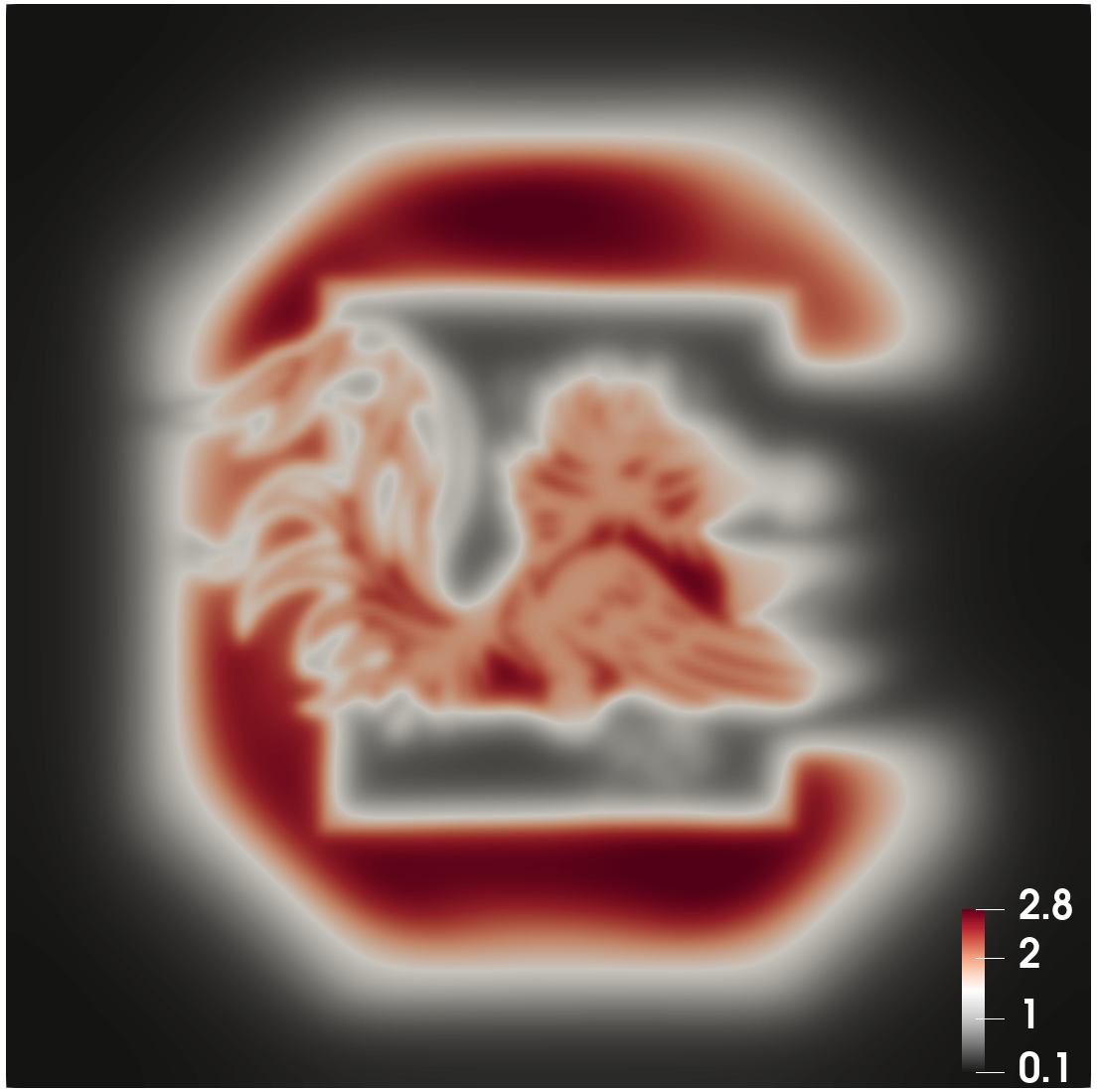}
\includegraphics[width=0.192\textwidth]{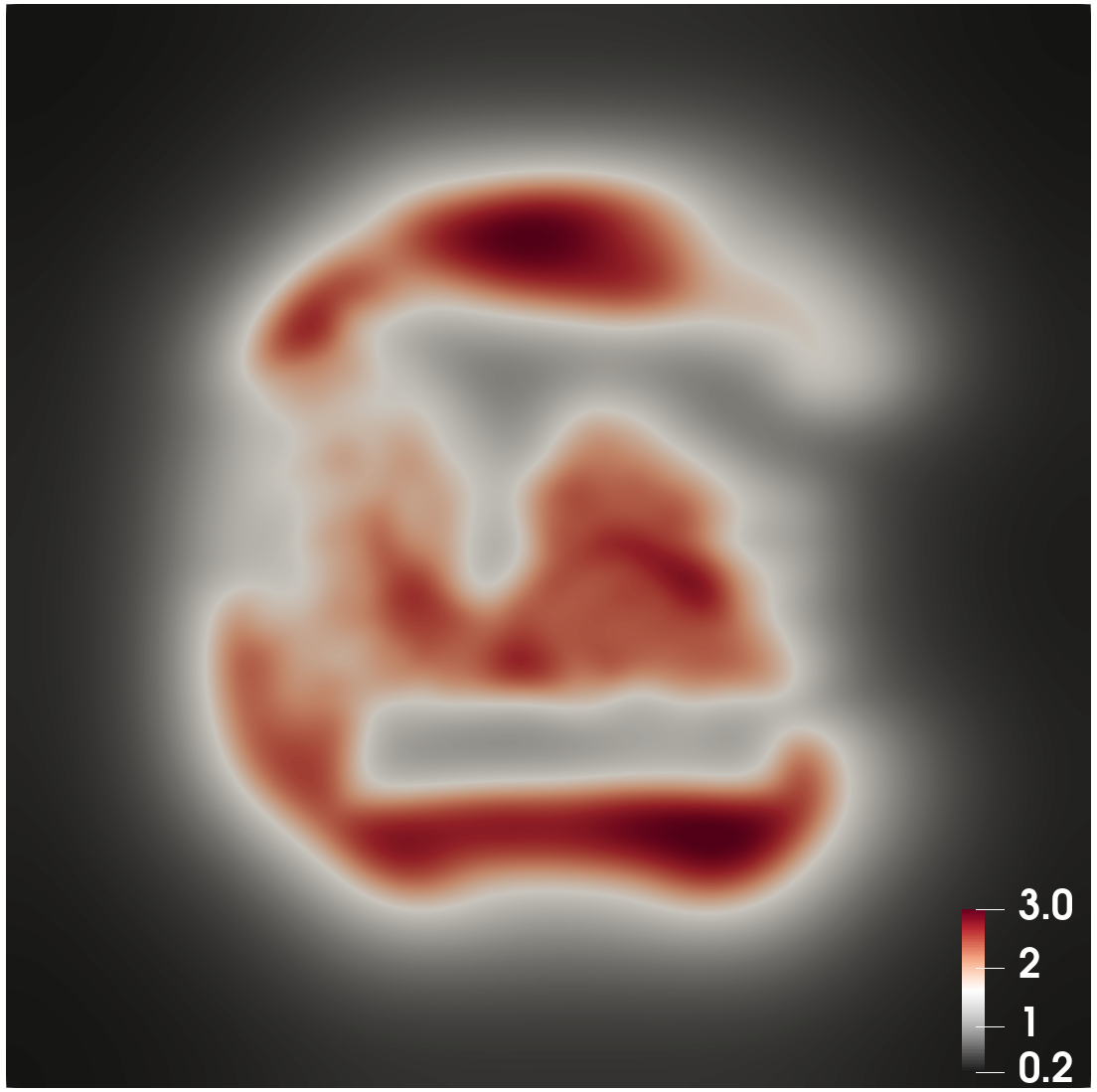}
\includegraphics[width=0.192\textwidth]{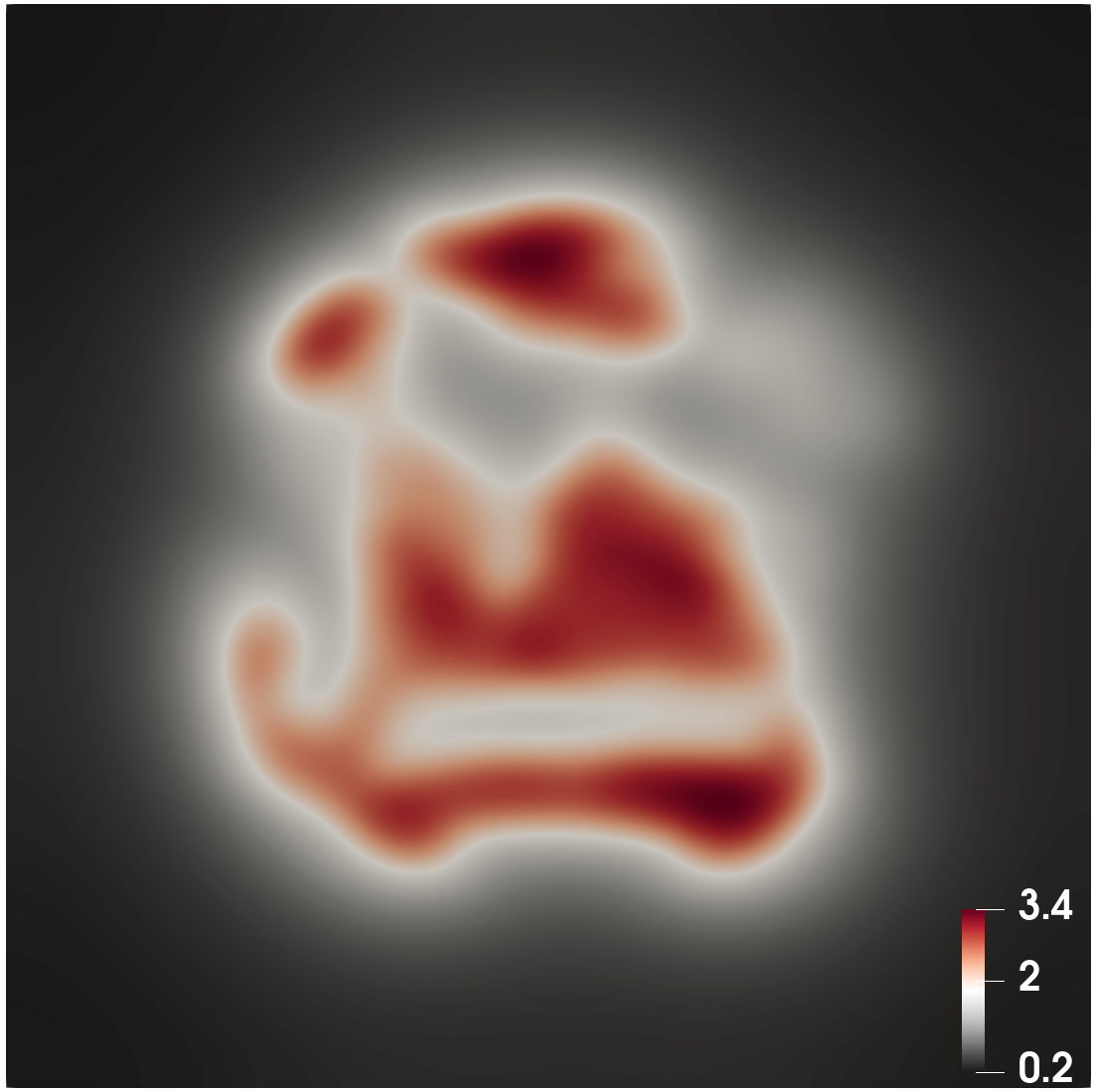}
\includegraphics[width=0.192\textwidth]{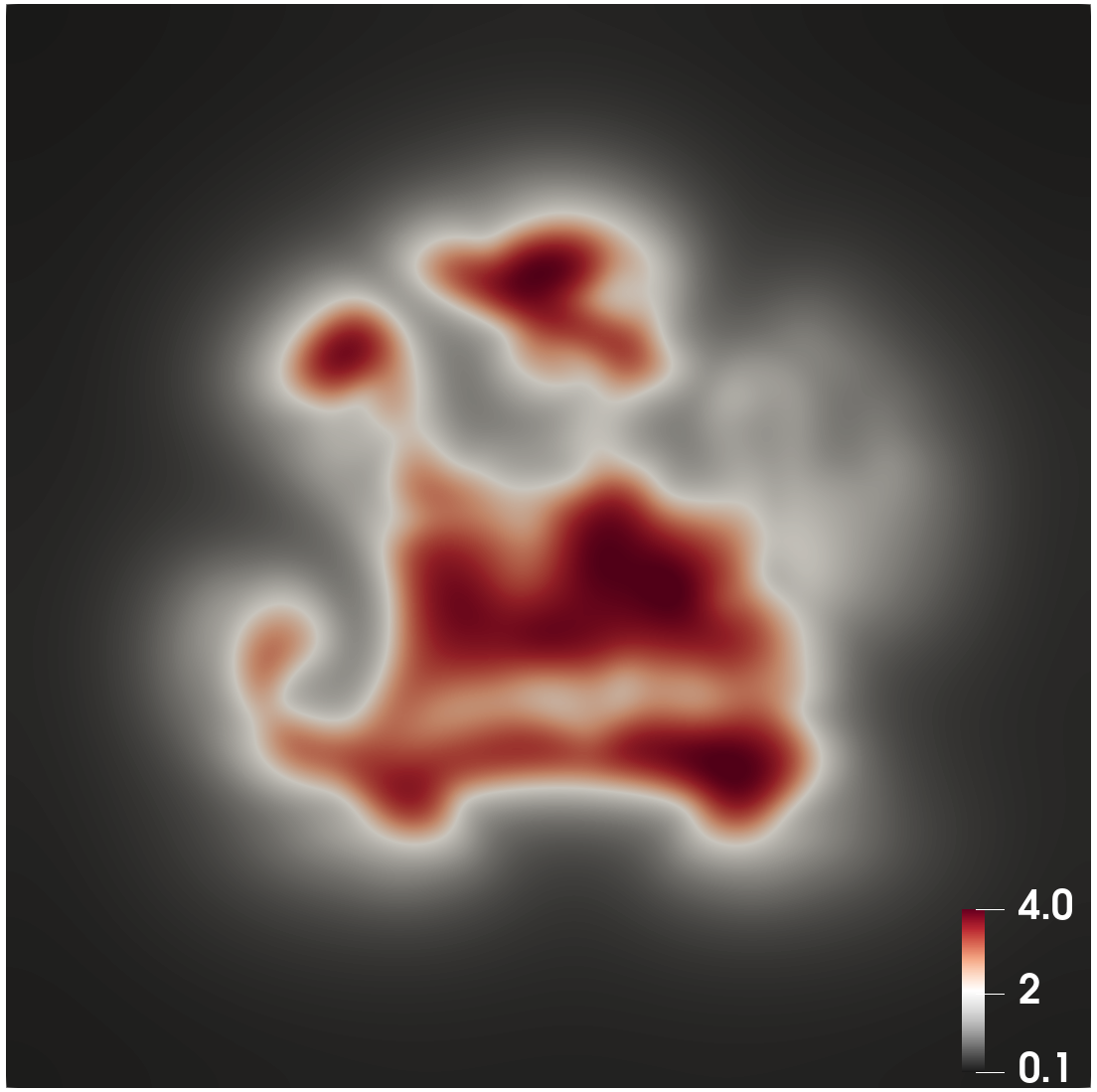}
\includegraphics[width=0.192\textwidth]{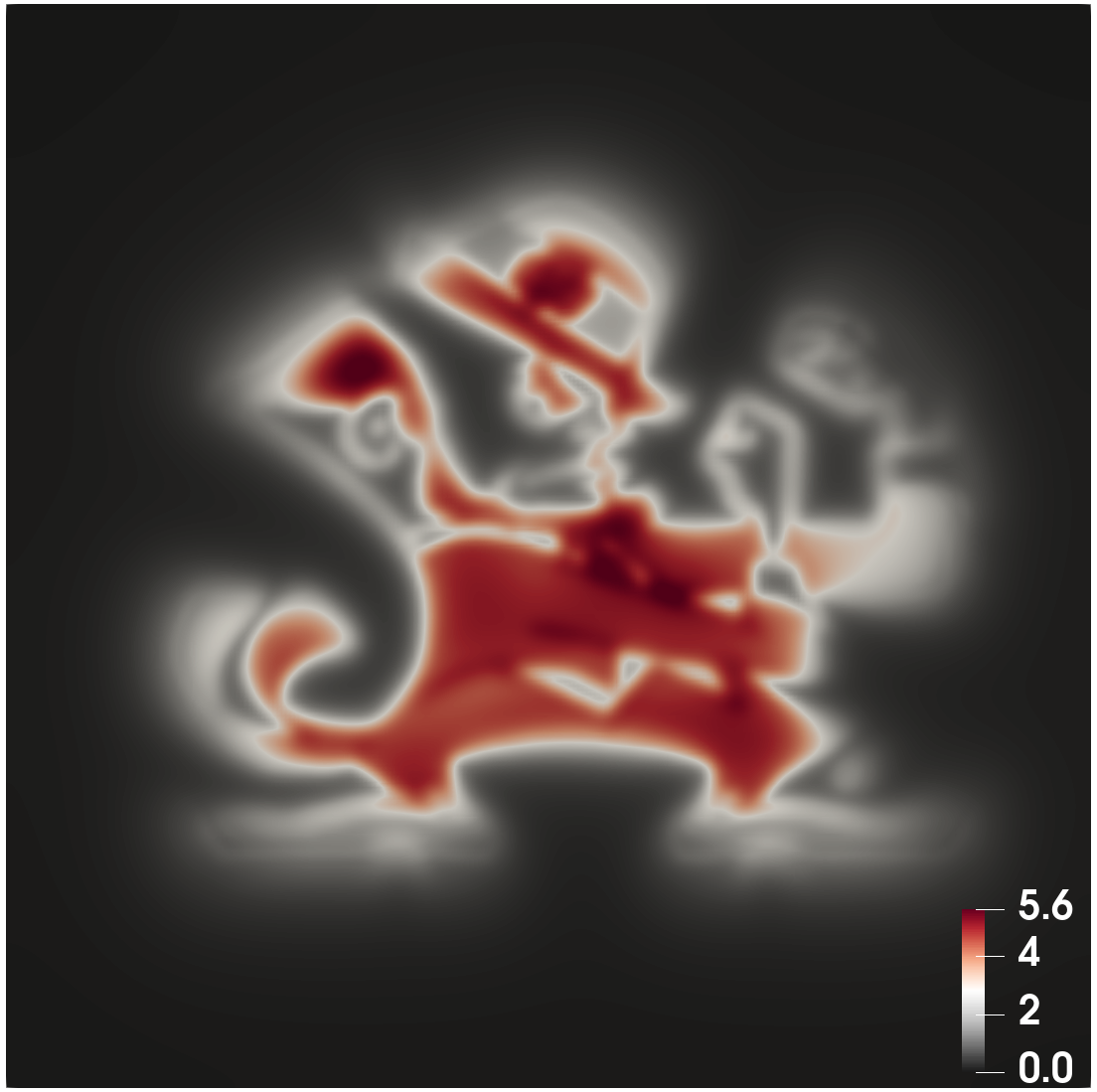}
}

\caption{Example \ref{ex4}. Initial density: USC. Terminal density: ND.
Snapshots of $\rho$ at 
$t=$ 0.1,0.3,0.5,0.7,0.9 (left to right).
}
\label{fig:den-case4C}
\end{figure}

Snapshots of the density contour at different times are shown in Figure~\ref{fig:den-case4A} for (i) ND$\rightarrow$ UCLA, in Figure~\ref{fig:den-case4B} for (ii) UCLA$\rightarrow$ USC,
and in Figure~\ref{fig:den-case4C} for (iii) USC$\rightarrow$ ND. 
We observe in these figures that Case 1 (OT) produce the most sharp results for the density evolution, and that  both interaction costs in Case 2/3 have a strong smoothing effect which blur the density profile, where Case 3 with $A(\rho)=0.01/\rho$ also leads to an everywhere positive density.

\section{Conclusion}
\label{sec:summary}
This paper applies high-order accurate finite element methods to compute optimal transport (OT) and mean field games (MFG). To our best knowledge, it is the first time to apply high order numerical methods in OT and MFGs. We verify the accuracy of algorithms through numerical examples. In future works, we shall investigate the numerical property of high-order accuracy FEM methods in OT and MFG-related dynamics. We expect they will have vast applications in computational physics, social science, biology modeling, pandemics control, and computer vision. We also expect to apply high order FEM in generalized mean field control formalisms to compute implicit-in-time fluid dynamics~\cite{li2022computational,li2021controlling,li2022controlling,liu2023primal}.

\bibliographystyle{amsplain}
\bibliography{references}
\end{document}